\documentclass{article}

\usepackage{arxiv}

\usepackage{fancyhdr,lipsum}
\makeatletter

\usepackage[cp1252]{inputenc}
\usepackage[T1]{fontenc}    
\usepackage{hyperref}       
\usepackage{url}            
\usepackage{booktabs}       
\usepackage{amsfonts}       
\usepackage{nicefrac}       
\usepackage{microtype}      
\usepackage{lipsum}
\usepackage{graphicx}

\usepackage[cp1252]{inputenc}
\usepackage{textcomp}
\usepackage{amsmath}
\usepackage{amssymb}
\usepackage{graphicx}
\usepackage{adjustbox}
\usepackage[abs]{overpic}
\usepackage{hyphenat}
\usepackage{inputenc}
\usepackage{amsthm}
\usepackage{booktabs}
\usepackage{bbold}
\usepackage{enumitem}
\setenumerate[1]{label = \arabic*.}
\setenumerate[2]{label*= \arabic*}
\setenumerate[3]{label*=.\arabic*}
\setenumerate[4]{label*=.\arabic*}
\setenumerate[4]{label*=.\arabic*}
\usepackage{xcolor}
\usepackage{multirow}
\usepackage[titletoc]{appendix}
\usepackage{anyfontsize}
\usepackage{framed}
\usepackage{footnote}
\usepackage{mathbbol}
\usepackage{subcaption}
\usepackage{placeins}

\usepackage{adjustbox}

\usepackage{algorithm}
\usepackage{algpseudocode}

\usepackage{url}
\usepackage{comment}
\usepackage{xcolor}

\usepackage{marginnote}

\usepackage{tabulary}
\usepackage{bbm} 
\newtheoremstyle{dotless}{}{}{\itshape}{}{\bfseries}{}{ }{}
\theoremstyle{dotless}

\makeatletter

\renewcommand{\mathbf}{\boldsymbol}

\graphicspath{{figures/}}

\usepackage{color}
\usepackage{soul}
\usepackage{amsthm}

\usepackage{units}

\newcommand{\ba}{\mbox{\boldmath{$a$}}}
\newcommand{\bA}{\mbox{\boldmath{$A$}}}
\newcommand{\bb}{\mbox{\boldmath{$b$}}}

\newcommand{\bd}{\mbox{\boldmath{$d$}}}

\newcommand{\bg}{\mbox{\boldmath{$g$}}}
\newcommand{\bI}{\mbox{\boldmath{$I$}}}
\newcommand{\bK}{\mbox{\boldmath{$K$}}}

\newcommand{\bn}{\mbox{\boldmath{$n$}}}
\newcommand{\bN}{\mbox{\boldmath{$N$}}}
\newcommand{\bp}{\mbox{\boldmath{$p$}}}

\newcommand{\bs}{\mbox{\boldmath{$s$}}}
\newcommand{\bt}{\mbox{\boldmath{$t$}}}

\newcommand{\bu}{\mbox{\boldmath{$u$}}}
\newcommand{\bU}{\mbox{\boldmath{$U$}}}

\newcommand{\bv}{\mbox{\boldmath{$v$}}}

\newcommand{\bx}{\mbox{\boldmath{$x$}}}

\newcommand{\bepsilon}{\mbox{\boldmath{$\varepsilon$}}}

\newcommand{\bsigma}{\mbox{\boldmath{$\sigma$}}}

\newfont{\twelvemsb}{msbm10 at 11.6pt}

\renewcommand{\div}{\mathop{\rm div}\nolimits}

\newcommand{\tr}{\mathop{\rm tr}}

\title{Quasi-optimal mesh generation for the virtual element method: A fully adaptive remeshing procedure}

\author{
	Daniel van Huyssteen \\
	Institute of Applied Mechanics, 
	Friedrich-Alexander-Universität Erlangen-Nürnberg
	Erlangen, 91058, Germany \\
	\texttt{daniel.van.huyssteen@fau.de} \\
	\And
	Felipe Lopez Rivarola \\
	Facultad de Ingeniería,
	Universidad de Buenos Aires,
	Buenos Aires, C1127AAR, Argentina\\
	\And
	Guillermo Etse \\
	Facultad de Ingeniería,
	Universidad de Buenos Aires,
	Buenos Aires, C1127AAR, Argentina\\
	\AND
	Paul Steinmann \\
	Institute of Applied Mechanics, 
	Friedrich-Alexander-Universität Erlangen-Nürnberg
	Erlangen, 91058, Germany \\
}

\begin{document}
	
	\maketitle
	
	
	\begin{abstract}
	The mesh flexibility offered by the virtual element method through the permission of arbitrary element geometries, and the seamless incorporation of `hanging' nodes, has made the method increasingly attractive in the context of adaptive remeshing.
	There exists a healthy literature concerning error estimation and adaptive refinement techniques for virtual elements while the topic of adaptive coarsening (i.e. de-refinement) is in its infancy.
	The creation of a quasi-optimal mesh is based on the principle of quasi-even error distribution over the elements which inherently relies on localized refinement and coarsening techniques. Thus, necessitating a fully adaptive remeshing procedure. 
	In this work a novel fully adaptive remeshing procedure for the virtual element method is presented. Additionally, novel procedures are proposed for the identification of elements qualifying for refinement or coarsening based on user-defined targets. Specifically, a target percentage error, target number of elements, or target number of nodes can be selected. 
	Numerical results demonstrate that the adaptive remeshing procedure can meet any prescribed target while creating a quasi-optimal mesh. 
	The proposed fully adaptive procedure is of particular interest in engineering applications requiring an efficient simulation of a given accuracy, or desiring a simulation with the maximum possible accuracy for a given computational constraint.
	\end{abstract}
	
	
	\newpage
	\noindent
	
	\section{Introduction}
	\label{S:Introduction}
	An optimal mesh represents the pinnacle for the numerical analyst or engineer. An optimal mesh satisfies a specified accuracy target while meeting some `mesh optimality criterion'. Through this combination of objectives an optimal mesh represents the most computationally efficient solution to a specific problem.
	
	In the context of the finite element method (FEM), the generation of optimal meshes has been studied since (at least) the early 1970s. These works typically considered a mesh to be optimal if it met a specified accuracy target with the mesh optimality criterion of minimizing the number of element or nodes in the mesh \cite{wriggers1970comparison,oliveira1971optimization}. A somewhat more recent approach has emerged in which a mesh is considered optimal if it meets the target accuracy with the optimality criterion of equal error distribution over all elements in the mesh \cite{Zhu1988,ONATE1993,Li1995}. In either case, the mesh optimality criterion results in the computation of an element resizing parameter that defines whether a particular element should be refined, coarsened (i.e. de-refined), or is optimally sized. 
	
	In recent times the attention paid to, and number of publications concerning, mesh optimization comprising combined refinement and coarsening processes has seemingly waned. However, there remains great interest in adaptive refinement techniques. The decreased attention in combined processes is likely due to ever increasing advances in computational power that permit very fine mesh simulations of requisite accuracy. These fine mesh simulations may be sufficiently accurate, however, they are inherently inefficient and wasteful of resources. Furthermore, the perceived implementation difficulties, and limited number of works investigating adaptive coarsening techniques may deter one from implementing a traditional mesh optimization procedure comprising combined refinement and coarsening processes. 
	
	Adaptive remeshing techniques for the FEM are already well-established and there exists a wide range of approaches to \textit{a-posteriori} error estimation \cite{babuvvska1978error,Babuska1979,Zienkiewicz1987,Zienkiewicz1992,Zienkiewicz1992a,Zienkiewicz1995} and a variety of tools/packages for the creation of updated meshes \cite{Arndt2022,Kirk2006}. 
	Performing localized refinement or coarsening of finite element meshes is non-trivial as significant manipulation of not only the elements being adapted but also of the surrounding elements is required to preserve the method's conformity. 
	In general, coarsening of finite element meshes is more complex than refinement. As such, most (possibly all) coarsening processes performed using finite elements only reverse previously performed refinement to return a mesh, or parts thereof, to an initially coarser state, see for example \cite{Park2012}. This is particularly problematic in the context of mesh optimization which necessitates the ability to locally refine or coarsen any given mesh which may not contain information about a previously coarser state. 
	
	The introduction of the virtual element method (VEM) gave rise to many new opportunities in the context of adaptive remeshing.
	The VEM is an extension of the FEM that permits arbitrary polygonal and polyhedral element geometries in two- and three-dimensions respectively \cite{VEIGA2012,Veiga2014}. A feature of the VEM of particular interest in the context of adaptive remeshing is the permission of arbitrarily many nodes along an element's edge. That is, nodes that would be considered `hanging' in a finite element context are trivially incorporated into the VEM formulation \cite{VEMContactWriggers2016,Wriggers2019}. 
	The geometric robustness of the VEM has been demonstrated with the method exhibiting optimal convergence behaviour in cases of challenging, including strongly non-convex, element geometries \cite{Sorgente2021,Sorgente2022,Huyssteen2020,Huyssteen2021}. Additionally, in cases of distorted, and possibly stretched, element geometries that could arise during adaptive remeshing (particularly during anisotropic remeshing) the VEM stabilization term can be easily tuned to improve the accuracy of the method \cite{ReddyvanHuyssteen2021,DvH_BDR_MeshQuality}.
	Furthermore, the robustness of the VEM under challenging numerical conditions, such as near-incompressibility and/or near-inextensibility, is increasingly well reported \cite{WriggersIsotropic2017,Wriggers2023,Tang2020,Reddy2019,Huyssteen2020,Huyssteen2021}. 
	The geometric flexibility and numerical robustness mean that the VEM is particularly well-suited to problems involving fully adaptive remeshing and mesh optimization.
	
	Adaptive remeshing techniques for the VEM is an area of rapidly growing interest. There are many works concerning \textit{a-posteriori} error estimation \cite{Cangiani2017,Veiga2015a,Berrone2017,Mora2017,Chi2019,Guo2019,Wei2023,NguyenThanh2018} and several approaches have been presented for localized refinement of the unstructured polygonal element geometries permitted by the method \cite{NguyenThanh2018,vanHuyssteen2022,Huyssteen2022,Berrone2021}. Furthermore, recent works have proposed adaptive coarsening (i.e. de-refinement) techniques for VEM meshes with attention paid to the identification of elements, or groups of elements, to coarsen and the presentation of algorithms for performing the coarsening \cite{Choi2020,Huyssteen2024}. 
	
	The geometric and numerical suitability of the VEM for problems involving mesh optimization, and the high degree of efficacy of even standalone adaptive refinement and coarsening procedures, strongly motivate the development of a mesh optimization procedure for the VEM. 
	The generation of a truly optimal mesh would require complex algorithms for the precise resizing of individual elements based on the resizing parameter. Furthermore, since error-estimation is approximate, and depends on the mesh, the resizing would need to be performed iteratively. It is unlikely that this iterative procedure would stabilize and yield a fully optimal mesh within a reasonable (or even finite) number of iterations.
	Thus, in this work novel procedures are proposed for the creation of `quasi-optimal' meshes of virtual elements. These procedures consider various targets. Specifically, a target accuracy, target number of elements, or target number of nodes can be set. Additionally, the mesh optimality criterion of `quasi-equal' error distribution is chosen due to its simplicity of implementation and practical suitability. The procedures comprise the novel combination of adaptive refinement and coarsening algorithms previously proposed and numerically studied by the authors, as well as novel algorithms for the identification of elements to refine or coarsen. Finally, in this work the accuracy and efficacy of the novel fully adaptive procedure with the novel element selection algorithms is measured through an approximation of the well-known energy error.

	The structure of the rest of this work is as follows. The governing equations of linear elasticity are set out in Section~\ref{sec:GovEq}. This is followed in Section~\ref{sec:VEM} by a description of the first-order virtual element method. 
	The procedures used to generate, refine, and coarsen meshes are presented in Section~\ref{sec:MeshGenerationRefinementAndCoarsening}. 
	This is followed, in Section~\ref{sec:ErrorEstimationAndPrediction}, by a description of the procedures used to compute local and global error estimations. Thereafter, the procedures for the selection of elements to resize via refinement and coarsening are presented in Section~\ref{sec:ElementSelectionForRefinementAndCoarsening} for various remeshing strategies. 
	Section~\ref{sec:Results} comprises a set of numerical results through which the performance of the various remeshing strategies is evaluated.
	Finally, the work concludes in Section~\ref{sec:Conclusion} with a discussion of the results.
	
	\section{Governing equations of linear elasticity} 
	\label{sec:GovEq}
\begin{figure}[ht!]
		\centering
		\includegraphics[width=0.25\textwidth]{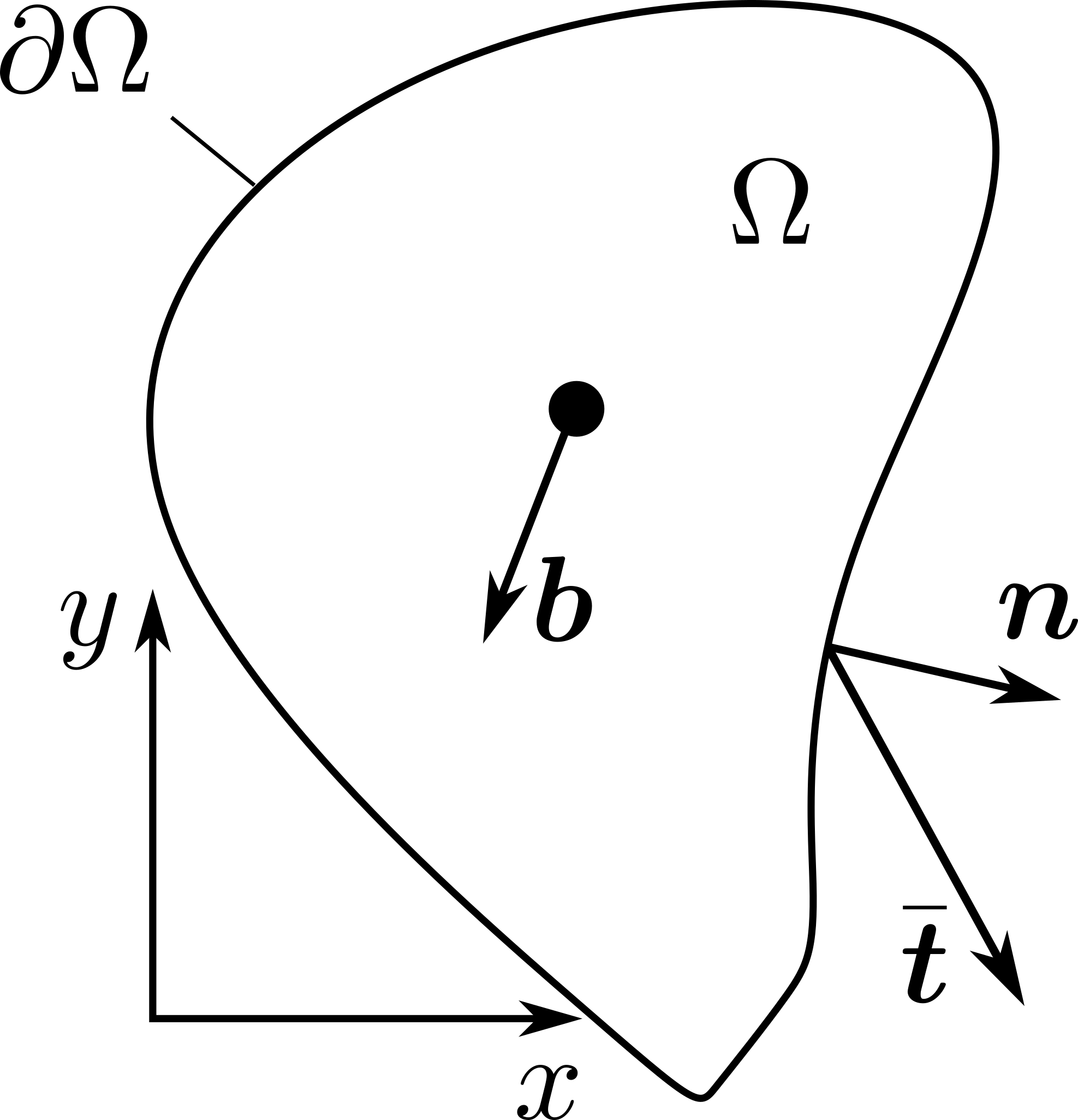}
		\caption{Arbitrary elastic body $\Omega$ with boundary ${\partial \Omega}$ subject to body ${\bb}$ force and traction ${\bar{\bt}}$.
			\label{fig:ElasticBody}}
	\end{figure}

	Consider an arbitrary elastic body occupying a plane, bounded, domain ${\Omega \subset \mathbb{R}^{2} }$ subject to a traction ${\bar{\bt}}$ and body force ${\bb}$ (see Figure~\ref{fig:ElasticBody}).
	The boundary ${\partial \Omega}$ has an outward facing normal denoted by $\bn$ and comprises a non-trivial Dirichlet part $\Gamma_{D}$ and a Neumann part $\Gamma_{N}$ such that ${\Gamma_{D} \cap \Gamma_{N} = \emptyset}$ and ${\overline{\Gamma_{D} \cup \Gamma_{N}}=\partial \Omega}$.
	
	\FloatBarrier
	\FloatBarrier
	
	In this work small displacements are assumed and the strain-displacement relation is given by 
	\begin{equation}
		\bepsilon \left( \bu \right) = \frac{1}{2} \left[\nabla\,\bu + \left[ \nabla\,\bu \right]^{T} \right] \, . \label{eqn:DisplacementStrain}
	\end{equation}
	Here the displacement is denoted by $\bu$, ${\bepsilon}$ is the symmetric infinitesimal strain tensor and ${\nabla \left( \bullet \right) = \frac{\partial \left( \bullet \right)_{i} }{\partial \, x_{j}} \, \boldsymbol{e}_{i} \otimes \boldsymbol{e}_{j} }$ is the gradient of a vector quantity. Additionally, linear elasticity is assumed and the stress-strain relation is given by
	\begin{equation}
		\bsigma = \mathbb{C} : \bepsilon \,. \label{eqn:StressStrain1}
	\end{equation}
	Here, ${\bsigma}$ is the Cauchy stress tensor and ${\mathbb{C}}$ is a fourth-order constitutive tensor. For a linear elastic and isotropic material (\ref{eqn:StressStrain1}) is given by
	\begin{eqnarray}
		\bsigma = \lambda \tr \left( \bepsilon \right)  \bI + 2\mu \, \bepsilon \, , \label{eqn:StressStrain2}
	\end{eqnarray}
	where ${\tr \left( \bullet \right)}$ denotes the trace, $\bI$ is the second-order identity tensor, and $\lambda$ and $\mu$ are the well-known Lam\'{e} parameters.
	
	For equilibrium it is required that 
	\begin{equation}
		\div \, \bsigma + \bb = \boldsymbol{0} \, , \label{eqn:Equilibrium}
	\end{equation}
	where ${\div \left( \bullet \right) = \frac{\partial \left( \bullet \right)_{ij} }{\partial \, x_{j}} \boldsymbol{e}_{i} }$ is the divergence of a tensor quantity.
	The Dirichlet and Neumann boundary conditions are given by 
	\begin{align}
		&\bu = \bg \quad \text{on } \Gamma_{D} \, , \text{ and} \label{eqn:DirichletBC} \\
		&\bsigma \cdot \bn = \bar{\bt} \quad \text{on } \Gamma_{N} \, , \label{eqn:NeumannBC}
	\end{align}
	respectively, with $\bg$ and $\bar{\bt}$ denoting prescribed displacements and tractions respectively.
	Equations (\ref{eqn:StressStrain2})-(\ref{eqn:NeumannBC}), together with the displacement-strain relationship (\ref{eqn:DisplacementStrain}), constitute the boundary-value problem for a linear elastic isotropic body.
	
	\subsection{Weak form}
	\label{subsec:WeakForm}
	The space of square-integrable functions on $\Omega$ is hereinafter denoted by ${\mathcal{L}^{2}\left(\Omega\right)}$. The Sobolev space of functions that, together with their first derivatives, are square-integrable on $\Omega$ is hereinafter denoted by ${\mathcal{H}^{1}\left(\Omega\right)}$. Additionally, the function space $\mathcal{V}$ is introduced and defined such that
	\begin{align}
		\mathcal{V} = \left[ \mathcal{H}^{1}_{D} \left(\Omega\right) \right]^{d} 
		=
		\left\{ \bv \, | \, v_{i} \in \mathcal{H}^{1}\left(\Omega\right), \, \bv = \boldsymbol{0} \,\, \text{on} \,\, \Gamma_{D}  \right\} \, 
	\end{align}
	where ${d=2}$ is the dimension.
	Furthermore, the function ${\bu_{g}\in \left[ \mathcal{H}^{1} \left(\Omega\right) \right]^{d} }$ is introduced satisfying (\ref{eqn:DirichletBC}) such that ${\bu_{g}|_{\Gamma_{D}}=\bg}$.
	
	The bilinear form ${a\left(\cdot,\cdot\right)}$, where ${a : \left[ \mathcal{H}^{1} \left(\Omega\right) \right]^{d}  \times \left[ \mathcal{H}^{1} \left(\Omega\right) \right]^{d} \rightarrow \mathbb{R}}$, and the linear functional ${\ell \left(\cdot\right)}$, where ${\ell : \left[ \mathcal{H}^{1} \left(\Omega\right) \right]^{d} \rightarrow \mathbb{R}}$, are defined respectively by
	\begin{equation}
		a\left(\bu, \, \bv \right) = \int_{\Omega} \bsigma \left( \bu \right) : \bepsilon \left( \bv \right) \, dx \, , \label{eqn:BilinearForm}
	\end{equation}
	and 
	\begin{equation}
		\ell \left( \bv \right) = \int_{\Omega} \bb \cdot \bv \, dx + \int_{\Gamma_{N}} \bar{\bt} \cdot \bv \, ds - a\left(\bu_{g},\, \bv \right) \, . \label{eqn:LinearFuntional} 
	\end{equation}
	
	The weak form of the problem is then: given ${\bb \in \left[ \mathcal{L}^{2}\left(\Omega\right) \right]^{d} }$ and ${ \bar{\bt} \in \left[ \mathcal{L}^{2}\left(\Gamma_{N}\right) \right]^{d} }$, find ${\bU \in \left[ \mathcal{H}^{1}\left(\Omega\right) \right]^{d} }$ such that 
	\begin{equation}
		\bU = \bu + \bu_{g} \, , \quad \bu \in \mathcal{V} \, ,
	\end{equation}
	and
	\begin{equation}
		a\left(\bu , \, \bv \right) = \ell \left( \bv \right) \, , \quad \forall \bv \in \mathcal{V} \, . \label{eqn:BilinearFormFinal}
	\end{equation}
	
	\section{The virtual element method}
	\label{sec:VEM}
	
\begin{figure}[ht!]
		\centering
		\includegraphics[width=0.33\textwidth]{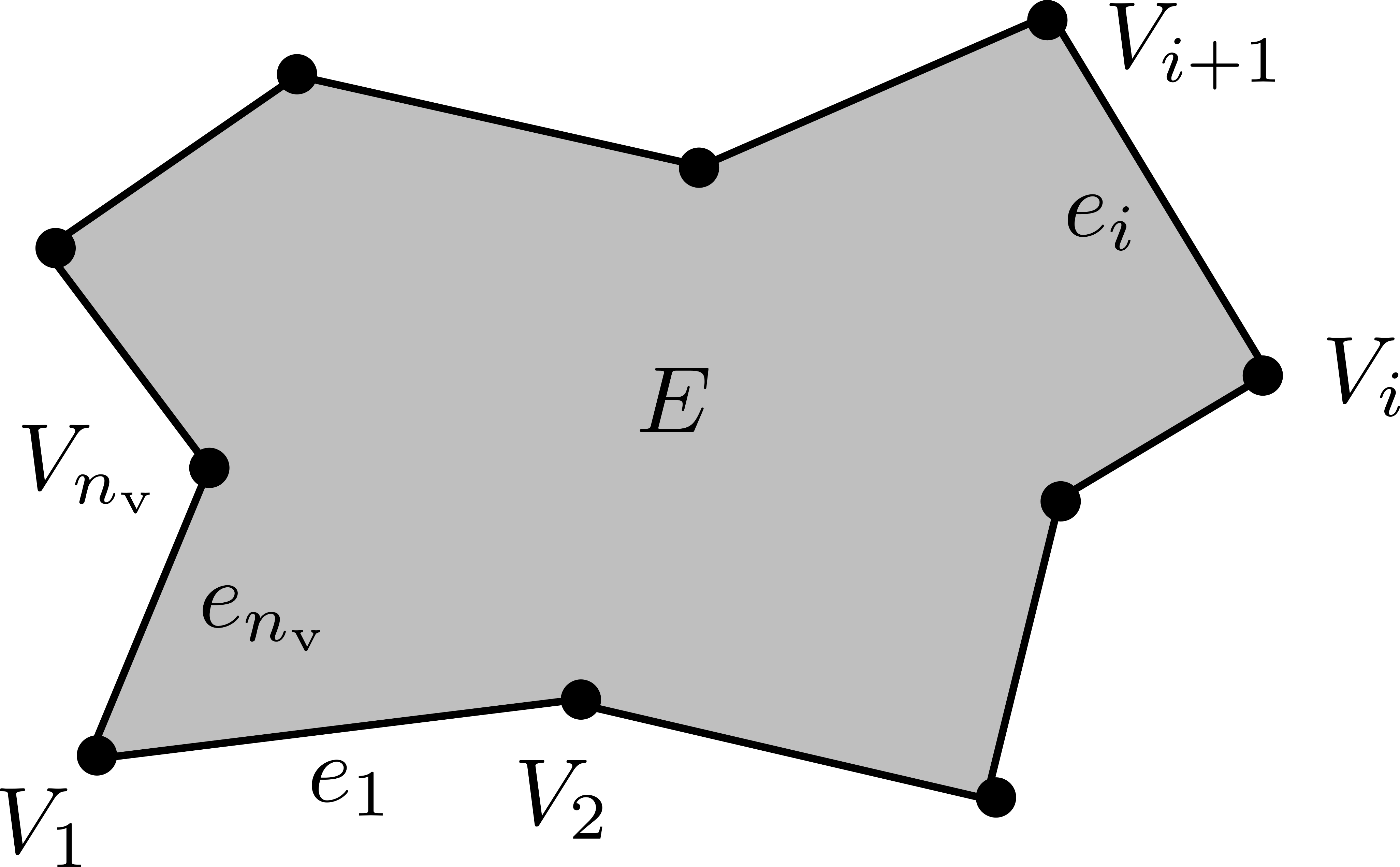}
		\caption{Sample virtual element $E$ with $n_{\rm v}$ vertices where edge $e_{i}$ connects vertices $V_{i}$ and $V_{i+1}$.
			\label{fig:SampleElement}}
	\end{figure}

	The domain $\Omega$ is partitioned into a mesh of non-overlapping arbitrary polygonal elements\footnote{If $\Omega$ is not polygonal the mesh will be an approximation of $\Omega$.} $E$ with $\overline{\cup E}=\overline{\Omega}$. Here $E$ denotes the element domain and $\partial E$ its boundary, with ${\overline{(\, \bullet \,)}}$ denoting the closure of a set.
	An example of a typical first-order element is depicted in Figure~\ref{fig:SampleElement} with edge $e_{i}$ connecting vertices $V_{i}$ and $V_{i+1}$. Here ${i=1,\dots,n_{\rm v}}$ with $n_{\rm v}$ denoting the total number of element vertices.
	
	\FloatBarrier
	\FloatBarrier
	A conforming approximation of order $k$ is constructed in a space ${\mathcal{V}^{h} \subset \mathcal{V}}$ where $\mathcal{V}^{h}$ is built-up element-wise and comprises vector valued functions $\bv_{h}$.
	The functions $\bv_{h}$ are those that are $\mathcal{C}^{0}$ continuous on the domain $\Omega$, are polynomials of degree ${\leq \, k}$ on element edges, and whose strain gradient divergence is a polynomial of degree ${\leq \, k-2}$ on an element (see \cite{Artioli2017}). 
	For the most general case of an approximation of arbitrary order $k$ the space $\mathcal{V}^{h}|_{E}$ is defined as 
	\begin{equation}
		\mathcal{V}^{h}|_{E} = \left\{ \bv_{h} \in \mathcal{V} \, | \, \bv_{h} \in \left[ \mathcal{C}^{0}(E) \right]^{2} \, , \, \nabla^{2} \, \bv_{h}  \in \mathcal{P}_{k-2} \text{ on } E \, , \, \bv_{h}|_{e} \in \mathcal{P}_{k}(e)  \right\} \,. \label{eqn:ArbitraryVEMSpace}
	\end{equation}
	Here ${\mathcal{P}_{k}(X)}$ is the space of polynomials of degree ${\leq \, k}$ on the set ${X \, \subset \, \mathbb{R}^{d} }$ with ${d=1,\,2}$ and ${\nabla^{2}=\nabla\cdot\nabla}$ is the Laplacian operator. 
	In this work a first-order, i.e. ${k=1}$, approximation is considered, thus (\ref{eqn:ArbitraryVEMSpace}) simplifies to 
	\begin{equation}
		\mathcal{V}^{h}|_{E} = \left\{ \bv_{h} \in \mathcal{V} \, | \, \bv_{h} \in \left[ \mathcal{C}^{0}(E) \right]^{2} \, , \, \nabla^{2} \, \bv_{h}  = \boldsymbol{0} \text{ on } E \, , \, \bv_{h}|_{e} \in \mathcal{P}_{1}(e)  \right\} \,. 
	\end{equation}
	
	All computations will be performed on element edges and it is convenient to write, for element $E$,
	\begin{equation}
		\bv_{h}|_{\partial E} = \bN \cdot \bd^{E} \,. \label{eqn:DisplacementTrace}
	\end{equation}
	Here, $\bN$ is a matrix of standard linear Lagrangian basis functions and $\bd^{E}$ is a ${{2}n_{\rm v} \times 1}$ vector of the degrees of freedom associated with $E$. The virtual basis functions are not known, nor required to be known on $E$; their traces, however, are known and are simple Lagrangian functions.
	
	The virtual element projection for a first-order formulation ${\Pi \, : \, \mathcal{V}^{h}|_{E} \rightarrow \mathcal{P}_{0}(E)  }$ is required to satisfy
	\begin{equation}
		\int_{E} \Pi \, \bv_{h} \cdot \bepsilon\left( \bp \right) \, dx = \int_{E} \bepsilon\left(\bv_{h}\right) \cdot \bepsilon\left( \bp \right) \, dx \quad \forall \bp \in \mathcal{P}_{1} \,, \label{eqn:Projection}
	\end{equation}
	where ${\Pi \, \bv_{h}}$ represents the ${\mathcal{L}^{2}}$ projection of the symmetric gradient of ${\bv_{h}}$ onto constants \cite{Artioli2017}. Since the projection is constant at element-level, after applying integration by parts to (\ref{eqn:Projection}), and considering (\ref{eqn:DisplacementTrace}), the components of the projection can be computed as
	\begin{align}
		\left(\Pi\,\bv_{h}\right)_{ij} &= \frac{1}{2}\frac{1}{|E|}  \sum_{e\in\partial E} \int_{e} \left[ N_{iA} \, d_{A}^{E} \, n_{j} + N_{jA} \, d_{A}^{E} \, n_{i}\right] ds \,, \label{eqn:ComputeProjection}
	\end{align}
	where summation is implied over repeated indices.
	
	The virtual element approximation of the bilinear form (\ref{eqn:BilinearForm}) is constructed by writing 
	\begin{align}
		a^{E}\left(\bu,\,\bv\right) :&= a\left(\bu,\,\bv\right)|_{E} 
		= \int_{E} \bepsilon\left(\bv_{h}\right) : \left[ \mathbb{C} : \bepsilon\left(\bu_{h}\right) \right] dx \, , \label{eqn:ElementBilinearForm}
	\end{align}
	where ${a^{E}\left(\cdot,\cdot\right)}$ is the contribution of element $E$ to the bilinear form ${a\left(\cdot,\cdot\right)}$. Consideration of (\ref{eqn:ComputeProjection}) allows (\ref{eqn:ElementBilinearForm}) to be written as (see \cite{Reddy2019})
	\begin{align}
		a^{E}\left(\bu_{h},\,\bv_{h}\right) 
		&= \underbrace{\int_{E} \Pi\,\bv_{h} : \left[ \mathbb{C} : \Pi\,\bu_{h} \right] dx }_{\text{Consistency term}} 
		+ \underbrace{\int_{E} \left[ \bepsilon\left( \bv_{h} \right) : \left[ \mathbb{C} : \bepsilon \left( \bu_{h} \right) \right] - \Pi\,\bv_{h} : \left[ \mathbb{C} : \Pi \, \bu_{h} \right] \right] dx }_{\text{Stabilization term}} \,, \label{eqn:ExpandedBilinearForm}
	\end{align}
	where the remainder term is discretized by means of a stabilization.
	
	\subsection{The consistency term} 
	\label{subsec:ConsistencyTerm}
	
	The projection (\ref{eqn:ComputeProjection}), and thus the consistency term, can be computed exactly yielding
	\begin{equation}
		a_{\rm c}^{E}\left(\bu_{h},\,\bv_{h}\right) \, = \, \int_{E} \Pi\,\bv_{h} : \left[ \mathbb{C} : \Pi\,\bu_{h} \right] dx \, = \, \widehat{\bd}^{E} \cdot \left[ \bK_{\rm c}^{E} \cdot \bd^{E} \right] \,.
	\end{equation}
	
	Here $\bK_{\rm c}^{E}$ is the consistency part of the stiffness matrix of element $E$ with ${\widehat{\bd}^{E}}$ and ${\bd^{E}}$ the degrees of freedom of $\bv_{h}$ and $\bu_{h}$ respectively that are associated with element $E$.
	
	\subsection{The stabilization term} 
	\label{subsec:Stab}
	The remainder term cannot be computed exactly and is approximated by means of a discrete stabilization term \cite{Gain2014,Veiga2015}.
	The approximation employed in this work is motivated by seeking to approximate the difference between the element degrees of freedom $\bd^{E}$ and the nodal values of a linear function that is closest to $\bd^{E}$ in some way (see \cite{Artioli2017,Reddy2019}). 
	The nodal values of the linear function are given by
	\begin{equation}
		\widetilde{\bd} = \boldsymbol{\mathcal{D}} \cdot \bs \,. \label{eqn:LinearApprox}
	\end{equation}
	Here $\bs$ is a vector of the degrees of freedom of the linear function and $\boldsymbol{\mathcal{D}}$ is a matrix relating $\widetilde{\bd}$ to $\bs$ with respect to a scaled monomial basis. For the full expression of $\boldsymbol{\mathcal{D}}$ see \cite{Artioli2017,Reddy2019}.
	After some manipulation (see, again, \cite{Reddy2019}) the stabilization term of the bilinear form can be approximated as
	\begin{equation}
		a_{\text{stab}}^{E}\left(\bu_{h},\,\bv_{h}\right) \, = \,
		\int_{E} \left[ \bepsilon\left( \bv_{h} \right) : \left[ \mathbb{C} : \bepsilon \left( \bu_{h} \right) \right] - \Pi\,\bv_{h} : \left[ \mathbb{C} : \Pi \, \bu_{h} \right]  \right] dx \, \approx \, \widehat{\bd}^{E} \cdot \left[ \bK_{\rm s}^{E} \cdot \bd^{E} \right]  \, ,
	\end{equation}
	where $\bK_{\rm s}^{E}$ is the stabilization part of the stiffness matrix of element $E$ and is defined as
	\begin{equation}
		\bK_{\rm s}^{E} = \mu \left[ \bI - \boldsymbol{\mathcal{D}} \cdot \left[ \boldsymbol{\mathcal{D}}^{T} \cdot \boldsymbol{\mathcal{D}}\right]^{-1} \cdot  \boldsymbol{\mathcal{D}}^{T} \right] \, .
	\end{equation}
	The total element stiffness matrix ${\bK^{E}}$ is then computed as the sum of the consistency and stabilization matrices.
	
	\section{Mesh generation, refinement and coarsening} 
	\label{sec:MeshGenerationRefinementAndCoarsening}
	In this section the procedures used to generate meshes, refine elements and coarsen patches of elements are described.
	
	\subsection{Mesh generation} 
	\label{subsec:MeshGeneration}
\begin{figure}[ht!]
		\centering
		\includegraphics[width=0.95\textwidth]{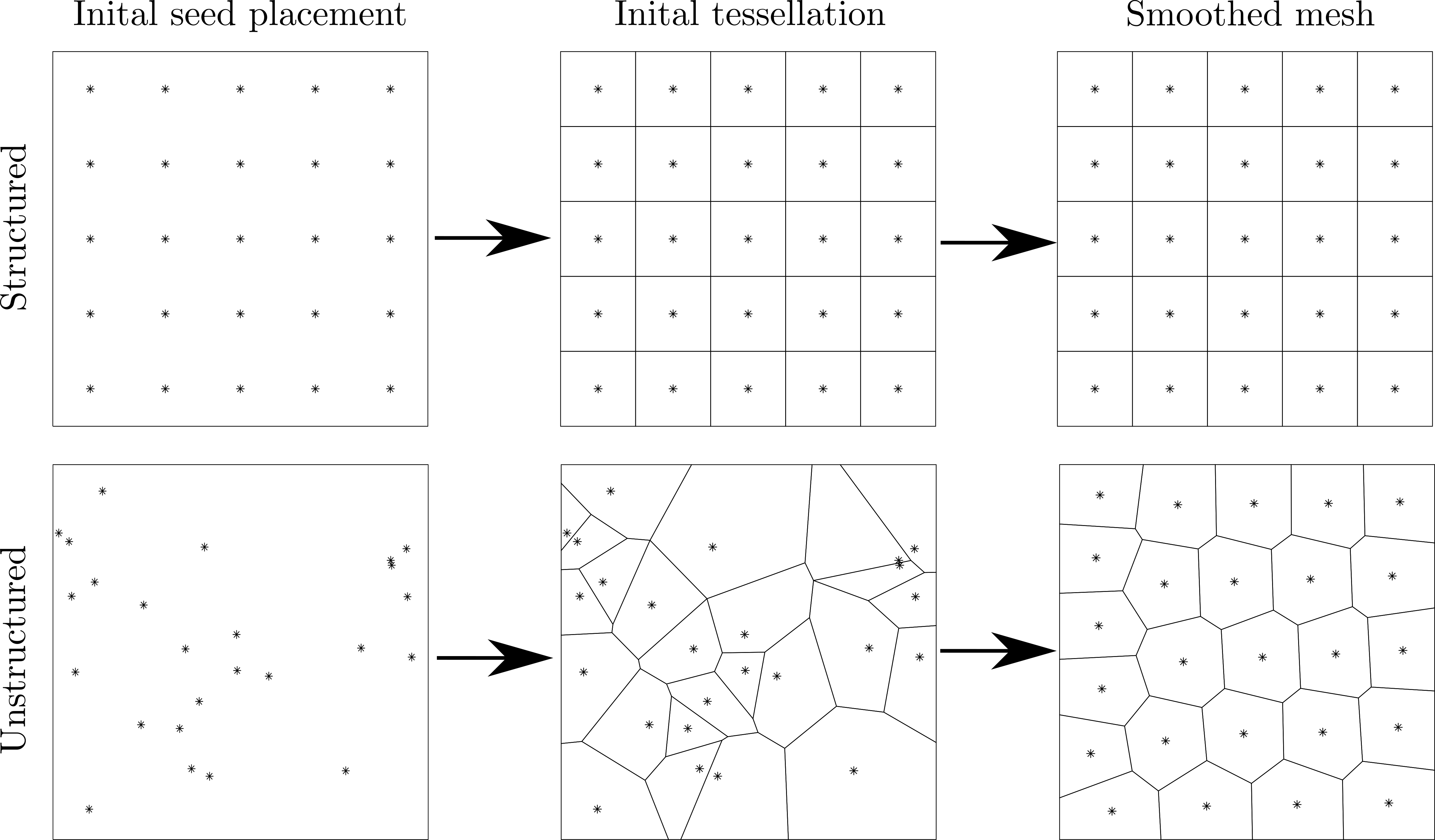}
		\caption{Mesh generation procedure for structured and unstructured/Voronoi meshes.
			\label{fig:MeshGeneration}}
	\end{figure}

	The mesh generation procedure used in this work is identical to that described in \cite{Huyssteen2022,Huyssteen2024} and is summarized briefly here for the sake of self-containment.
	All meshes are created by Voronoi tessellation of a set of seed points. Seed points will be generated in both structured and unstructured sets to create structured and unstructured meshes respectively. 
	In the case of structured meshes seeds points are placed to form a structured grid, while in the case of unstructured/Voronoi meshes seeds are placed arbitrarily within the problem domain. Hereinafter the terms `unstructured' and `Voronoi' meshes will be used interchangeably to refer to meshes created from arbitrarily placed seed points.
	An initial Voronoi tessellation of the seed points is created using PolyMesher \cite{PolyMesher}. Then, a smoothing algorithm in PolyMesher is used to iteratively modify the locations of the seed points to create a mesh in which all elements have approximately equal areas. Clearly, in the case of structured meshes the smoothing step is trivial. 
	The mesh generation procedure is illustrated in Figure~\ref{fig:MeshGeneration} where the top and bottom rows depict the generation of structured and unstructured/Voronoi meshes respectively.  
	
	\FloatBarrier
	\FloatBarrier
	
	\subsection{Mesh refinement} 
	\label{subsec:MeshRefinement}
\begin{figure}[ht!]
		\centering
		\includegraphics[width=0.95\textwidth]{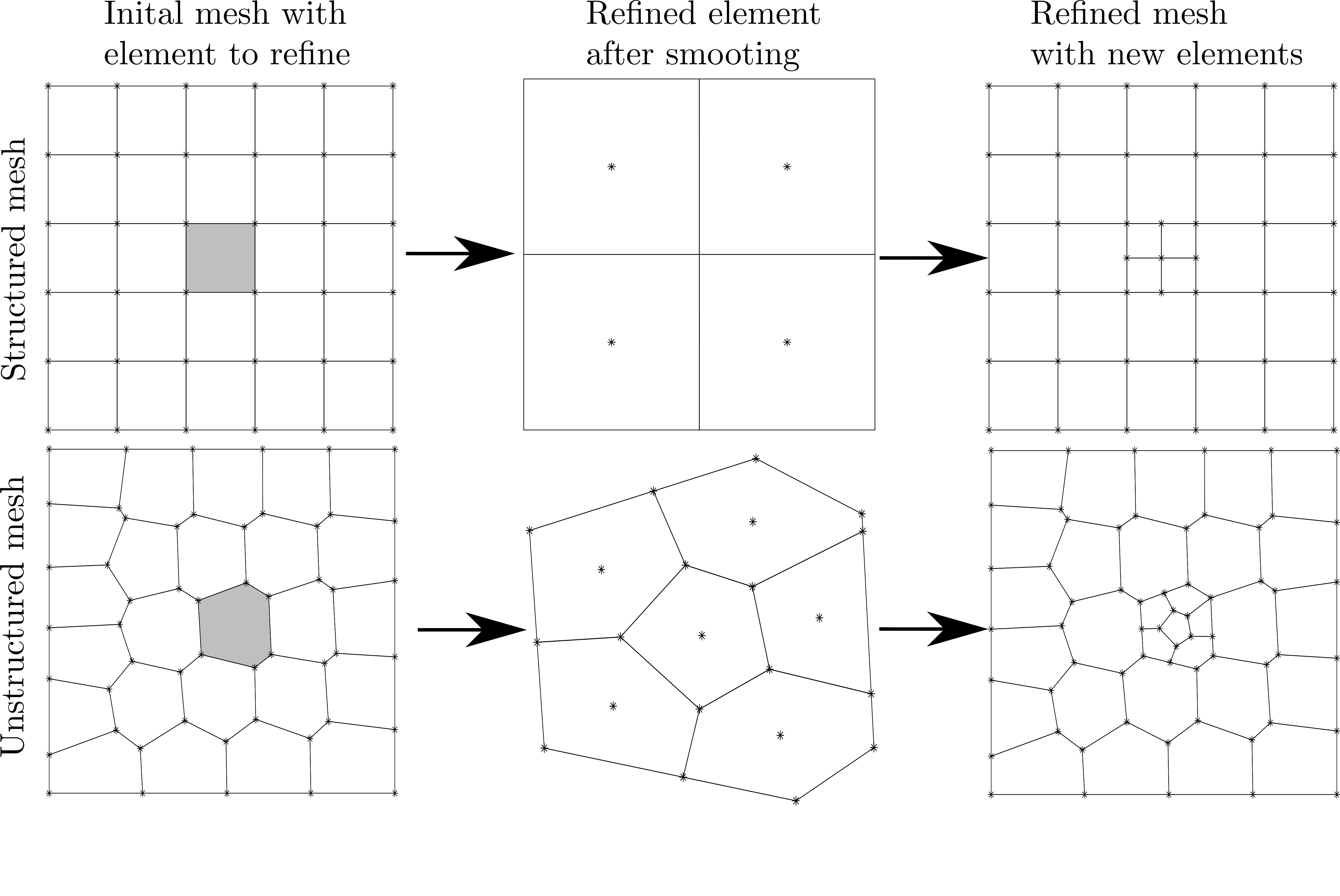}
		\caption{Refinement procedure for structured and unstructured/Voronoi elements.
			\label{fig:MeshRefinement}}
	\end{figure}

	The mesh refinement procedure used in this work is identical to that described in \cite{Huyssteen2022} and is summarized briefly here for the sake of self-containment.
	Once an element has been marked for refinement the process is performed using a modified version of PolyMesher \cite{PolyMesher}. An overview of the element refinement procedure is illustrated in Figure~\ref{fig:MeshRefinement} for structured and unstructured/Voronoi meshes.
	The element marked for refinement is indicated in grey within the initial mesh. 
	Refinement is performed by subdividing a marked element into smaller elements via Voronoi tessellation of a set of seed points, similar to the mesh generation process.
	For simplicity the number of seed points is chosen to be equal to the number of nodes of the element. In the case of structured meshes the seeds are placed in a structured grid, while in the case of unstructured/Voronoi elements the seeds are placed randomly within the element. An initial Voronoi tessellation of the seeds is created and then smoothed using PolyMesher.
	After smoothing, a procedure is used to `optimize' the positions of the newly created nodes that lie on the edges of the original element (see \cite{Huyssteen2022}). The smoothed and optimized elements are depicted in the right-hand column of Figure~\ref{fig:MeshRefinement}. 
	
	\FloatBarrier
	\FloatBarrier
	
	\subsection{Mesh coarsening} 
	\label{subsec:MeshCoarsening}
\begin{figure}[ht!]
		\centering
		\includegraphics[width=0.95\textwidth]{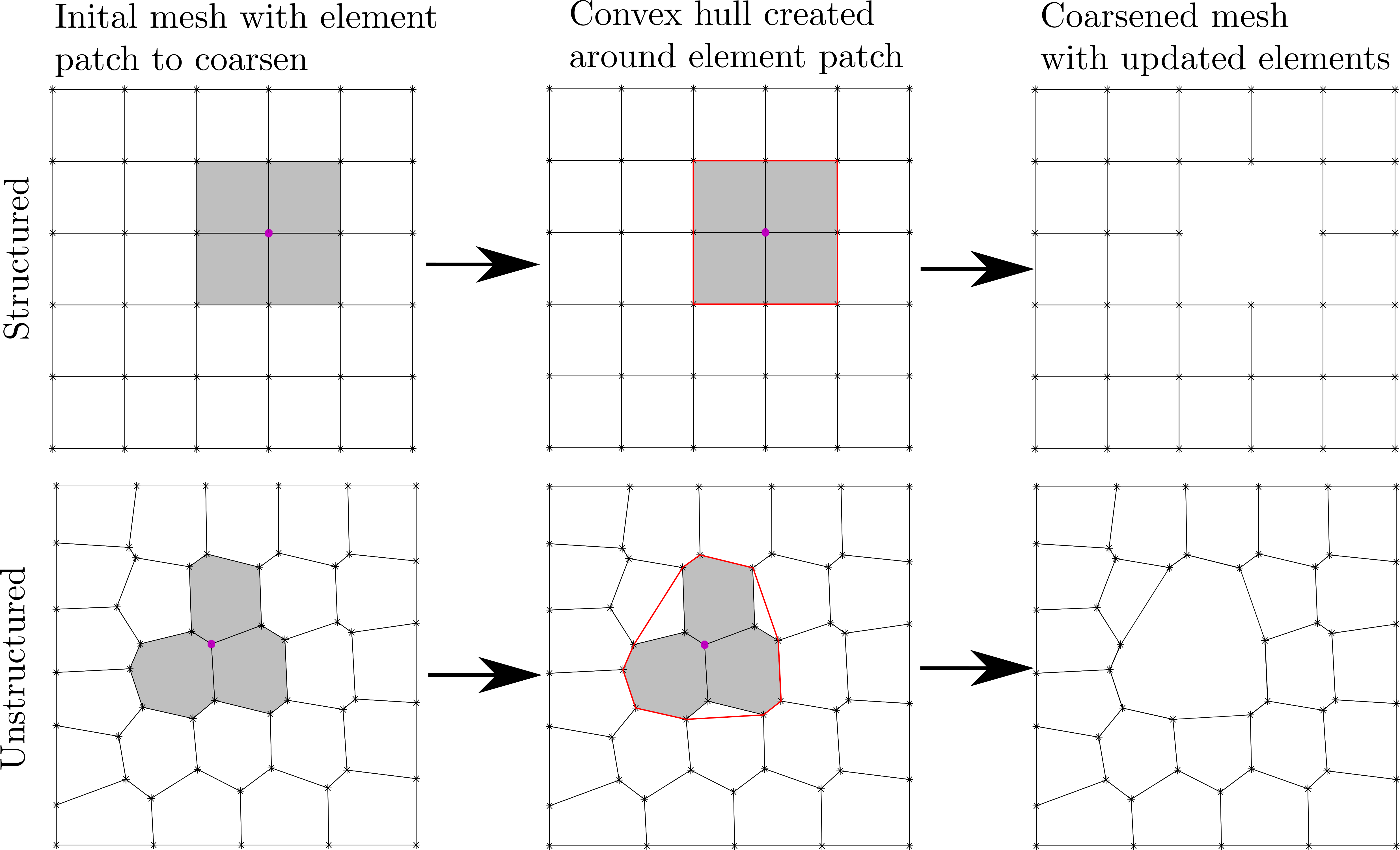}
		\caption{Coarsening procedure for structured and unstructured/Voronoi elements.
			\label{fig:MeshCoarsening}}
	\end{figure}

	The mesh coarsening procedure used in this work is identical to that described in \cite{Huyssteen2024} and is summarized briefly here for the sake of self-containment and illustrated in Figure~\ref{fig:MeshCoarsening}.
	The patch/group of elements to be coarsened/combined is indicated in grey. The geometry of the coarsened element is created by constructing a convex hull around the patch of elements as indicated in red. The geometries of the elements in the patch are modified to coincide with the convex hull using the edge straightening procedure proposed in \cite{Huyssteen2024}. Once the geometries of the marked elements, and the surrounding elements, have been modified the marked elements are deleted and one new element is created using the geometry of the convex hull.
	
	\FloatBarrier
	\FloatBarrier
	
	\section{Error estimation and prediction}
	\label{sec:ErrorEstimationAndPrediction}
	In this section the procedures for calculating the approximate local (element-level) and global errors are presented along with the procedure for the prediction of error after coarsening.
	In this work error is measured through the well-known energy error norm \cite{Zienkiewicz1987}. 
	
	\subsection{Global error estimation}
	\label{subsec:GlobalErrorEstimation}
	The global error in the $\mathcal{H}^{1}$ semi-norm, i.e. the energy error norm, is defined as 
	\begin{equation}
		\|e\|_{\mathcal{H}^{1}} =  
		\left[\frac{1}{2} \int_{\Omega} 
		\left[ \bsigma^{\text{ex}} - \bsigma^{h} \right]^{T}
		\mathbb{D}^{-1}  
		\left[ \bsigma^{\text{ex}} - \bsigma^{h} \right] d\Omega
		\right]^{0.5} \, , \label{eqn:EnergyExact}
	\end{equation}
	where $\bsigma^{\text{ex}}$ is the exact/analytical stress solution and $\mathbb{D}$ is the constitutive matrix. 
	In practical applications the exact stress is typically unknown and is replaced with an approximation $\bsigma^{\ast}$(see Section~\ref{subsec:SuperConvergentPatchRecovery}). 
	A relative energy error $\|e\|_{\text{rel}}$ is introduced and defined as the ratio of the energy error $\|e\|_{\mathcal{H}^{1}}$ to the elastic energy of the deformed body $\|U\|$ such that
	\begin{equation}
		\|e\|_{\text{rel}} = \frac{\|e\|_{\mathcal{H}^{1}}}{\|U\|} \, ,
	\end{equation}
	where the elastic energy is computed as
	\begin{equation}
		\|U\| = \left[\frac{1}{2} \int_{\Omega} 
		\left[ \bsigma^{\text{ex}} \right]^{T}
		\mathbb{D}^{-1}  
		\left[ \bsigma^{\text{ex}} \right] d\Omega
		\right]^{0.5} \,. 
	\end{equation}
	The global error $\|e\|_{\mathcal{H}^{1}}$ and global energy $\|U\|$ can be computed as a sum of element-level contributions given by 
	\begin{equation}
		\|e\|_{\mathcal{H}^{1}} = \left[\frac{1}{2} \sum_{i=1}^{n_{\text{el}}} e_i \right]^{0.5}
		\quad \text{and \quad} 
		\|U\| = \left[\frac{1}{2} \sum_{i=1}^{n_{\text{el}}} U_i \right]^{0.5} \label{eqn:GlobalFromLocal}
	\end{equation}
	respectively. Here $n_{\text{el}}$ is the total number of elements in the domain with $e_i$ and $U_i$ respectively the element-level error and energy contributions.
	
	\subsection{Local (element-level) error estimation}
	\label{subsec:LocalErrorEstimation}
	The element-level error contribution to the global energy error is computed, approximately, for the $i$-th element as 
	\begin{equation}
		e_i \approx  
		\frac{|E_{i}|}{n_{\text{v}}^{i}} \sum_{j=1}^{n_{\text{v}}^{i}} 
		\left[
		\left[ \bsigma^{\ast}\left(\bx_{j}\right) - \bsigma^{h}\left(\bx_{j}\right) \right]^{T} \mathbb{D}^{-1}    
		\left[ \bsigma^{\ast}\left(\bx_{j}\right) - \bsigma^{h}\left(\bx_{j}\right) \right]
		\right] \,,
	\end{equation}
	where $\bsigma^{\ast}$ is an approximation of the exact stress $\bsigma^{\text{ex}}$. Additionally, the area of the $i$-th element is denoted by $|E_{i}|$ and $n_{\text{v}}^{i}$ denotes the number of nodes/vertices.
	Furthermore, analogously to (\ref{eqn:GlobalFromLocal}), the energy error on a single element can be approximated as 
	\begin{equation}
		\|e_i\|_{\mathcal{H}^{1}} \approx  \left[ \frac{1}{2} \, e_i \right]^{0.5} \, .
	\end{equation}
	Similarly, the element-level energy contribution is approximated as 
	\begin{equation}
		U_i \approx  
		\frac{|E_{i}|}{n_{\text{v}}^{i}} \sum_{j=1}^{n_{\text{v}}^{i}} 
		\left[
		\left[ \bsigma^{\ast}\left(\bx_{j}\right) \right]^{T} \mathbb{D}^{-1}    
		\left[ \bsigma^{\ast}\left(\bx_{j}\right) \right]
		\right] \,.
	\end{equation}
	
	\subsection{Local (common node) error prediction}
	\label{subsec:LocalErrorPrediction}
	An energy error prediction is introduced that aims to predict how much coarsening a particular patch of elements would increase the local and global approximations of the energy error (see~\cite{Huyssteen2024}). 
	Here, a patch refers to all of the elements connected to a particular node. Thus, every node has an associated patch of elements and the energy error prediction is computed for each node. 
	The energy error prediction $\|e_{\text{p}_i}\|_{\mathcal{H}^1}$ is approximated over patch $i$ as 
	\begin{equation}
		\|e_{\text{p}_i}\|_{\mathcal{H}^1} \approx
		\left[ \frac{1}{2} 
		\frac{|E_{\text{p}_i}|}{n_{\text{v}}^{\text{p}_i}} \sum_{j=1}^{n_{\text{v}}^{\text{p}_i}} 
		\left[
		\left[ \bsigma^{\ast}\left(\bx_{j}\right) - \bar{\bsigma}^{h}_{\text{p}_i}\left(\bx_{j}\right) \right]^{T} \mathbb{D}^{-1}    
		\left[ \bsigma^{\ast}\left(\bx_{j}\right) - \bar{\bsigma}^{h}_{\text{p}_i}\left(\bx_{j}\right) \right]
		\right] \right]^{0.5} \, .
		\label{eqn:EnergyErrorIndicator}
	\end{equation}
	Here $|E_{\text{p}_i}|$ denotes the area of patch $i$ and $n_{\text{v}}^{\text{p}_i}$ is the number of unique vertices/nodes associated with the patch. Additionally, $\bar{\bsigma}^{h}_{\text{p}_i}$ denotes the `predicted´ stress over the coarsened patch computed as the weighted average of the element stresses on the patch.
	
	\subsection{Super-convergent patch recovery}
	\label{subsec:SuperConvergentPatchRecovery}
	In practical applications the exact stress $\bsigma^{\text{ex}}$ is typically unknown and is replaced with an approximation $\bsigma^{\ast}$. A simple and effective approach in the VEM context is to compute $\bsigma^{\ast}$ at only the nodal positions using a patch-based recovery technique based on super-convergent sampling points (see \cite{NguyenThanh2018,Huyssteen2024}). 
	
\begin{figure}[ht!]
		\centering
		\includegraphics[width=0.8\textwidth]{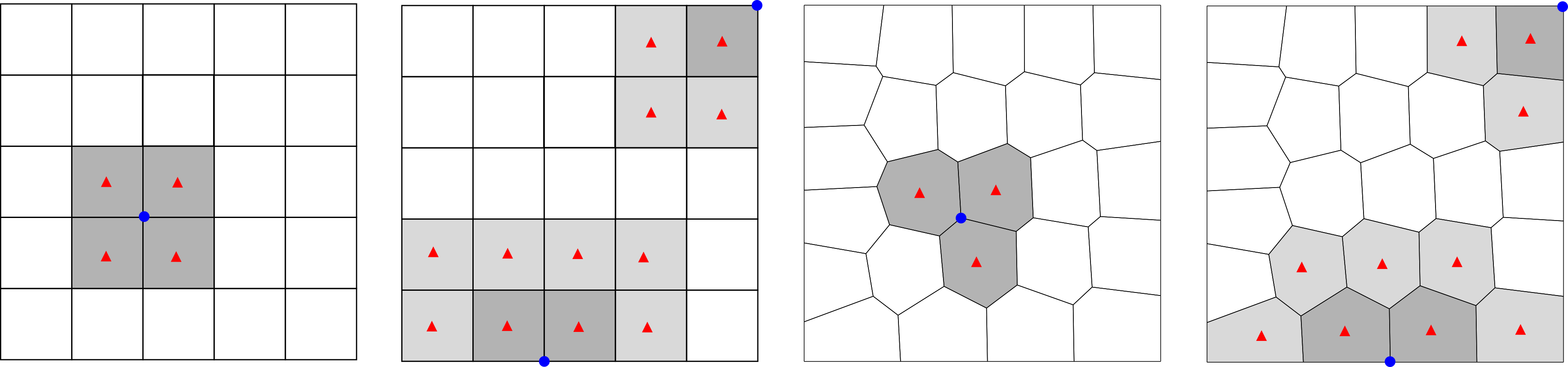}
		\caption{Standard (dark grey) and extended (light grey) patches of elements connected to the node at which to compute the super-convergent stress (blue circle) with super-convergent sampling points (red triangles).
			\label{fig:SuperConvergentStress}}
	\end{figure}

	In this work a low-order VEM is considered where the approximation of the stress field is piece-wise constant. Thus, the approximation of $\bsigma^{\ast}$ should be piece-wise linear and can computed at each node via a least-squares best fit over a patch of elements. The super-convergent stress at a node is computed by considering the patch of elements connected to the node. The location of the centroids of the elements in the patch are treated as the super-convergent sampling points and the element-level stresses are assigned as the degrees of freedom of the sampling points. Since a linear fit is required, at least three sampling points are needed to determine a unique fit. Thus, in cases where a node is connected to less than three elements the patch is enlarged to increase the number of sampling points. Specifically, the patch is enlarged to include elements that are connected to any of the elements in the original patch. For clarity, a few examples of element patches and sampling points are depicted in Figure~\ref{fig:SuperConvergentStress}. Here, the node at which the super-convergent stress is to be computed is indicated as a blue circle, the elements in the patch connected to the node are indicated in dark grey, and (if applicable) the elements included in the enlarged patch are indicated in light grey. Additionally, the locations of the sampling points are indicated as red triangles.
	
	\FloatBarrier
	\FloatBarrier 
	
	The super-convergent stress component ${\sigma^{\ast}_{i}}$ computed over a specific patch is given by
	\begin{equation}
		\sigma^{\ast}_{i} = \bp \left(x, \, y\right) \, \ba_{i} =
		\begin{bmatrix}
			1 & x & y
		\end{bmatrix}
		\begin{bmatrix}
			a_{i}^{1} \\ a_{i}^{2} \\ a_{i}^{3}
		\end{bmatrix}
	\end{equation}
	where $\ba_{i}$ are the degrees of freedom of the super-convergent stress component. The degrees of freedom are computed as 
	\begin{equation}
		\ba_{i} = \bA^{-1} \bb_{i}
	\end{equation}
	where
	\begin{equation}
		\bA = \sum_{k=1}^{n_{\text{sp}}} \bp\left(x_k,\,y_k\right)^{T} \bp\left(x_k,\,y_k\right)
		\quad \text{and} \quad 
		\bb_{i} = \sum_{k=1}^{n_{\text{sp}}} \bp\left(x_k,\,y_k\right)^{T} \sigma^{h}_{i}\left(x_k,\,y_k\right)
	\end{equation}
	respectively. Here $n_{\text{sp}}$ is the number of sampling points, $x_k$ and $y_k$ are the coordinates of the sampling points, and $\sigma^{h}_{i}$ is the stress component at the sampling point (computed via (\ref{eqn:ComputeProjection})).
	
	\section{Procedure for the selection elements to refine and coarsen} 
	\label{sec:ElementSelectionForRefinementAndCoarsening}
	In this section the procedures proposed for the identification of elements qualifying for refinement, and element patches qualifying for coarsening, are presented for various remeshing targets. Specifically, procedures are presented for a target global error, a target number of elements, and a target number of nodes. The error-based target is suited to applications in which the engineer is performing analysis/simulation work with a specified requisite accuracy. The adaptive procedure will then generate a mesh that meets the accuracy target with a quasi-minimal computational load.
	The resource-based targets (i.e., the element and node targets) are suited to applications in which the engineer has a specific computational constraint. The adaptive procedure will then generate a mesh that meets the resource target with a quasi-minimal error. 
	A restriction on the number of elements permitted is a common constraint. However, the computational load, in terms of array sizes and memory allocation, is directly related to the number of degrees of freedom of the system. Thus, motivating the additional presentation of a node-based target.  
	
	\subsection{Target error} 
	\label{subsec:ElementSelectionTargetError}
	The global relative error target $\|e\|_{\text{rel}}^{\text{targ}}$ is set by the user based on their requirements and would typically fall in the range of ${1 - 10 \%}$. From $\|e\|_{\text{rel}}^{\text{targ}}$ and the elastic energy a specific global error target $\|e\|_{\mathcal{H}^{1}}^{\text{targ}}$ is computed as 
	\begin{equation}
		\|e\|_{\mathcal{H}^{1}}^{\text{targ}} = \|\%\, e\|_{\text{rel}}^{\text{targ}} \cdot \|U\| \,.
	\end{equation}
	Assuming an optimal mesh with even error distribution, and considering (\ref{eqn:GlobalFromLocal}), a target element-level error contribution is computed as 
	\begin{equation}
		e_{\text{targ}} = \frac{2\, \left( \|e\|_{\mathcal{H}^{1}}^{\text{targ}} \right)^2}{n_{\text{el}}} \,. \label{eqn:ErrorTarget}
	\end{equation}
	Thus, if ${e_i = e_{\text{targ}} \, \forall i \in [1,\, n_{\text{el}}]}$ the mesh would be optimal and the specified error target would be satisfied. Finally, a corresponding target element level energy error is computed as
	\begin{equation}
		\|e_\text{loc}\|_{\mathcal{H}^{1}}^{\text{targ}} = \left[ \frac{1}{2} \, e_{\text{targ}} \right]^{0.5} \,.
	\end{equation}
	Here the subscript \textit{loc} is introduced for clarity to distinguish the local ${\|e_\text{loc}\|_{\mathcal{H}^{1}}^{\text{targ}}}$ and global ${\|e\|_{\mathcal{H}^{1}}^{\text{targ}}}$ error targets.
	
	Since the objective of this work is to create a quasi-optimal mesh, with quasi-even error distribution, an allowable target error range is introduced. The bounds of this range are based on well-known convergence behaviours. Since a first-order VEM is considered it is expected that under uniform refinement the method would exhibit ${\mathcal{O}(h^1)}$ convergence. That is, if every element is refined the global error $\|e\|_{\mathcal{H}^{1}}$ should decrease by a factor of half. Subsequently, it is expected that if a single element $i$ is refined its local error $\|e_i\|_{\mathcal{H}^{1}}$ should decrease by a factor of a quarter. Therefore, the upper and lower bounds of the allowable element-level error range are chosen to be ${ \lceil e_{\text{loc}}^{\text{targ}} \rceil = 2 \, \|e_\text{loc}\|_{\mathcal{H}^{1}}^{\text{targ}}}$ and ${ \lfloor e_{\text{loc}}^{\text{targ}} \rfloor = 0.5 \, \|e_\text{loc}\|_{\mathcal{H}^{1}}^{\text{targ}}}$ respectively. These bounds cover the error range spanned by one refinement, or one coarsening, iteration of a single element. That is, if an element has a local error equivalent to ${ \lceil e_{\text{loc}}^{\text{targ}} \rceil }$ and is refined its `children'/`successor' elements would each have an error of ${\lfloor e_{\text{loc}}^{\text{targ}} \rfloor}$. 
	Based on these error bounds an element is marked for refinement if ${\|e_i\|_{\mathcal{H}^{1}} > \lceil e_{\text{loc}}^{\text{targ}} \rceil}$ and an element patch is marked for coarsening if ${ \|e_{\text{p}_i}\|_{\mathcal{H}^1} < \lceil e_{\text{loc}}^{\text{targ}} \rceil }$. In addition to the error-prediction criterion, an element patch can only be coarsened if it meets a geometric eligibility criterion. In short, an element patch is eligible for coarsening if the geometry of the coarsened patch does not modify the geometry of the problem domain (for details see \cite{Huyssteen2024}).
	
\begin{figure}[ht!]
		\centering
		\includegraphics[width=0.9\textwidth]{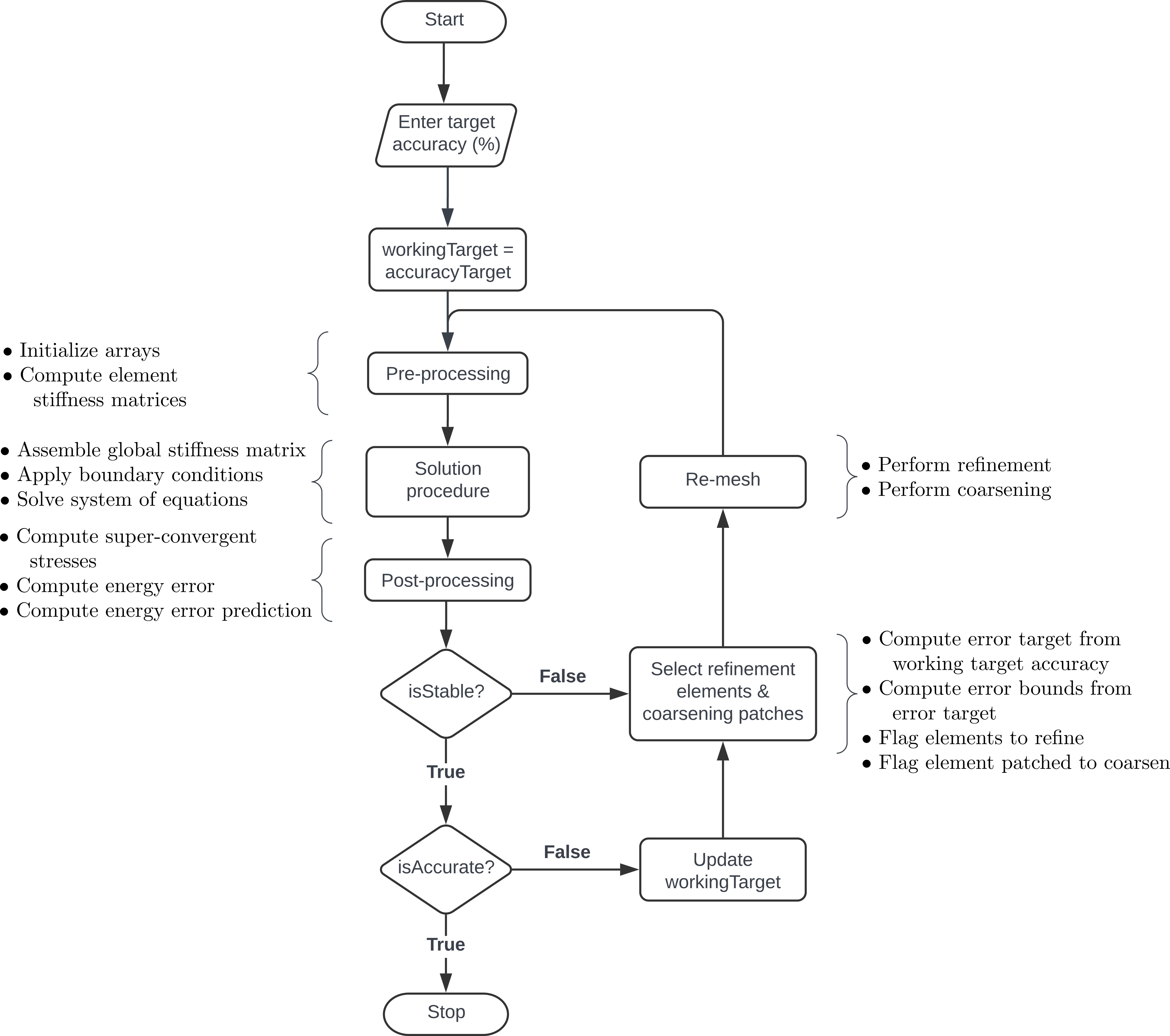}
		\captionof{figure}{Flow diagram for target global error. \label{flow:TargetError}}
	\end{figure}

	An overview of the adaptive procedure for a specified target error is presented in Figure~\ref{flow:TargetError}.
	The user selects a target accuracy, e.g. ${\|e\|_{\text{rel}}^{\text{targ}} = 3\%}$ and this value is set as the working target. The pre-processing, solution procedure, and post processing steps are all performed in a similar manner to a typical finite or virtual element program. A query is made to check if the system is stable. For the system to be stable the global number of nodes and global error must not deviate by more than 1\% for structured meshes and 2\% for Voronoi meshes for at least three successive iterations/loops. If the system is not stable then elements are marked for refinement and element patches are marked for coarsening using the procedure described previously. If the system is stable it is checked for accuracy. The solution is considered sufficiently accurate if the approximation of the global error is within 1\% of the specified accuracy target. I.e., if the target accuracy is ${\|e\|_{\text{rel}}^{\text{targ}} = 3\%}$ then the solution is accurate if ${ 2.97\% \leq \|e\|_{\text{rel}} \leq 3.03\% }$. In the case of Voronoi meshes a 2\% deviation from the target accuracy is permitted.
	In rare cases the global error of a stable system is not sufficiently accurate. In these cases an updated working target accuracy is computed from which updated error bounds are determined. The updated working target is computed by subtracting half of the current discrepancy. For example, if the approximate global accuracy is 3.3\% and the current working target is 3\%, the updated working target will be 2.85\%. Conversely, if the approximate global accuracy is 2.7\% and the current working target is 3\%, the updated working target will be 3.15\%. It has been found through experimentation that the computation of an updated working target is most common in cases of larger target accuracies, typically for ${\|e\|_{\text{rel}}^{\text{targ}} > 8\%}$.
	
	\FloatBarrier
	\FloatBarrier

	\subsection{Target number of elements} 
	\label{subsec:ElementSelectionTargetNumElements}
\begin{figure}[ht!]
		\centering
		\includegraphics[width=0.75\textwidth]{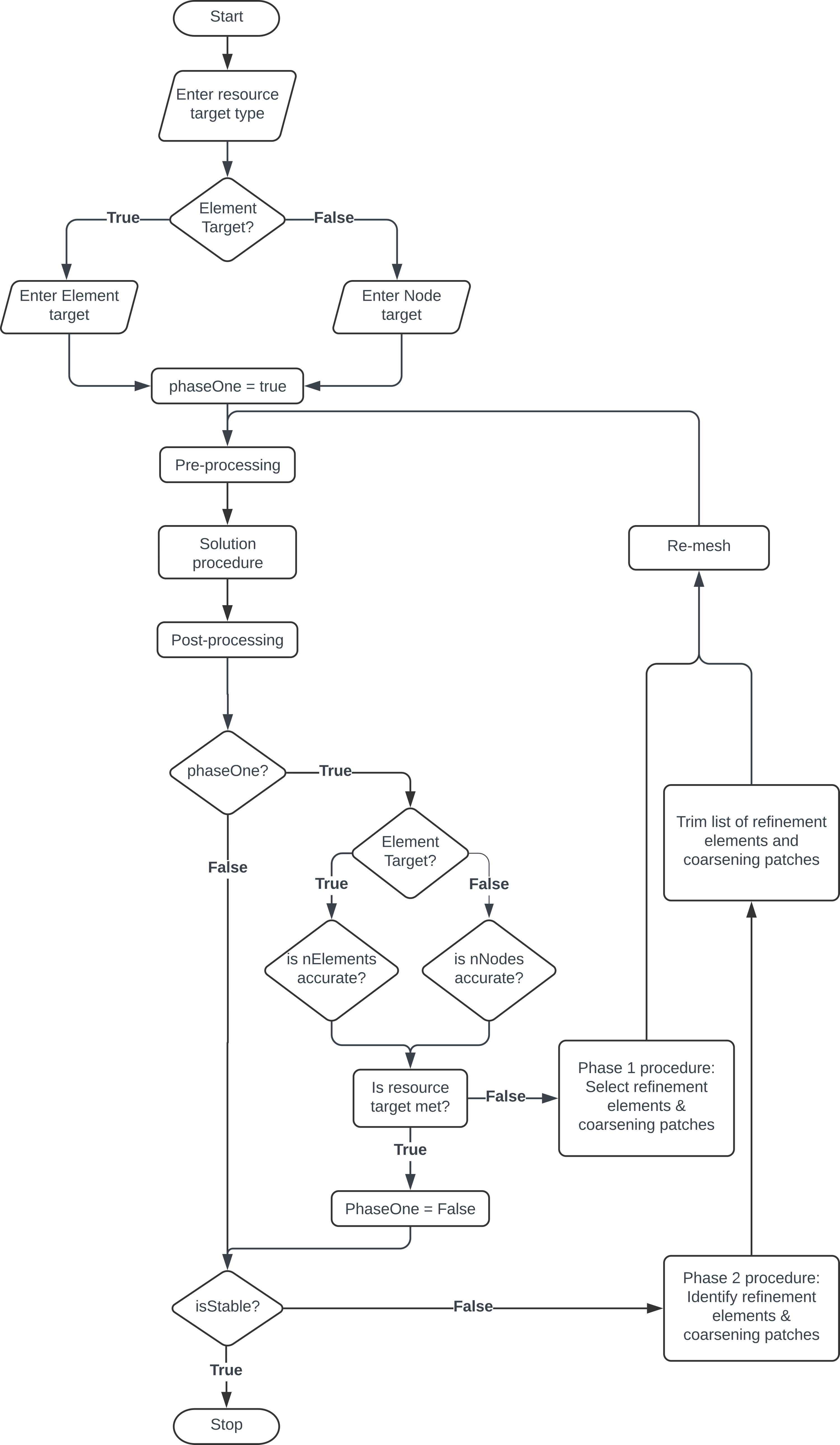}
		\captionof{figure}{Flow diagram for target number of elements or nodes. \label{flow:TargetElementsOrNodes}}
	\end{figure}

	An overview of the program flow for the adaptive procedure for a target number of elements is presented in Figure~\ref{flow:TargetElementsOrNodes}.
	Due to the similarities in the procedures for a target number of elements or target number of nodes (see Section~\ref{subsec:ElementSelectionTargetNumNodes}), as well as for brevity, the program flow in Figure~\ref{flow:TargetElementsOrNodes} is expanded to cover both types of resource-based objective.
	The procedure begins with the user inputting the desired type of resource target (i.e., element target or node target). Then, the specific target is entered. In this case the user inputs the target number of elements ${n_{\text{el}}^{\text{targ}}}$. Thereafter the adaptive procedure comprises two distinct phases. The first phase's objective is to meet the specified element target and the second phase's objective is to keep the number of elements approximately constant while optimizing the error distribution. 
	
	During the first phase of the procedure refinement or coarsening is performed based on simple mesh assumptions and two parameters are introduced. The number of elements that a parent is sub-divided into during refinement is denoted by ${n_{\text{refine}}}$. The number of elements that are grouped together to form one new element during coarsening is denoted by ${n_{\text{coarsen}}}$.
	In the case of structured meshes it is known that refining an element will create four `children' elements ${n_{\text{refine}}^{\text{struct}} = 4}$ and coarsening one element patch will, most often, combine four smaller elements into one larger element ${n_{\text{coarsen}}^{\text{struct}} = 4}$. 
	Therefore, if the current number of elements ${n_{\text{el}}}$ is less than the target number of elements and ${ n_{\text{el}}^{\text{targ}} / n_{\text{el}} \geq n_{\text{refine}} }$ then all elements are marked for refinement. 
	Alternatively, if ${ n_{\text{el}} < n_{\text{el}}^{\text{targ}} }$ and 
	${ n_{\text{el}}^{\text{targ}} / n_{\text{el}} < n_{\text{refine}} }$ the number of elements to mark for refinement is 
	${ n_{\text{el}}^{\text{ref}} = \left[ n_{\text{el}}^{\text{targ}} - n_{\text{el}} \right] / \left[ n_{\text{refine}} - 1 \right] }$ and elements are marked based on their local energy error approximation ${ \|e_i\|_{\mathcal{H}^{1}} }$ in descending order.
	If the current number of elements is greater than the target number of elements and ${ n_{\text{el}} / n_{\text{el}}^{\text{targ}} \geq n_{\text{coarsen}} }$ then all eligible element patches are marked for coarsening.
	Alternatively, if ${ n_{\text{el}} > n_{\text{el}}^{\text{targ}} }$ and ${ n_{\text{el}} / n_{\text{el}}^{\text{targ}} < n_{\text{coarsen}} }$ then the number of element patches to mark for coarsening is
	${ n_{\text{patch}}^{\text{coarsen}} = \left[ n_{\text{el}} - n_{\text{el}}^{\text{targ}} \right] / \left[ n_{\text{coarsen}} - 1 \right] }$ and element patches are marked for coarsening based on their energy error prediction 
	${ \|e_{\text{p}_i}\|_{\mathcal{H}^1} }$ in ascending order. 
	In the case of Voronoi meshes the procedure is the same albeit with ${n_{\text{refine}}^{\text{vrn}} = 5}$ and ${n_{\text{coarsen}}^{\text{vrn}} = 3}$.
	The procedure is repeated iteratively until the resource usage is sufficiently accurate. In the case of a target number of elements the resource usage is sufficiently accurate if the current number of elements ${n_{\text{el}}}$ is within 1\% of the target number of elements ${ n_{\text{el}}^{\text{targ}}  }$ Similarly, in the case of a target number of nodes the resource usage is sufficiently accurate if the current number of nodes ${n_{\text{v}}}$ is within 1\% of the target number of nodes ${ n_{\text{v}}^{\text{targ}}  }$. Once the resource usage is sufficiently accurate the first phase is complete.
	
	During the second phase elements are marked for refinement or coarsening based on their local errors in a similar manner to that described in Section~\ref{subsec:ElementSelectionTargetError}. From the current energy error and number of elements an element-level error target is computed in the same manner as (\ref{eqn:ErrorTarget}). Thereafter, and as described in Section~\ref{subsec:ElementSelectionTargetError}, upper and lower error bounds are computed, elements are identified for refinement, and element patches are identified for coarsening. Before the refinement and coarsening can be performed consideration must be made to keep the number of elements approximately constant. Specifically, the number of elements added by refinement must equal the number of elements removed by coarsening. Therefore, the lists of elements identified for refinement and element patches identified for coarsening must be trimmed. The number of elements added to the system if all identified elements are refined is computed as 
	${n_{\text{add}} = \left[ n_{\text{refine}} - 1 \right] n_{\text{el}}^{\text{refine}} }$ where ${n_{\text{el}}^{\text{refine}} }$ denotes the number of elements identified for refinement. Similarly, the number of elements removed from the system if all identified element patches are coarsened is computed as 
	${n_{\text{rem}} = \left[ n_{\text{coarsen}} - 1 \right] n_{\text{patch}}^{\text{coarsen}} }$ where ${n_{\text{patch}}^{\text{coarsen}} }$ denotes the number of element patches identified for coarsening. The total number of elements to modify is then computed as ${n_{\text{mod}} = \text{min}\left(n_{\text{add}},\, n_{\text{rem}} \right) }$.
	Updated numbers of elements to refine and element patches to coarsen are then computed as 
	${ n_{\text{el}}^{\text{refine}} = n_{\text{mod}} / \left[ n_{\text{refine}} - 1 \right] }$ and 
	${ n_{\text{patch}}^{\text{coarsen}} = n_{\text{mod}} / \left[ n_{\text{coarsen}} - 1 \right] }$ respectively. Elements are then marked for refinement based on their local energy error approximation ${ \|e_i\|_{\mathcal{H}^{1}} }$ in descending order and element patches are marked for coarsening based on their energy error prediction 
	${ \|e_{\text{p}_i}\|_{\mathcal{H}^1} }$ in ascending order. 
	The procedure is repeated iteratively until the same stability criteria as described in Section~\ref{subsec:ElementSelectionTargetError} are met.
	
	\FloatBarrier
	\FloatBarrier
	
	\subsection{Target number of nodes} 
	\label{subsec:ElementSelectionTargetNumNodes}
	Since the adaptive remeshing process involves refining and coarsening elements the number of nodes in a mesh cannot be directly controlled, rather it is a consequence of the number of elements in, and type of, the mesh. As such, the procedure for meeting a target number of nodes ${n_{\text{v}}^{\text{targ}}}$ is based upon approximate relations between the number of elements and the number of nodes for a given mesh type. 
	
	While generating the results for error target and element target based computations a database was built-up comprising the number of elements and nodes in a mesh. This database comprises entries for a broad range of error and element targets, uniform and adapted meshes, and a variety of problem types not presented in this work (for brevity). From the data an approximate exponential relationship was determined for the number of nodes per element given by ${r_{n/e} = a\, \left( n_{\text{v}} \right)^b}$. The parameters of the relationship are presented in Table~\ref{tab:Parameters}.
	
	\FloatBarrier
	\begin{table}[ht!]
		\centering 
		\begin{tabular}{|c|c|c|}
			\hline
			Mesh type & $n_{\text{v}} \leq 1000$  & $n_{\text{v}} > 1000$    \\ \hline
			Structured & \begin{tabular}[c]{@{}c@{}}$a = 2.2763$\\ $b = -0.102$ \end{tabular} 
			& \begin{tabular}[c]{@{}c@{}}$a = 1.5032$\\ $b = -0.04$\end{tabular} \\ \hline
			Voronoi & \begin{tabular}[c]{@{}c@{}}$a = 2.2225$\\ $b = -0.054$\end{tabular} 
			& \begin{tabular}[c]{@{}c@{}}$a = 2.0871$\\ $b = -0.044$\end{tabular} \\ \hline
		\end{tabular}
		\captionof{table}{Parameters for approximate relationship between number of nodes per element and number of nodes in the mesh. \label{tab:Parameters}}
	\end{table}
	\FloatBarrier
	
	Using the approximate relationship between the number of nodes per element and the number of nodes in a mesh the procedure for a target number of elements presented in Section~\ref{subsec:ElementSelectionTargetNumElements} can be trivially modified for a target number of nodes. Since the modification is trivial, presentation of the procedure for a target number of nodes is omitted for brevity. However, an overview of the program flow is presented in Figure~\ref{flow:TargetElementsOrNodes}.
	
	\section{Numerical Results}
	\label{sec:Results}
	
	In this section numerical results are presented to demonstrate the efficacy of the proposed adaptive remeshing and quasi-optimization procedures for various targets.
	The efficacy is evaluated in the $\mathcal{H}^{1}$ semi-norm, i.e. the energy error norm, as described in Section~\ref{sec:ErrorEstimationAndPrediction}.
	
	In the examples that follow the material is isotropic with a Young's modulus of ${E=1~\rm{Pa}}$, a Poisson's ratio of ${\nu=0.3}$, and the shear modulus is computed as ${\mu = E /2 \left[1+\nu\right]}$. 
	Additionally, example problems with quite large deformations are presented. While the material parameters and large deformations may not be realistic for the linear elastic material model and small strain theory used, they are useful to demonstrate the behaviour of the various adaptive remeshing procedures, and are helpful in providing an intuition of where meshes should be refined or coarsened. Furthermore, since the material is linear elastic and small strain theory is used, larger deformations have the same effect as magnifying smaller deformations, which is useful when studying the nature of the mesh adaptation.
	
	\subsection{L-shaped domain}
	\label{subsec:LDomain}
\begin{figure}[ht!]
		\centering
		\begin{subfigure}[t]{0.295\textwidth}
			\centering
			\includegraphics[width=0.95\textwidth]{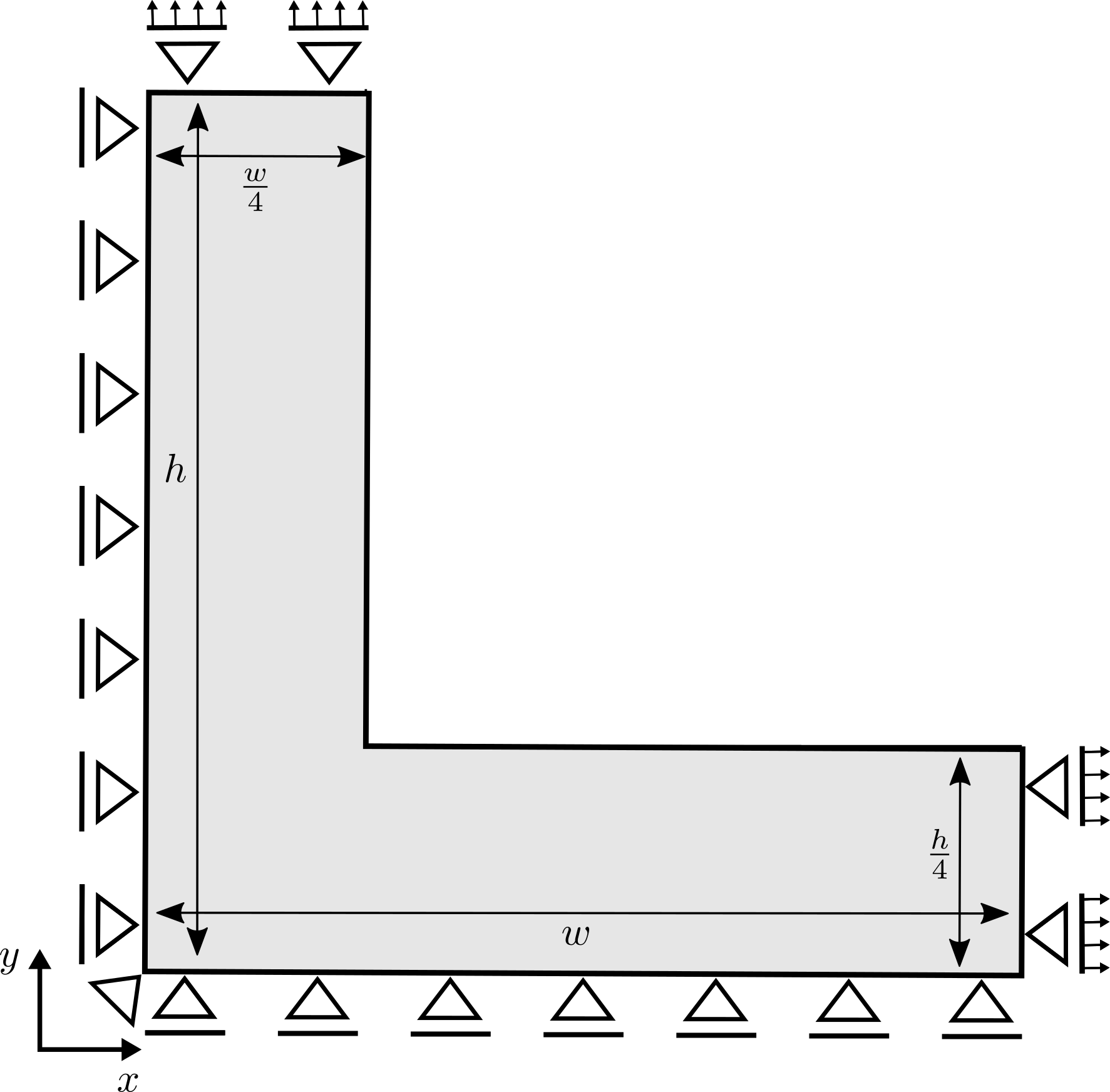}
			\caption{Problem geometry}
		\end{subfigure}%
		\begin{subfigure}[t]{0.345\textwidth}
			\centering
			\includegraphics[width=0.85\textwidth]{{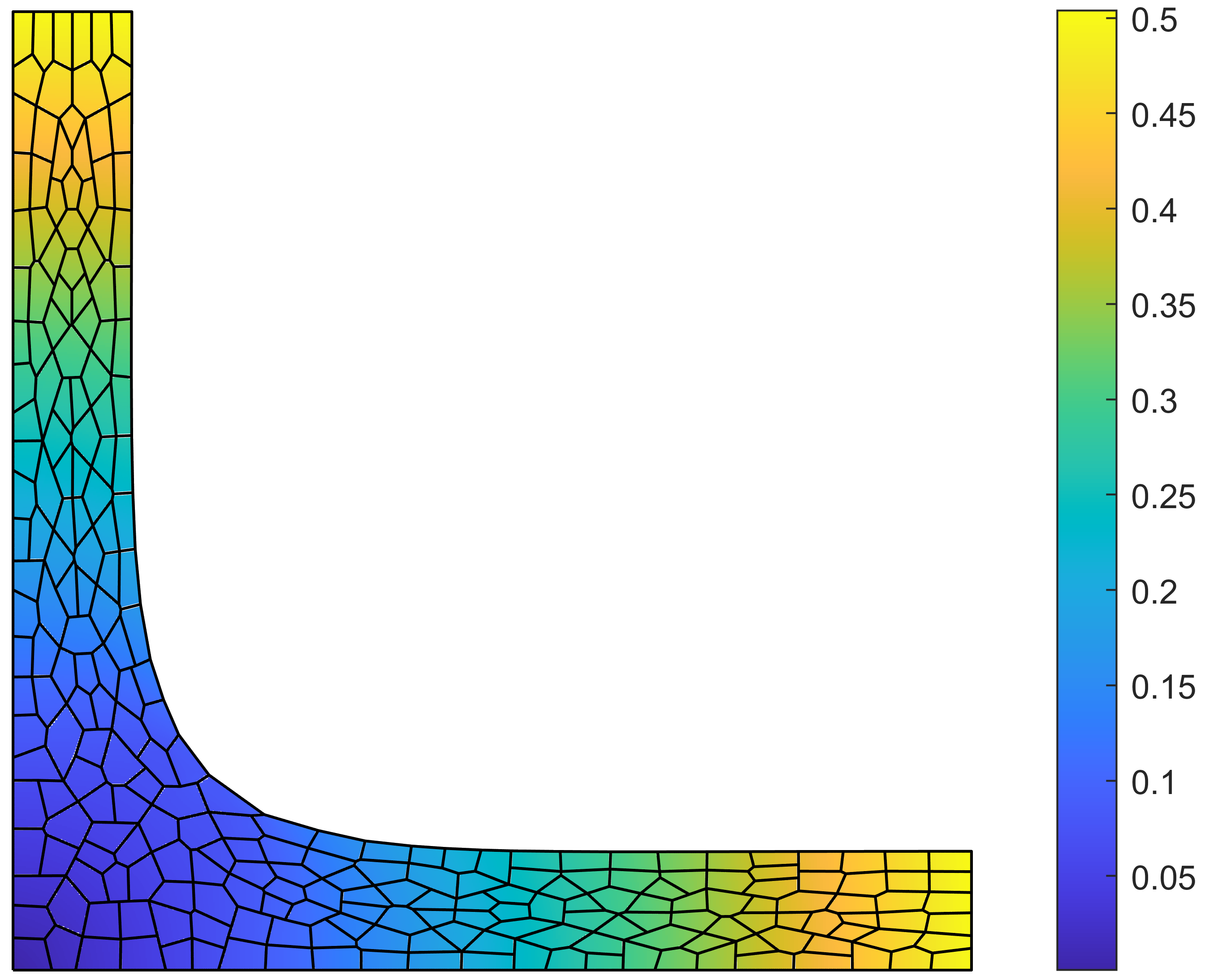}}
			\caption{Deformed configuration}
		\end{subfigure}
		\begin{subfigure}[t]{0.345\textwidth}
			\centering
			\includegraphics[width=0.85\textwidth]{{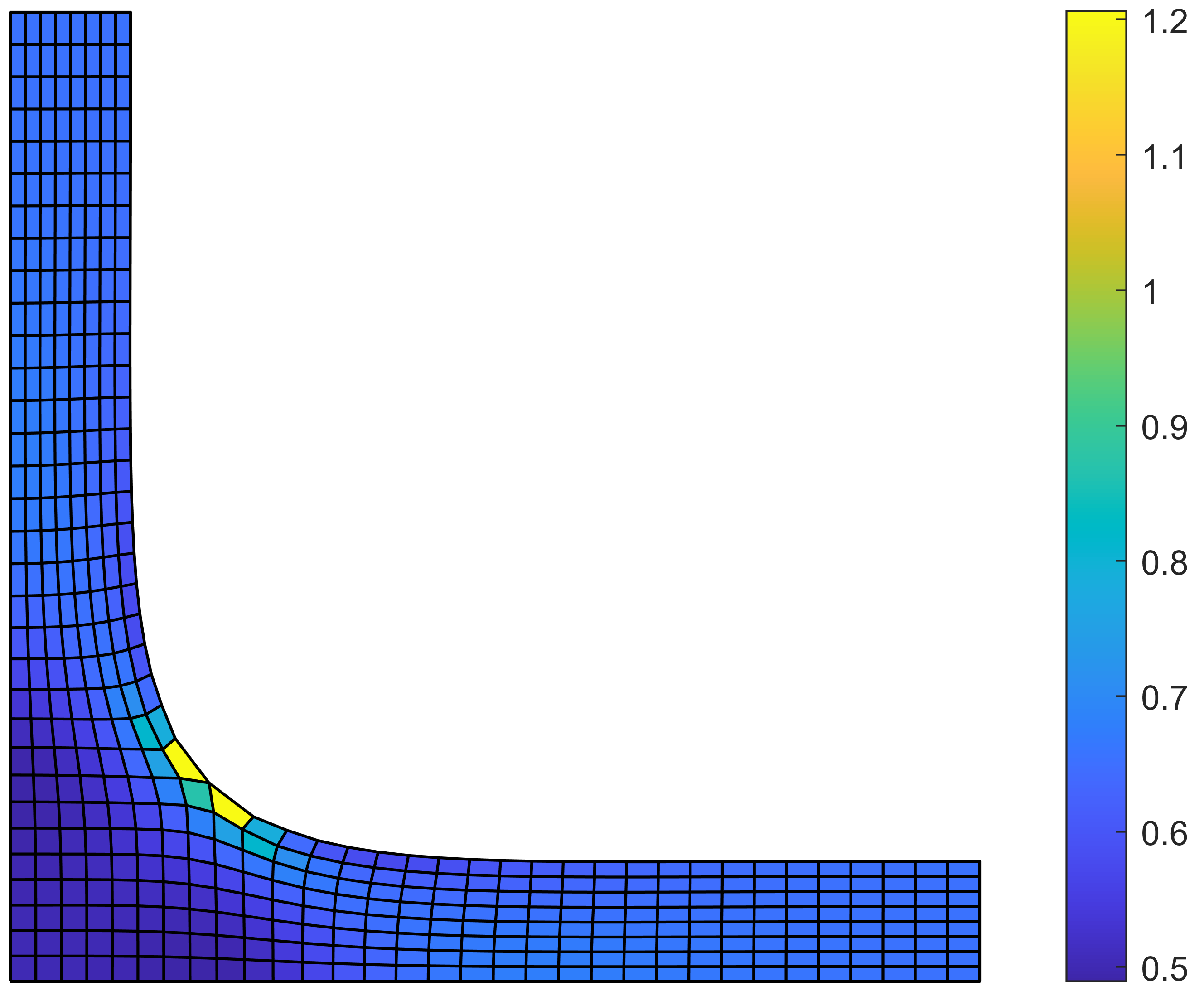}}
			\caption{Von Mises stress}
		\end{subfigure}
		\caption{L-shaped domain (a) geometry, (b) deformed configuration, and (c) von Mises stress distribution (log-scale).
			\label{fig:LShapedDomainGeometry}}
	\end{figure}

	The L-shaped domain problem comprises a domain of width ${w=1~\rm{m}}$ and height ${h=1~\rm{m}}$ where the horizontal and vertical thickness of the L are ${\frac{w}{4}}$ and ${\frac{h}{4}}$ respectively. The bottom and left-hand edges of the domain are constrained vertically and horizontally respectively, with the bottom left corner fully constrained. The upper and right-hand edges are subject to prescribed displacements of ${\bar{u}_{y}=0.5~\rm{m}}$ and ${\bar{u}_{x}=0.5~\rm{m}}$ respectively, with the displacements of the edges unconstrained in the $x$- and $y$-directions respectively (see Figure~\ref{fig:LShapedDomainGeometry}(a)). 
	The deformed configuration of the body is illustrated in Figures~\ref{fig:LShapedDomainGeometry}(b) and (c) with the displacement magnitude $|\bu|$ and von Mises stress respectively plotted on the colour axis.
	From these figures it is clear that the more complex (i.e., more heterogeneous) deformations, high stresses, and high stress gradients are localized to the internal corner of the L where the domain's geometry induces a stress singularity. 
	Conversely, the deformation throughout the rest of the domain is much simpler (i.e., more homogeneous) and the stresses are much smoother and lower. As such, the L-shaped domain problem represents a thorough test and ideal application for an adaptive remeshing procedure.
	
	\FloatBarrier
	\FloatBarrier
	
	\subsubsection{Target error}
	\label{subsubsec:L_Domain_TargetError}
	
\begin{figure}[ht!]
		\centering
		\begin{subfigure}[t]{0.33\textwidth}
			\centering
			\includegraphics[width=0.95\textwidth,height=0.95\textwidth]{{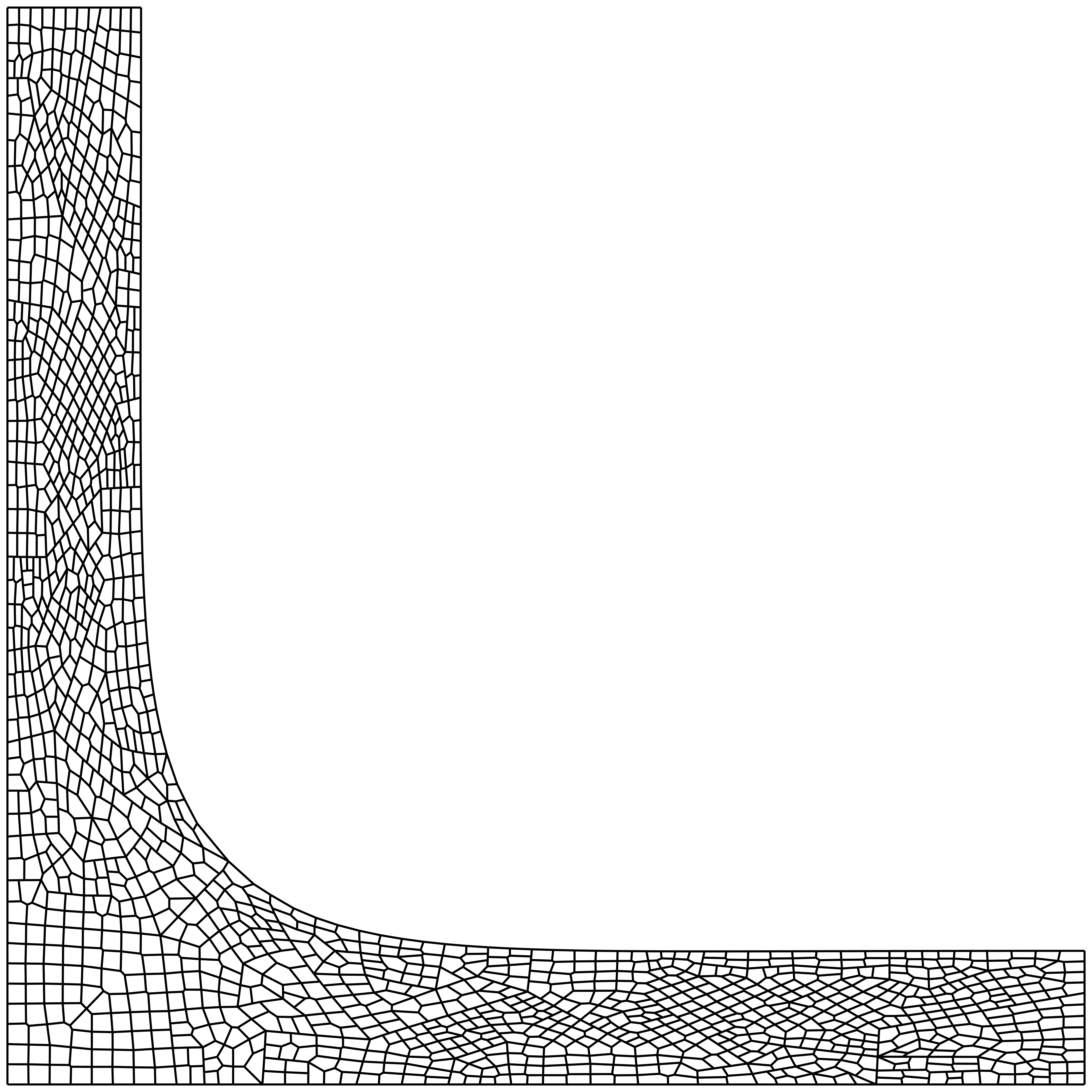}}
			\caption{Step 1: Initial uniform mesh}
		\end{subfigure}%
		\begin{subfigure}[t]{0.33\textwidth}
			\centering
			\includegraphics[width=0.95\textwidth,height=0.95\textwidth]{{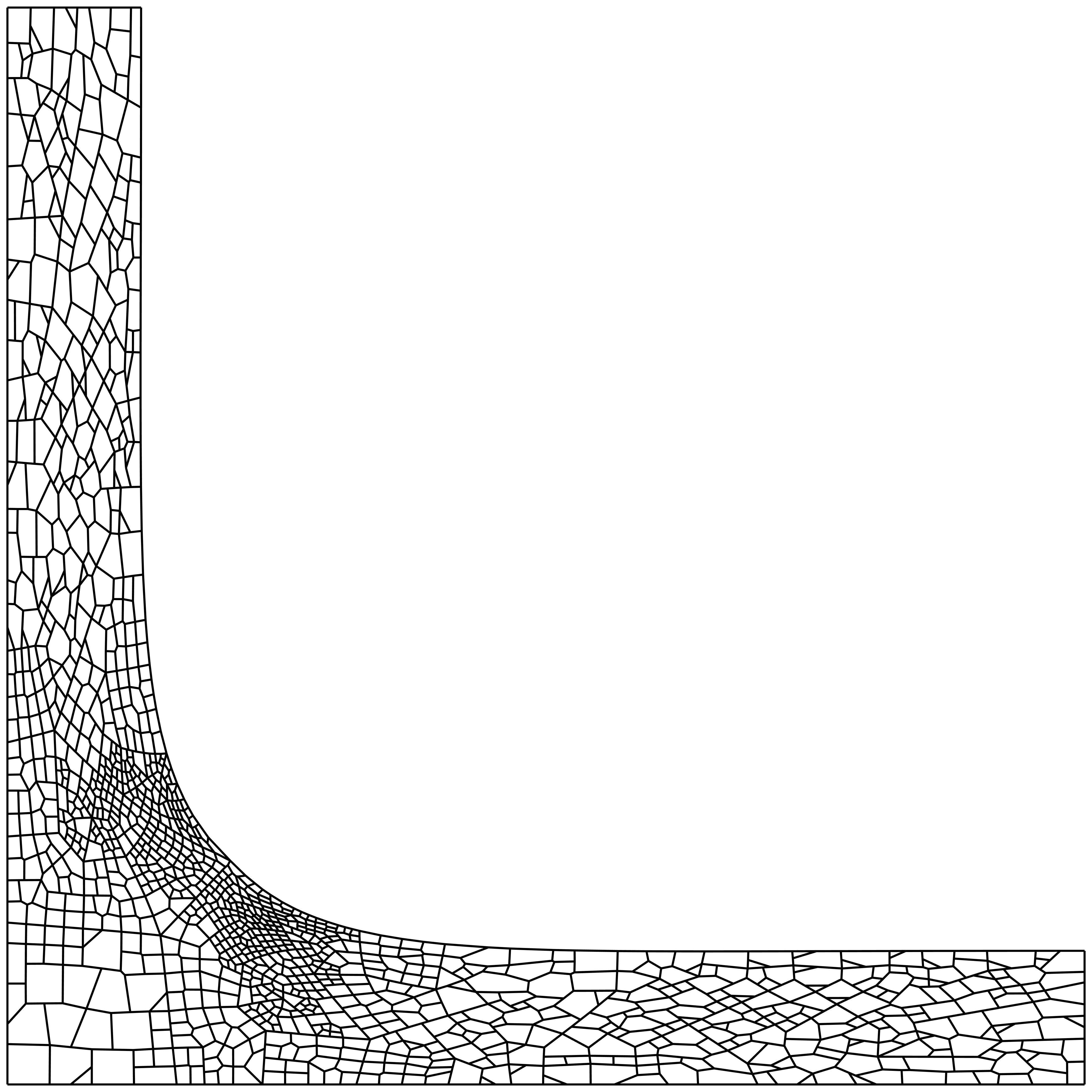}}
			\caption{Step 2}
		\end{subfigure}%
		\begin{subfigure}[t]{0.33\textwidth}
			\centering
			\includegraphics[width=0.95\textwidth,height=0.95\textwidth]{{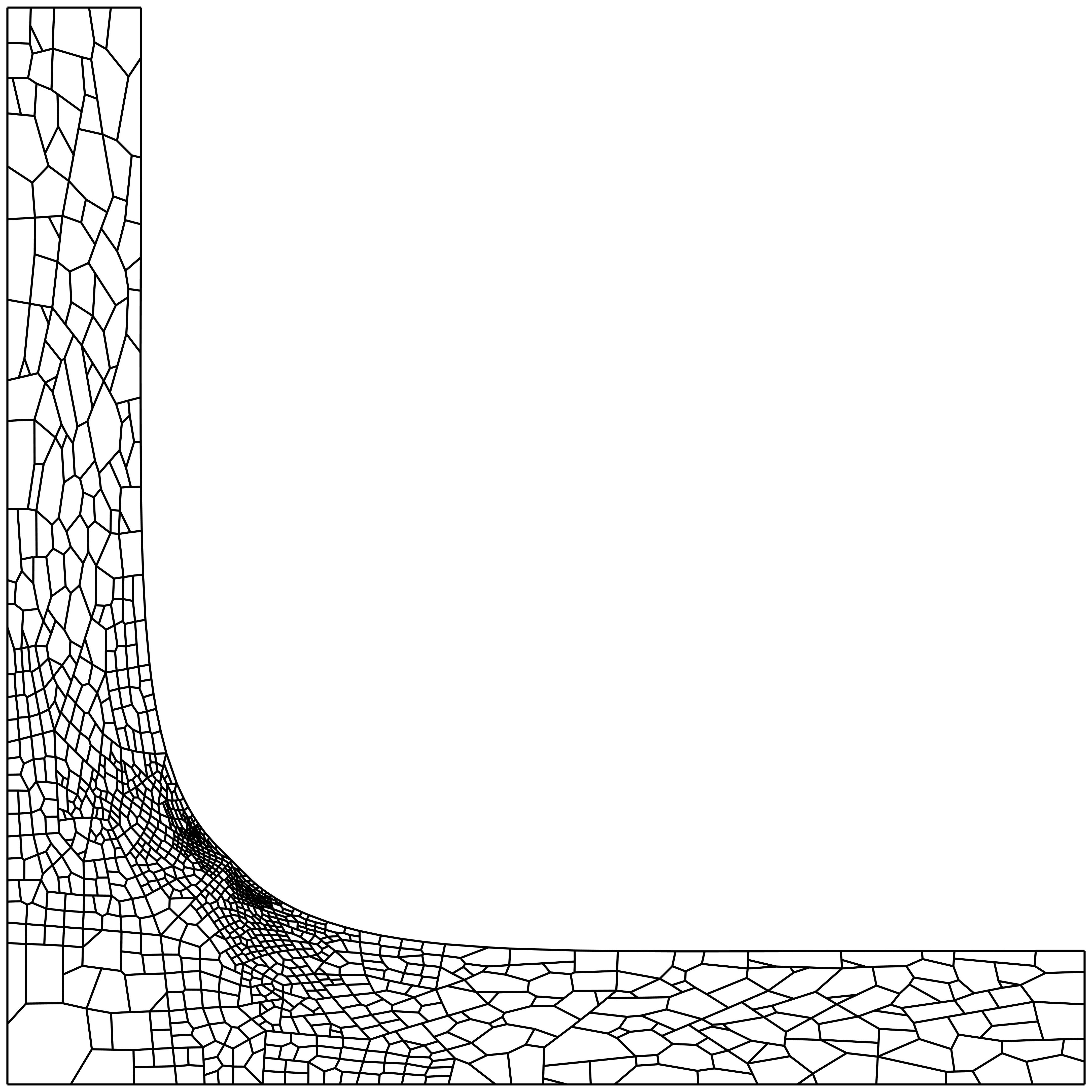}}
			\caption{Step 3}
		\end{subfigure}
		\vskip \baselineskip 
		\begin{subfigure}[t]{0.33\textwidth}
			\centering
			\includegraphics[width=0.95\textwidth,height=0.95\textwidth]{{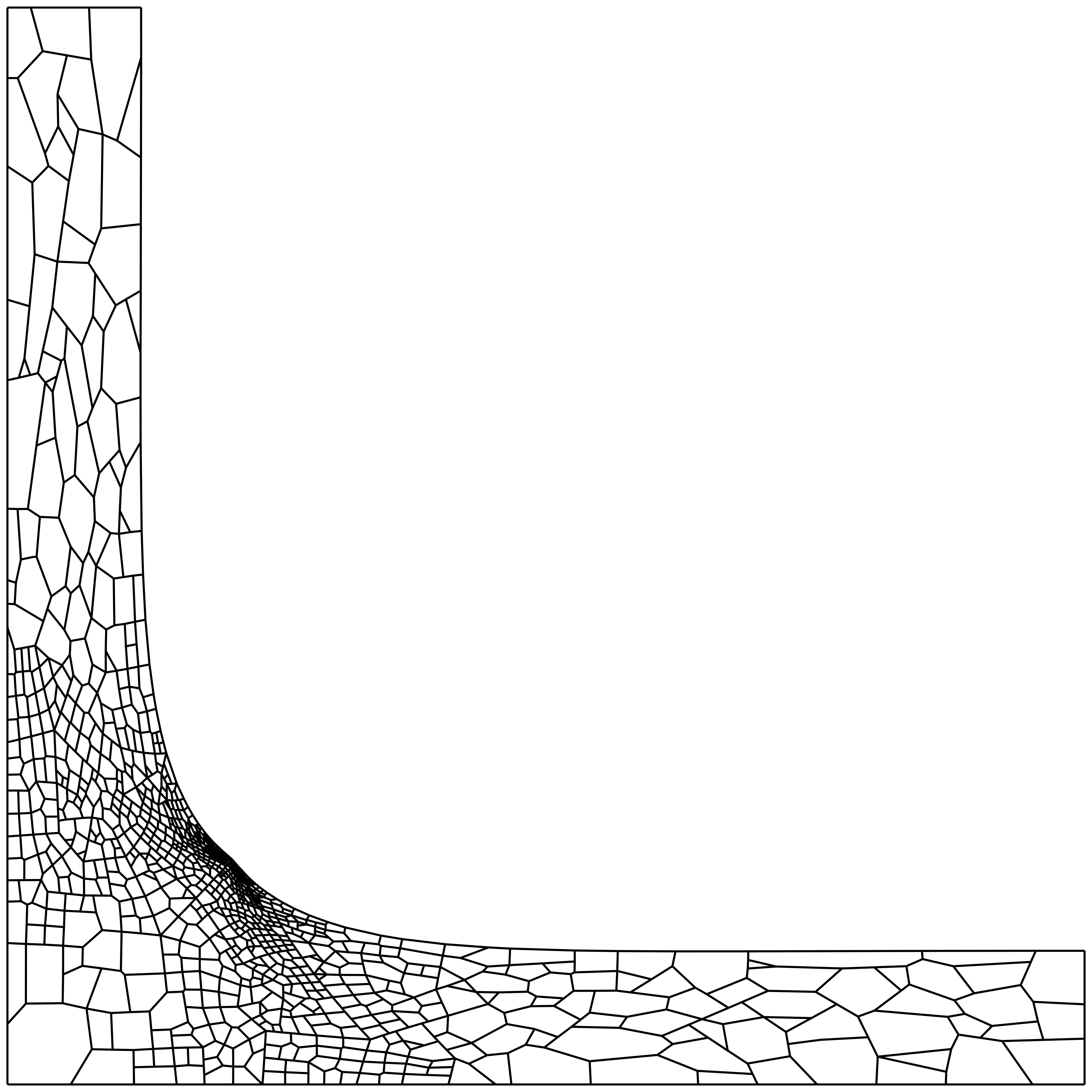}}
			\caption{Step 4}
		\end{subfigure}%
		\begin{subfigure}[t]{0.33\textwidth}
			\centering
			\includegraphics[width=0.95\textwidth,height=0.95\textwidth]{{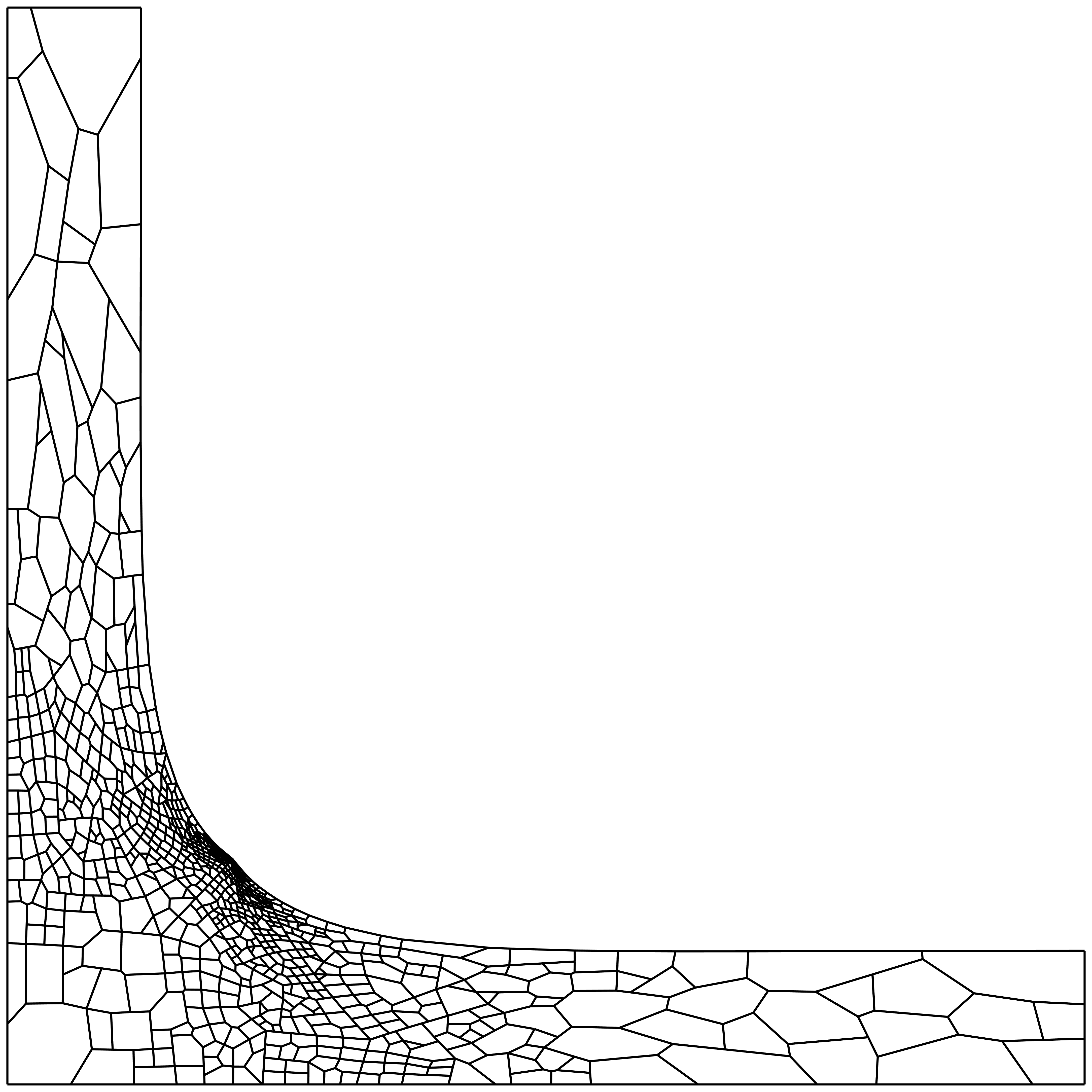}}
			\caption{Step 5}
		\end{subfigure}%
		\begin{subfigure}[t]{0.33\textwidth}
			\centering
			\includegraphics[width=0.95\textwidth,height=0.95\textwidth]{{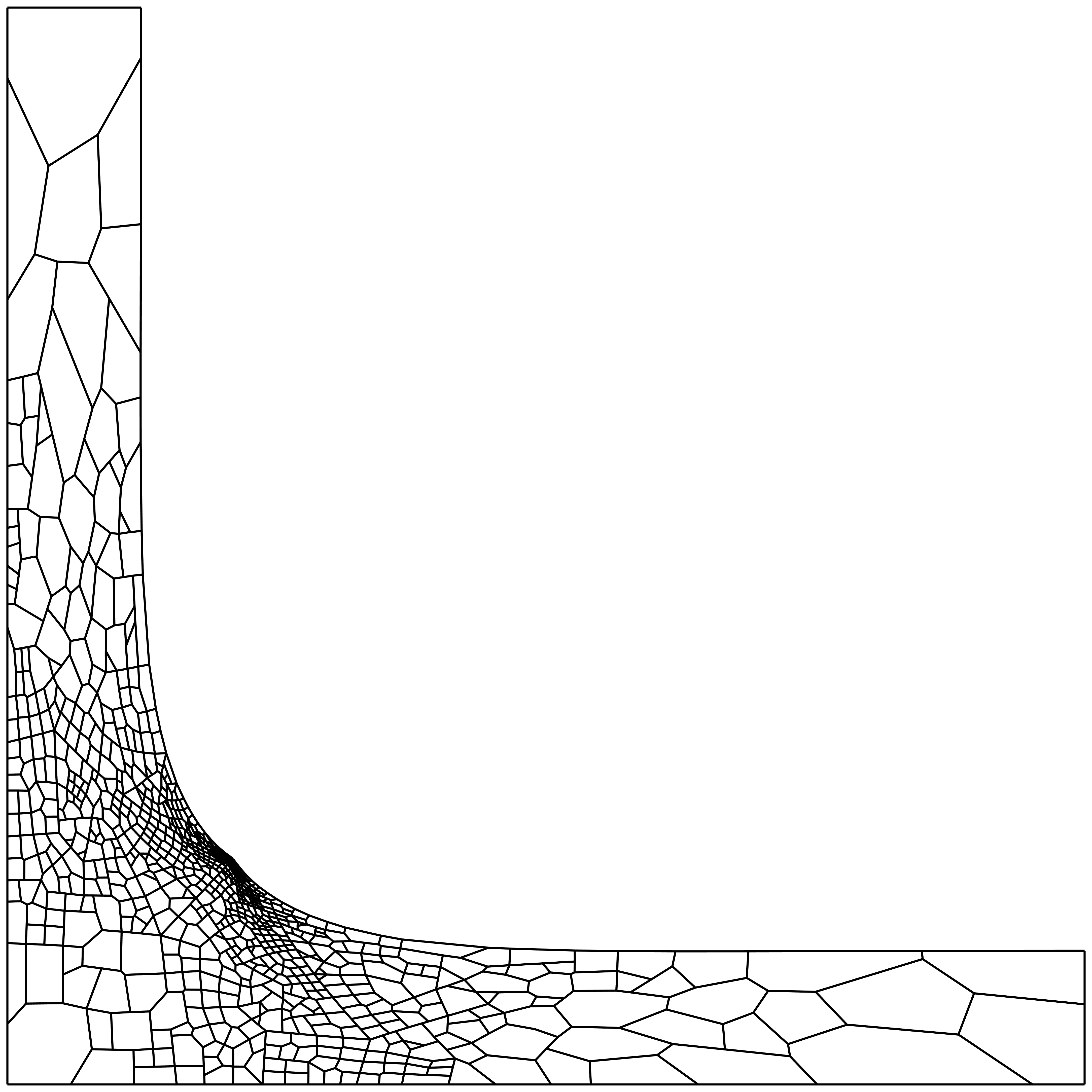}}
			\caption{Step 6: Final adapted mesh}
		\end{subfigure}
		\caption{Mesh evolution with the fully adaptive remeshing procedure for the L-shaped domain problem from an initial uniform Voronoi mesh with ${\|e\|_{\text{rel}}^{\text{targ}} = 3\%}$.
			\label{fig:LDomain_TargetError_MeshEvolution}}
	\end{figure}

	The mesh evolution during the adaptive remeshing process for the L-shaped domain problem is depicted in Figure~\ref{fig:LDomain_TargetError_MeshEvolution} for an initially uniform Voronoi mesh with an error target of ${\|e\|_{\text{rel}}^{\text{targ}} = 3\%}$.
	The adaptive procedure generates a very sensible and intuitive mesh evolution for this problem.
	The areas of the domain with the highest stresses and stress gradients are increasingly refined while the regions of the domain experiencing simpler (quasi-homogeneous) deformations and lower stresses are coarsened.
	Furthermore, the strongest concentration of elements and the highest refinement level is generated in the internal corner of the L and coincident with the stress singularity. 
	
	\FloatBarrier
	\FloatBarrier
	
\begin{figure}[ht!]
		\centering
		\begin{subfigure}[t]{0.33\textwidth}
			\centering
			\includegraphics[width=0.95\textwidth,height=0.95\textwidth]{{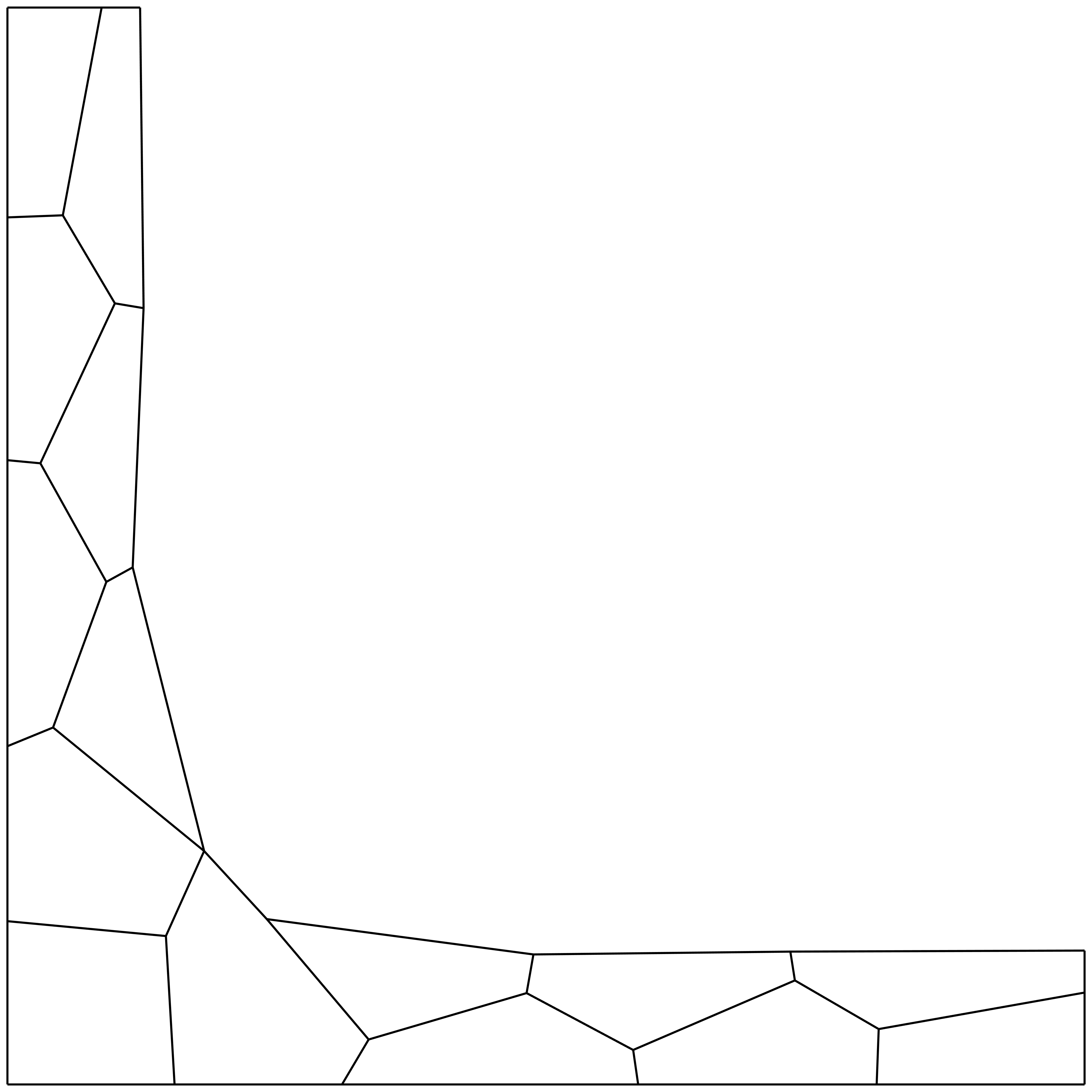}}
			\caption{Coarse mesh: Initial mesh}
		\end{subfigure}%
		\begin{subfigure}[t]{0.33\textwidth}
			\centering
			\includegraphics[width=0.95\textwidth,height=0.95\textwidth]{{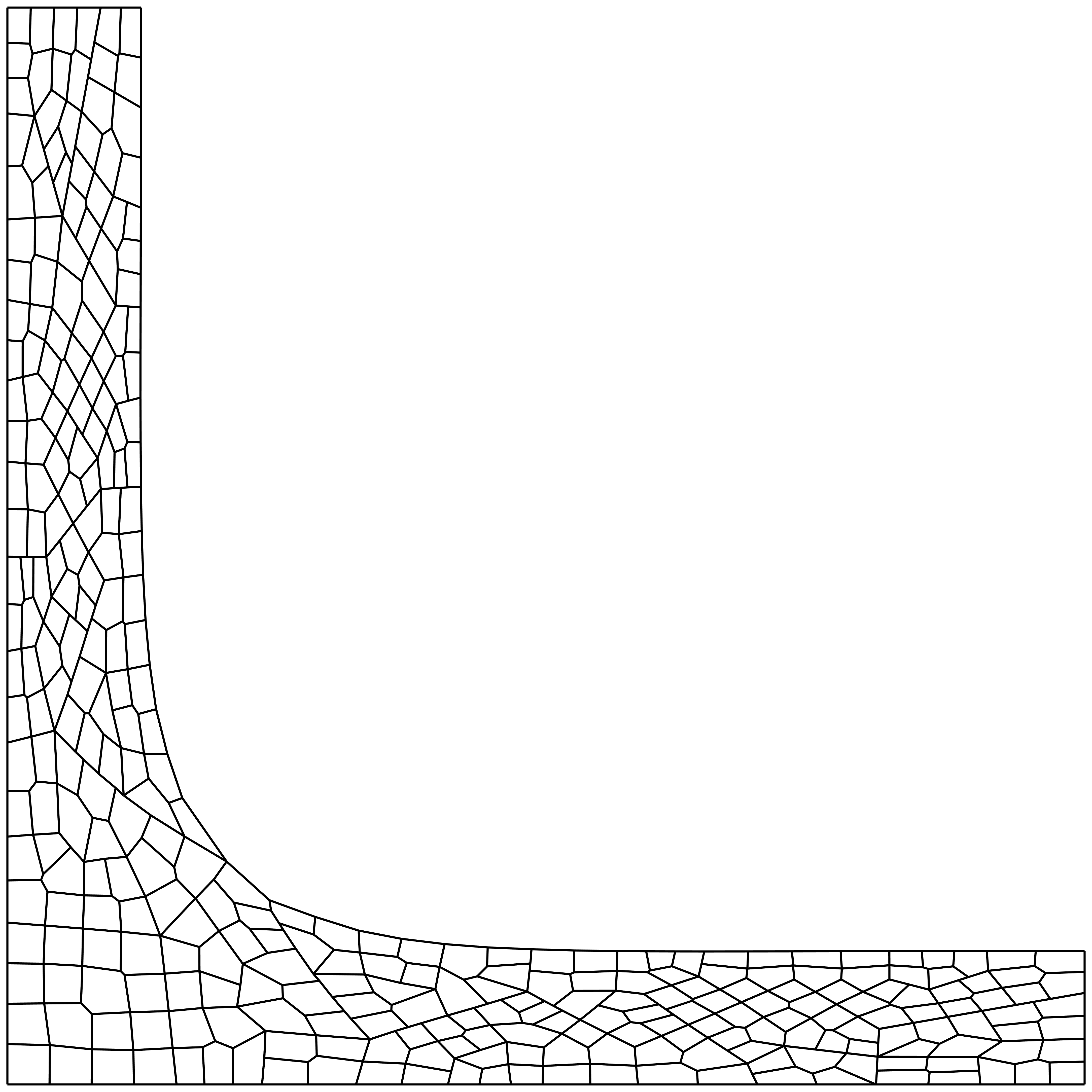}}
			\caption{Intermediate mesh: Initial mesh}
		\end{subfigure}%
		\begin{subfigure}[t]{0.33\textwidth}
			\centering
			\includegraphics[width=0.95\textwidth,height=0.95\textwidth]{{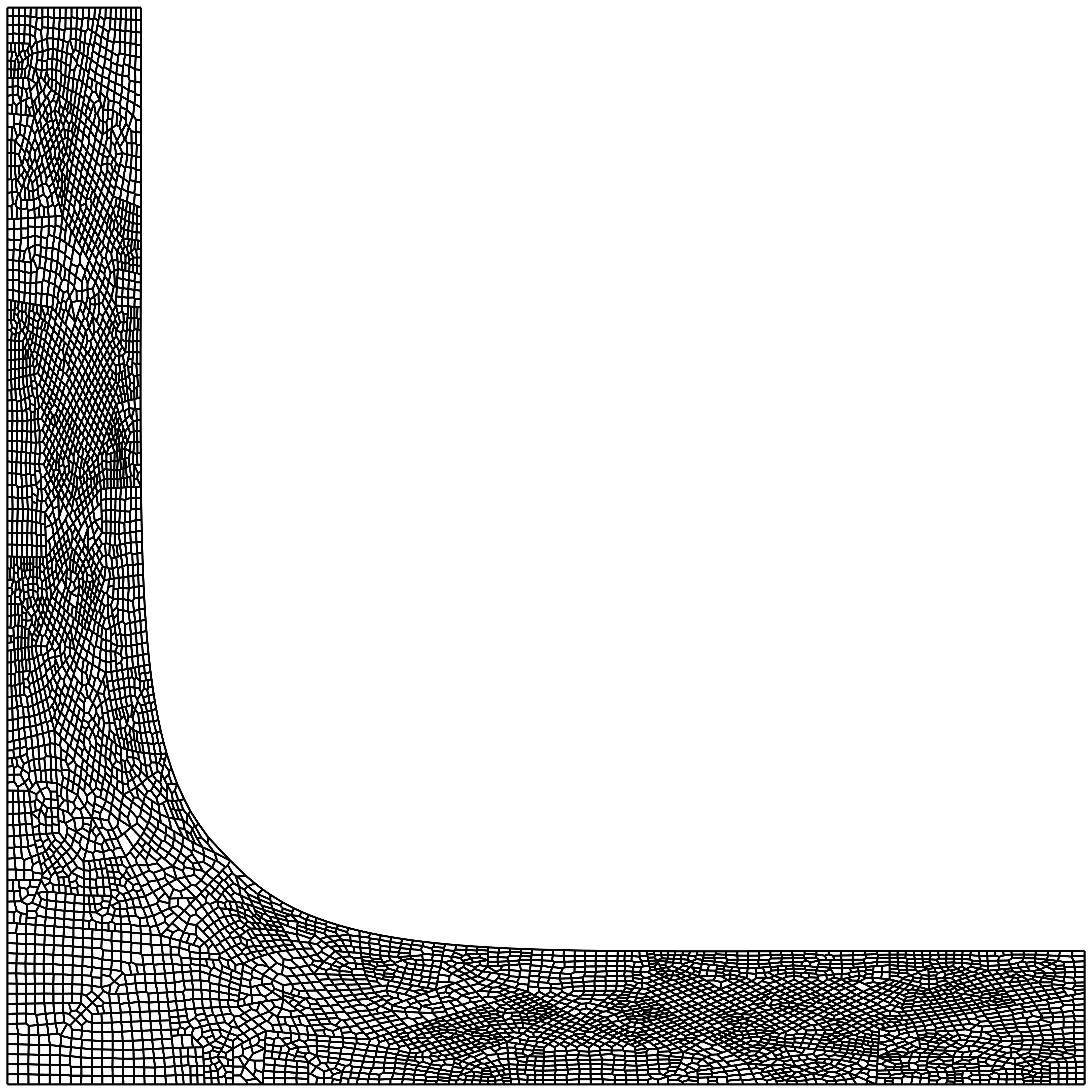}}
			\caption{Fine mesh: Initial mesh}
		\end{subfigure}
		\vskip \baselineskip 
		\begin{subfigure}[t]{0.33\textwidth}
			\centering
			\includegraphics[width=0.95\textwidth,height=0.95\textwidth]{{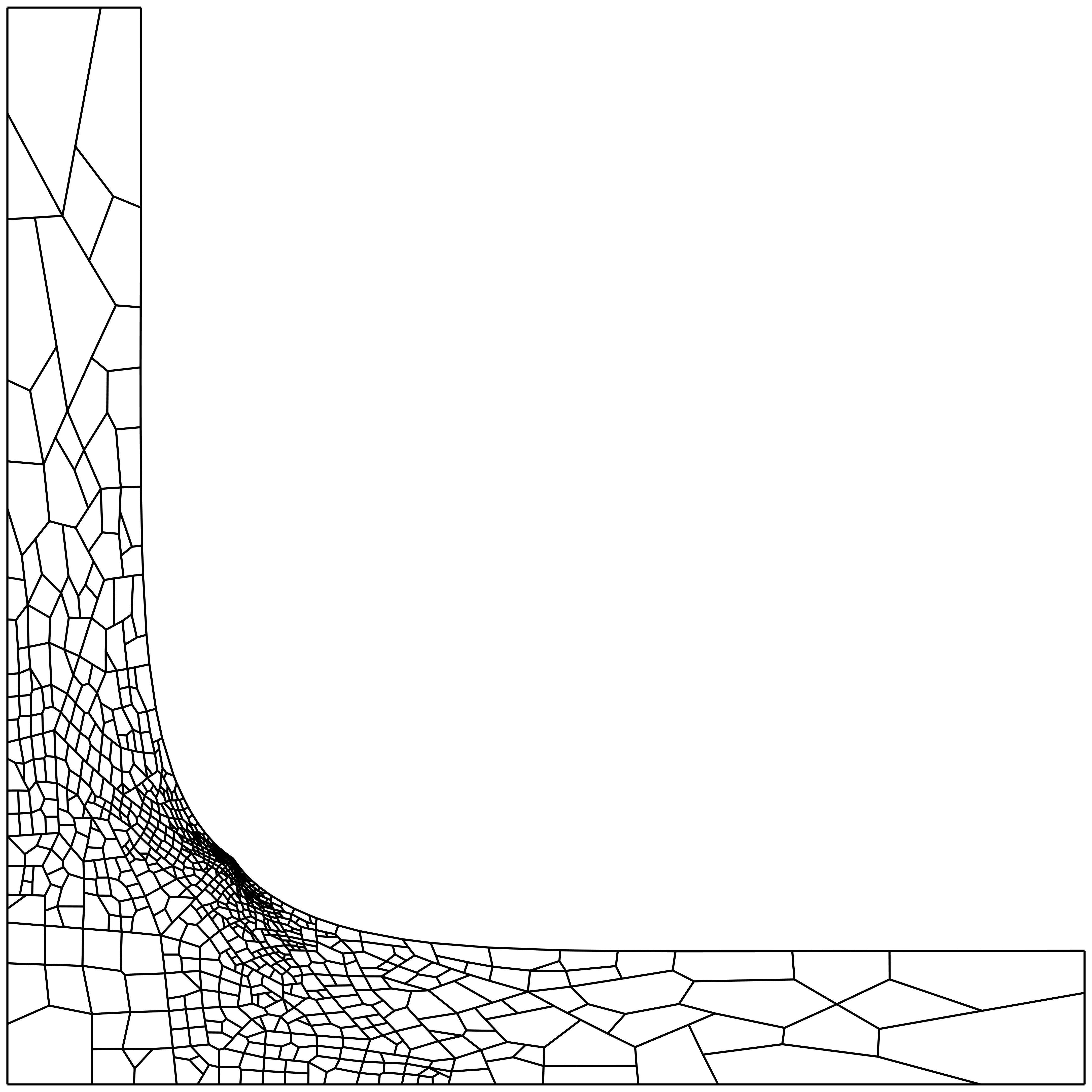}}
			\caption{Coarse mesh: Final mesh}
		\end{subfigure}%
		\begin{subfigure}[t]{0.33\textwidth}
			\centering
			\includegraphics[width=0.95\textwidth,height=0.95\textwidth]{{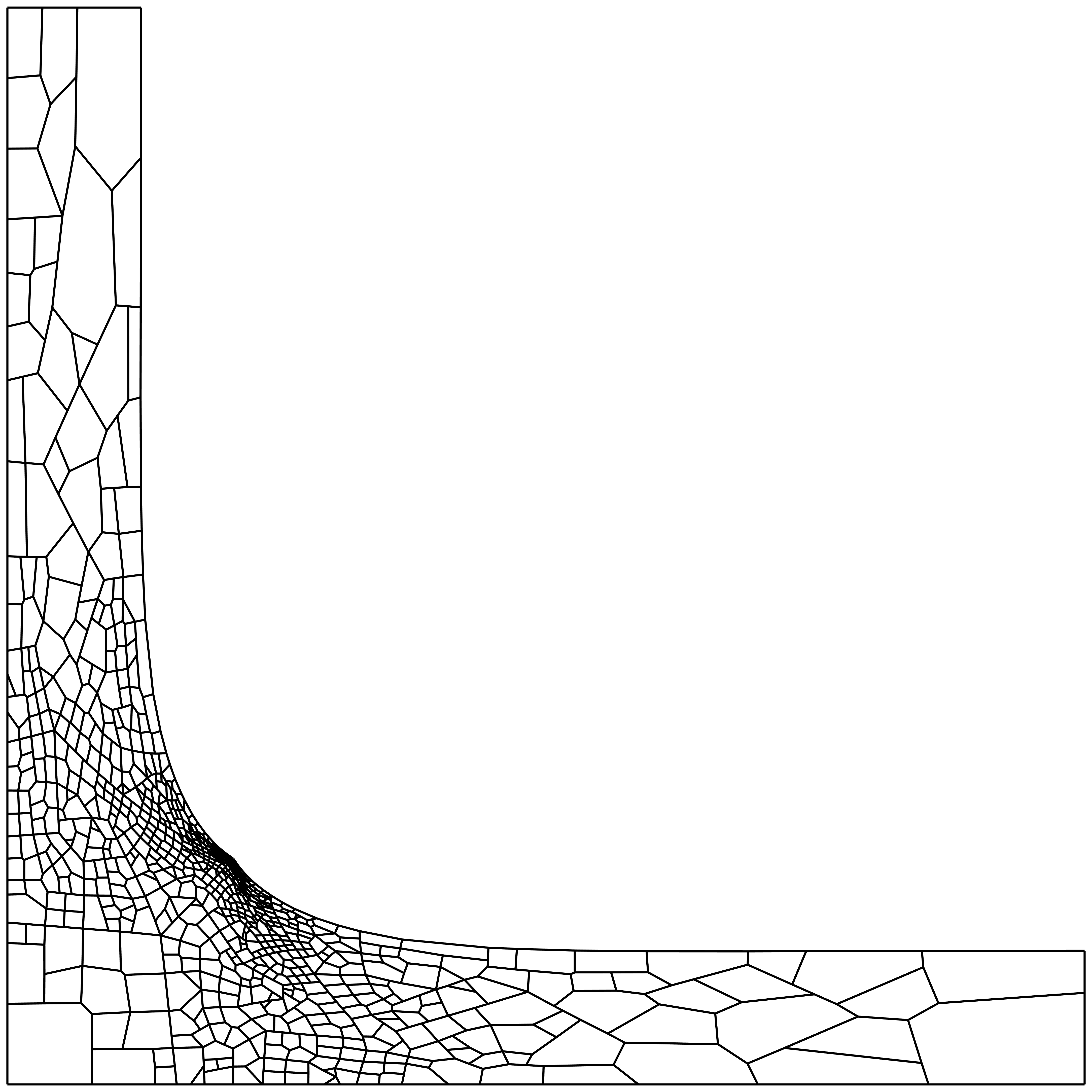}}
			\caption{Intermediate mesh: Final mesh}
		\end{subfigure}%
		\begin{subfigure}[t]{0.33\textwidth}
			\centering
			\includegraphics[width=0.95\textwidth,height=0.95\textwidth]{{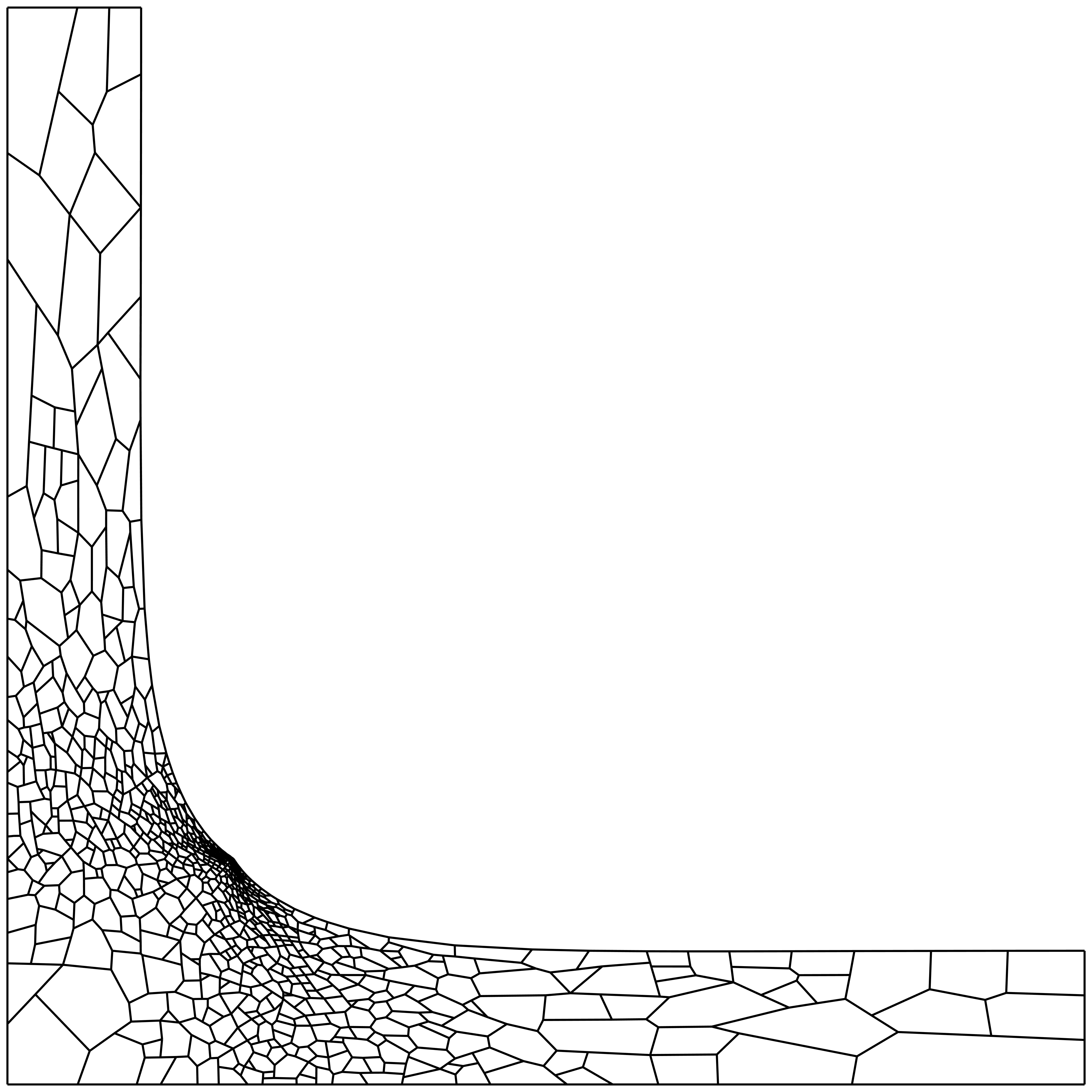}}
			\caption{Fine mesh: Final mesh}
		\end{subfigure}
		\caption{Initial and final adapted meshes generated by the fully adaptive remeshing procedure for the L-shaped domain problem from initial uniform Voronoi meshes of various densities with ${\|e\|_{\text{rel}}^{\text{targ}} = 3\%}$.
			\label{fig:LDomain_TargetError_MeshEvolution_VariousInitialMeshes}}
	\end{figure}

	The results of the adaptive remeshing process for the L-shaped domain problem are depicted in Figure~\ref{fig:LDomain_TargetError_MeshEvolution_VariousInitialMeshes} for initially uniform Voronoi meshes of various densities with an error target of ${\|e\|_{\text{rel}}^{\text{targ}} = 3\%}$. The top row of figures depicts the initial meshes while the bottom row depicts the final adapted meshes after the error target and termination criteria have been met.
	The final adapted meshes exhibit the same sensible and intuitive element distribution as observed in Figure~\ref{fig:LDomain_TargetError_MeshEvolution}. Most notably, the final adapted meshes are almost identical for all initial uniform meshes considered. Thus, the output, or final result, of the adaptive remeshing procedure is independent of the initial mesh and depends only on the specified error target. Furthermore, the cases of the `Intermediate' and `Fine' initial meshes demonstrate that the fully adaptive remeshing procedure is able to perform coarsening from any mesh and does not require knowledge of a previously coarser state. This ability distinguishes the proposed fully adaptive procedure from other procedures surveyed in the literature. 
	
	\FloatBarrier
	\FloatBarrier
	
\begin{figure}[ht!]
		\centering
		\begin{subfigure}[t]{0.495\textwidth}
			\centering
			\includegraphics[width=0.95\textwidth]{{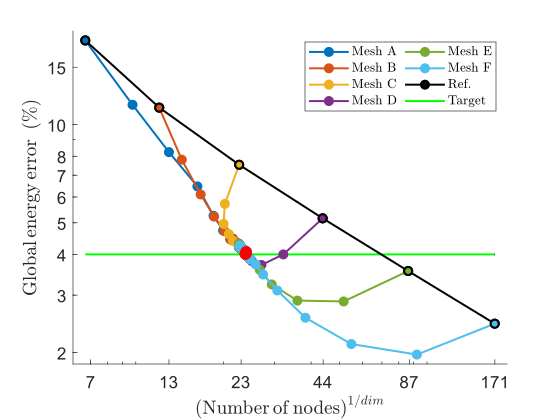}}
			\caption{$\|e\|_{\text{targ}}^{\text{rel}} = 4\%$}
		\end{subfigure}%
		\begin{subfigure}[t]{0.495\textwidth}
			\centering
			\includegraphics[width=0.95\textwidth]{{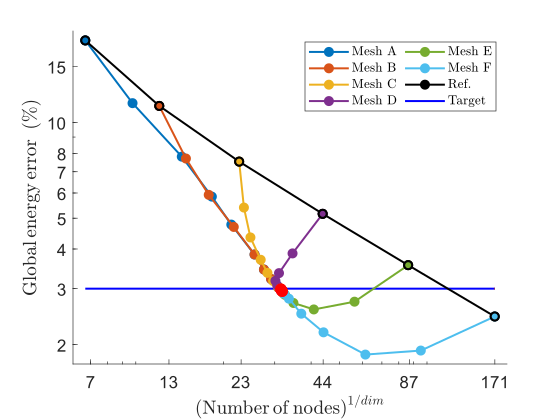}}
			\caption{$\|e\|_{\text{targ}}^{\text{rel}} = 3\%$}
		\end{subfigure}
		\vskip \baselineskip 
		\begin{subfigure}[t]{0.495\textwidth}
			\centering
			\includegraphics[width=0.95\textwidth]{{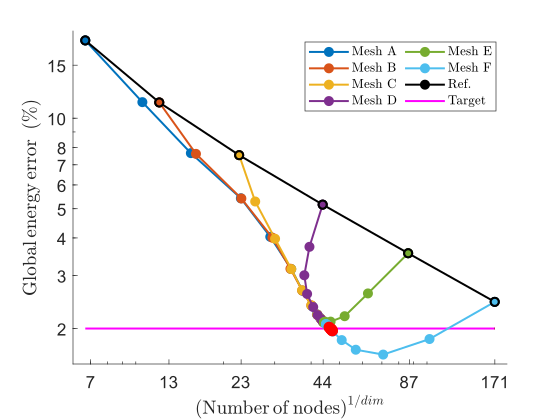}}
			\caption{$\|e\|_{\text{targ}}^{\text{rel}} = 2\%$}
		\end{subfigure}%
		\begin{subfigure}[t]{0.495\textwidth}
			\centering
			\includegraphics[width=0.95\textwidth]{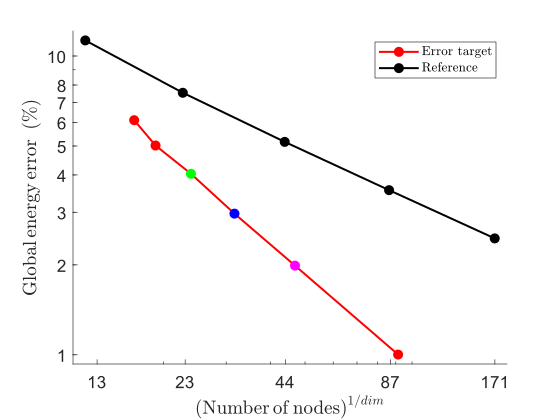}
			\caption{Various error targets}
		\end{subfigure}
		\caption{Energy error vs number of nodes for the L-shaped domain problem on structured meshes of various initial densities with (a) ${\|e\|_{\text{rel}}^{\text{targ}} = 4\%}$, (b) ${\|e\|_{\text{rel}}^{\text{targ}} = 3\%}$, (c) ${\|e\|_{\text{rel}}^{\text{targ}} = 2\%}$, and (d) average final/adapted mesh error for various error targets. 
			\label{fig:LDomain_TargetError_ErrorConvergence_Structured}}
	\end{figure}

	The convergence behaviour of the energy error approximation vs the number of nodes in the mesh is depicted on a logarithmic scale in Figures~\ref{fig:LDomain_TargetError_ErrorConvergence_Structured}(a)-(c) for the L-shaped domain problem on structured meshes. The convergence behaviour is plotted for cases of several initially uniform structured meshes of varying density (denoted by `Meshes A-F') for various error targets. For readability purposes the outline of the marker denoting the first step, or initial mesh, is indicated in black and a red marker is used to indicate the final adapted mesh result for each of the initial meshes. Additionally, the black `Reference' curve corresponds to the standard convergence behaviour under uniform refinement, i.e. all elements are refined.
	Furthermore, for each error target considered (not all have been shown here) an average final mesh result is computed, i.e. the average position of the red markers for each error target. These averaged results are plotted in Figure~\ref{fig:LDomain_TargetError_ErrorConvergence_Structured}(d) along with the reference uniform convergence curve. 
	Where applicable, the markers in Figure~\ref{fig:LDomain_TargetError_ErrorConvergence_Structured}(d) are colored to match their corresponding targets depicted in Figures~\ref{fig:LDomain_TargetError_ErrorConvergence_Structured}(a)-(c).
	From Figures~\ref{fig:LDomain_TargetError_ErrorConvergence_Structured}(a)-(c) it is clear that the fully adaptive procedure is able to meet the specified global error targets from any initial mesh. Furthermore, the final adapted meshes contain an almost identical number of nodes for a specific error target. This, again, demonstrates that the performance of the fully adaptive remeshing procedure is independent of the initial mesh.
	From Figure~\ref{fig:LDomain_TargetError_ErrorConvergence_Structured}(d) it is clear that the outputs of the fully adaptive procedure for various error targets exhibits a linear convergence rate. Since the procedure aims to generate a quasi-optimal mesh, it is expected that this is the (approximately) optimal convergence rate for this problem.
	
	\FloatBarrier
	\FloatBarrier
	
\begin{figure}[ht!]
		\centering
		\begin{subfigure}[t]{0.495\textwidth}
			\centering
			\includegraphics[width=0.95\textwidth]{{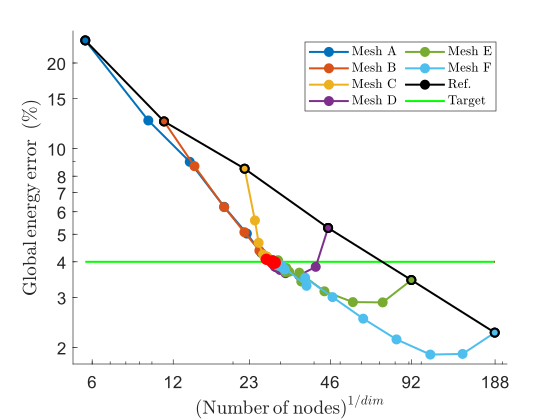}}
			\caption{$\|e\|_{\text{targ}}^{\text{rel}} = 4\%$}
		\end{subfigure}%
		\begin{subfigure}[t]{0.495\textwidth}
			\centering
			\includegraphics[width=0.95\textwidth]{{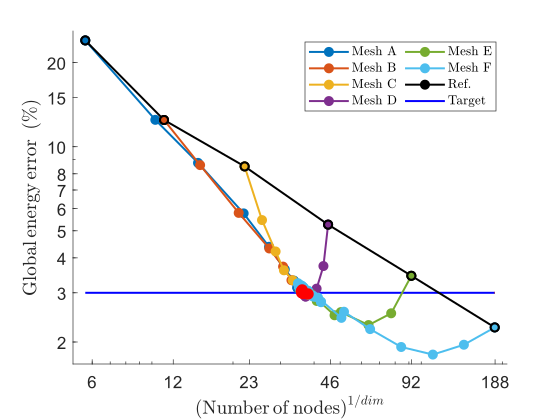}}
			\caption{$\|e\|_{\text{targ}}^{\text{rel}} = 3\%$}
		\end{subfigure}
		\vskip \baselineskip 
		\begin{subfigure}[t]{0.495\textwidth}
			\centering
			\includegraphics[width=0.95\textwidth]{{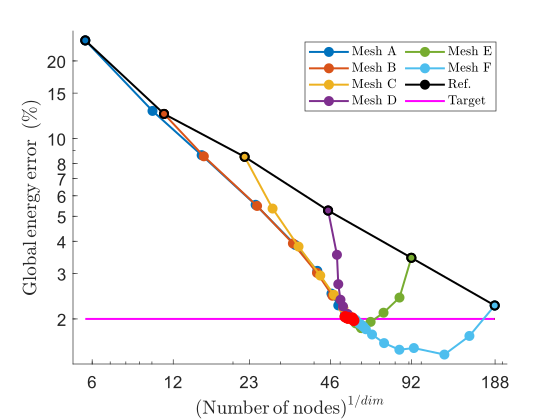}}
			\caption{$\|e\|_{\text{targ}}^{\text{rel}} = 2\%$}
		\end{subfigure}%
		\begin{subfigure}[t]{0.495\textwidth}
			\centering
			\includegraphics[width=0.95\textwidth]{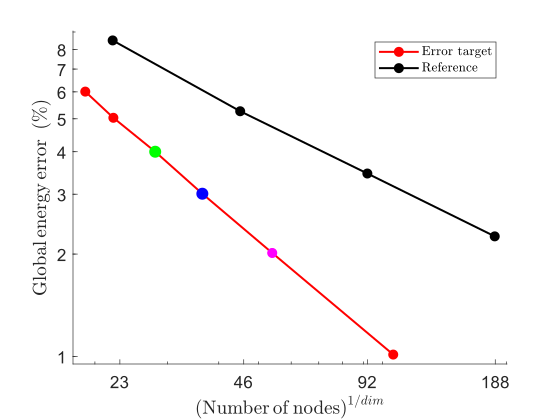}
			\caption{Various error targets}
		\end{subfigure}
		\caption{Energy error vs number of nodes for the L-shaped domain problem on Voronoi meshes of various initial densities with (a) ${\|e\|_{\text{rel}}^{\text{targ}} = 4\%}$, (b) ${\|e\|_{\text{rel}}^{\text{targ}} = 3\%}$, (c) ${\|e\|_{\text{rel}}^{\text{targ}} = 2\%}$, and (d) average final/adapted mesh error for various error targets. 
			\label{fig:LDomain_TargetError_ErrorConvergence_Voronoi}}
	\end{figure}

	The convergence behaviour of the energy error approximation vs the number of nodes in the mesh is depicted on a logarithmic scale in Figures~\ref{fig:LDomain_TargetError_ErrorConvergence_Voronoi}(a)-(c) for the L-shaped domain problem on Voronoi meshes. The convergence behaviour is plotted for cases of several initially uniform Voronoi meshes of varying density for various error targets. 
	Additionally, the averaged final adapted mesh result (red markers) for all error targets considered is plotted in Figure~\ref{fig:LDomain_TargetError_ErrorConvergence_Voronoi}(d).
	Where applicable, the markers in Figure~\ref{fig:LDomain_TargetError_ErrorConvergence_Voronoi}(d) are colored to match their corresponding targets depicted in Figures~\ref{fig:LDomain_TargetError_ErrorConvergence_Voronoi}(a)-(c).
	The behaviours exhibited in Figure~\ref{fig:LDomain_TargetError_ErrorConvergence_Voronoi} are almost identical to those observed in Figure~\ref{fig:LDomain_TargetError_ErrorConvergence_Structured} for structured meshes. The only discernible difference between the two sets of figures is that in the case of Voronoi meshes there are very small differences in the number of nodes of the final adapted mesh results (as indicated by the positions of the red markers). These differences are a result of the inherent randomness in Voronoi meshes and the randomness involved in the refinement of Voronoi elements and are not a pathology of the adaptive procedure.
	Thus, from Figures~\ref{fig:LDomain_TargetError_ErrorConvergence_Structured} and \ref{fig:LDomain_TargetError_ErrorConvergence_Voronoi} it is clear that the fully adaptive procedure can meet any specified error target on both structured and Voronoi meshes.
	
	\FloatBarrier
	\FloatBarrier
	
\begin{figure}[ht!]
		\centering
		\begin{subfigure}[t]{0.495\textwidth}
			\centering
			\includegraphics[width=0.95\textwidth]{{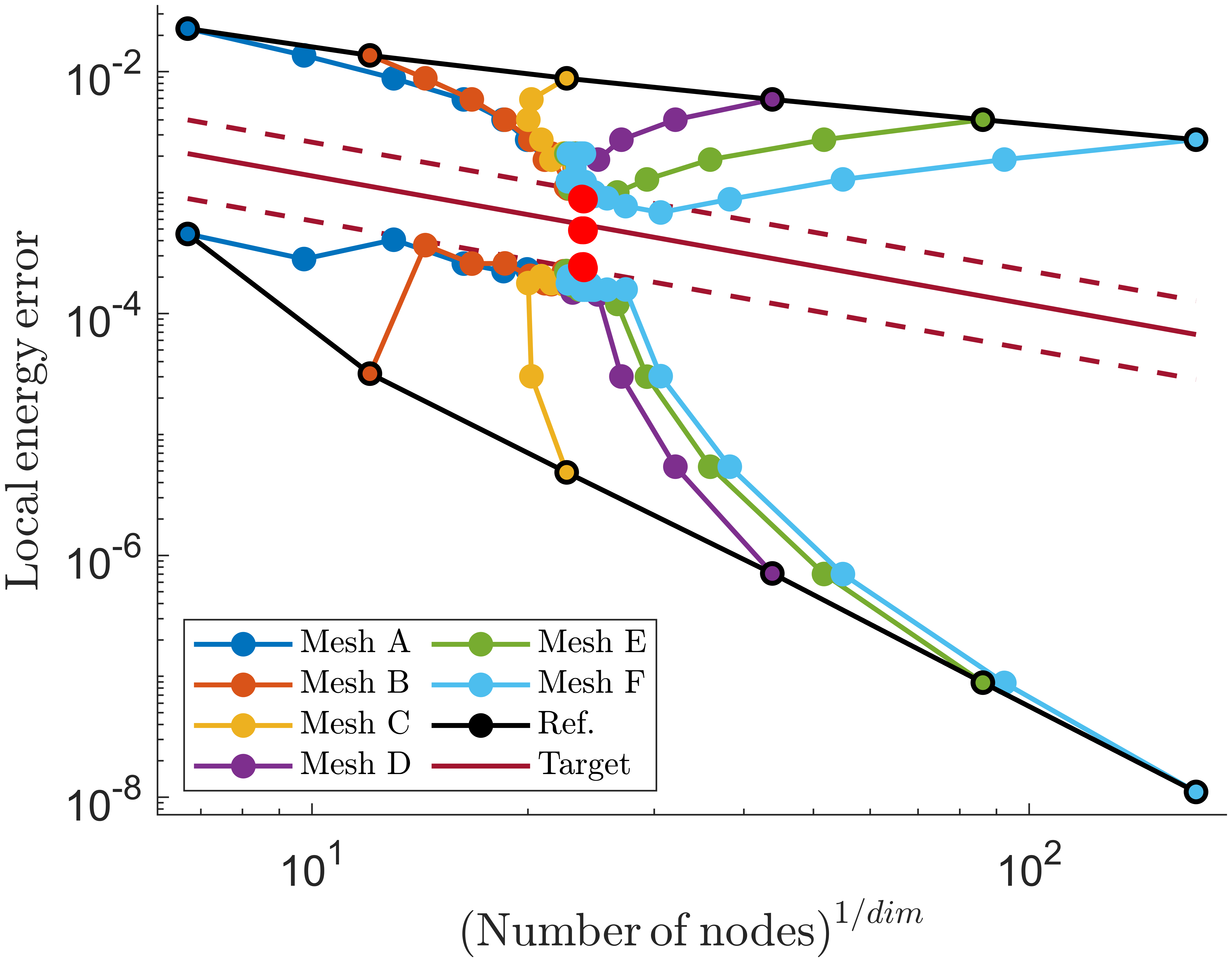}}
			\caption{$\|e\|_{\text{targ}}^{\text{rel}} = 4\%$ - Max and min local error}
			\vspace*{-3mm}
		\end{subfigure}%
		\begin{subfigure}[t]{0.495\textwidth}
			\centering
			\includegraphics[width=0.95\textwidth]{{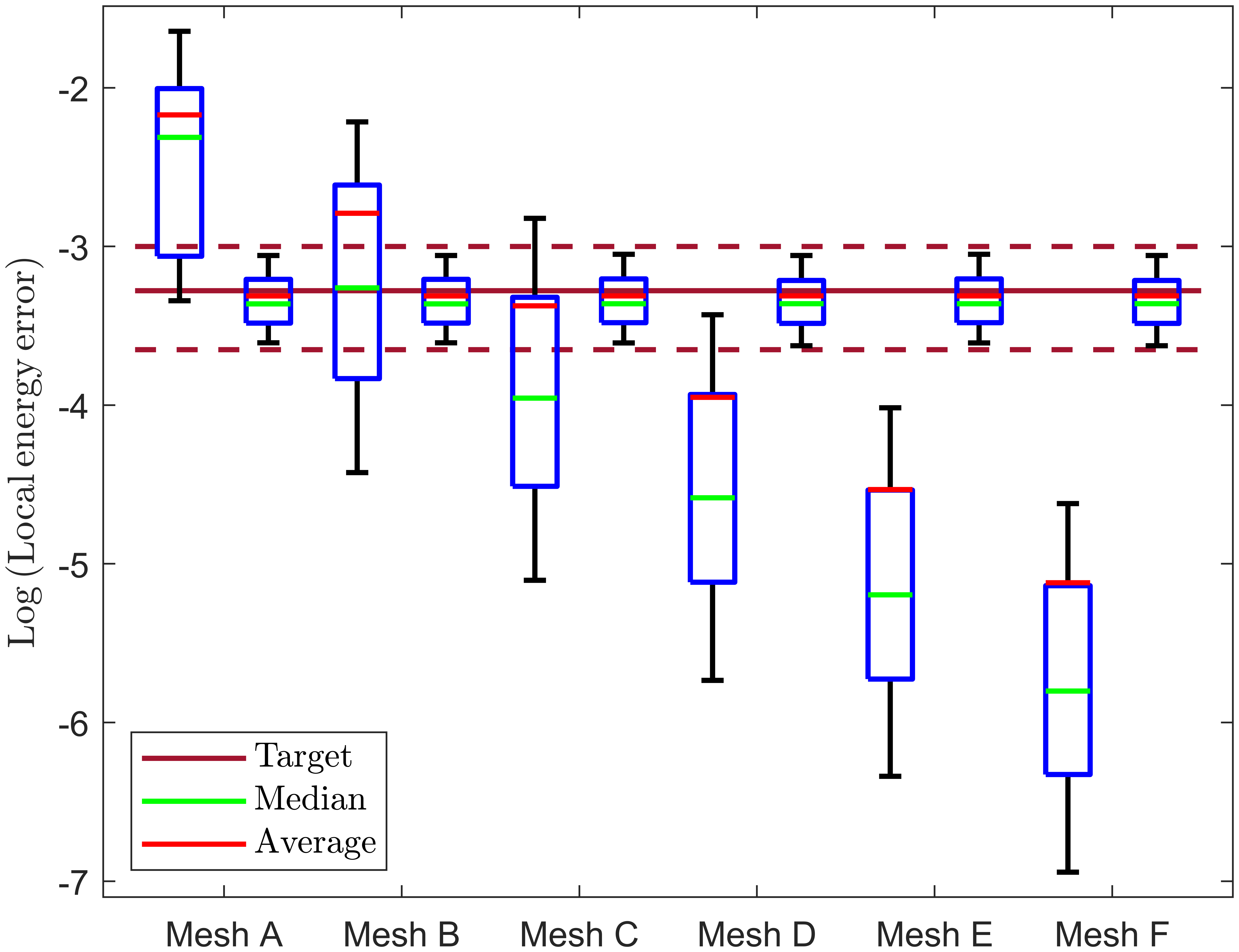}}
			\caption{$\|e\|_{\text{targ}}^{\text{rel}} = 4\%$ - Local error distribution}
			\vspace*{-3mm}
		\end{subfigure}
		\vskip \baselineskip 
		\begin{subfigure}[t]{0.495\textwidth}
			\centering
			\includegraphics[width=0.95\textwidth]{{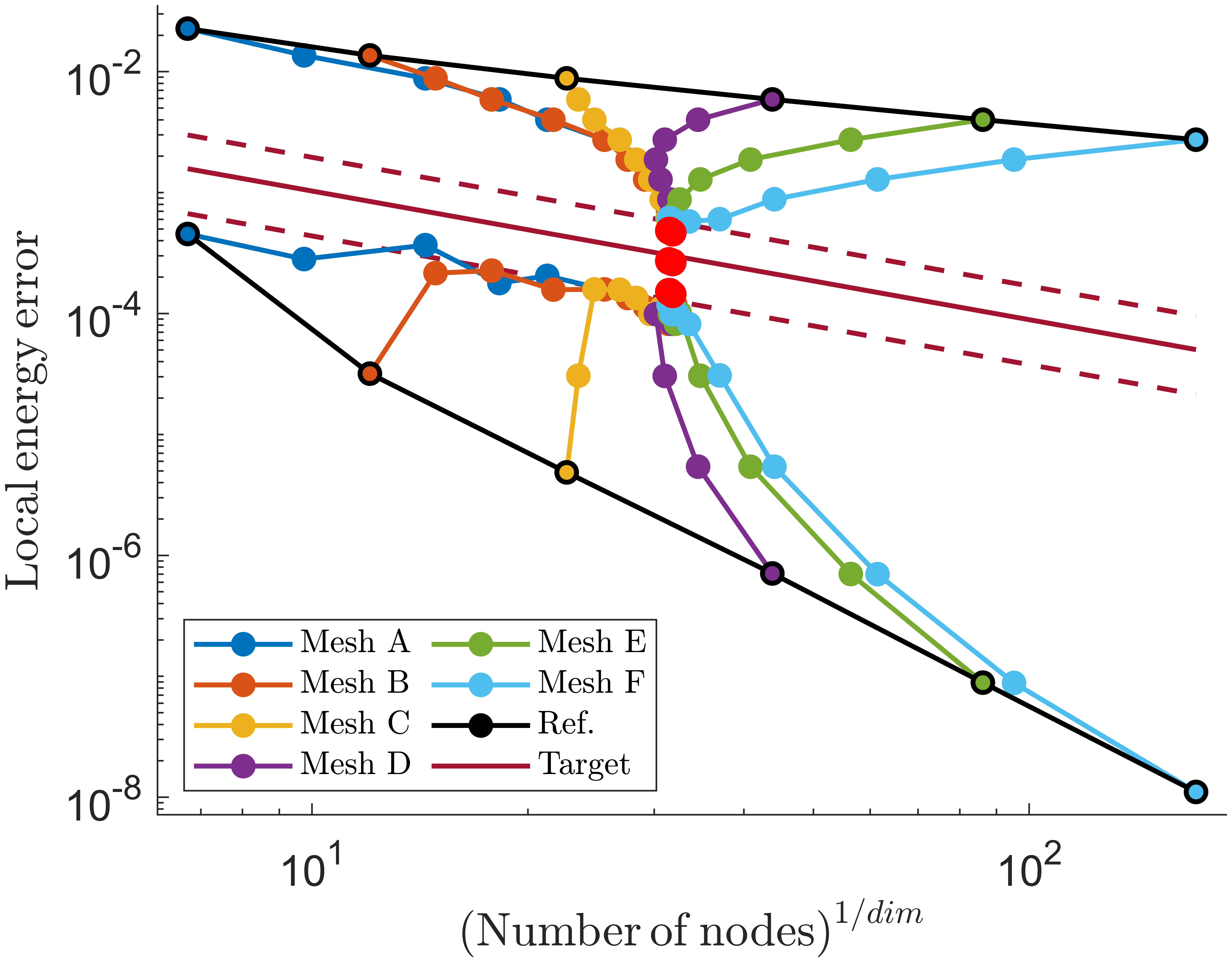}}
			\caption{$\|e\|_{\text{targ}}^{\text{rel}} = 3\%$ - Max and min local error}
			\vspace*{-3mm}
		\end{subfigure}%
		\begin{subfigure}[t]{0.495\textwidth}
			\centering
			\includegraphics[width=0.95\textwidth]{{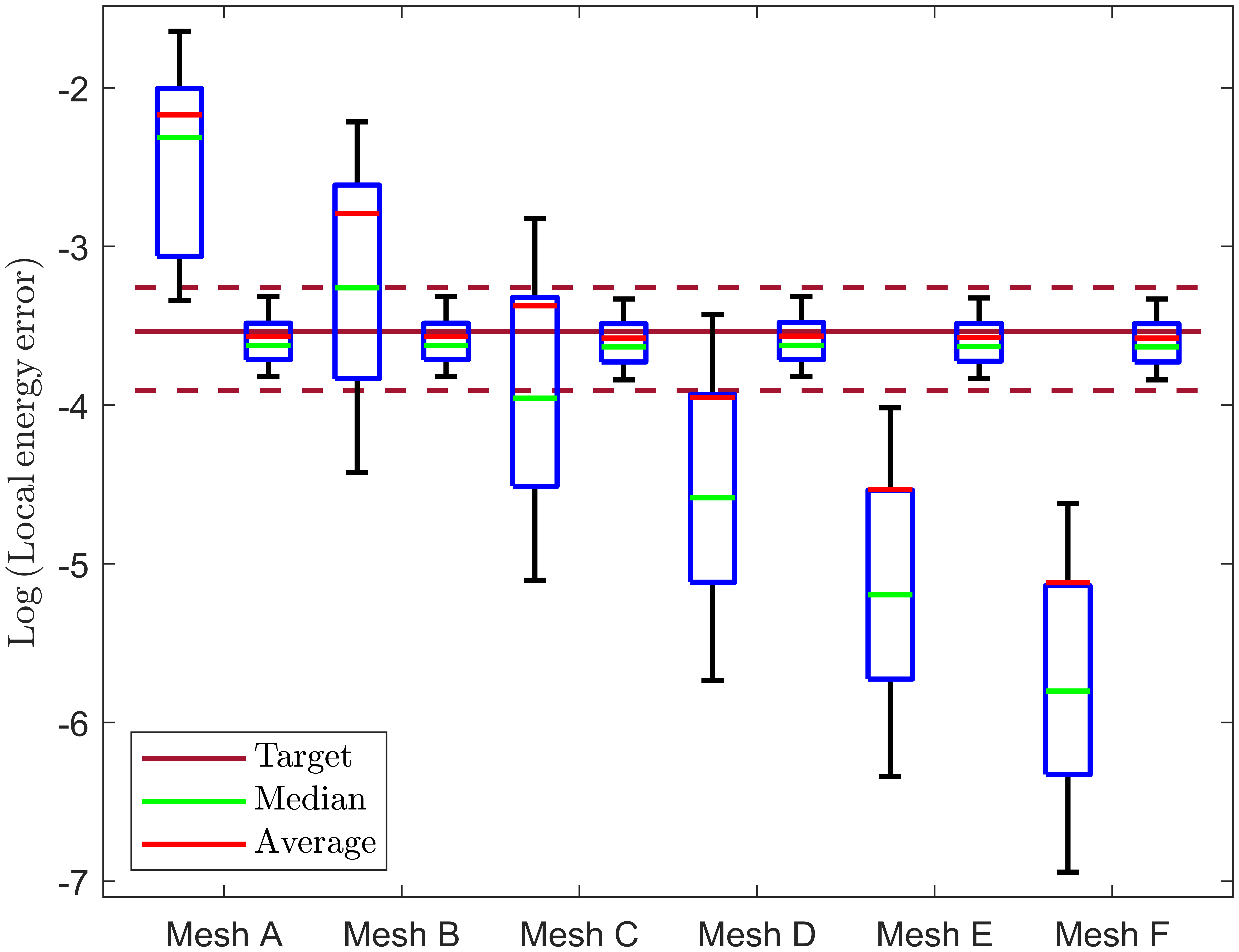}}
			\caption{$\|e\|_{\text{targ}}^{\text{rel}} = 3\%$ - Local error distribution}
			\vspace*{-3mm}
		\end{subfigure}
		\vskip \baselineskip 
		\begin{subfigure}[t]{0.495\textwidth}
			\centering
			\includegraphics[width=0.95\textwidth]{{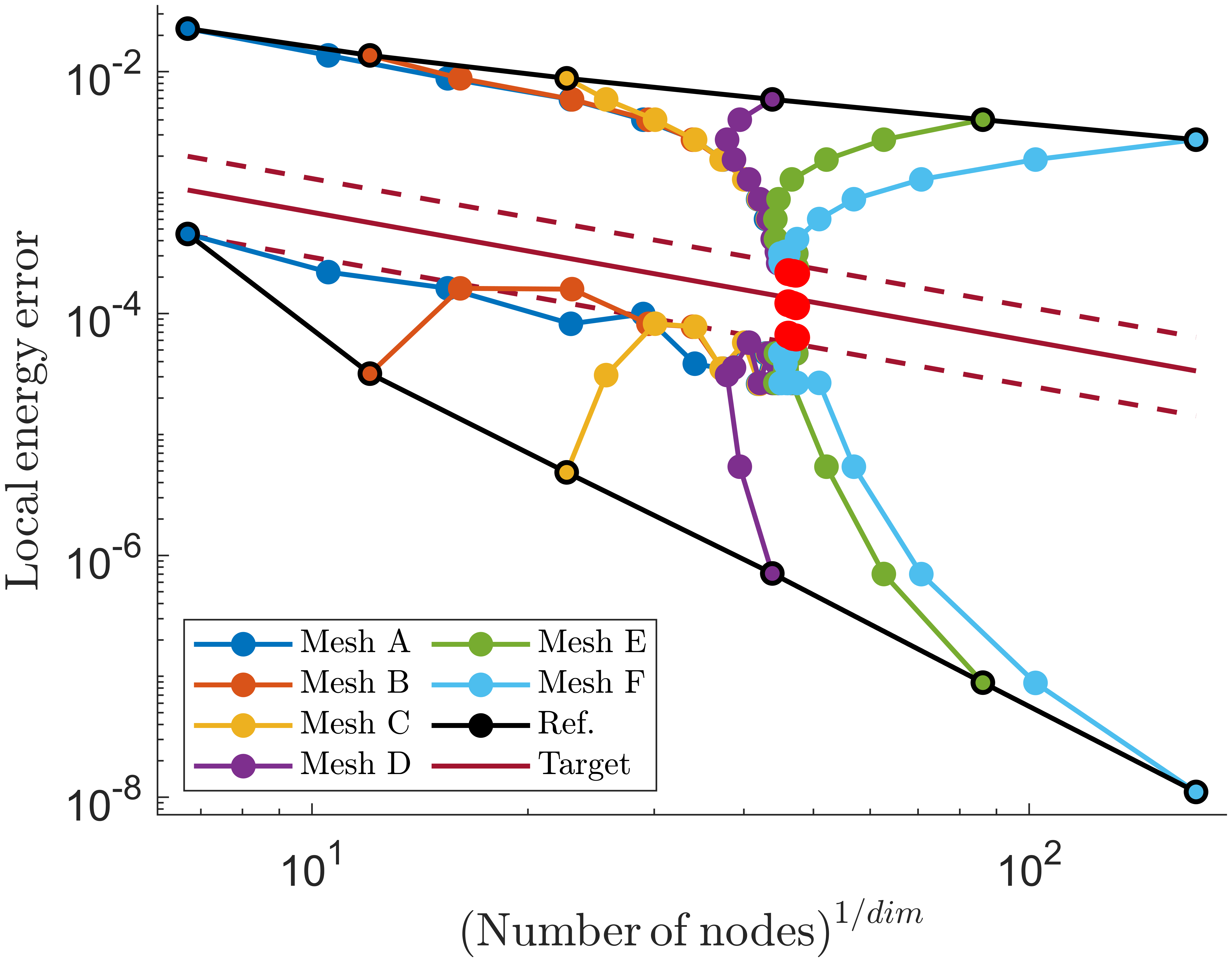}}
			\caption{$\|e\|_{\text{targ}}^{\text{rel}} = 2\%$ - Max and min local error}
			\vspace*{-3mm}
		\end{subfigure}%
		\begin{subfigure}[t]{0.495\textwidth}
			\centering
			\includegraphics[width=0.95\textwidth]{{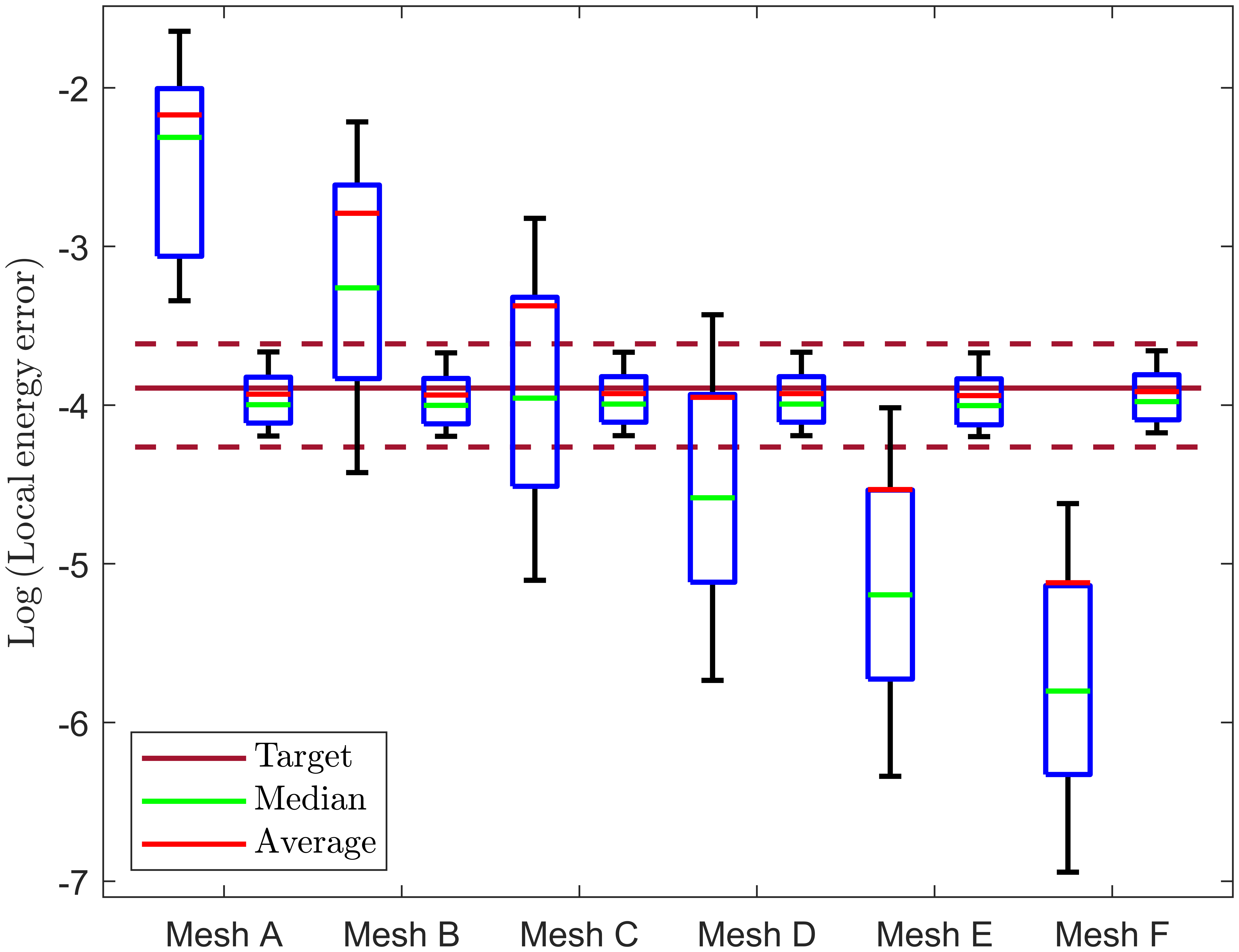}}
			\caption{$\|e\|_{\text{targ}}^{\text{rel}} = 2\%$ - Local error distribution}
			\vspace*{-3mm}
		\end{subfigure}
		\caption{Max and min energy error (left) and box and whisker plot of energy error distribution (right) for the L-shaped domain problem on structured meshes of various initial densities for a range of error targets.
			\label{fig:LDomain_TargetError_ErrorDistribution_Structured}}
	\end{figure}

	The error distribution during the mesh adaptation process for the L-shaped domain problem is depicted in Figure~\ref{fig:LDomain_TargetError_ErrorDistribution_Structured} for various error targets on structured meshes.
	
	The left column of figures depicts the evolution of the true/absolute maximum (top curves) and minimum (bottom curves) local element errors from initial meshes of varying density. Additionally, the optimal target element error is indicated by a solid maroon line as a function of mesh density. The dashed maroon lines indicate the target upper and lower element error bounds as described in Section~\ref{subsec:ElementSelectionTargetError}. Finally, three sets of red markers indicate the final adapted mesh result for each of the initial meshes. The central markers represent the average element error over all the elements in the mesh. The upper and lower markers respectively represent the 5\% trimmed maximum and minimum element errors. A slightly trimmed maximum and minimum are used to improve the readability of the graph while still accurately representing the underlying data.
	For all error targets considered the maximum and minimum errors respectively converge to the prescribed upper and lower error bounds. Furthermore, the average element-level error (as indicated by the central red markers) closely meets the element-level target for all considered global error targets.
	
	The right column of figures illustrates the nature of the distribution of the local element-level error over all of the elements in a mesh through a classical box and whisker plot. Here, the median and quartiles are computed in the standard way (i.e., from the full set of element-level data) while the maximum and minimum whiskers correspond to the 5\% trimmed data (i.e., equivalent to the red markers). Additionally, the average error of all elements in the mesh is computed and indicated on the figure. While the average is not typically considered in a box and whisker plot it is helpful in understanding the spread of the data.
	In these figures pairs of results correspond to the error evolution for a particular mesh. For example, the first (left-most) data corresponds to the initial uniform mesh error distribution for `Mesh A' and the second data corresponds to the error distribution of the corresponding final adapted mesh. 
	For all error targets considered the final adapted mesh error distributions fit within the specified upper and lower error bounds. Additionally, the upper and lower quartiles indicate that the majority of the element-level errors are very close to the element-level targets. 
	For the cases of the initial meshes the average and median error differ significantly which is a classical indicator of inequality within the dataset. Conversely, in the cases of the final adapted meshes the average and median are almost identical, thus, indicating the equality of the data and further demonstrating the narrow distribution of the element-level error. 
	
	The results presented in Figure~\ref{fig:LDomain_TargetError_ErrorConvergence_Structured} demonstrated that the fully adaptive procedure was able to meet all specified error targets on structured meshes. This is indicated by the red markers denoting the final adapted mesh lying exactly on the target error line. 
	The results presented in Figure~\ref{fig:LDomain_TargetError_ErrorDistribution_Structured} for structured meshes demonstrated that the average element-level error almost exactly met the element-level target as the red markers denoting the average error strongly overlap the solid maroon target line. Furthermore, the element-level errors were satisfactorily equal as they all fell within the specified target error range.
	Thus, the fully adaptive procedure successfully generated quasi-optimal meshes for the specified target errors on structured meshes.
	
	\FloatBarrier
	\FloatBarrier
	
\begin{figure}[ht!]
		\centering
		\begin{subfigure}[t]{0.495\textwidth}
			\centering
			\includegraphics[width=0.95\textwidth]{{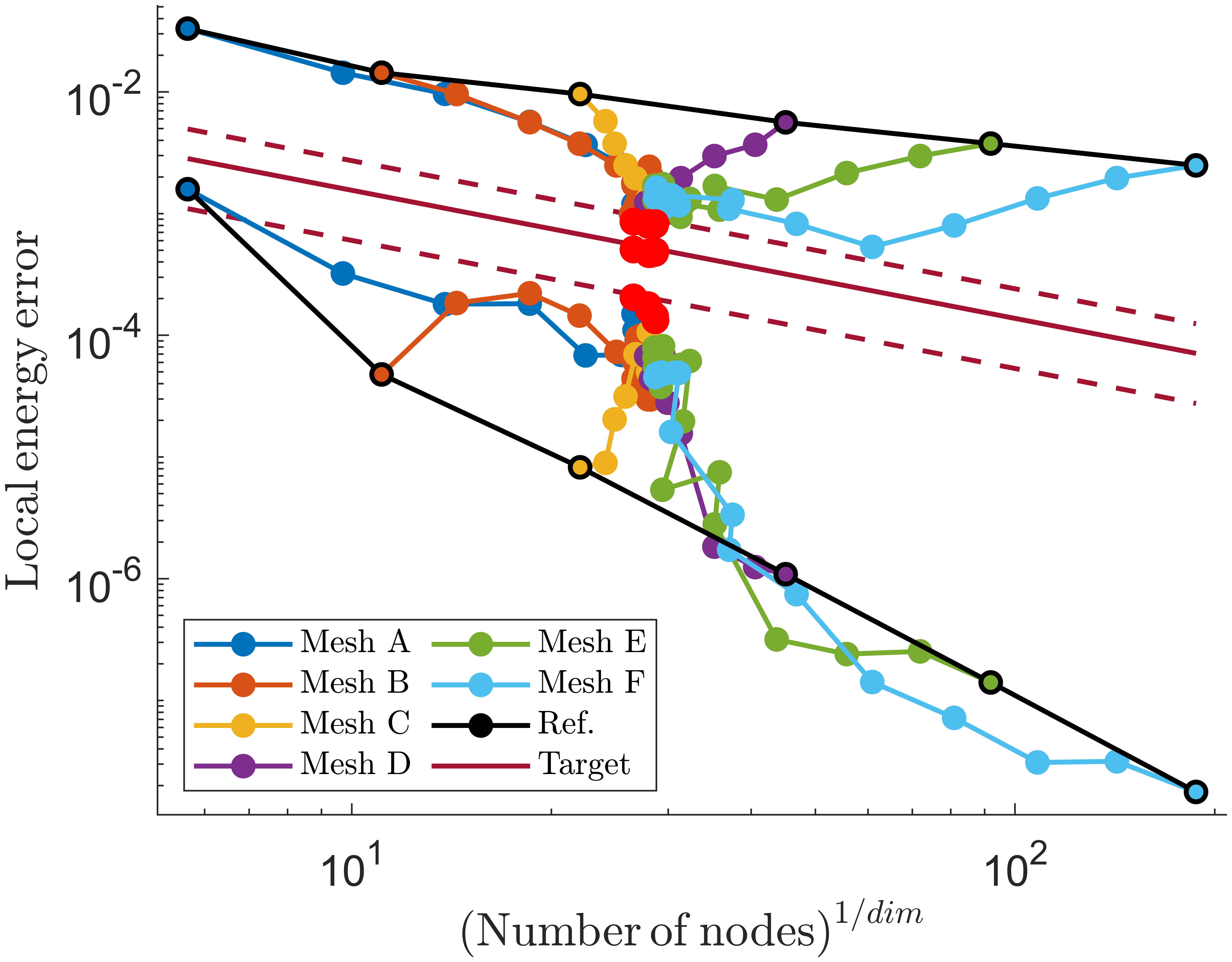}}
			\caption{$\|e\|_{\text{targ}}^{\text{rel}} = 4\%$ - Max and min local error}
			\vspace*{-3mm}
		\end{subfigure}%
		\begin{subfigure}[t]{0.495\textwidth}
			\centering
			\includegraphics[width=0.95\textwidth]{{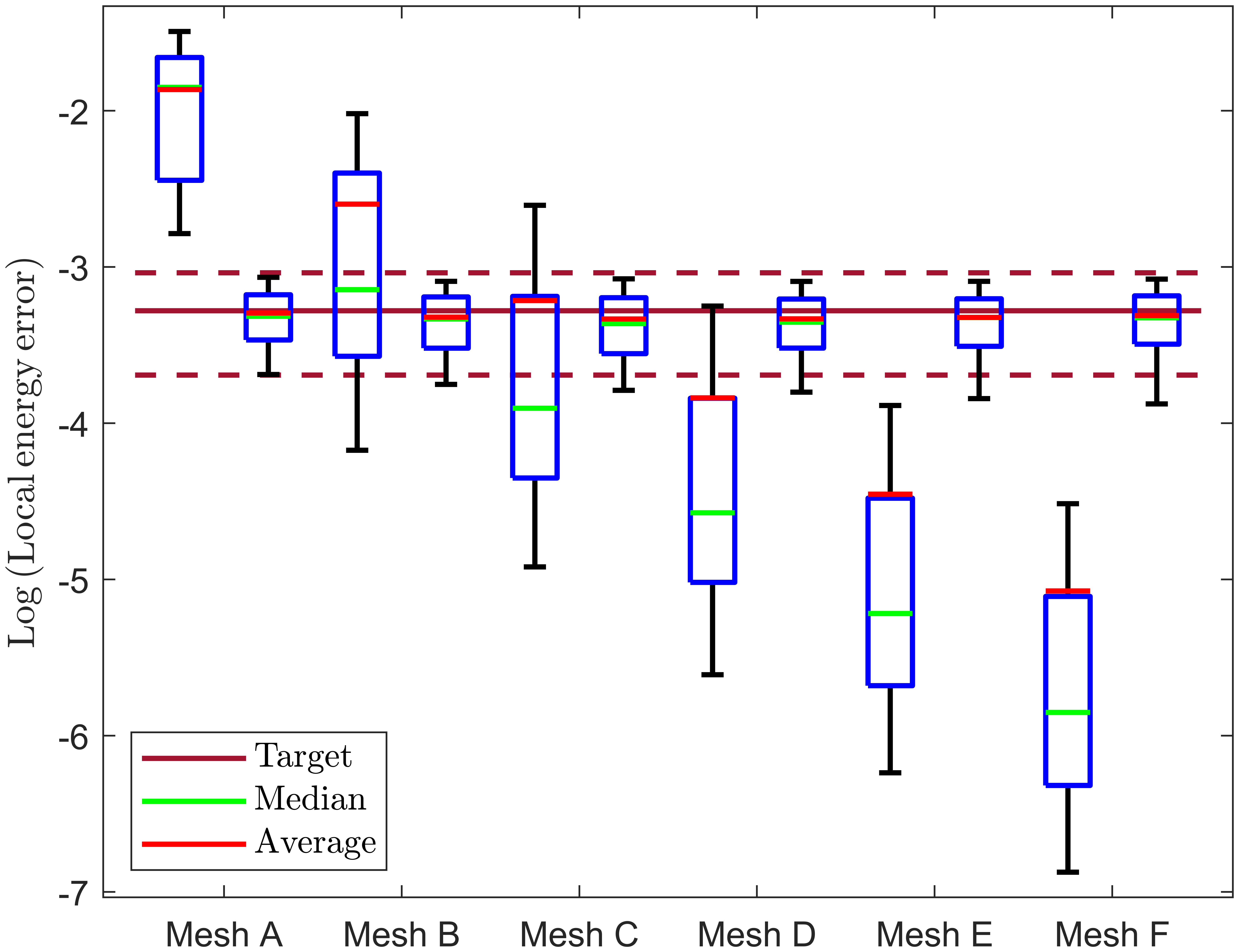}}
			\caption{$\|e\|_{\text{targ}}^{\text{rel}} = 4\%$ - Local error distribution}
			\vspace*{-3mm}
		\end{subfigure}
		\vskip \baselineskip 
		\begin{subfigure}[t]{0.495\textwidth}
			\centering
			\includegraphics[width=0.95\textwidth]{{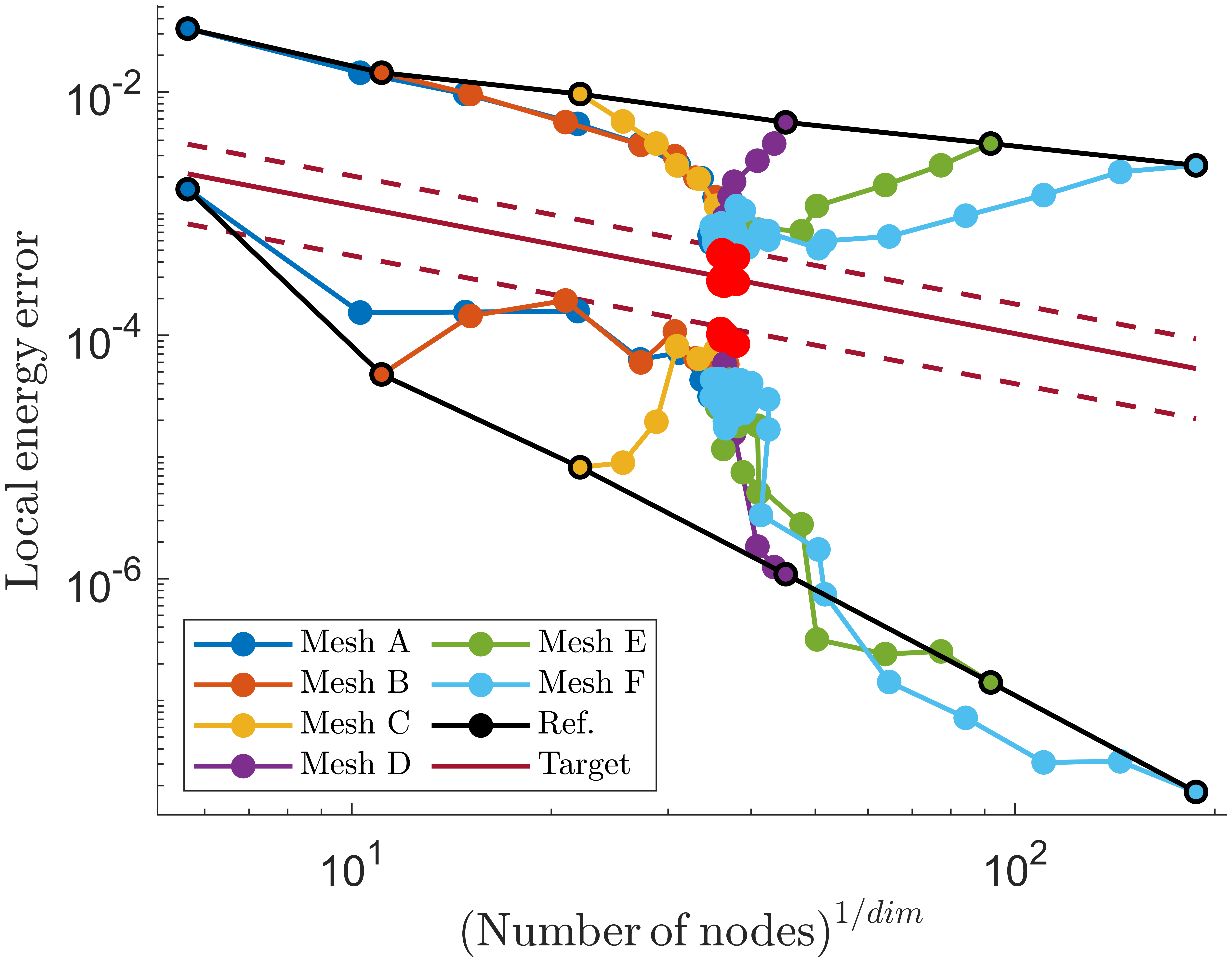}}
			\caption{$\|e\|_{\text{targ}}^{\text{rel}} = 3\%$ - Max and min local error}
			\vspace*{-3mm}
		\end{subfigure}%
		\begin{subfigure}[t]{0.495\textwidth}
			\centering
			\includegraphics[width=0.95\textwidth]{{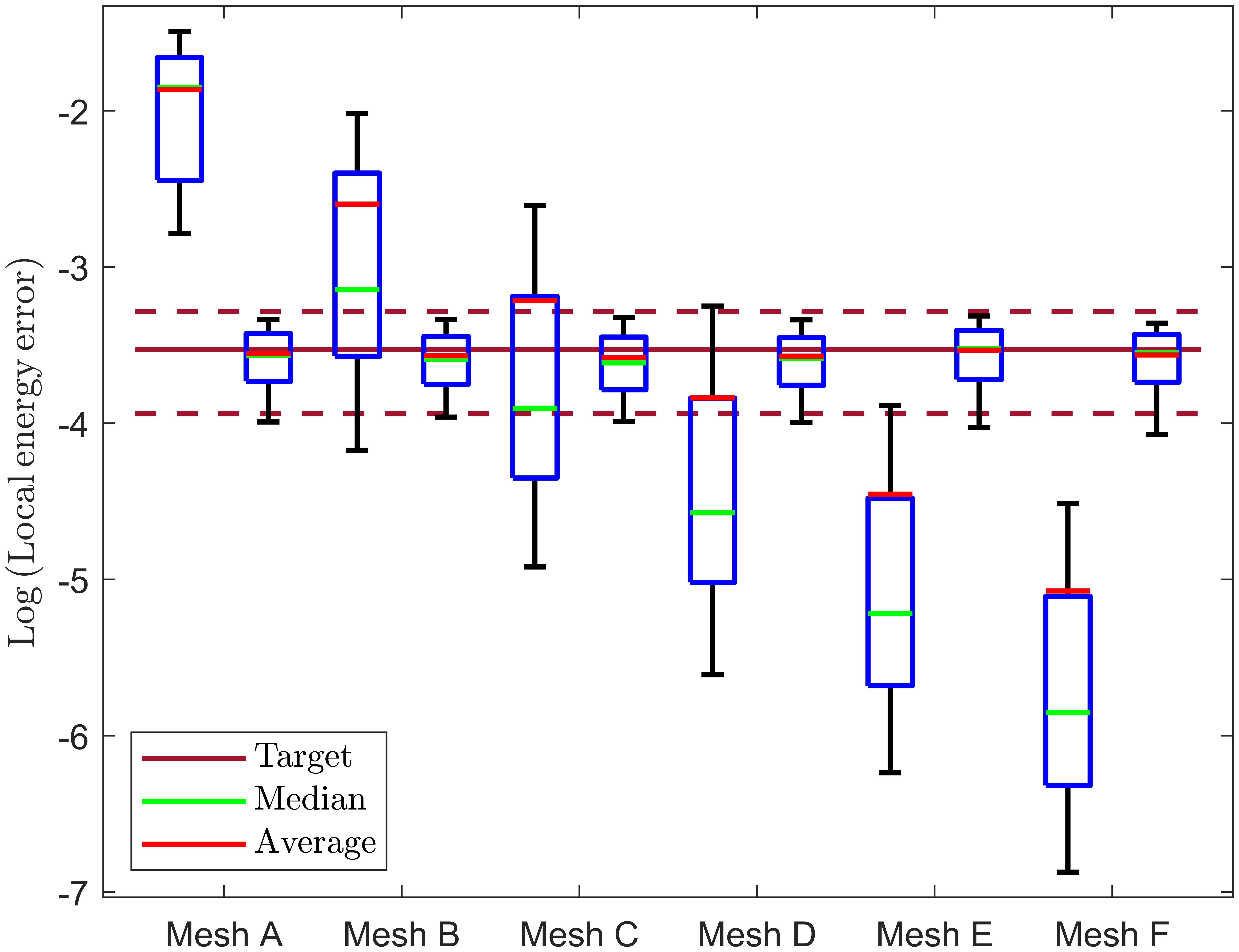}}
			\caption{$\|e\|_{\text{targ}}^{\text{rel}} = 3\%$ - Local error distribution}
			\vspace*{-3mm}
		\end{subfigure}
		\vskip \baselineskip 
		\begin{subfigure}[t]{0.495\textwidth}
			\centering
			\includegraphics[width=0.95\textwidth]{{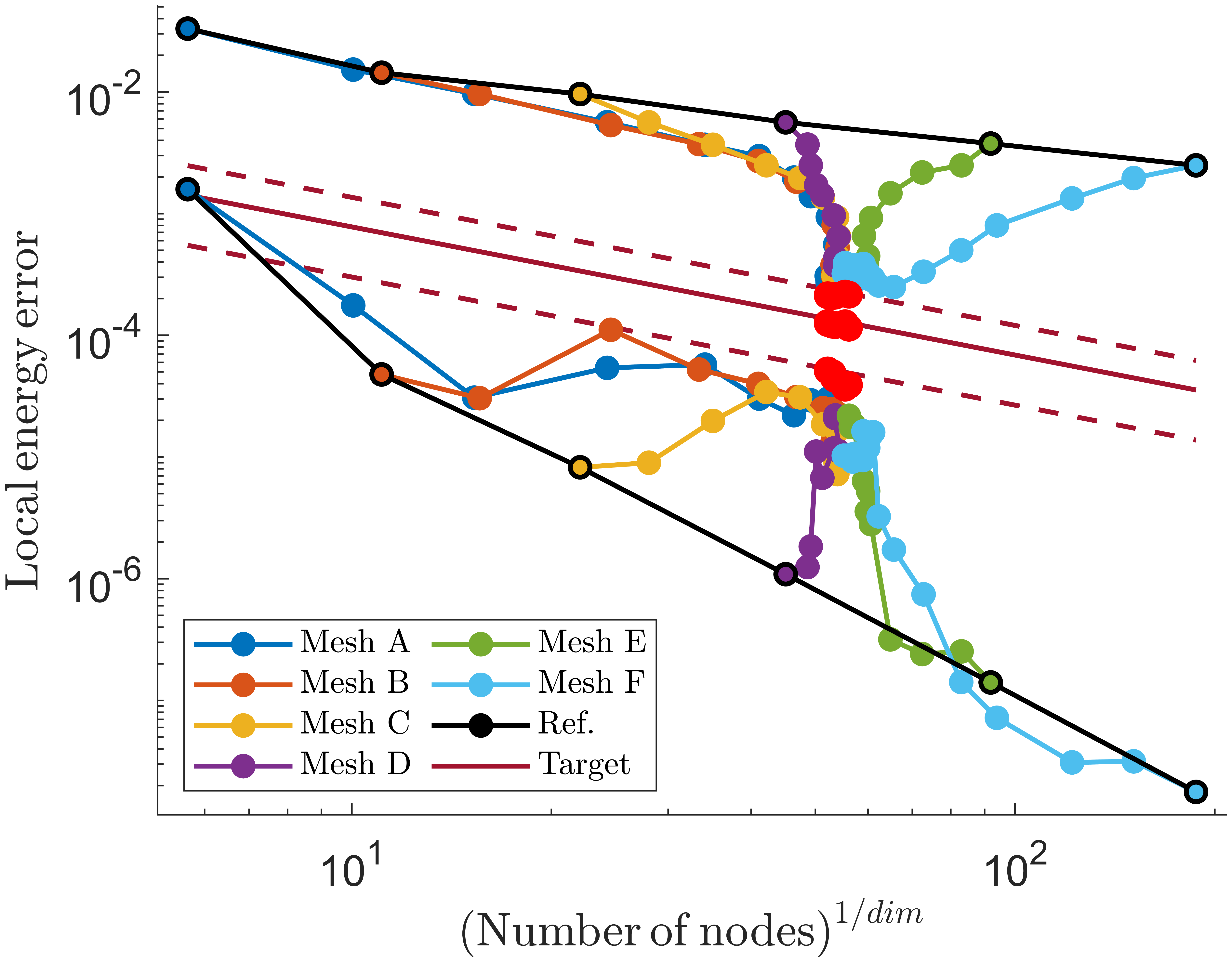}}
			\caption{$\|e\|_{\text{targ}}^{\text{rel}} = 2\%$ - Max and min local error}
			\vspace*{-3mm}
		\end{subfigure}%
		\begin{subfigure}[t]{0.495\textwidth}
			\centering
			\includegraphics[width=0.95\textwidth]{{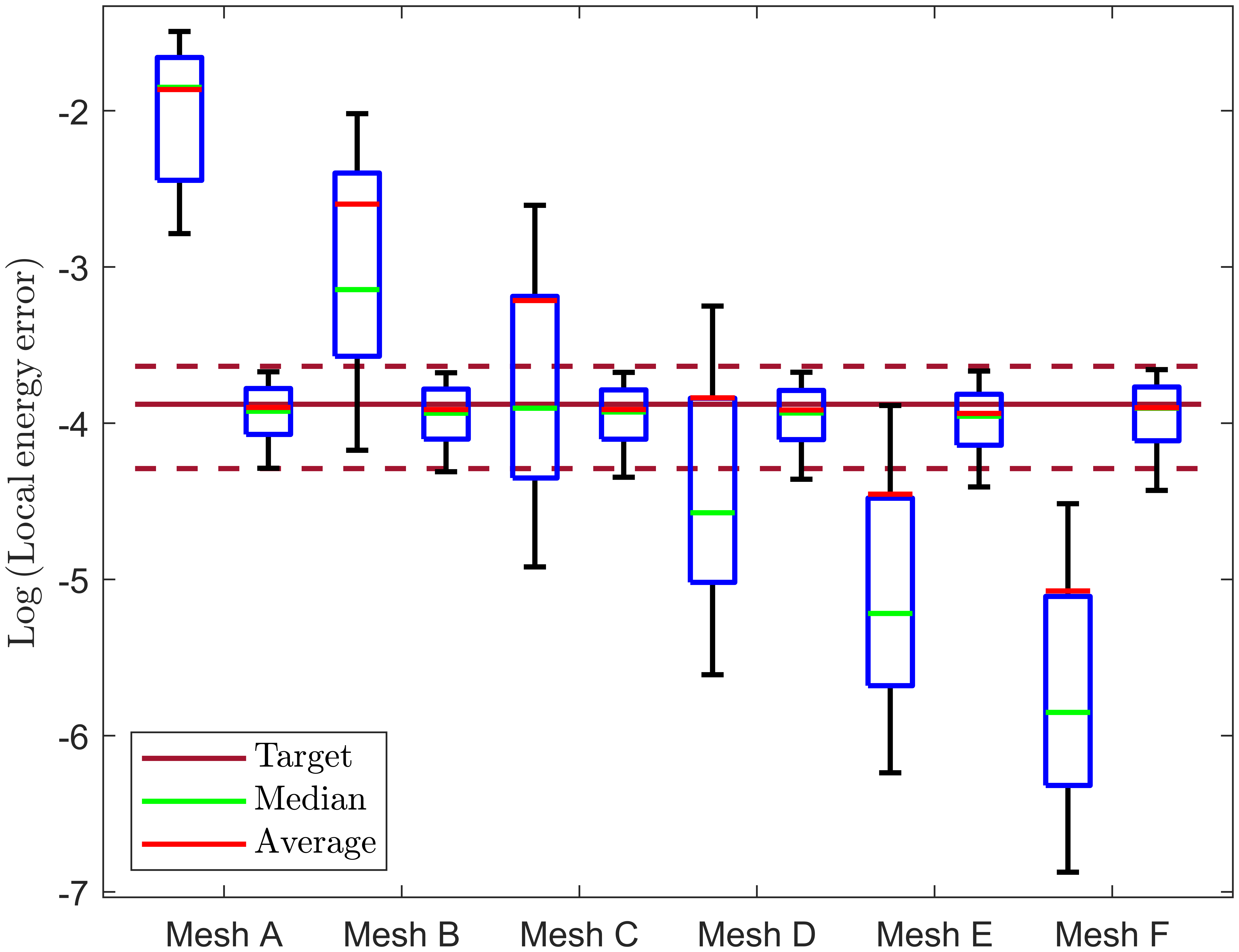}}
			\caption{$\|e\|_{\text{targ}}^{\text{rel}} = 2\%$ - Local error distribution}
			\vspace*{-3mm}
		\end{subfigure}
		\caption{Max and min energy error (left) and box and whisker plot of energy error distribution (right) for the L-shaped domain problem on Voronoi meshes of various initial densities for a range of error targets.
			\label{fig:LDomain_TargetError_ErrorDistribution_Voronoi}}
	\end{figure}

	The error distribution during the mesh adaptation process for the L-shaped domain problem is depicted in Figure~\ref{fig:LDomain_TargetError_ErrorDistribution_Voronoi} for various error targets on Voronoi meshes.
	The results presented here are very similar to those presented in Figure~\ref{fig:LDomain_TargetError_ErrorDistribution_Structured} for structured meshes. 
	The only noteworthy difference between the two sets of results is that in the case of Voronoi meshes the minimum element-level error occasionally falls slightly below the prescribed lower bound. However, the difference between the lower bound and the minimum error is small and is not indicative of a pathology in the fully adaptive procedure. The difference is most likely a result of the choice of termination criteria presented in Section~\ref{subsec:ElementSelectionTargetError}. The termination criteria is based only on the stability of the global error and of the number of elements in the mesh between remeshing iterations. Thus, the criteria does not consider if all of the element-level errors lie within the prescribed error bounds and it is possible that, if continued, all errors would fall within the bounds after several more remeshing iterations. 
	The presented termination criteria were chosen for efficiency reasons. A termination criteria based on all element-level errors falling within the error bounds was investigated but was found to be less efficient and less practical. It was found that a larger allowable error band had to be prescribed and that the number of remeshing iterations required before termination was far greater (approximately $50\%$ more) than in the case of the presented method, particularly in the case of Voronoi meshes. As such, the minimum element-level error occasionally falling below the lower bound was deemed an acceptable consequence of improved efficiency. Furthermore, the position of the lower quartiles in the box and whisker plots strongly indicate that the vast majority of element errors are above the lower bound.
	
	Similarly to the case of structured meshes, the results presented in Figure~\ref{fig:LDomain_TargetError_ErrorConvergence_Voronoi} demonstrated that the fully adaptive procedure was able to meet all specified error targets on Voronoi meshes. 
	Additionally, the results presented in Figure~\ref{fig:LDomain_TargetError_ErrorDistribution_Voronoi} for Voronoi meshes demonstrated that the average element-level error almost exactly met the element-level target.
	Finally, the element-level errors were satisfactorily equal as they almost all fell within the specified target error range.
	Thus, the fully adaptive procedure successfully generated quasi-optimal meshes for the specified target errors on Voronoi meshes.
	
	\FloatBarrier
	\FloatBarrier
	
	\subsubsection{Target number of elements}
	\label{subsubsec:L_Domain_TargetNumElements}
	
\begin{figure}[ht!]
		\centering
		\begin{subfigure}[t]{0.33\textwidth}
			\centering
			\includegraphics[width=0.95\textwidth,height=0.95\textwidth]{{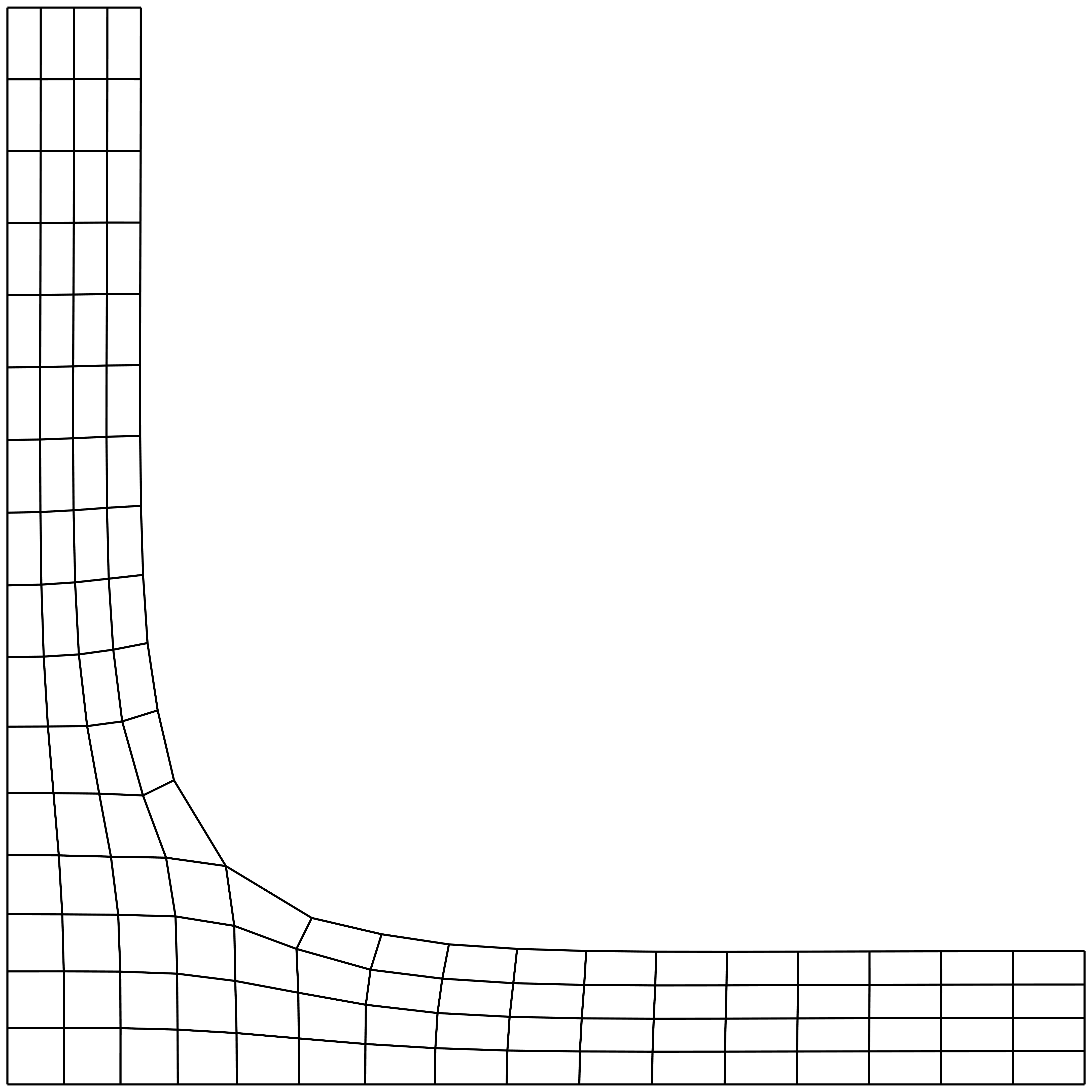}}
			\caption{Step 1: Initial uniform mesh}
		\end{subfigure}%
		\begin{subfigure}[t]{0.33\textwidth}
			\centering
			\includegraphics[width=0.95\textwidth,height=0.95\textwidth]{{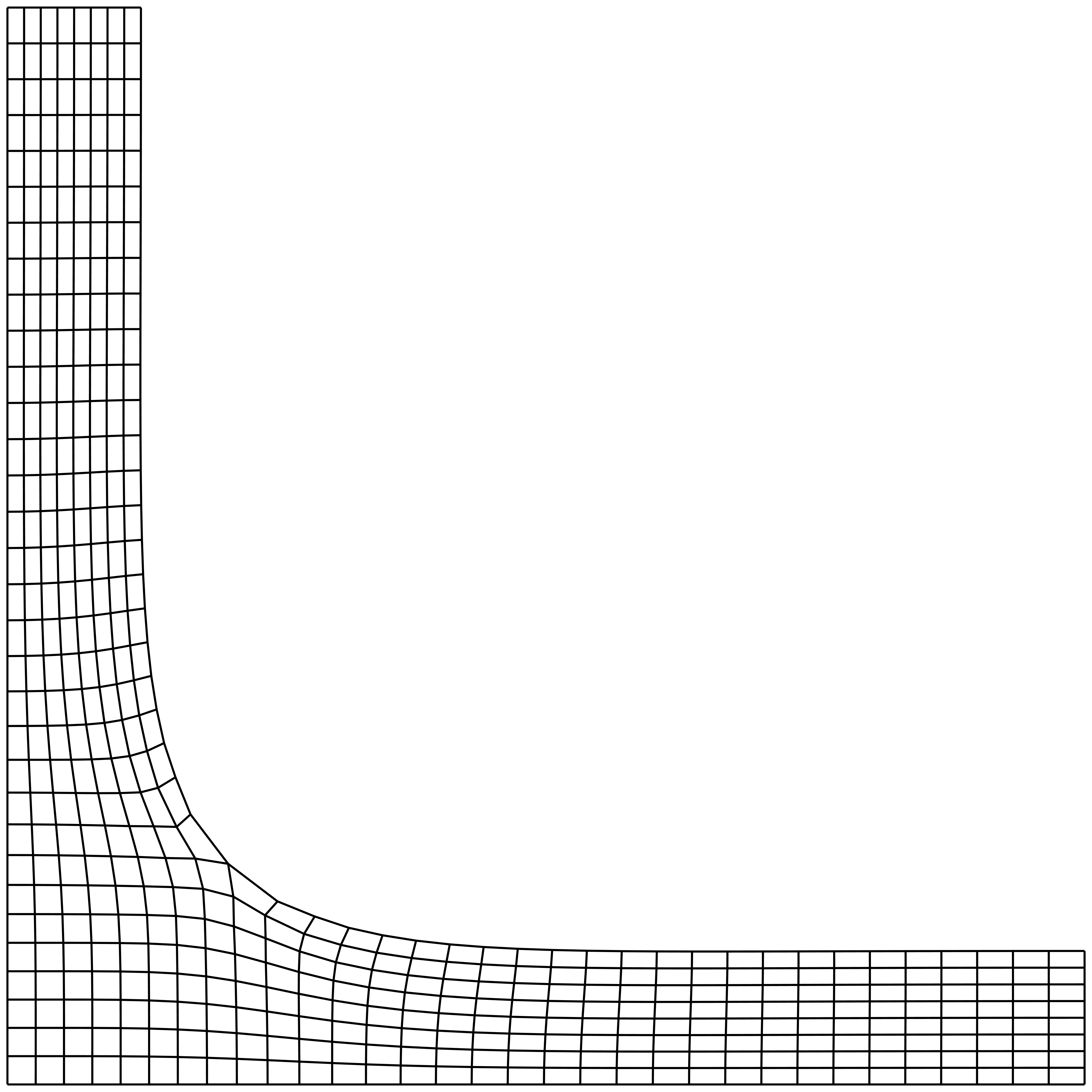}}
			\caption{Step 2}
		\end{subfigure}%
		\begin{subfigure}[t]{0.33\textwidth}
			\centering
			\includegraphics[width=0.95\textwidth,height=0.95\textwidth]{{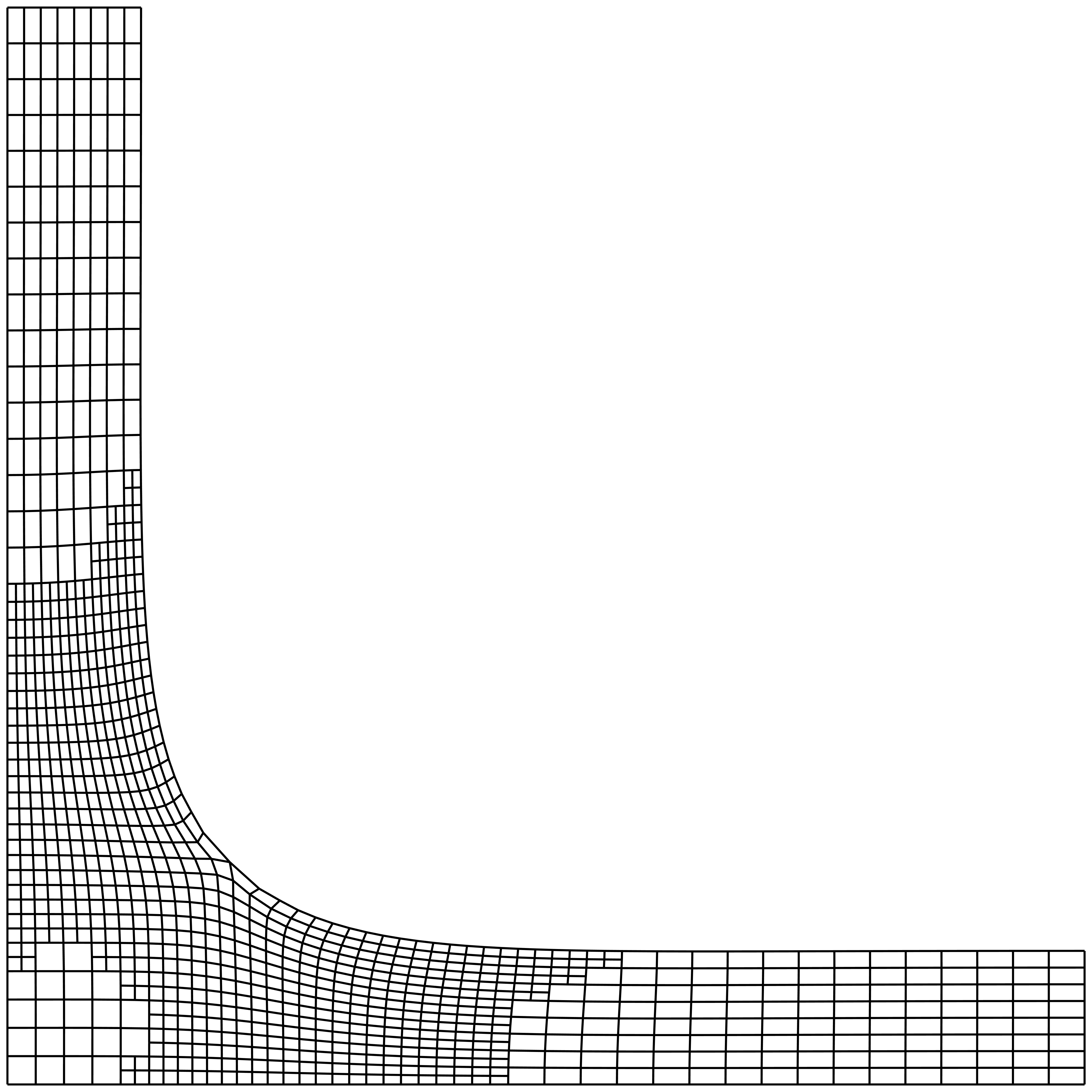}}
			\caption{Step 3}
		\end{subfigure}
		\vskip \baselineskip 
		\begin{subfigure}[t]{0.33\textwidth}
			\centering
			\includegraphics[width=0.95\textwidth,height=0.95\textwidth]{{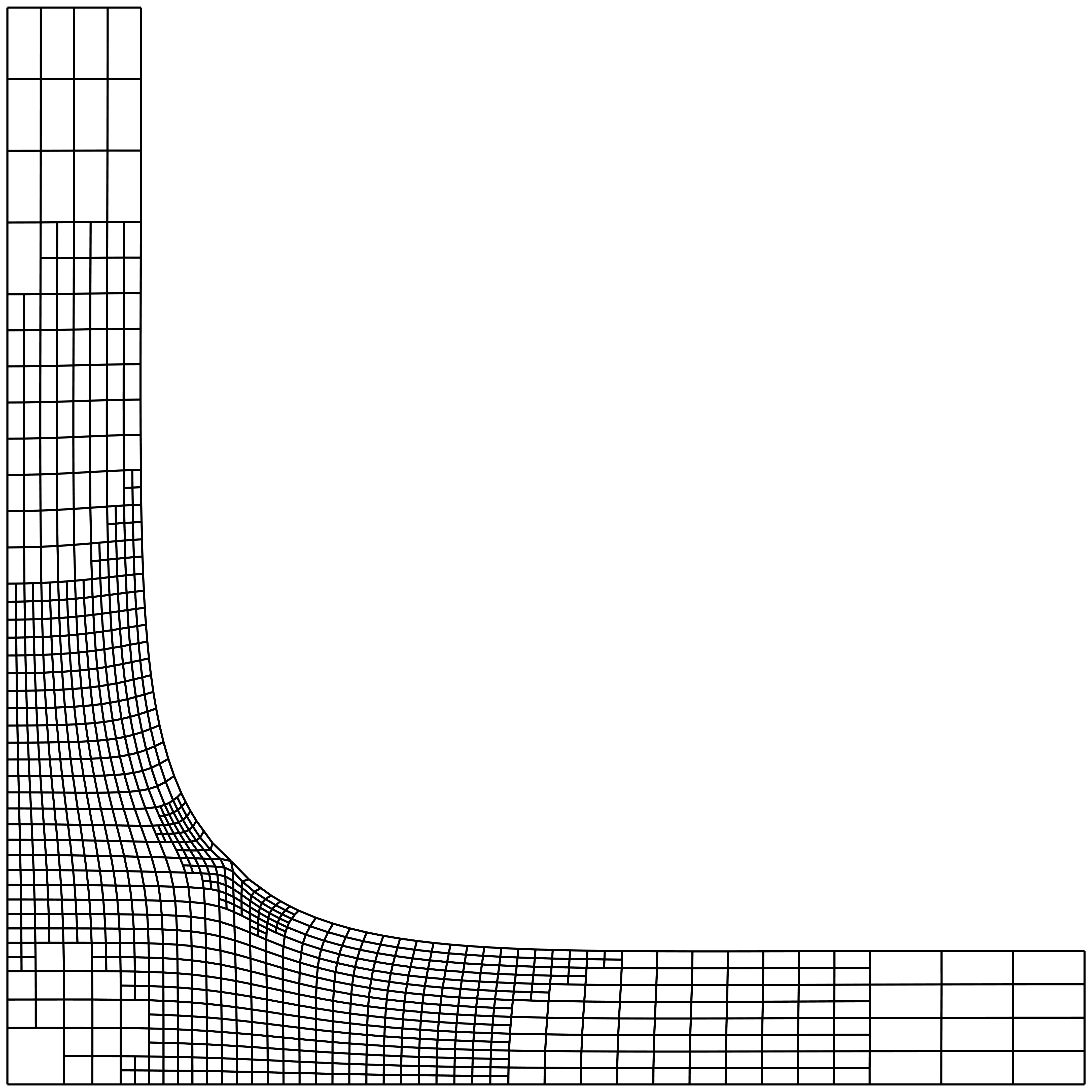}}
			\caption{Step 4}
		\end{subfigure}%
		\begin{subfigure}[t]{0.33\textwidth}
			\centering
			\includegraphics[width=0.95\textwidth,height=0.95\textwidth]{{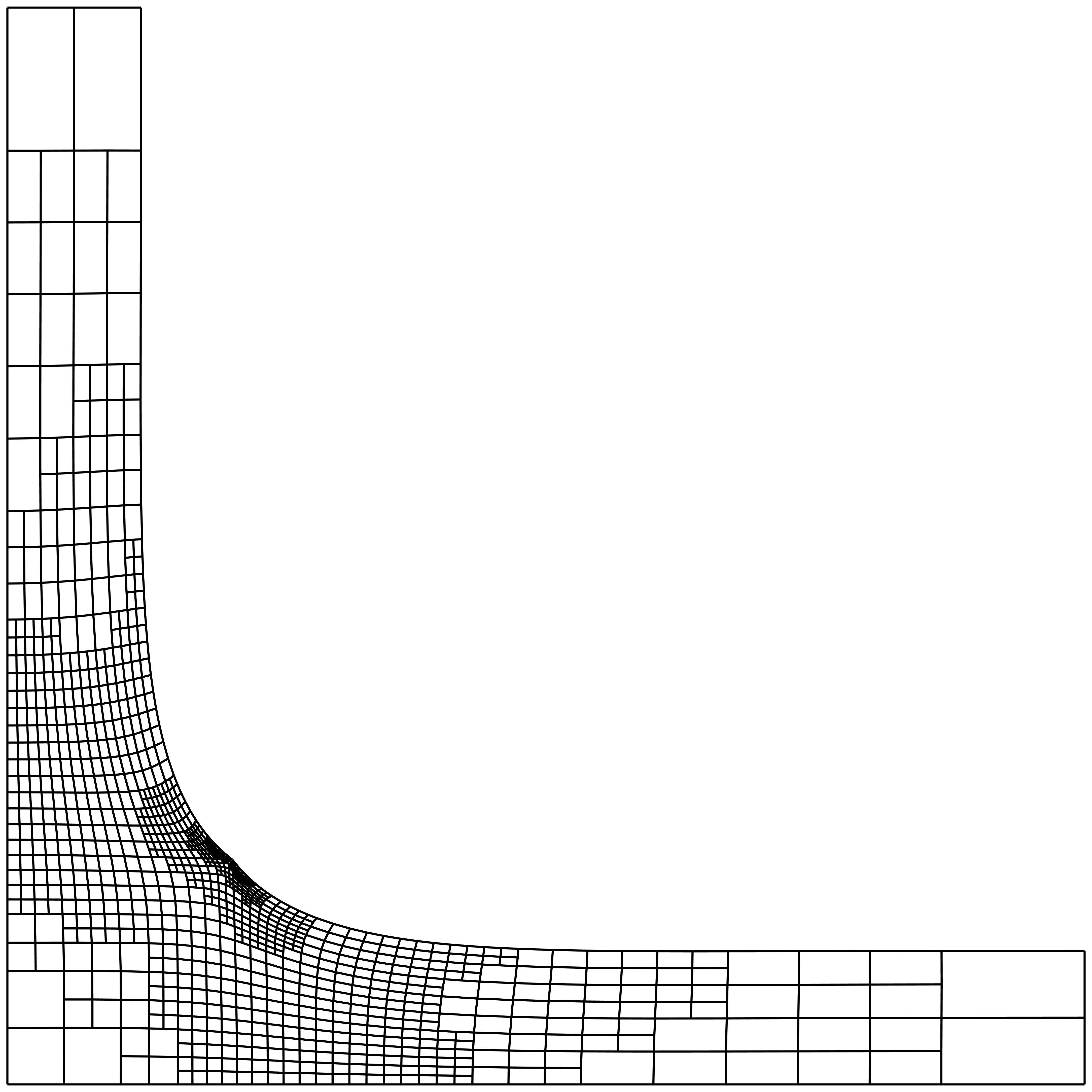}}
			\caption{Step 6}
		\end{subfigure}%
		\begin{subfigure}[t]{0.33\textwidth}
			\centering
			\includegraphics[width=0.95\textwidth,height=0.95\textwidth]{{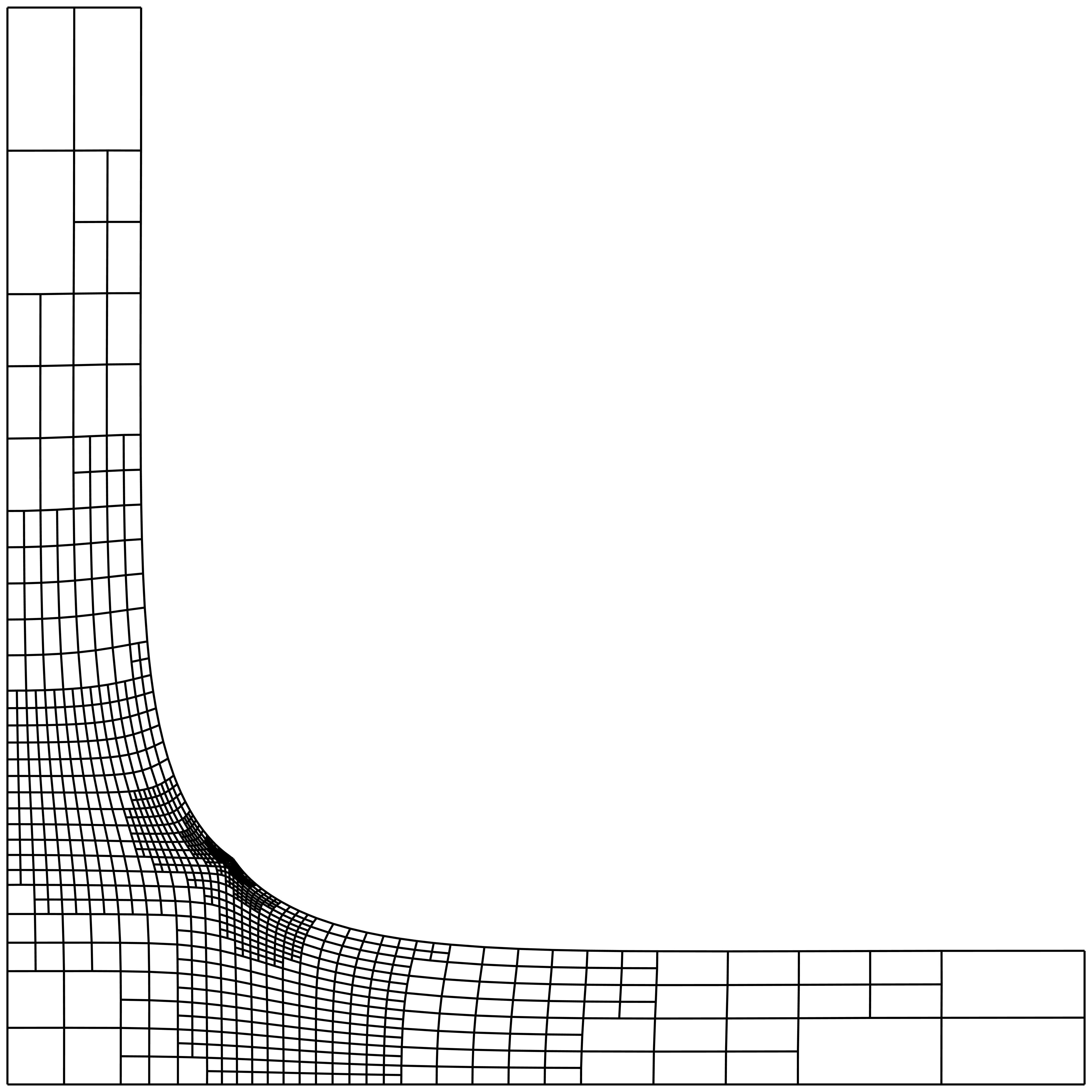}}
			\caption{Step 9: Final adapted mesh}
		\end{subfigure}
		\caption{Mesh evolution with the fully adaptive procedure for the L-shaped domain problem from an initial uniform structured mesh with ${n_{\text{el}}^{\text{targ}} = 1000}$.
			\label{fig:LDomain_TargetElements_MeshEvolution}}
	\end{figure}

	The mesh evolution during the fully adaptive remeshing process for the L-shaped domain problem is depicted in Figure~\ref{fig:LDomain_TargetElements_MeshEvolution} for an initially uniform structured mesh with an element target of ${n_{\text{el}}^{\text{targ}} = 1000}$.
	The differences in the remeshing steps between those presented here and those presented in Figure~\ref{fig:LDomain_TargetError_MeshEvolution} for an error target are immediately apparent.
	In Figure~\ref{fig:LDomain_TargetError_MeshEvolution} both refinement and coarsening are performed from the first remeshing iteration. Conversely, here only refinement is performed during the first phase of the adaptive procedure until the element target is met (this is the expected behaviour, see Section~\ref{subsec:ElementSelectionTargetNumElements}). It is clear that remeshing step 3 satisfies the error target and signifies the end of the first phase of the procedure as from step 4 onwards refinement and coarsening are performed simultaneously.
	During the second phase of the procedure refinement and coarsening are performed while keeping the number of elements approximately constant. During this phase the expected distribution of elements is achieved with the areas of the domain with the highest stresses and stress gradients becoming increasingly refined while the regions of the domain experiencing simpler deformations and lower stresses are coarsened.
	
	\FloatBarrier
	\FloatBarrier
	
\begin{figure}[ht!]
		\centering
		\begin{subfigure}[t]{0.33\textwidth}
			\centering
			\includegraphics[width=0.95\textwidth,height=0.95\textwidth]{{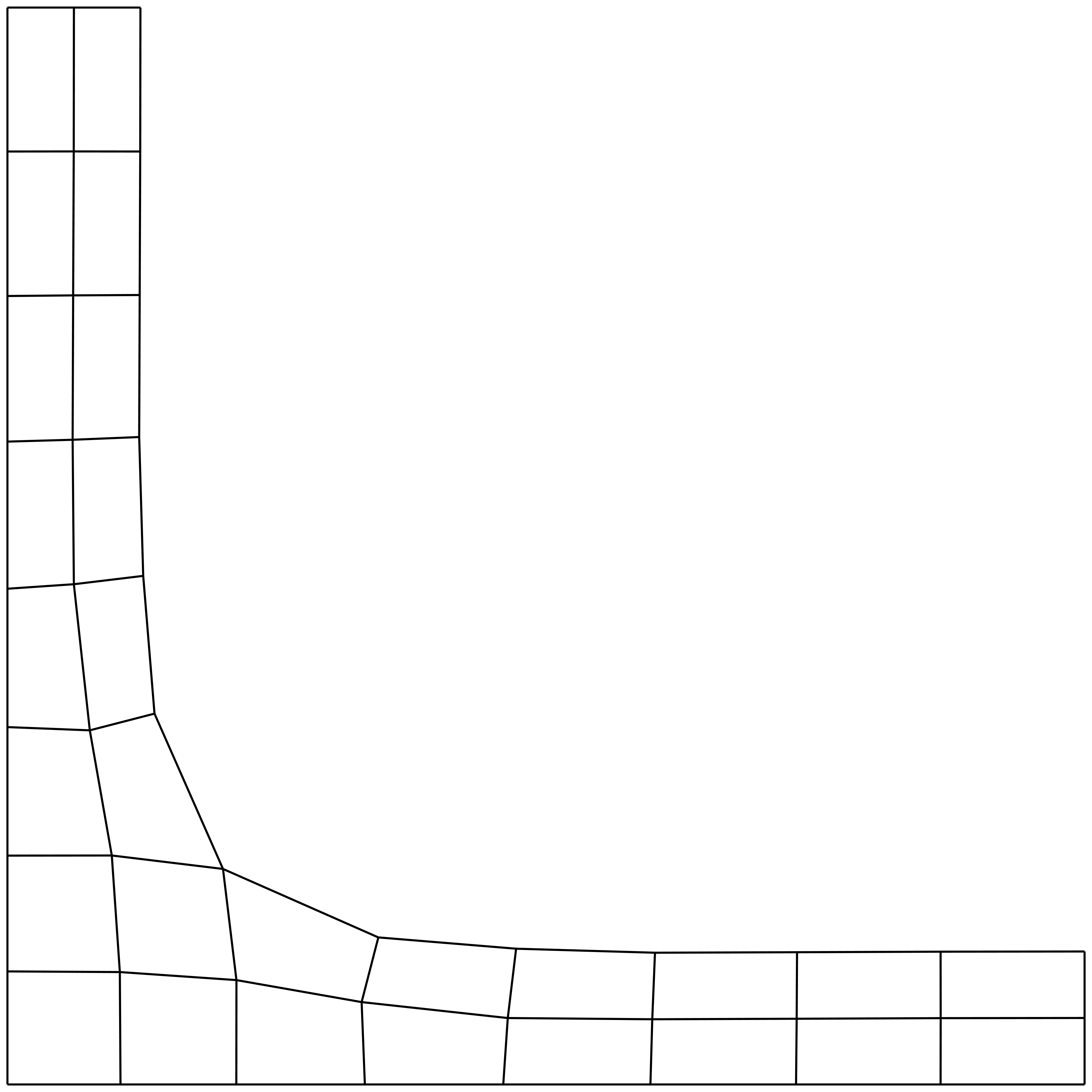}}
			\caption{Coarse mesh: Initial mesh}
		\end{subfigure}%
		\begin{subfigure}[t]{0.33\textwidth}
			\centering
			\includegraphics[width=0.95\textwidth,height=0.95\textwidth]{{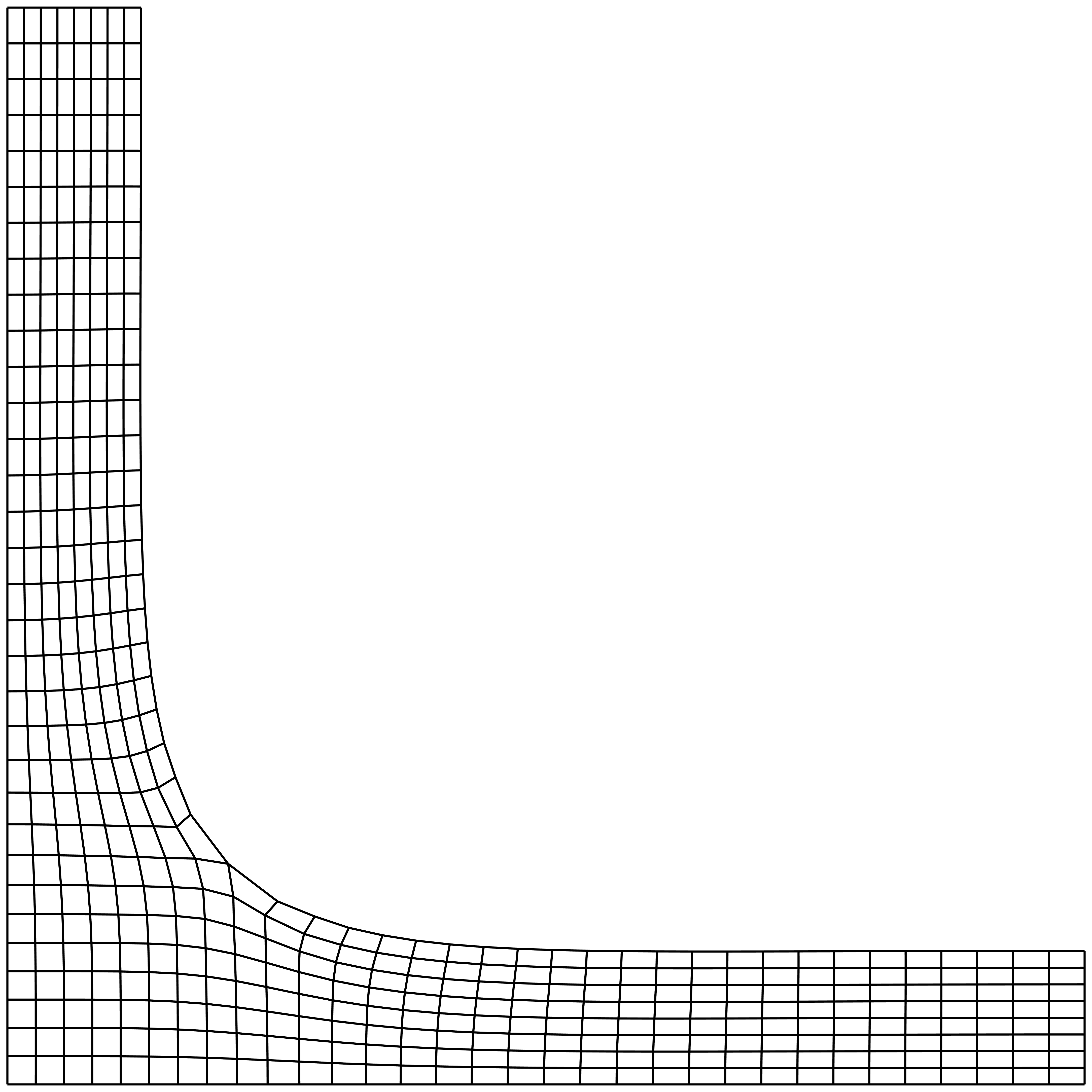}}
			\caption{Intermediate mesh: Initial mesh}
		\end{subfigure}%
		\begin{subfigure}[t]{0.33\textwidth}
			\centering
			\includegraphics[width=0.95\textwidth,height=0.95\textwidth]{{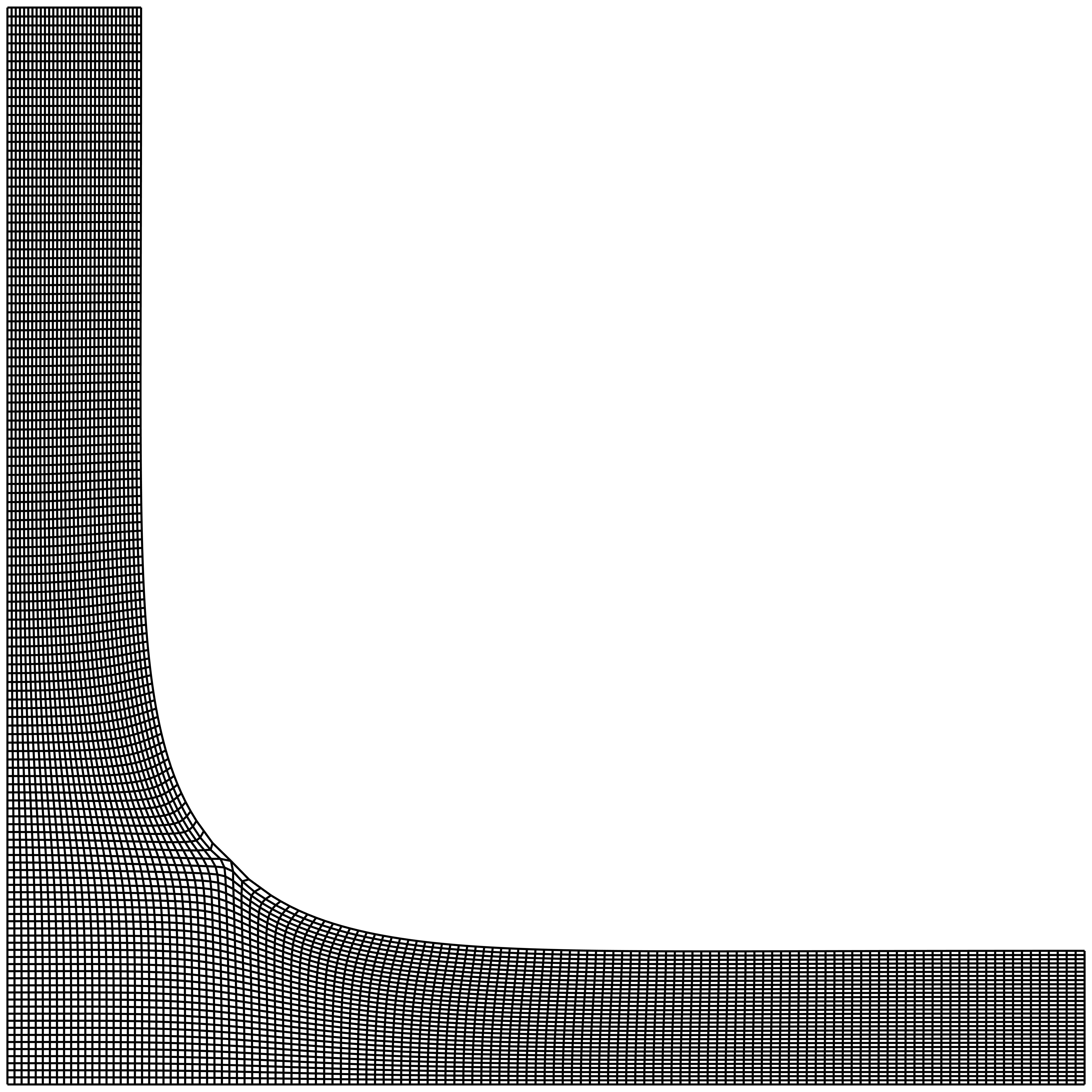}}
			\caption{Fine mesh: Initial mesh}
		\end{subfigure}
		\vskip \baselineskip 
		\begin{subfigure}[t]{0.33\textwidth}
			\centering
			\includegraphics[width=0.95\textwidth,height=0.95\textwidth]{{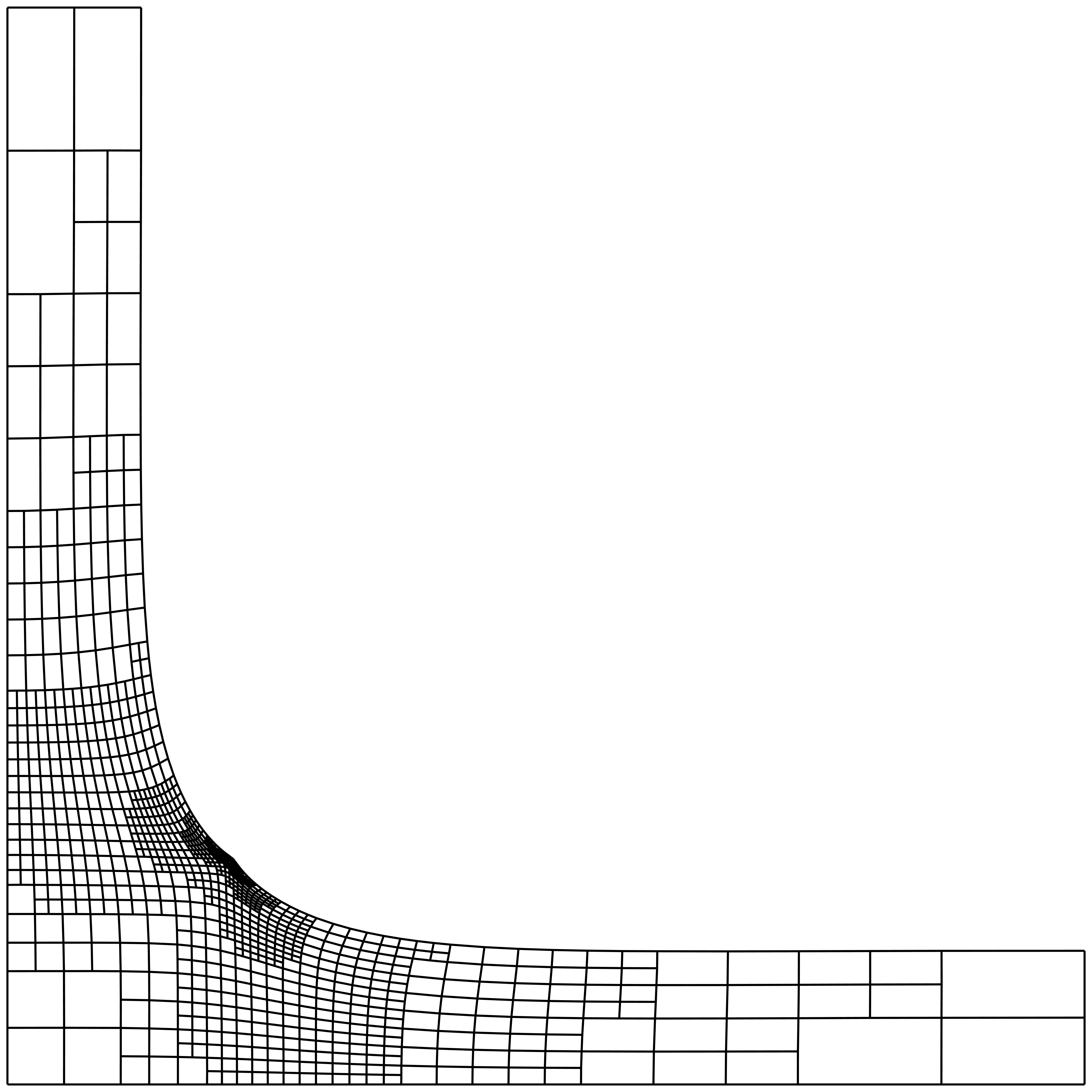}}
			\caption{Coarse mesh: Final mesh}
		\end{subfigure}%
		\begin{subfigure}[t]{0.33\textwidth}
			\centering
			\includegraphics[width=0.95\textwidth,height=0.95\textwidth]{{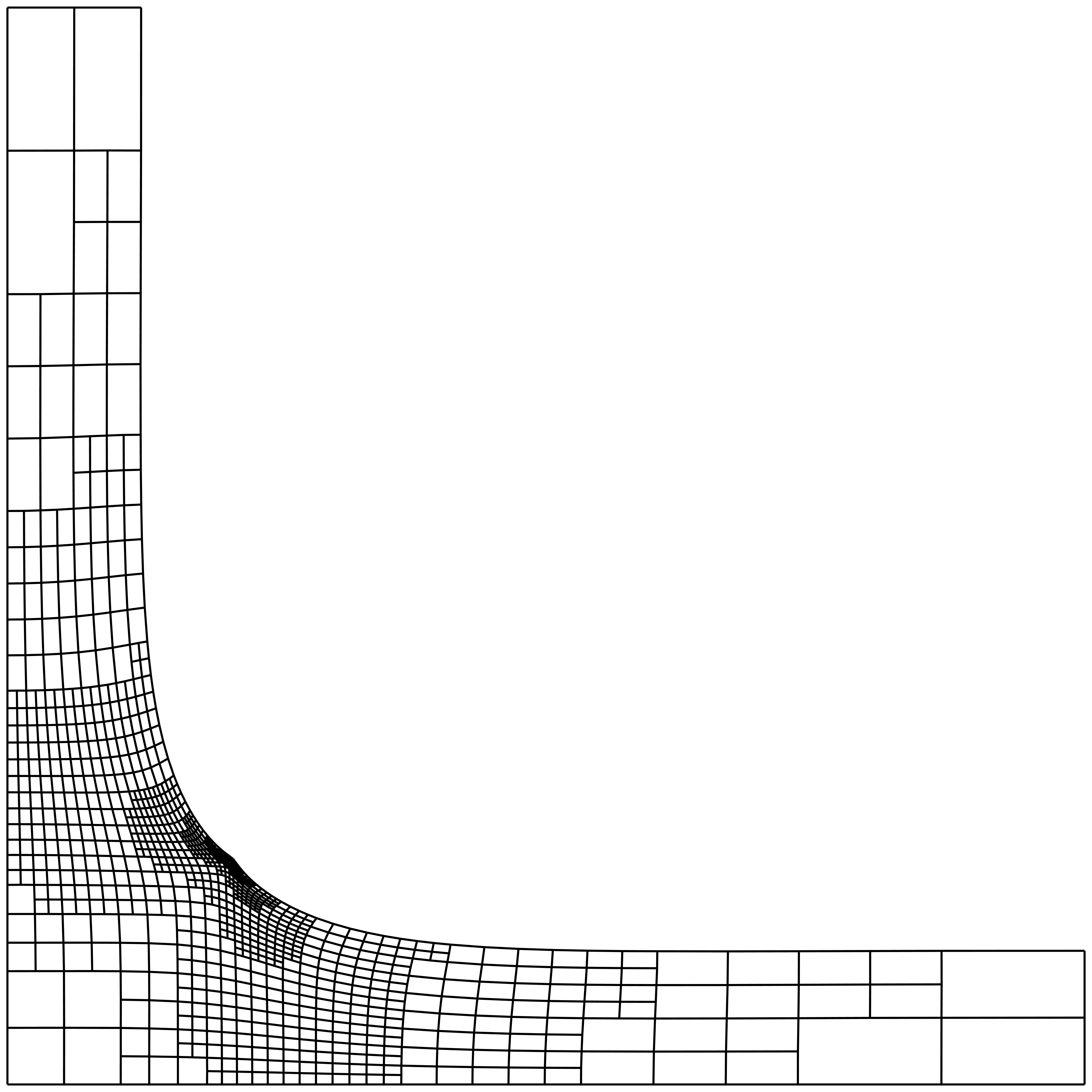}}
			\caption{Intermediate mesh: Final mesh}
		\end{subfigure}%
		\begin{subfigure}[t]{0.33\textwidth}
			\centering
			\includegraphics[width=0.95\textwidth,height=0.95\textwidth]{{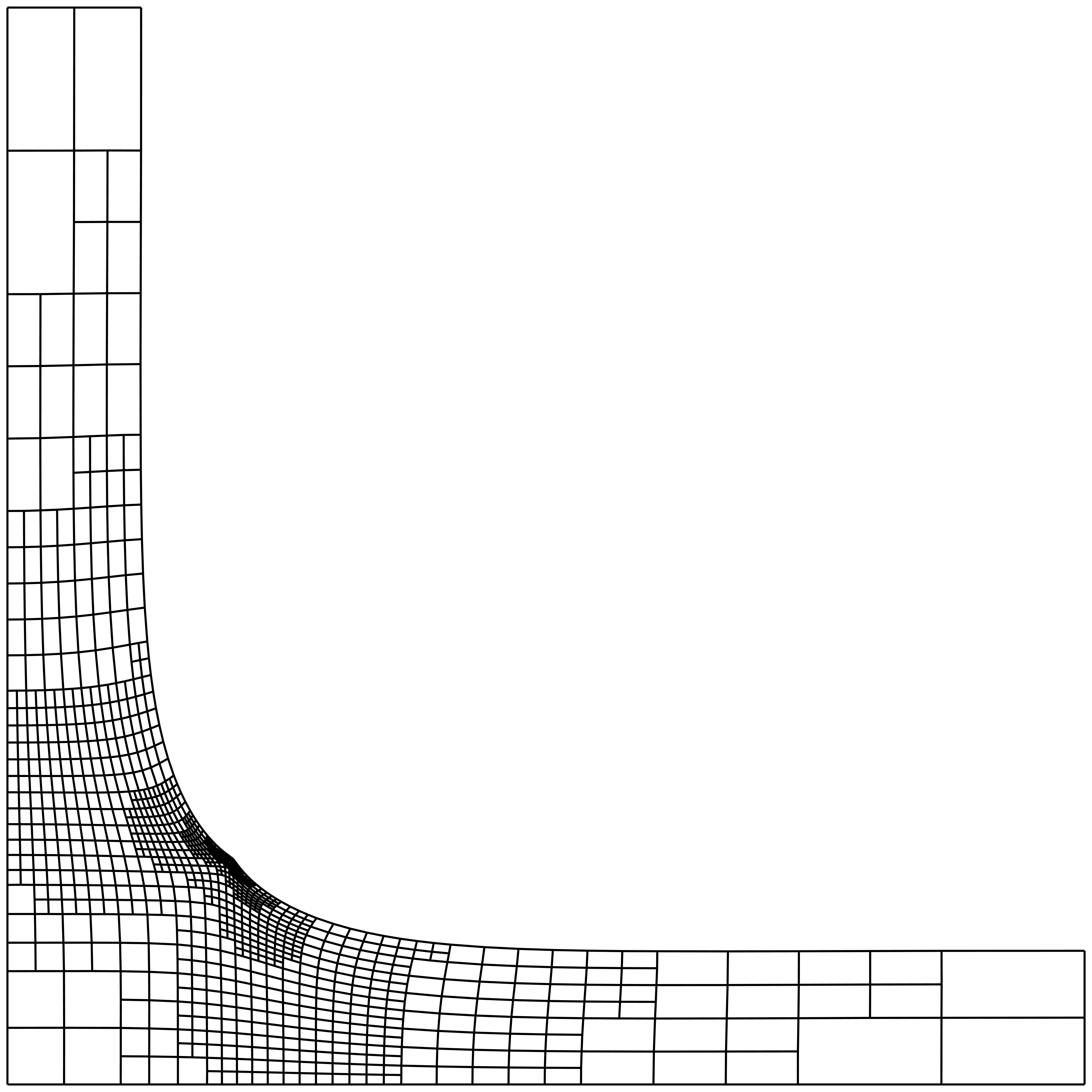}}
			\caption{Fine mesh: Final mesh}
		\end{subfigure}
		\caption{Initial and final adapted meshes using the fully adaptive remeshing procedure for the L-shaped domain problem from initial uniform structured meshes of various densities with ${n_{\text{el}}^{\text{targ}} = 1000}$.
			\label{fig:LDomain_TargetElements_MeshEvolution_VariousInitialMeshes}}
	\end{figure}

	The results of the adaptive remeshing process for the L-shaped domain problem are depicted in Figure~\ref{fig:LDomain_TargetElements_MeshEvolution_VariousInitialMeshes} for initially uniform structured meshes of varying density with an element target of ${n_{\text{el}}^{\text{targ}} = 1000}$. The top row of figures depicts the initial meshes while the bottom row depicts the final adapted meshes after the element target and termination criteria have been met.
	The final adapted meshes exhibit the same sensible and intuitive element distribution as observed in Figure~\ref{fig:LDomain_TargetElements_MeshEvolution}. Notably, the final adapted meshes are, again, almost identical for all initial uniform meshes considered. Thus, the output, or final result, of the fully adaptive remeshing procedure is independent of the initial mesh.
	
	\FloatBarrier
	\FloatBarrier
	
\begin{figure}[ht!]
		\centering
		\begin{subfigure}[t]{0.495\textwidth}
			\centering
			\includegraphics[width=0.95\textwidth]{{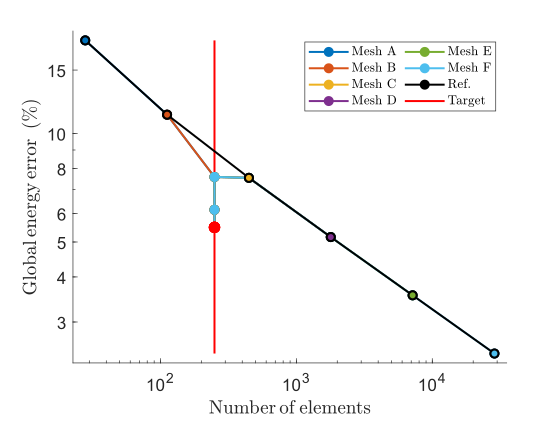}}
			\caption{$n_{\text{el}}^{\text{targ}} = 250$}
		\end{subfigure}%
		\begin{subfigure}[t]{0.495\textwidth}
			\centering
			\includegraphics[width=0.95\textwidth]{{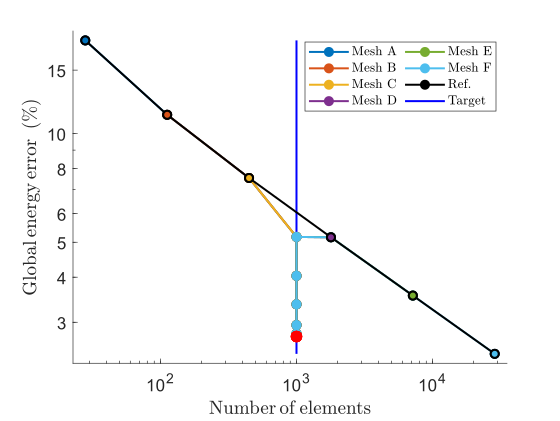}}
			\caption{$n_{\text{el}}^{\text{targ}} = 1000$}
		\end{subfigure}
		\vskip \baselineskip 
		\begin{subfigure}[t]{0.495\textwidth}
			\centering
			\includegraphics[width=0.95\textwidth]{{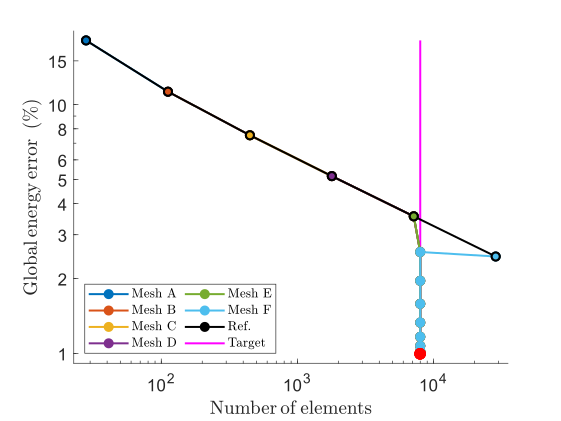}}
			\caption{$n_{\text{el}}^{\text{targ}} = 8000$}
		\end{subfigure}%
		\begin{subfigure}[t]{0.495\textwidth}
			\centering
			\includegraphics[width=0.95\textwidth]{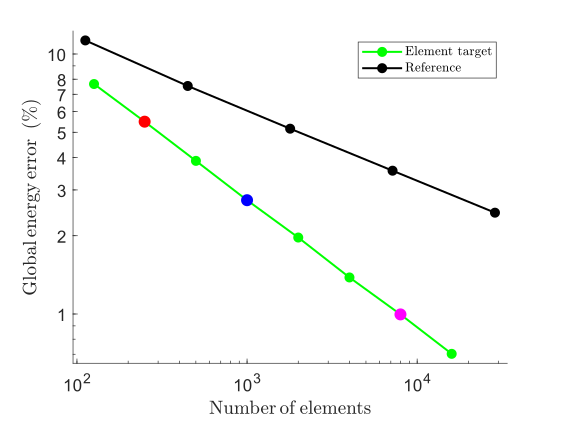}
			\caption{Various element targets}
		\end{subfigure}
		\caption{Energy error vs number of elements for the L-shaped domain problem on structured meshes of various initial densities with (a) ${n_{\text{el}}^{\text{targ}} = 250}$, (b) ${n_{\text{el}}^{\text{targ}} = 1000}$, (c) ${n_{\text{el}}^{\text{targ}} = 8000}$, and (d) average final/adapted mesh error for various element targets.  
			\label{fig:LDomain_TargetElements_ErrorConvergence_Structured}}
	\end{figure}

	The convergence behaviour of the energy error approximation vs the number of elements in the mesh is depicted on a logarithmic scale in Figures~\ref{fig:LDomain_TargetElements_ErrorConvergence_Structured}(a)-(c) for the L-shaped domain problem. The convergence behaviour is plotted for cases of several initially uniform structured meshes of varying density for various element targets. 
	Additionally, the averaged final adapted mesh result (red markers) for all element targets considered is plotted in Figure~\ref{fig:LDomain_TargetElements_ErrorConvergence_Structured}(d).
	Where applicable, the markers in Figure~\ref{fig:LDomain_TargetElements_ErrorConvergence_Structured}(d) are colored to match their corresponding targets depicted in Figures~\ref{fig:LDomain_TargetElements_ErrorConvergence_Structured}(a)-(c).
	From Figures~\ref{fig:LDomain_TargetElements_ErrorConvergence_Structured}(a)-(c) it is clear that the fully adaptive procedure is able to meet the specified element targets from any initial mesh.
	During the first phase of the adaptive procedure uniform refinement or coarsening is performed until the point at which a uniform refinement or coarsening would overshoot the element target. During this phase the convergence curves follow the uniform reference refinement curve. Thereafter, a partial refinement or coarsening is performed such that the error target is met. This step is the first time the curves of the adaptive procedure deviate from the reference curve and signifies the end of the first phase. 
	During the second phase the adaptive curves move vertically as the number of elements its held constant while selective refinement and coarsening is performed to evenly distribute the element-level error. This process has the effect of reducing the global error and is performed until the termination criteria is met (see Section~\ref{subsec:ElementSelectionTargetNumElements}).
	The final adapted meshes have an almost identical error for a specific element target. This, again, demonstrates that the performance of the fully adaptive remeshing procedure is independent of the initial mesh.
	From Figure~\ref{fig:LDomain_TargetElements_ErrorConvergence_Structured}(d) it is clear that the outputs of the fully adaptive procedure for various element targets exhibits a linear convergence rate. Since the procedure aims to generate a quasi-optimal mesh, it is again expected that this is the (approximately) optimal convergence rate for this problem.
	
	\FloatBarrier
	\FloatBarrier
	
\begin{figure}[ht!]
		\centering
		\begin{subfigure}[t]{0.495\textwidth}
			\centering
			\includegraphics[width=0.95\textwidth]{{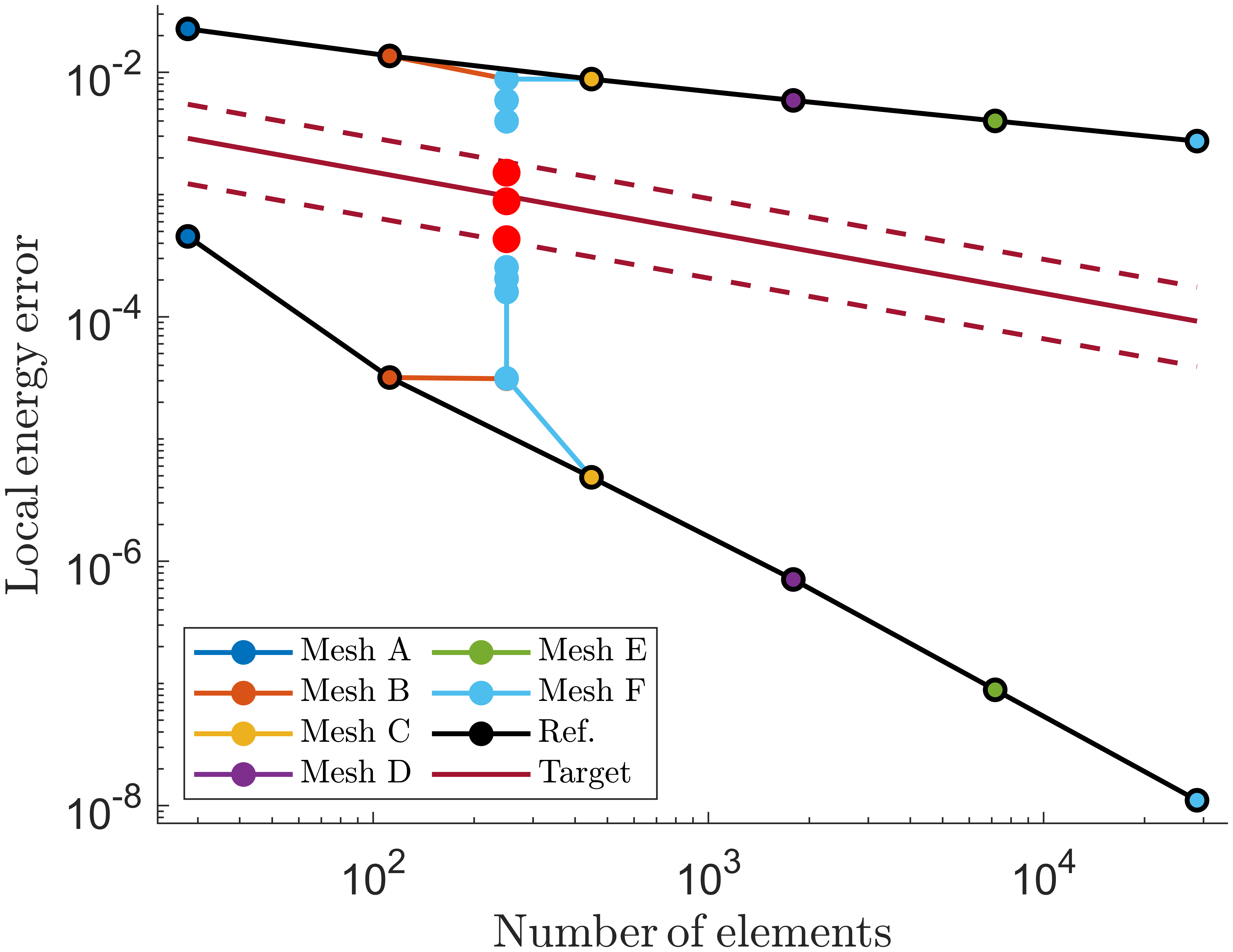}}
			\caption{$n_{\text{el}}^{\text{targ}} = 250$ - Max and min local error}
			\vspace*{-3mm}
		\end{subfigure}%
		\begin{subfigure}[t]{0.495\textwidth}
			\centering
			\includegraphics[width=0.95\textwidth]{{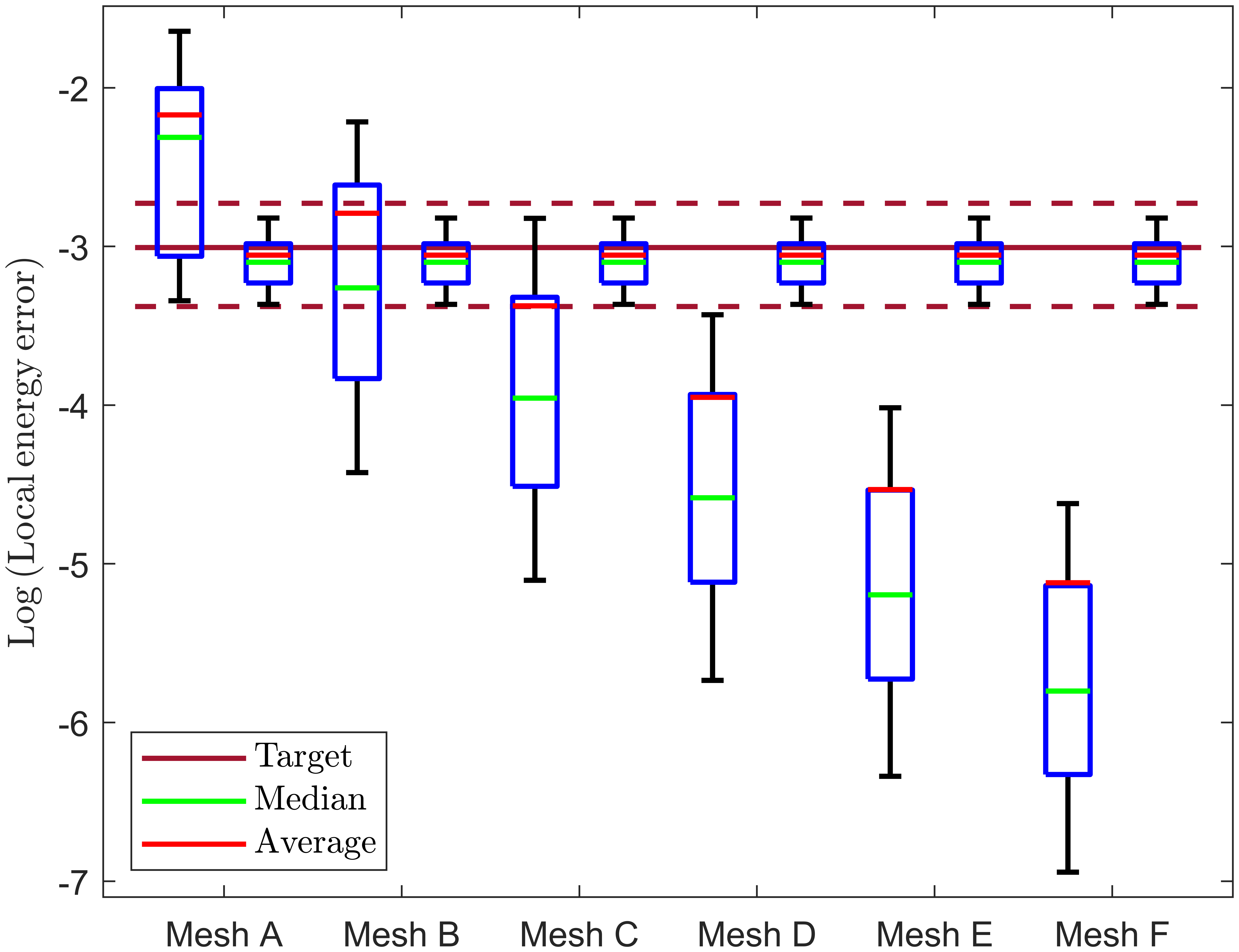}}
			\caption{$n_{\text{el}}^{\text{targ}} = 250$ - Local error distribution}
			\vspace*{-3mm}
		\end{subfigure}
		\vskip \baselineskip 
		\begin{subfigure}[t]{0.495\textwidth}
			\centering
			\includegraphics[width=0.95\textwidth]{{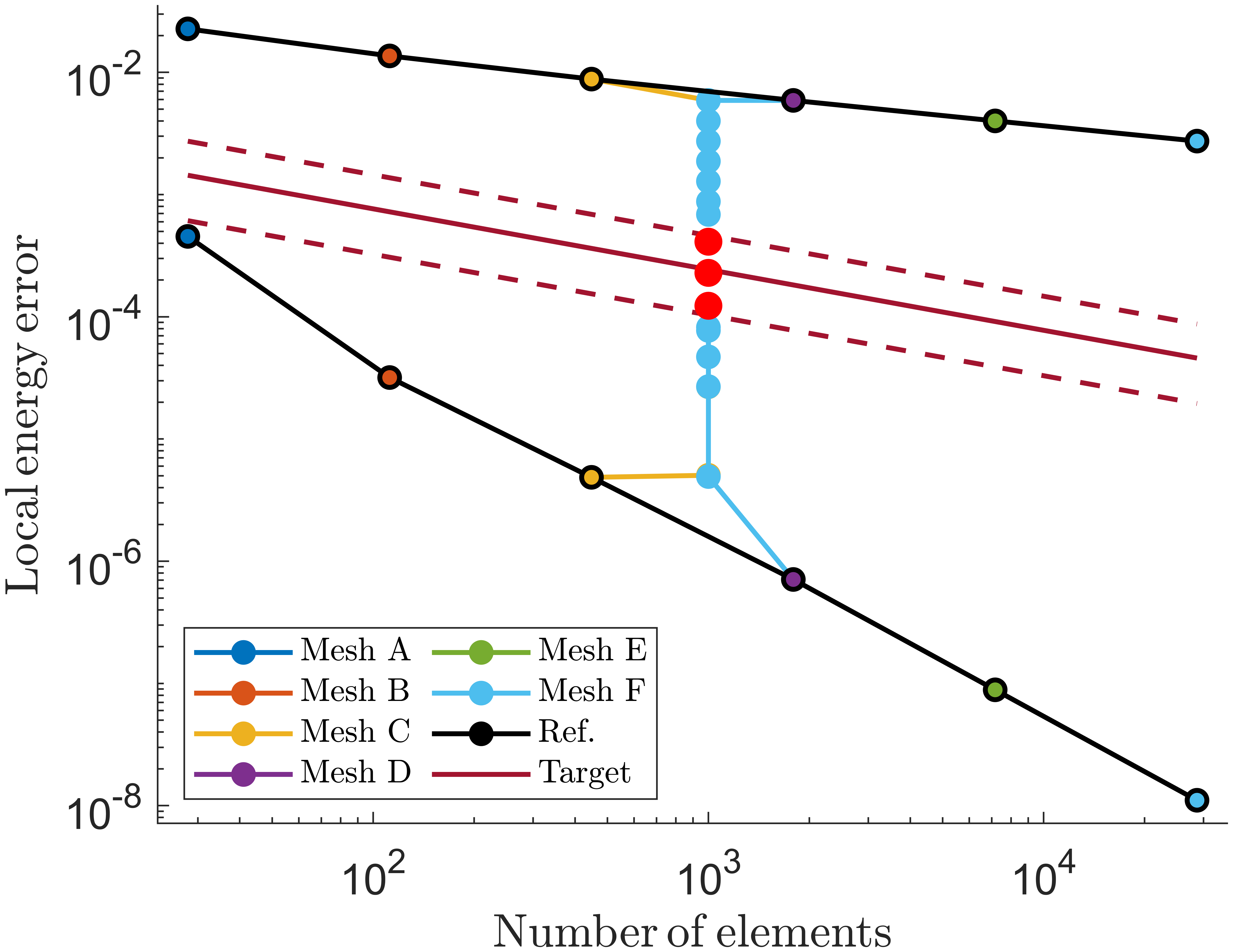}}
			\caption{$n_{\text{el}}^{\text{targ}} = 1000$ - Max and min local error}
			\vspace*{-3mm}
		\end{subfigure}%
		\begin{subfigure}[t]{0.495\textwidth}
			\centering
			\includegraphics[width=0.95\textwidth]{{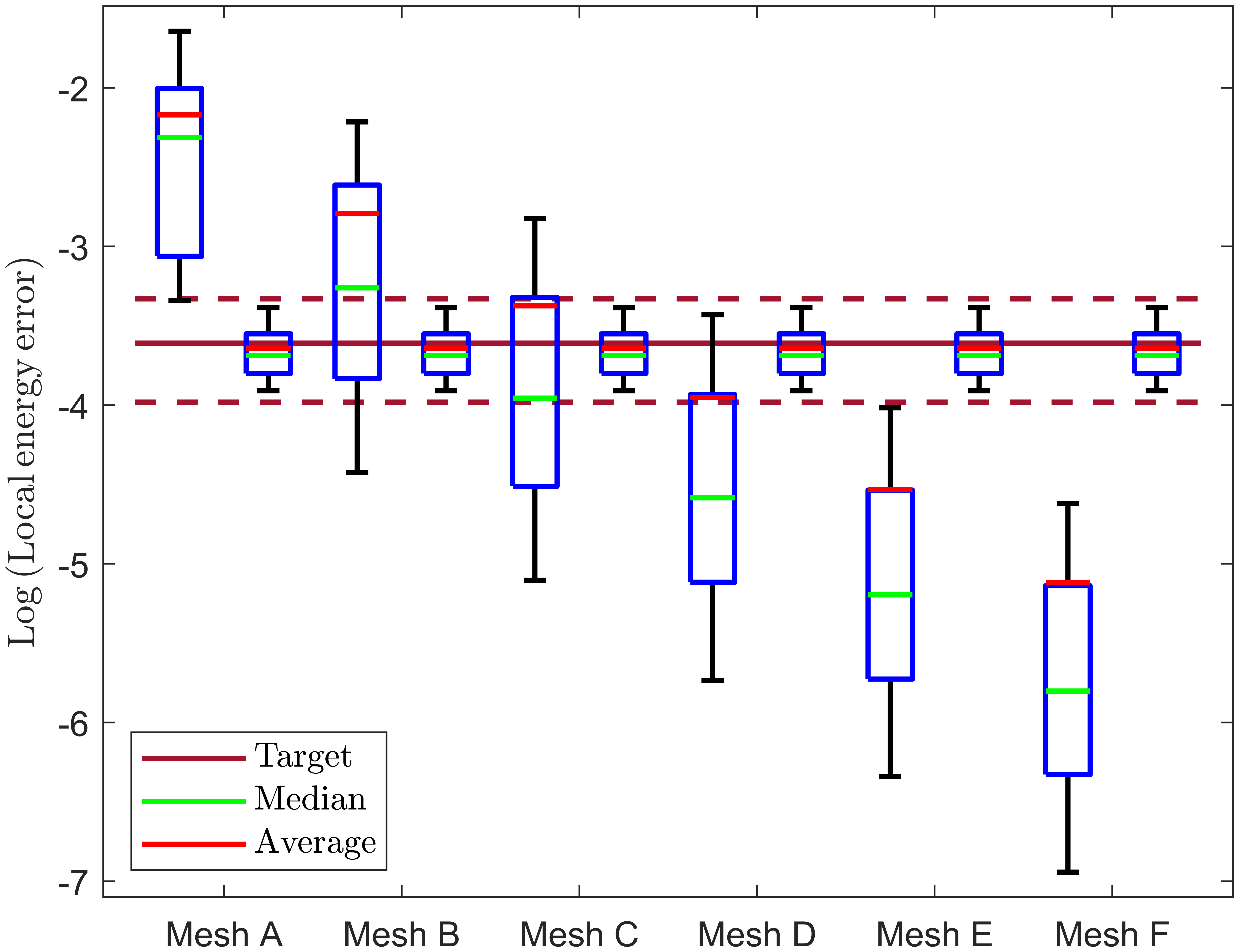}}
			\caption{$n_{\text{el}}^{\text{targ}} = 1000$ - Local error distribution}
			\vspace*{-3mm}
		\end{subfigure}
		\vskip \baselineskip 
		\begin{subfigure}[t]{0.495\textwidth}
			\centering
			\includegraphics[width=0.95\textwidth]{{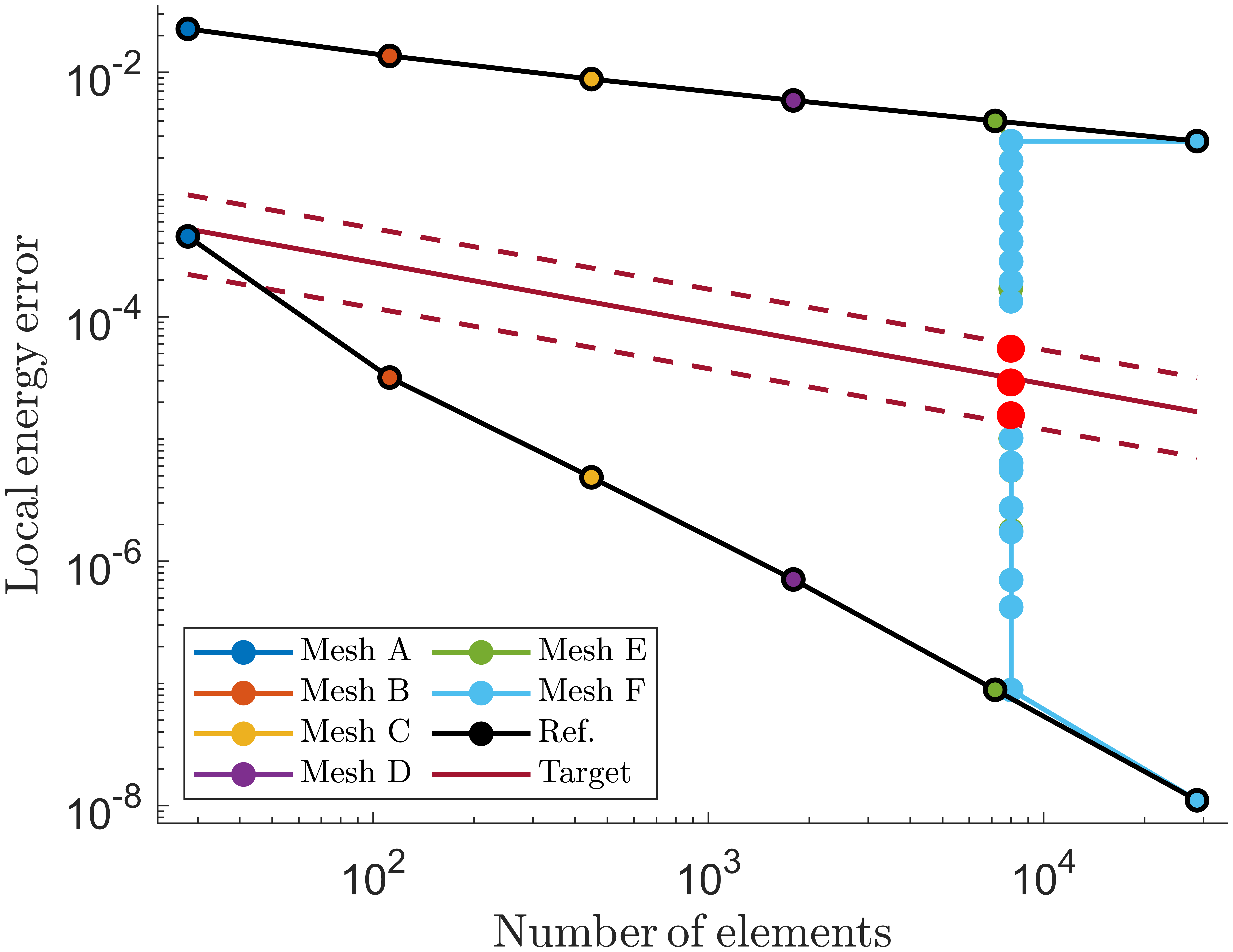}}
			\caption{$n_{\text{el}}^{\text{targ}} = 8000$ - Max and min local error}
			\vspace*{-3mm}
		\end{subfigure}%
		\begin{subfigure}[t]{0.495\textwidth}
			\centering
			\includegraphics[width=0.95\textwidth]{{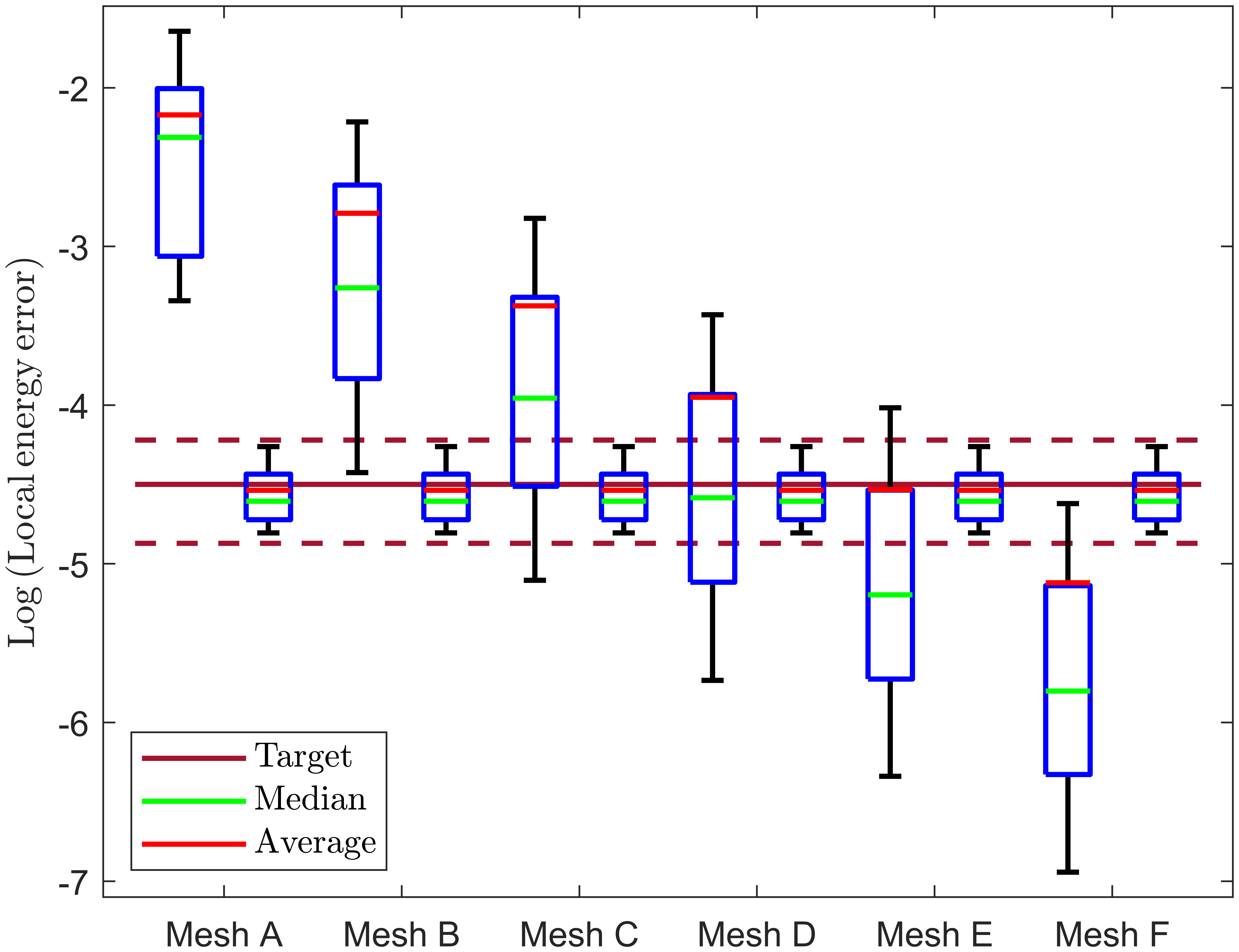}}
			\caption{$n_{\text{el}}^{\text{targ}} = 8000$ - Local error distribution}
			\vspace*{-3mm}
		\end{subfigure}
		\caption{Max and min energy error (left) and box and whisker plot of energy error distribution (right) for the L-shaped domain problem on structured meshes of various initial densities for a range of element targets.
			\label{fig:LDomain_TargetElements_ErrorDistribution_Structured}}
	\end{figure}

	The error distribution during the mesh adaptation process for the L-shaped domain problem is depicted in Figure~\ref{fig:LDomain_TargetElements_ErrorDistribution_Structured} for various element targets on structured meshes.
	The nature of the convergence exhibited by the left column of figures is significantly different to that exhibited in Figure~\ref{fig:LDomain_TargetError_ErrorConvergence_Structured} for a target error.
	The difference is a result of the two distinct adaptive phases for the case of a target number of nodes.
	However, the error distributions of the final adapted meshes, as indicated by the red markers and the box and whisker plots, are qualitatively very similar to those of Figure~\ref{fig:LDomain_TargetError_ErrorConvergence_Structured} for a target error.
	Specifically, the maximum and minimum errors fall within the upper and lower bounds, the upper and lower quartiles indicate a narrow distribution of error around the average, and the average error is almost identical to the target error. Additionally, the average and median errors are almost identical which, again, emphasises the narrow distribution of the element-level errors.
	
	The results presented in Figure~\ref{fig:LDomain_TargetElements_ErrorConvergence_Structured} demonstrated that the fully adaptive procedure was able to meet all specified element targets on structured meshes. This is indicated by the red markers denoting the final adapted mesh lying exactly on the target element line. 
	The results presented in Figure~\ref{fig:LDomain_TargetElements_ErrorDistribution_Structured} for structured meshes demonstrated that the average element-level error almost exactly met the element-level target as the red markers denoting the average error strongly overlap the solid maroon target line. Furthermore, the element-level errors were satisfactorily equal as they all fell within the specified target error range.
	Thus, the proposed fully adaptive procedure successfully generated quasi-optimal meshes for the specified element targets on structured meshes.
	
	\FloatBarrier
	\FloatBarrier
	
	\subsubsection{Target number of nodes}
	\label{subsubsec:L_Domain_TargetNumNodes}
	
\begin{figure}[ht!]
		\centering
		\begin{subfigure}[t]{0.33\textwidth}
			\centering
			\includegraphics[width=0.95\textwidth,height=0.95\textwidth]{{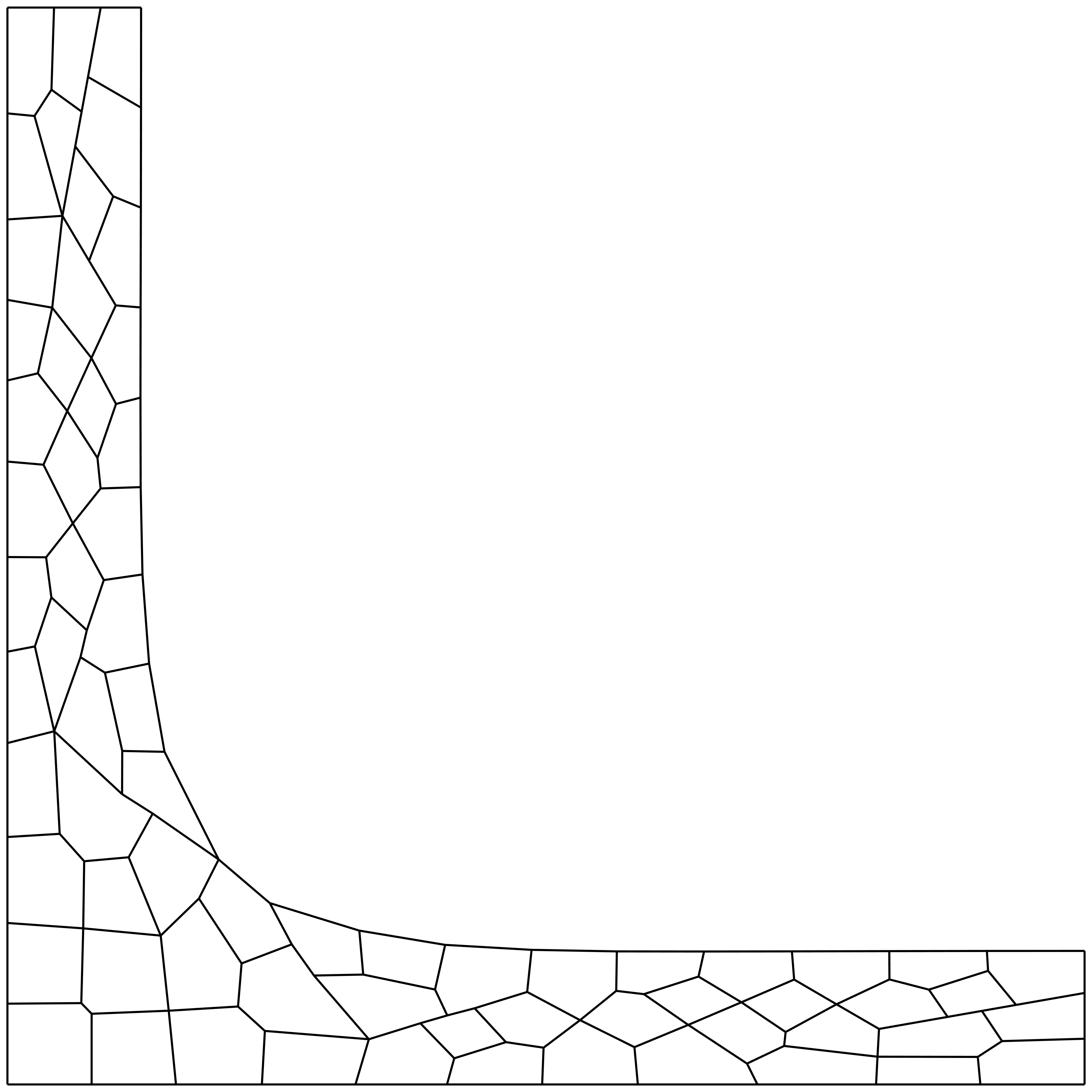}}
			\caption{Step 1: Initial uniform mesh}
		\end{subfigure}%
		\begin{subfigure}[t]{0.33\textwidth}
			\centering
			\includegraphics[width=0.95\textwidth,height=0.95\textwidth]{{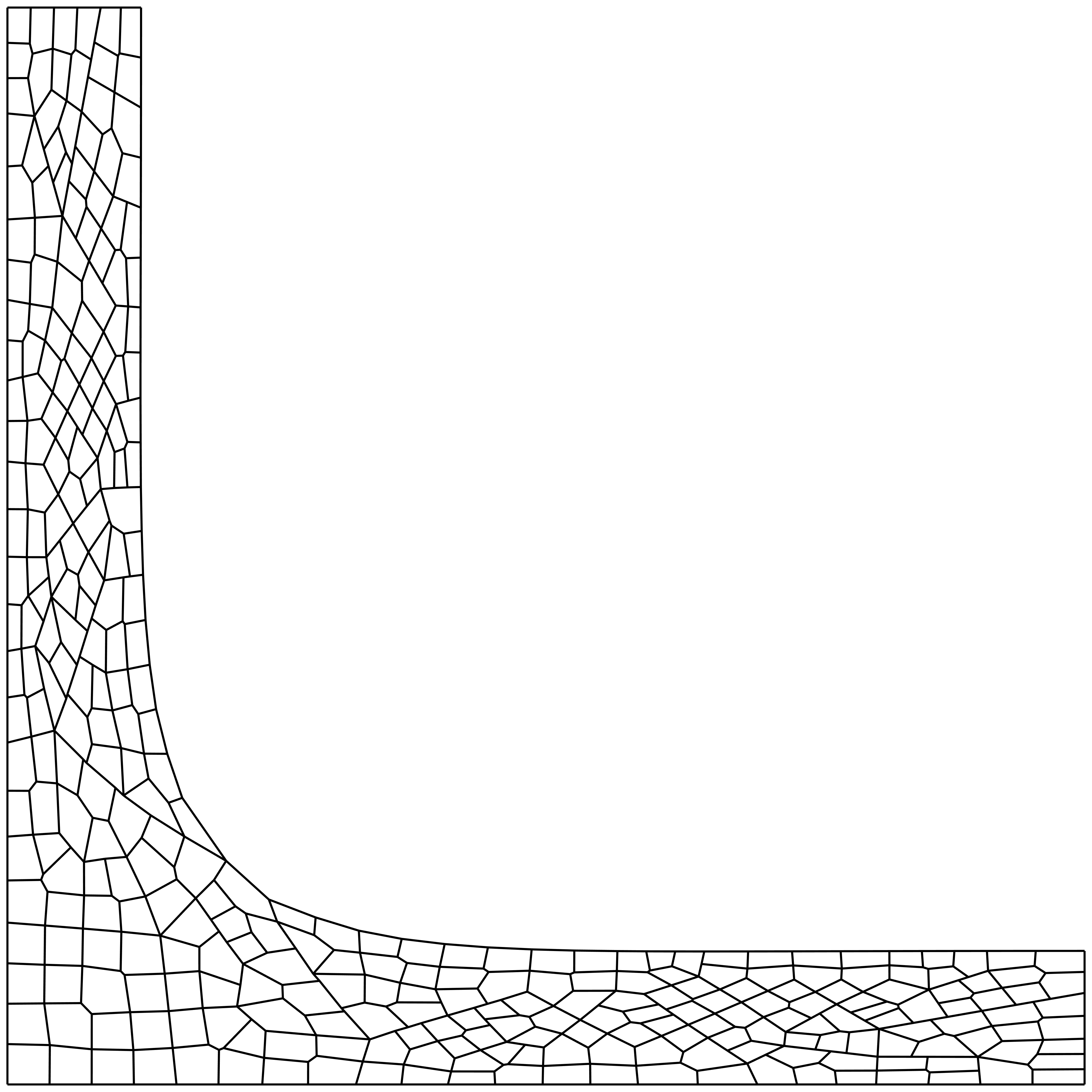}}
			\caption{Step 2}
		\end{subfigure}%
		\begin{subfigure}[t]{0.33\textwidth}
			\centering
			\includegraphics[width=0.95\textwidth,height=0.95\textwidth]{{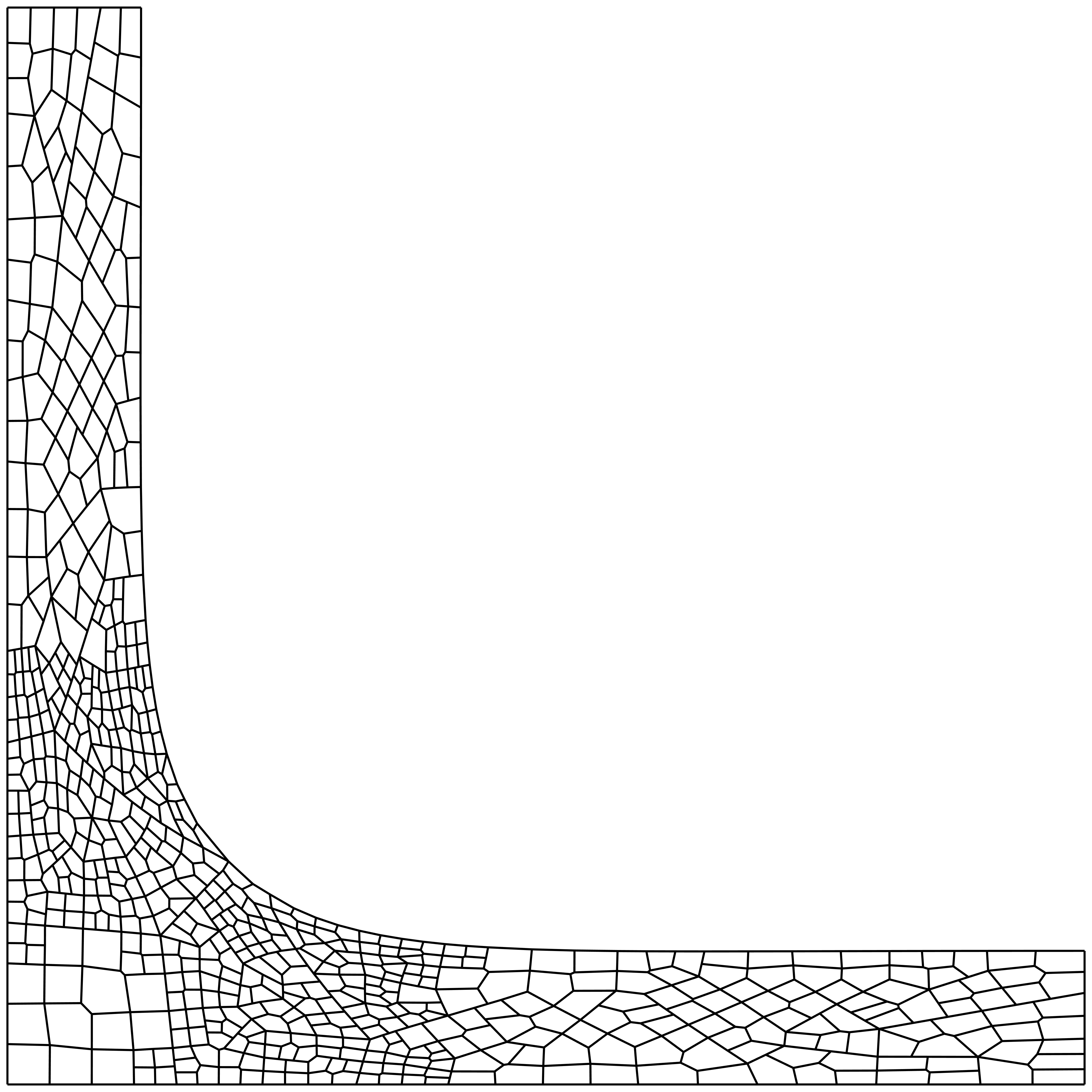}}
			\caption{Step 3}
		\end{subfigure}
		\vskip \baselineskip 
		\begin{subfigure}[t]{0.33\textwidth}
			\centering
			\includegraphics[width=0.95\textwidth,height=0.95\textwidth]{{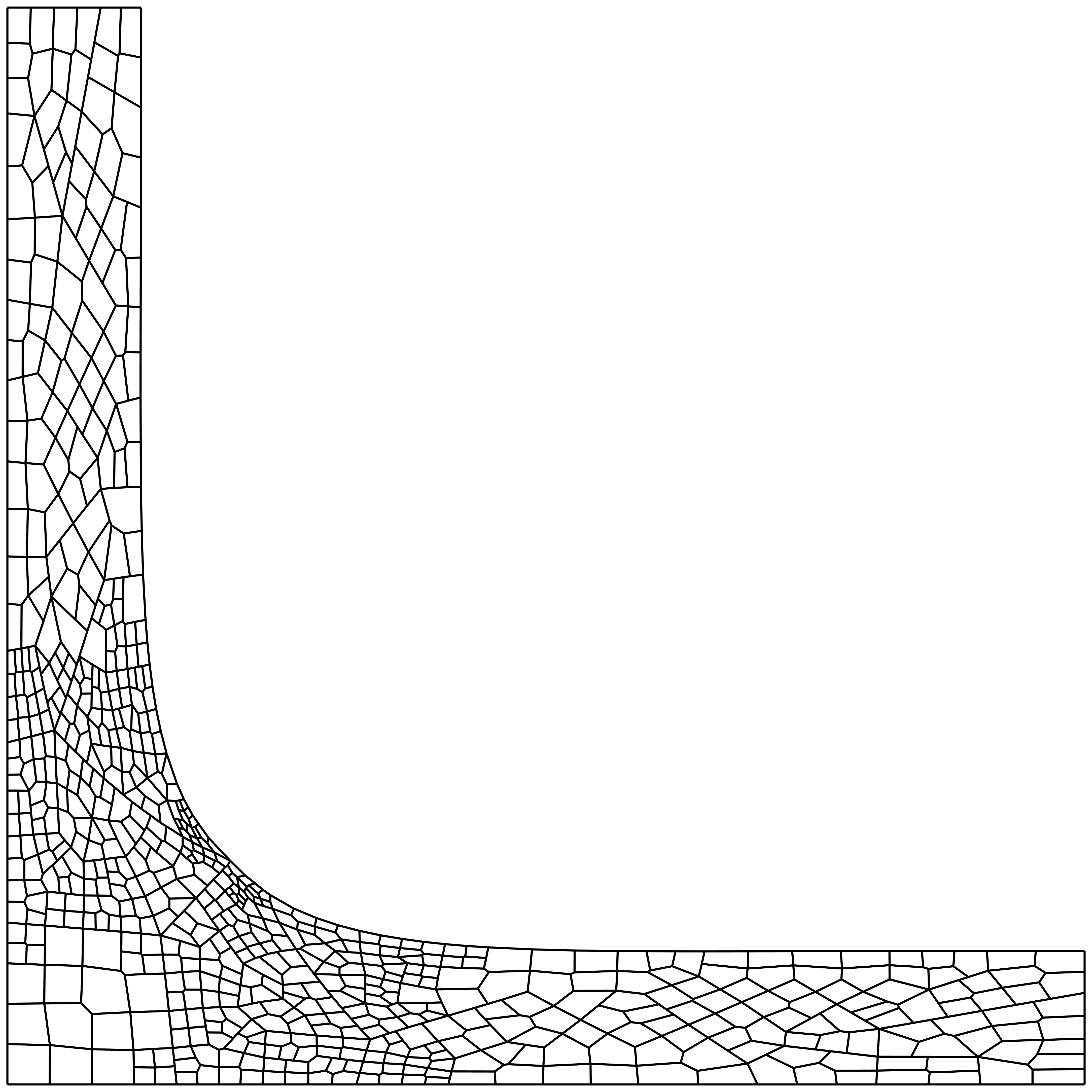}}
			\caption{Step 4}
		\end{subfigure}%
		\begin{subfigure}[t]{0.33\textwidth}
			\centering
			\includegraphics[width=0.95\textwidth,height=0.95\textwidth]{{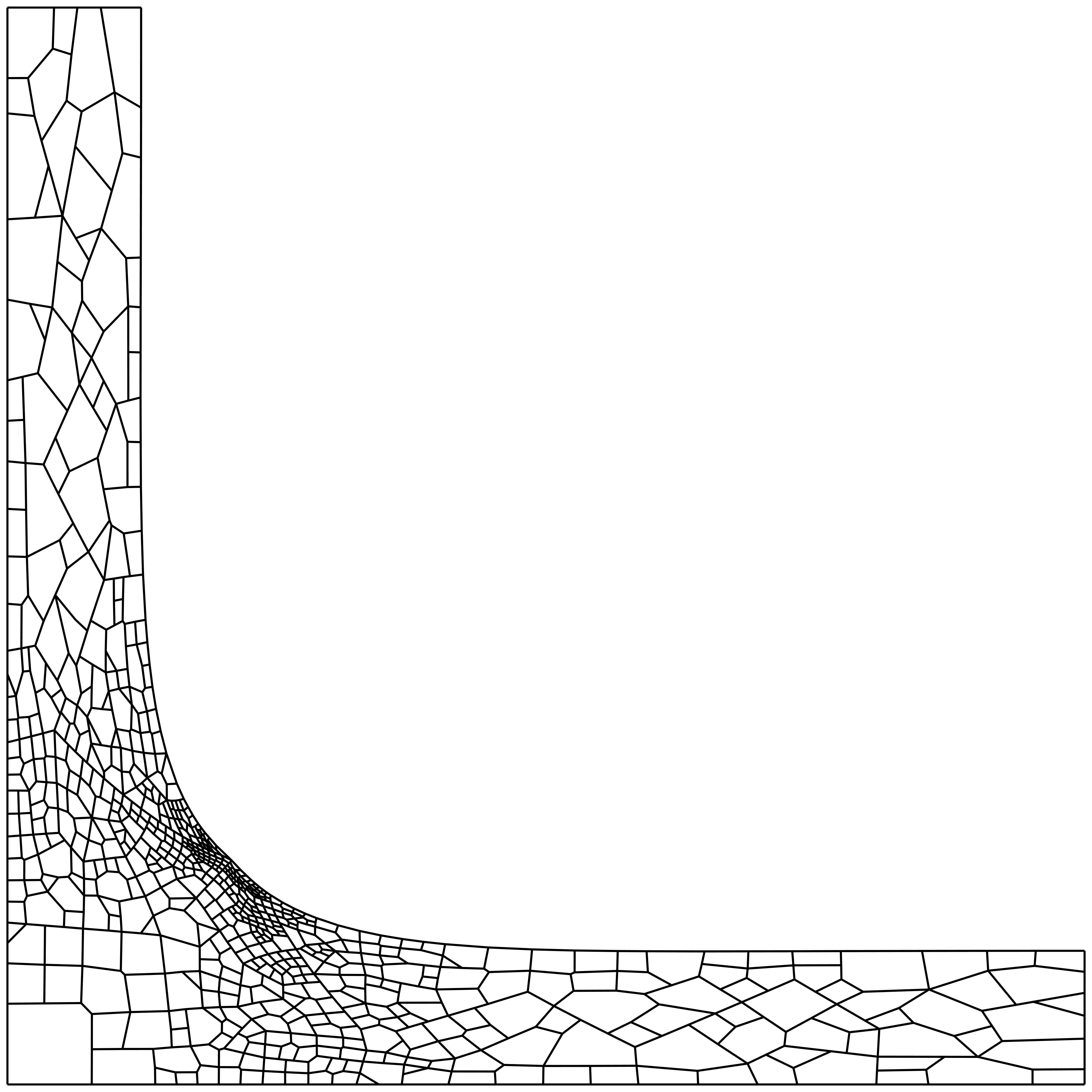}}
			\caption{Step 5}
		\end{subfigure}%
		\begin{subfigure}[t]{0.33\textwidth}
			\centering
			\includegraphics[width=0.95\textwidth,height=0.95\textwidth]{{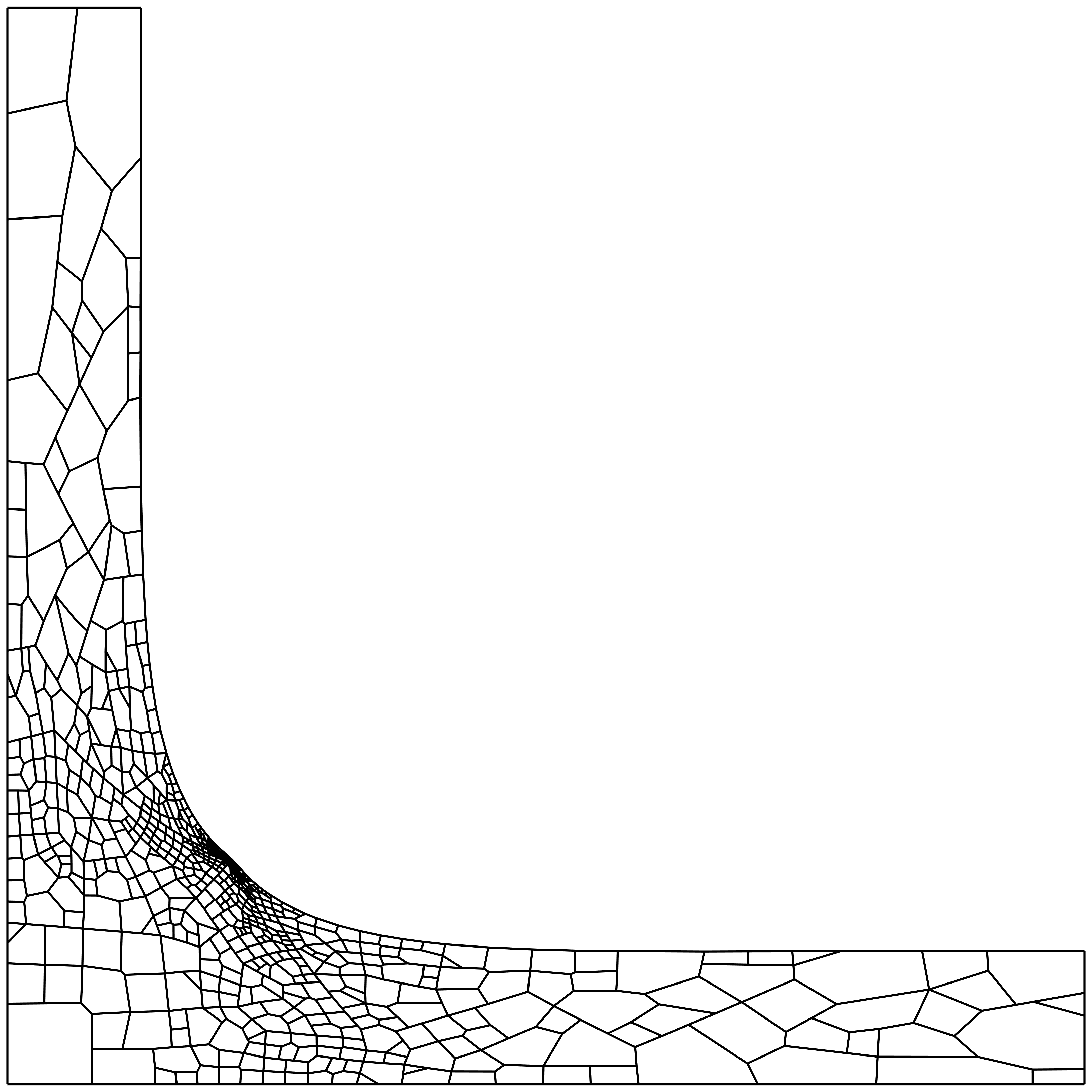}}
			\caption{Step 6: Final adapted mesh}
		\end{subfigure}
		\caption{Mesh evolution with the fully adaptive remeshing procedure for the L-shaped domain problem from an initial uniform Voronoi mesh with ${n_{\text{v}}^{\text{targ}} = 1000}$.
			\label{fig:LDomain_TargetNodes_MeshEvolution}}
	\end{figure}

	The mesh evolution during the fully adaptive remeshing process for the L-shaped domain problem is depicted in Figure~\ref{fig:LDomain_TargetNodes_MeshEvolution} for an initially uniform Voronoi mesh with a node target of ${n_{\text{v}}^{\text{targ}} = 1000}$.
	The first four steps correspond to the first phase of the adaptive remeshing procedure in which the target number of nodes is met. Thereafter, refinement and coarsening are performed simultaneously while keeping the number of nodes approximately constant. During this phase the expected distribution of elements is, again, achieved with the areas of the domain with the highest stresses and stress gradients becoming increasingly refined while the regions of the domain experiencing simpler (i.e., more uniform/homogeneous) deformations and lower stresses are coarsened.
	
	\FloatBarrier
	\FloatBarrier
	
\begin{figure}[ht!]
		\centering
		\begin{subfigure}[t]{0.33\textwidth}
			\centering
			\includegraphics[width=0.95\textwidth,height=0.95\textwidth]{{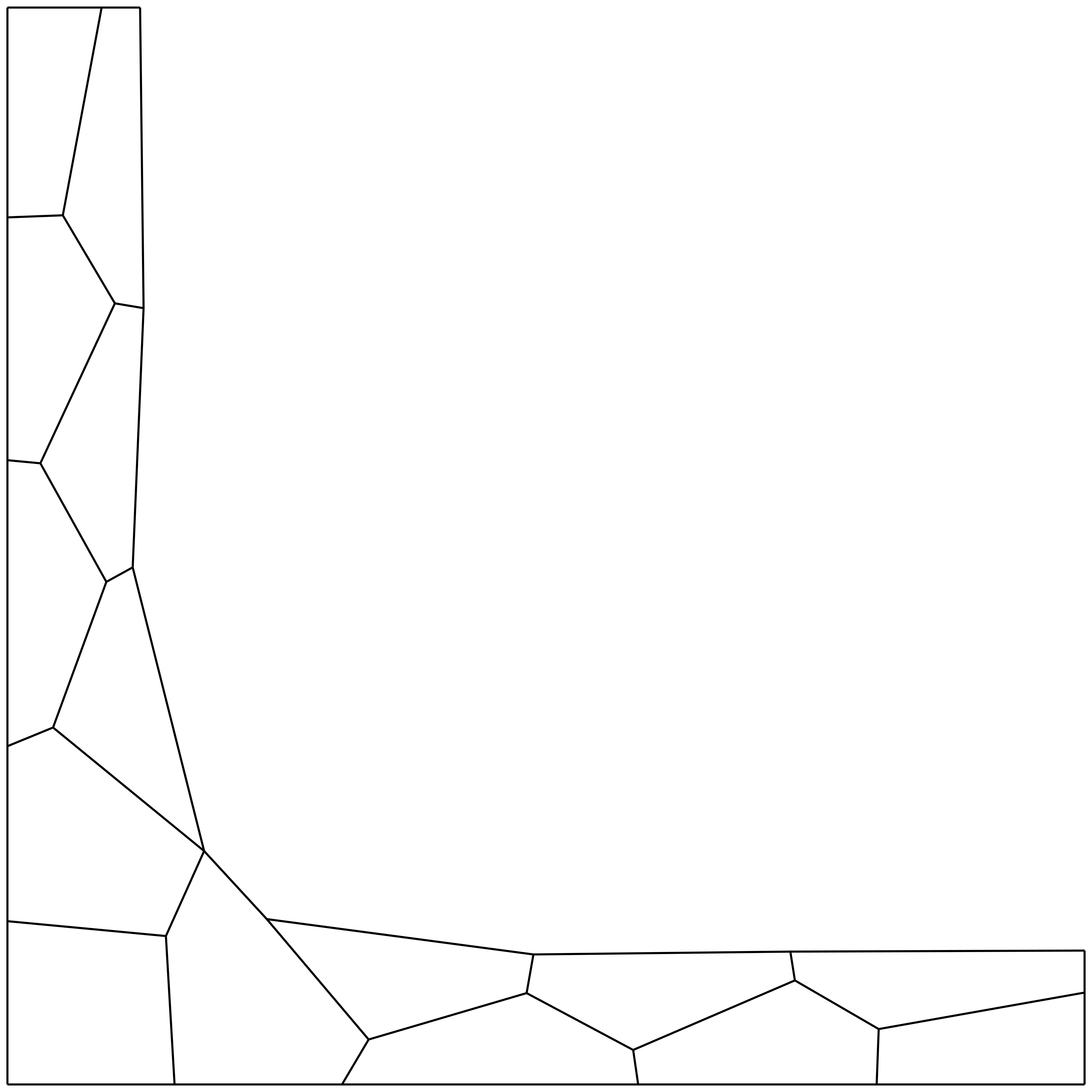}}
			\caption{Coarse mesh: Initial mesh}
		\end{subfigure}%
		\begin{subfigure}[t]{0.33\textwidth}
			\centering
			\includegraphics[width=0.95\textwidth,height=0.95\textwidth]{{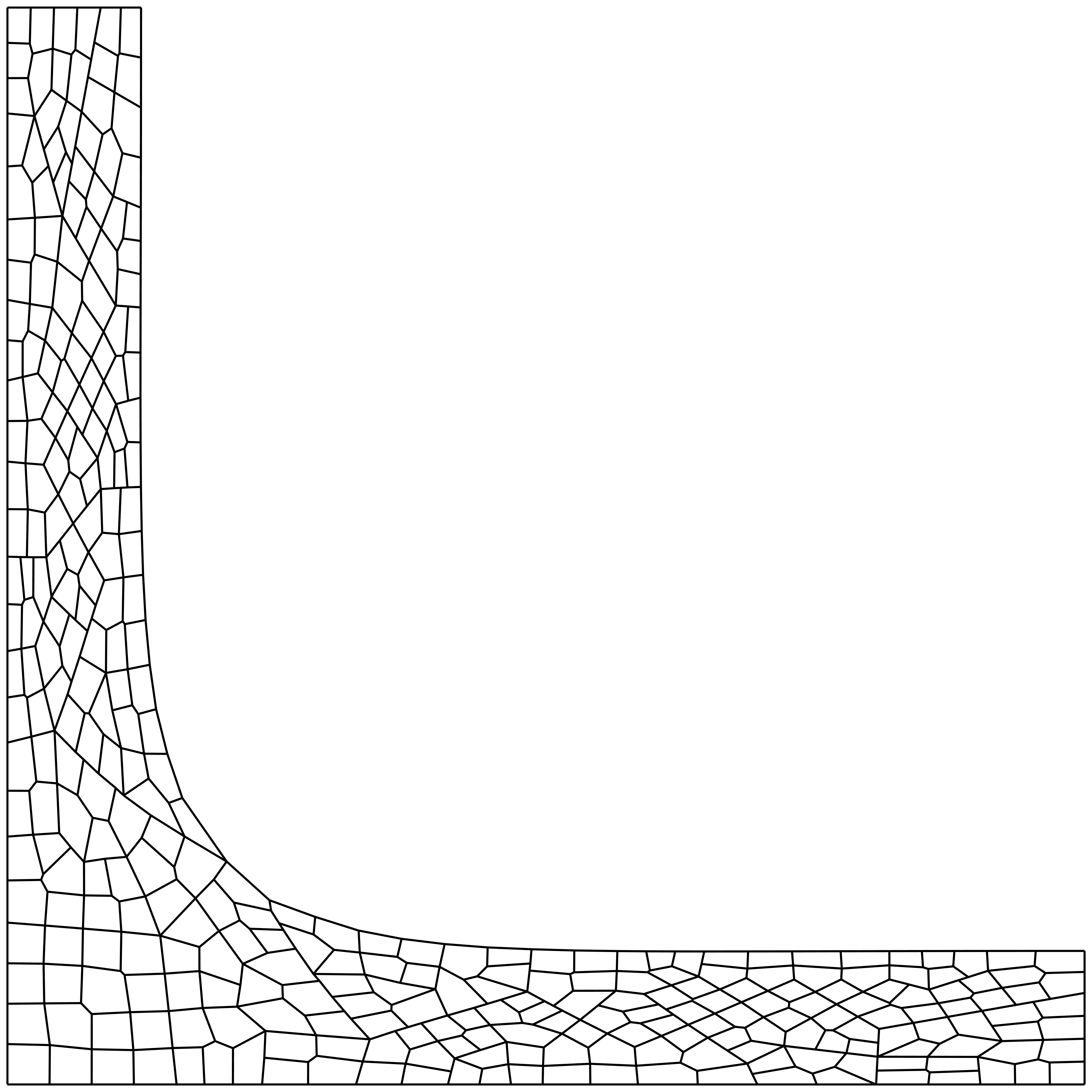}}
			\caption{Intermediate mesh: Initial mesh}
		\end{subfigure}%
		\begin{subfigure}[t]{0.33\textwidth}
			\centering
			\includegraphics[width=0.95\textwidth,height=0.95\textwidth]{{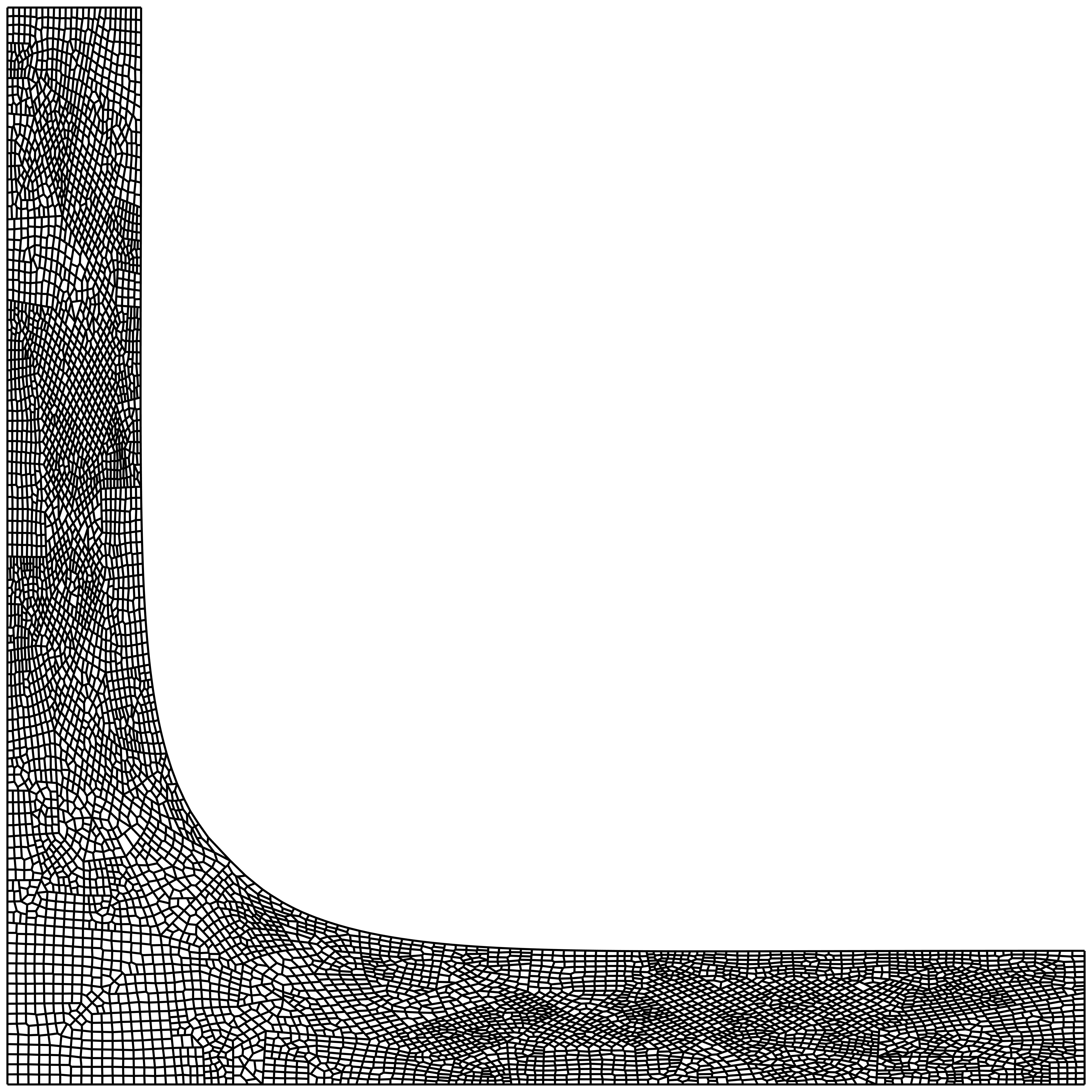}}
			\caption{Fine mesh: Initial mesh}
		\end{subfigure}
		\vskip \baselineskip 
		\begin{subfigure}[t]{0.33\textwidth}
			\centering
			\includegraphics[width=0.95\textwidth,height=0.95\textwidth]{{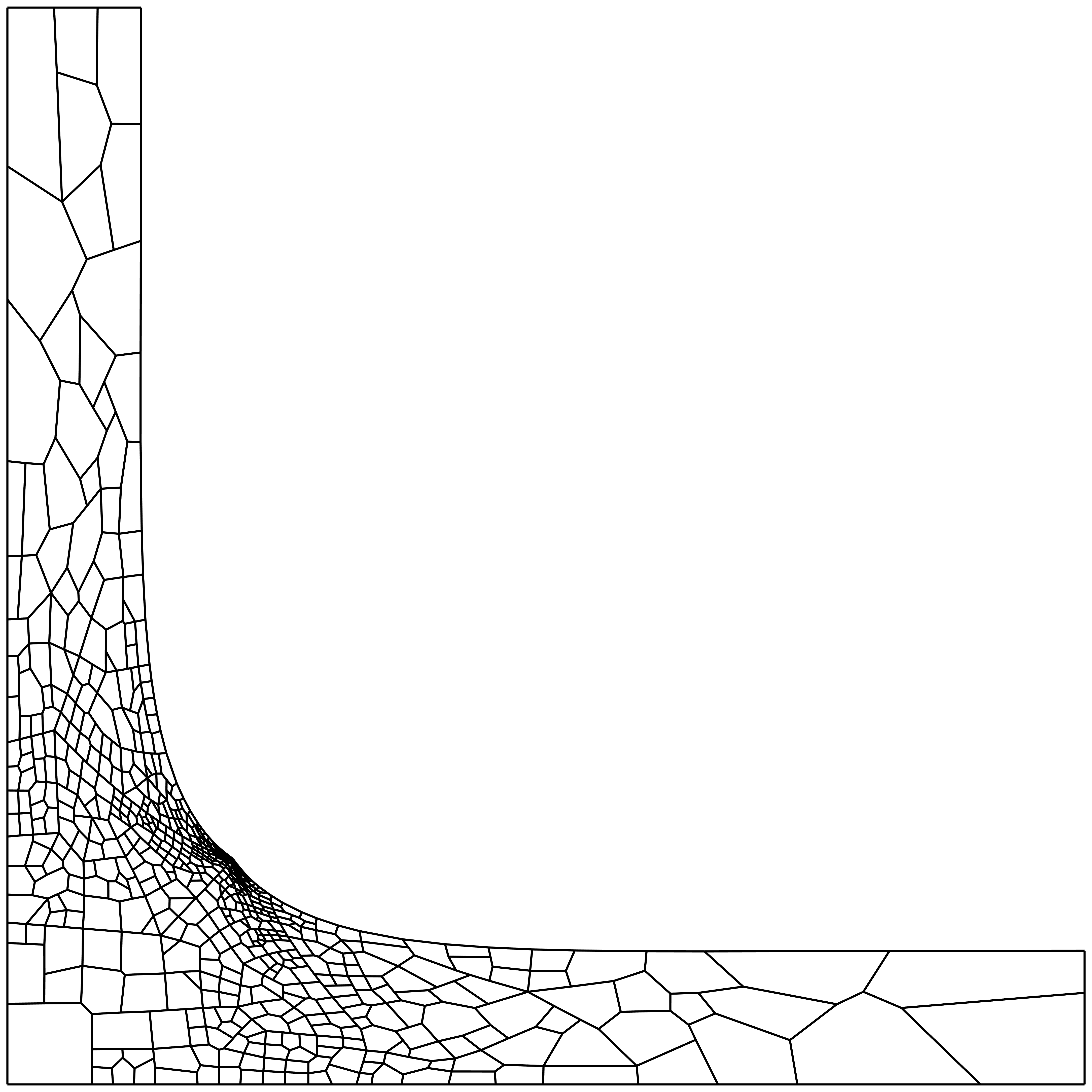}}
			\caption{Coarse mesh: Final mesh}
		\end{subfigure}%
		\begin{subfigure}[t]{0.33\textwidth}
			\centering
			\includegraphics[width=0.95\textwidth,height=0.95\textwidth]{{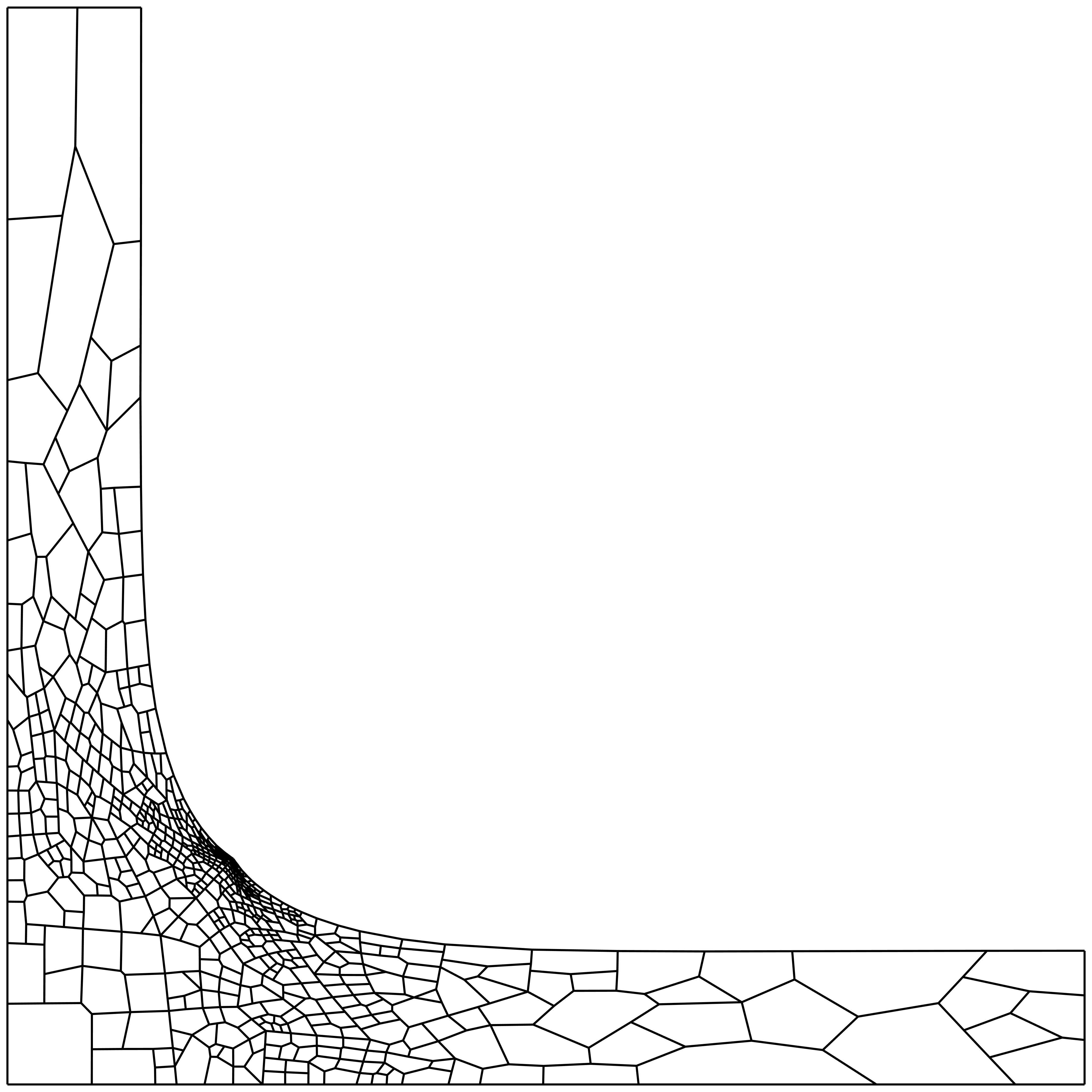}}
			\caption{Intermediate mesh: Final mesh}
		\end{subfigure}%
		\begin{subfigure}[t]{0.33\textwidth}
			\centering
			\includegraphics[width=0.95\textwidth,height=0.95\textwidth]{{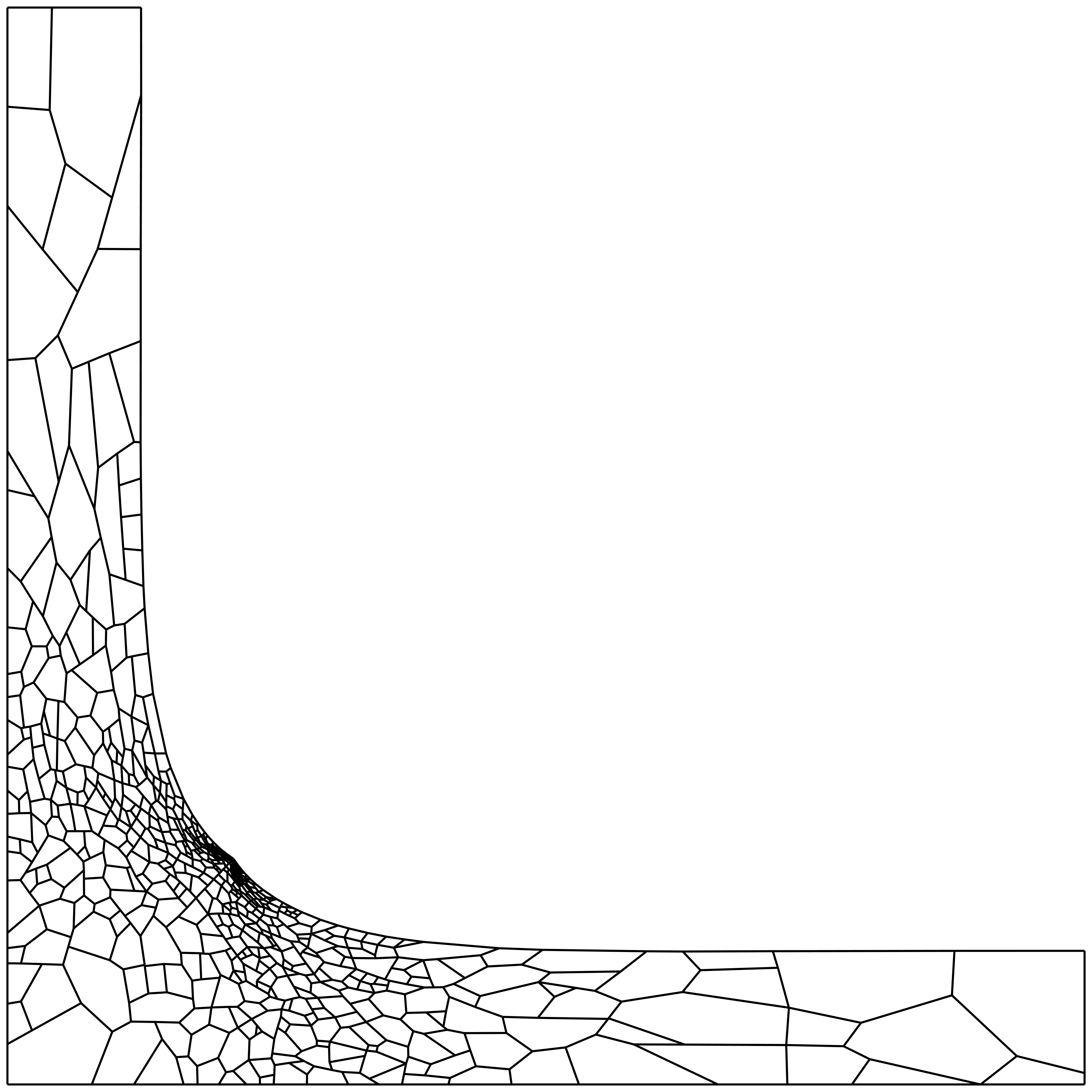}}
			\caption{Fine mesh: Final mesh}
		\end{subfigure}
		\caption{Initial and final adapted meshes using the fully adaptive remeshing procedure for the L-shaped domain problem from initial uniform Voronoi meshes of various densities with ${n_{\text{v}}^{\text{targ}} = 1000}$.
			\label{fig:LDomain_TargetNodes_MeshEvolution_VariousInitialMeshes}}
	\end{figure}

	The results of the fully adaptive remeshing process for the L-shaped domain problem are depicted in Figure~\ref{fig:LDomain_TargetNodes_MeshEvolution_VariousInitialMeshes} for initially uniform Voronoi meshes of varying density with a node target of ${n_{\text{v}}^{\text{targ}} = 1000}$. The top row of figures depicts the initial meshes while the bottom row depicts the final adapted meshes after the element target and termination criteria have been met.
	The final adapted meshes exhibit the same sensible and intuitive element distribution as observed in Figures~\ref{fig:LDomain_TargetElements_MeshEvolution} and \ref{fig:LDomain_TargetElements_MeshEvolution_VariousInitialMeshes}. The final adapted meshes are, again, almost identical for all initial uniform meshes considered. Thus, demonstrating that the final result of the fully adaptive remeshing procedure is independent of the initial mesh.
	
	\FloatBarrier
	\FloatBarrier
	
\begin{figure}[ht!]
		\centering
		\begin{subfigure}[t]{0.495\textwidth}
			\centering
			\includegraphics[width=0.95\textwidth]{{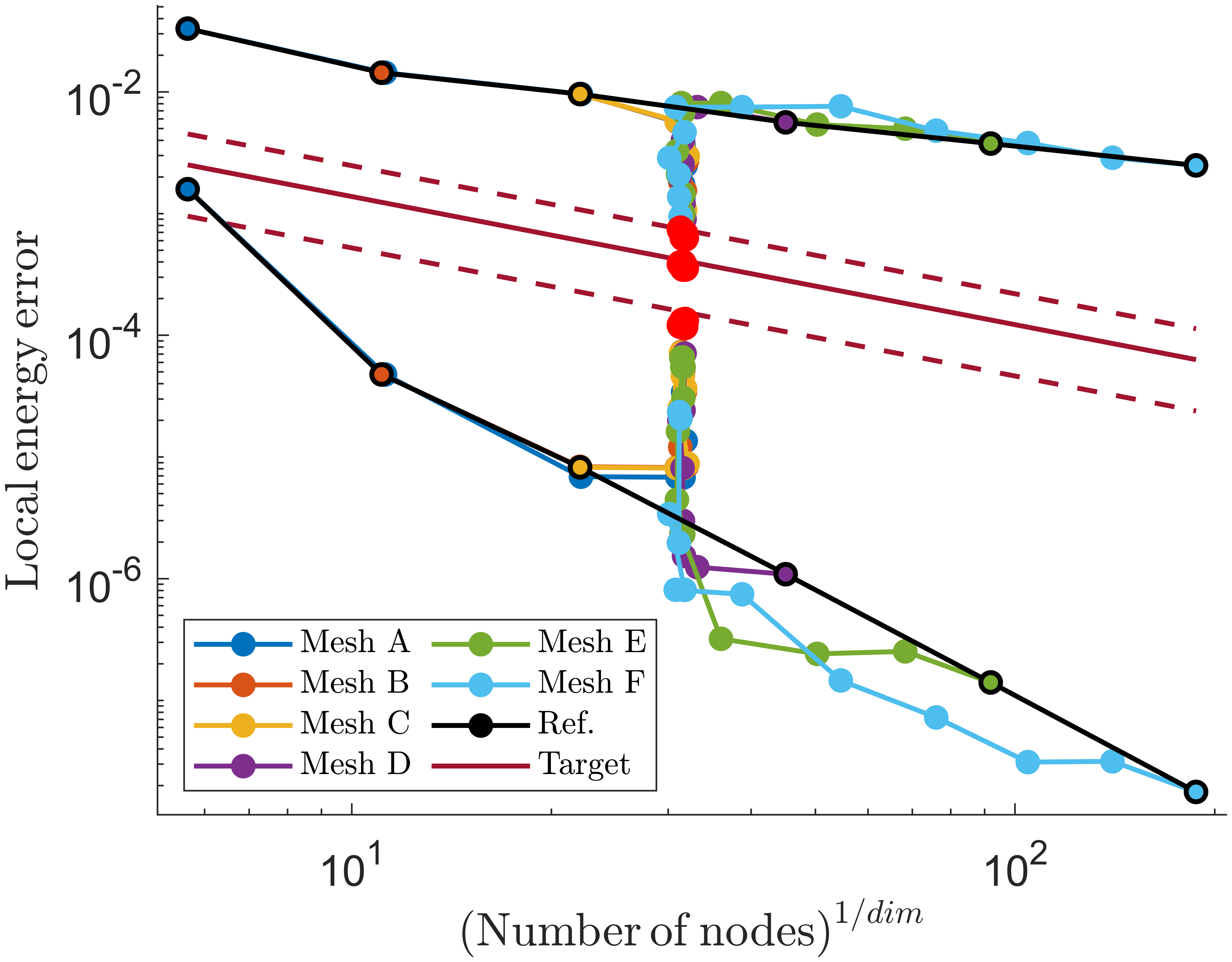}}
			\caption{$n_{\text{v}}^{\text{targ}} = 1000$ - Max and min local error}
			\vspace*{-3mm}
		\end{subfigure}%
		\begin{subfigure}[t]{0.495\textwidth}
			\centering
			\includegraphics[width=0.95\textwidth]{{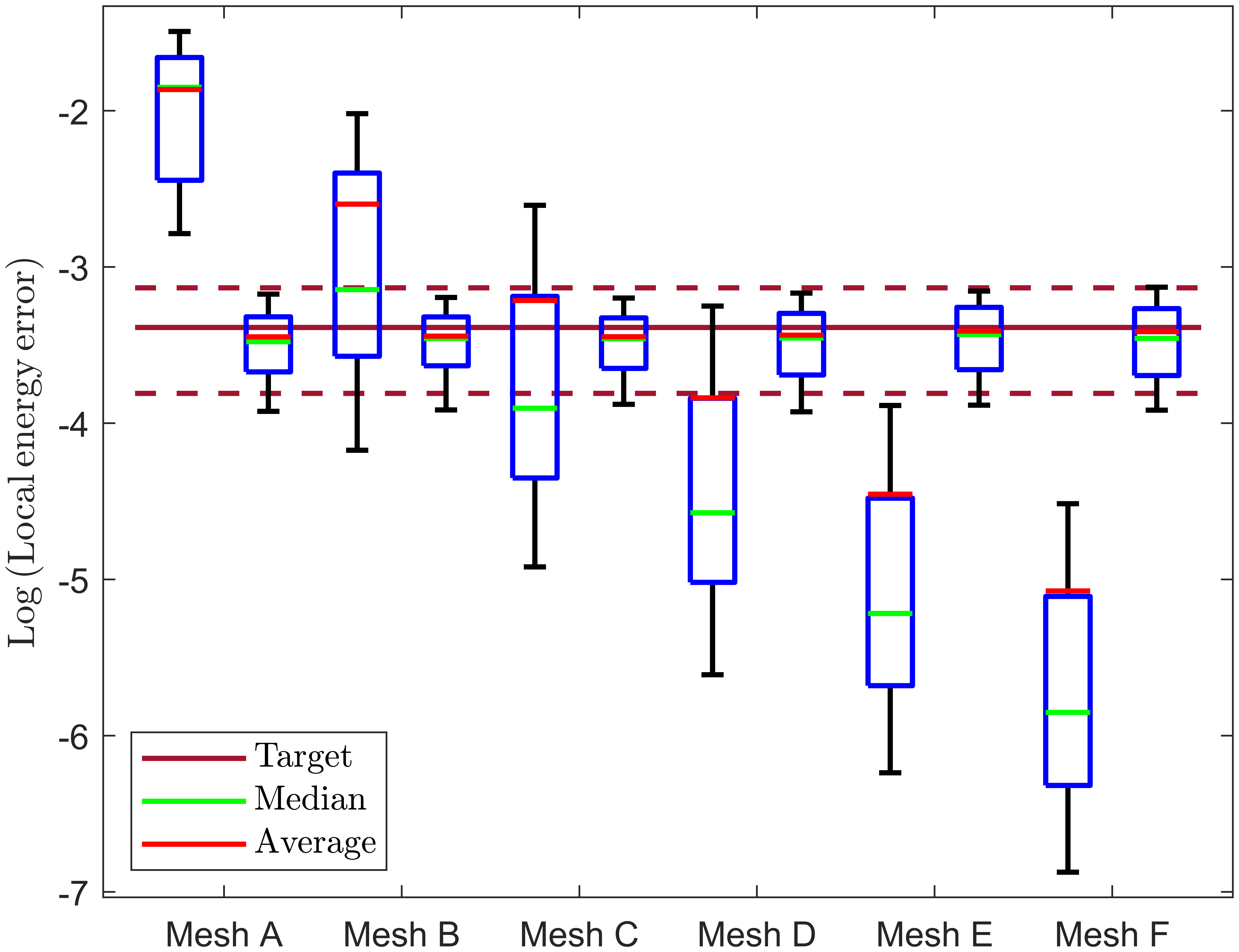}}
			\caption{$n_{\text{v}}^{\text{targ}} = 1000$ - Local error distribution}
			\vspace*{-3mm}
		\end{subfigure}
		\vskip \baselineskip 
		\begin{subfigure}[t]{0.495\textwidth}
			\centering
			\includegraphics[width=0.95\textwidth]{{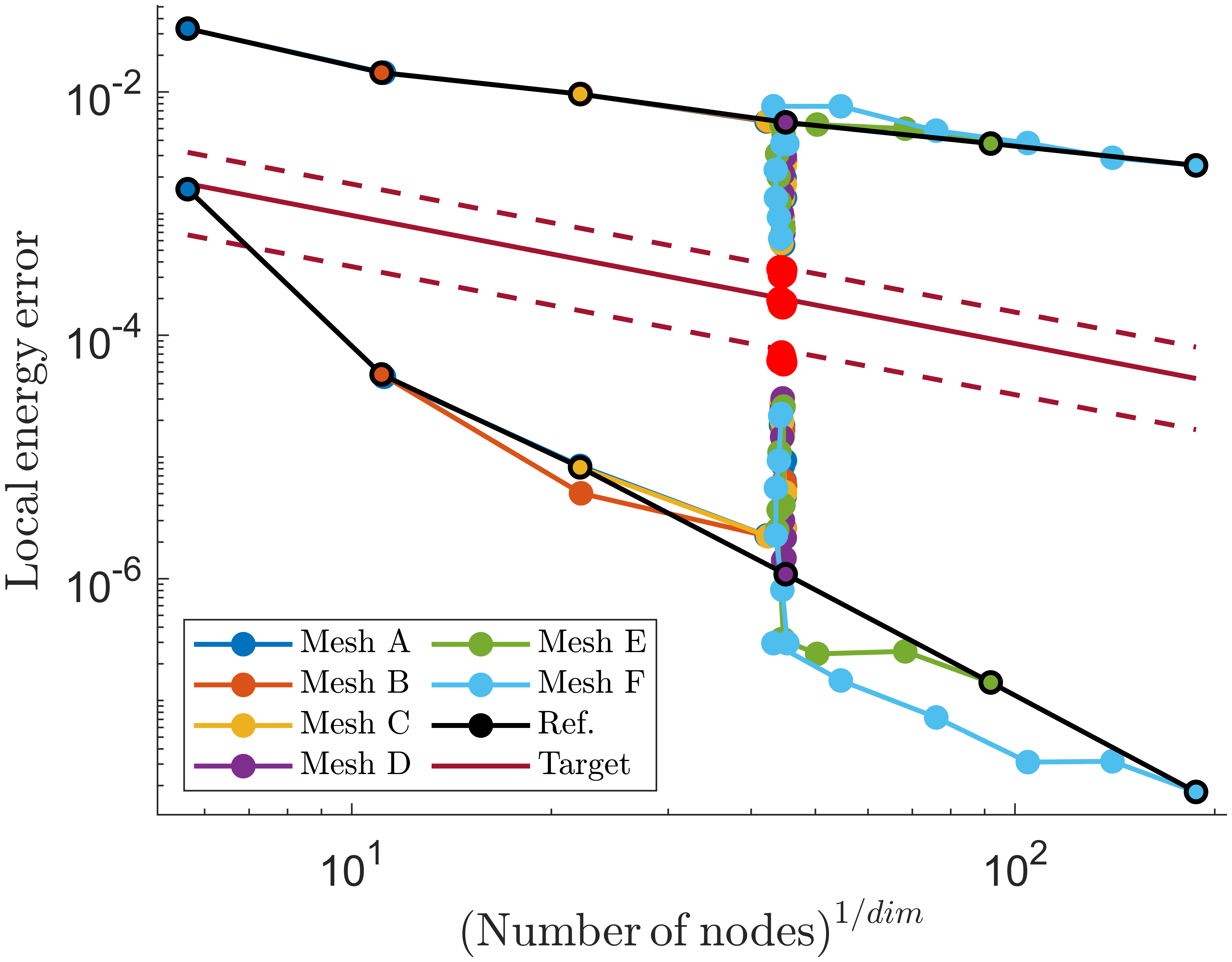}}
			\caption{$n_{\text{v}}^{\text{targ}} = 2000$ - Max and min local error}
			\vspace*{-3mm}
		\end{subfigure}%
		\begin{subfigure}[t]{0.495\textwidth}
			\centering
			\includegraphics[width=0.95\textwidth]{{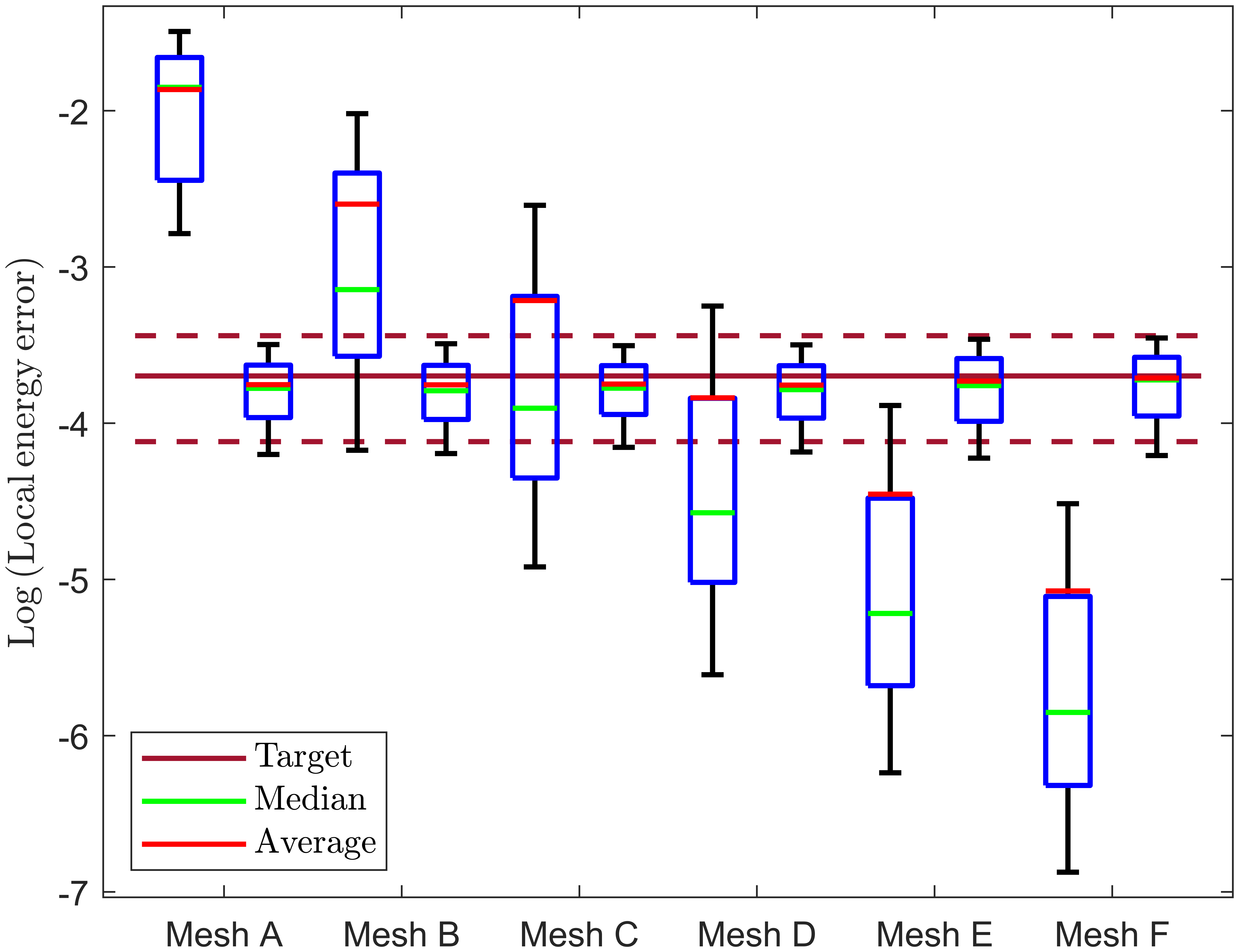}}
			\caption{$n_{\text{v}}^{\text{targ}} = 2000$ - Local error distribution}
			\vspace*{-3mm}
		\end{subfigure}
		\vskip \baselineskip 
		\begin{subfigure}[t]{0.495\textwidth}
			\centering
			\includegraphics[width=0.95\textwidth]{{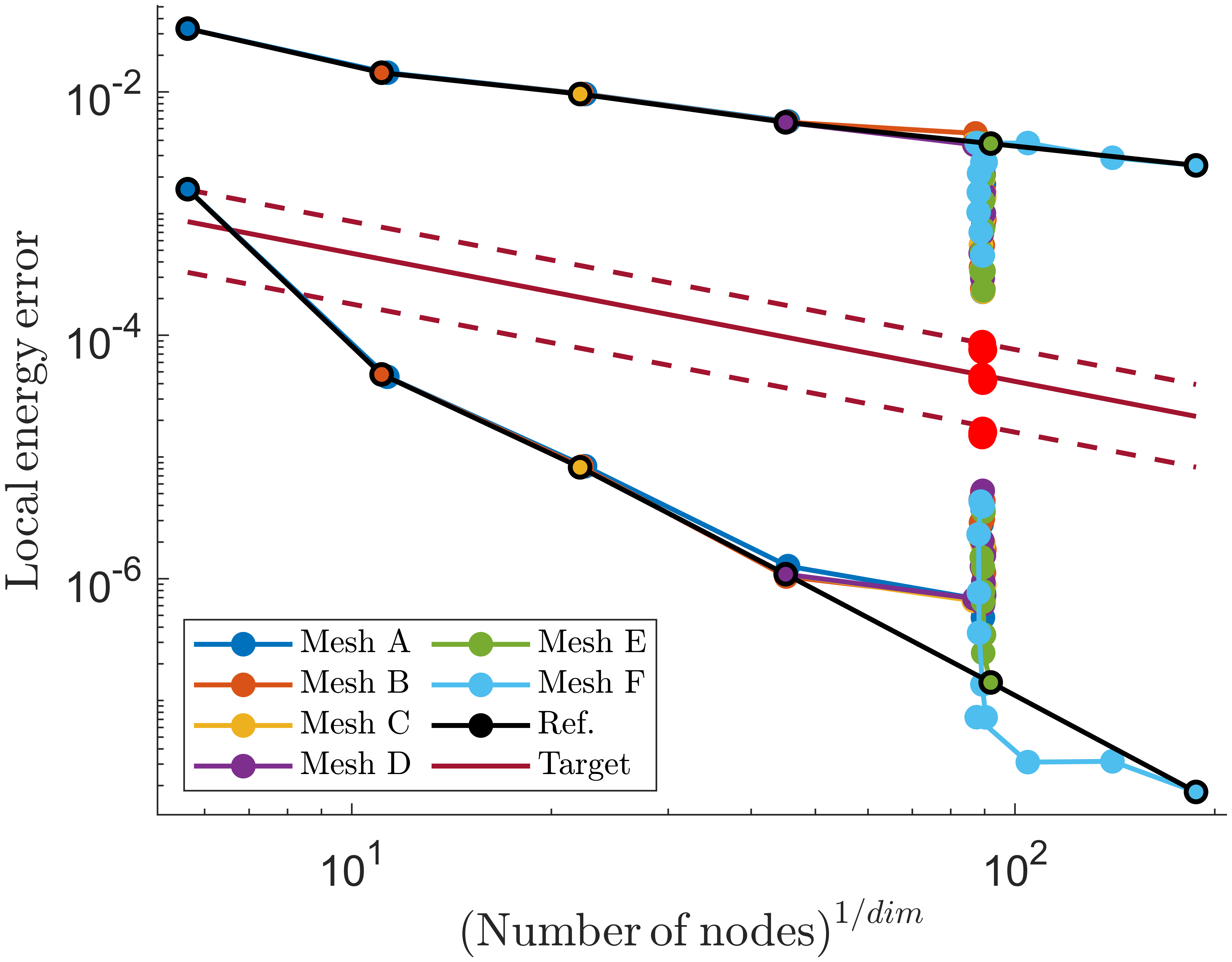}}
			\caption{$n_{\text{v}}^{\text{targ}} = 8000$ - Max and min local error}
			\vspace*{-3mm}
		\end{subfigure}%
		\begin{subfigure}[t]{0.495\textwidth}
			\centering
			\includegraphics[width=0.95\textwidth]{{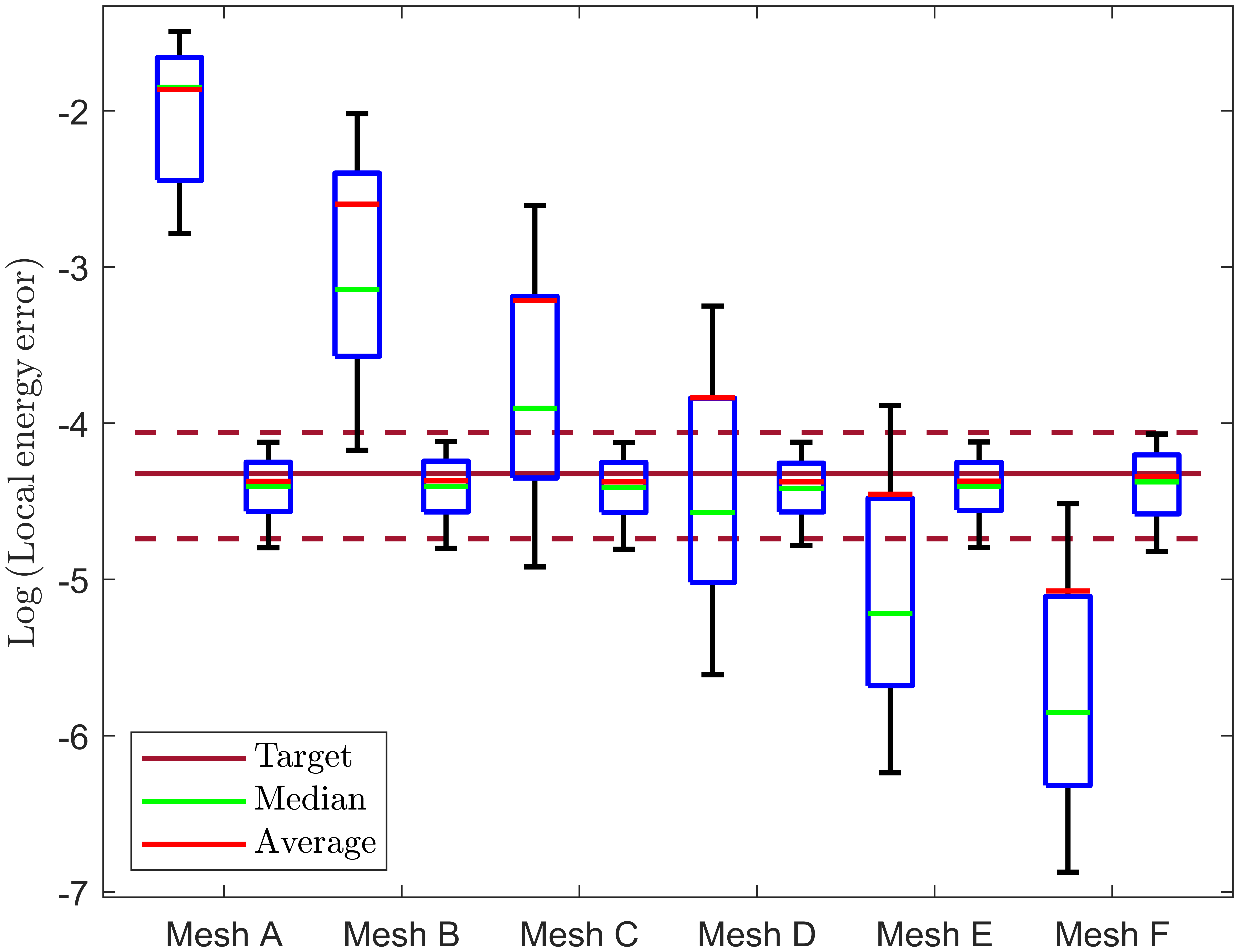}}
			\caption{$n_{\text{v}}^{\text{targ}} = 8000$ - Local error distribution}
			\vspace*{-3mm}
		\end{subfigure}
		\caption{Max and min energy error (left) and box and whisker plot of energy error distribution (right) for the L-shaped domain problem on Voronoi meshes of various initial densities for a range of node targets.
			\label{fig:LDomain_TargetNodes_ErrorDistribution_Voronoi}}
	\end{figure}

	The convergence behaviour of the energy error approximation vs the number of nodes in the mesh is depicted on a logarithmic scale in Figures~\ref{fig:LDomain_TargetNodes_ErrorDistribution_Voronoi}(a)-(c) for the L-shaped domain problem. The convergence behaviour is plotted for cases of several initially uniform Voronoi meshes of varying density for various node targets. 
	Additionally, the averaged final adapted mesh result (red markers) for all node targets considered is plotted in Figure~\ref{fig:LDomain_TargetNodes_ErrorDistribution_Voronoi}(d).
	Where applicable, the markers in Figure~\ref{fig:LDomain_TargetNodes_ErrorDistribution_Voronoi}(d) are colored to match their corresponding targets depicted in Figures~\ref{fig:LDomain_TargetNodes_ErrorDistribution_Voronoi}(a)-(c).
	From Figures~\ref{fig:LDomain_TargetNodes_ErrorDistribution_Voronoi}(a)-(c) it is clear that the fully adaptive procedure is able to meet the specified node targets from any initial mesh. Furthermore, the final adapted meshes have an almost identical error for a specific element target. Thus, demonstrating that the performance of the fully adaptive remeshing procedure is independent of the initial mesh.
	From Figure~\ref{fig:LDomain_TargetNodes_ErrorDistribution_Voronoi}(d) it is clear that the outputs of the adaptive procedure for various node targets exhibits a linear convergence rate which, again, is expected to be the (approximately) optimal convergence rate for this problem.
	
	\FloatBarrier
	\begin{figure}[ht!]
		\centering
		\begin{subfigure}[t]{0.495\textwidth}
			\centering
			\includegraphics[width=0.95\textwidth]{{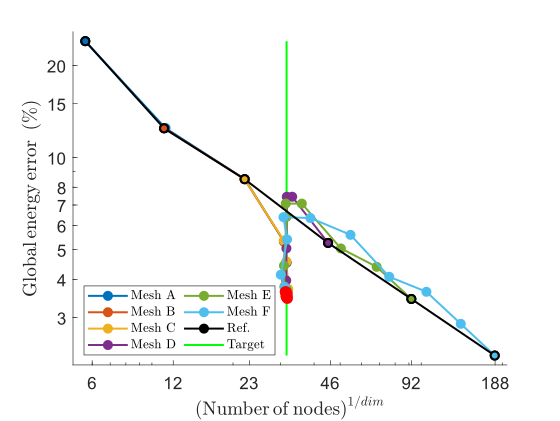}}
			\caption{$n_{\text{v}}^{\text{targ}} = 1000$}
		\end{subfigure}%
		\begin{subfigure}[t]{0.495\textwidth}
			\centering
			\includegraphics[width=0.95\textwidth]{{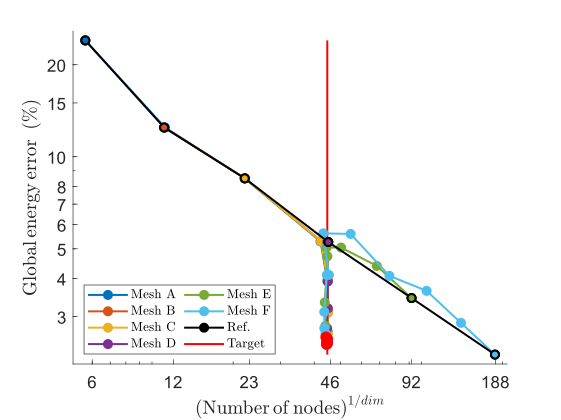}}
			\caption{$n_{\text{v}}^{\text{targ}} = 2000$}
		\end{subfigure}
		\vskip \baselineskip 
		\begin{subfigure}[t]{0.495\textwidth}
			\centering
			\includegraphics[width=0.95\textwidth]{{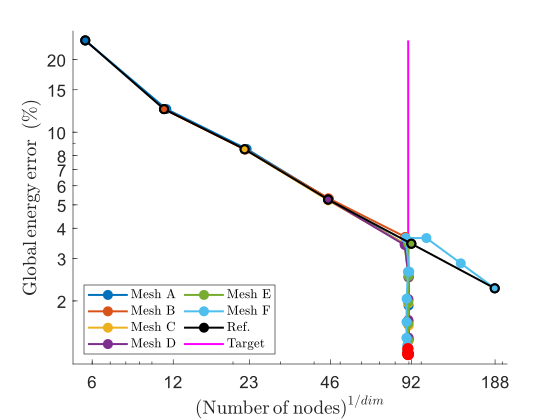}}
			\caption{$n_{\text{v}}^{\text{targ}} = 8000$}
		\end{subfigure}%
		\begin{subfigure}[t]{0.495\textwidth}
			\centering
			\includegraphics[width=0.95\textwidth]{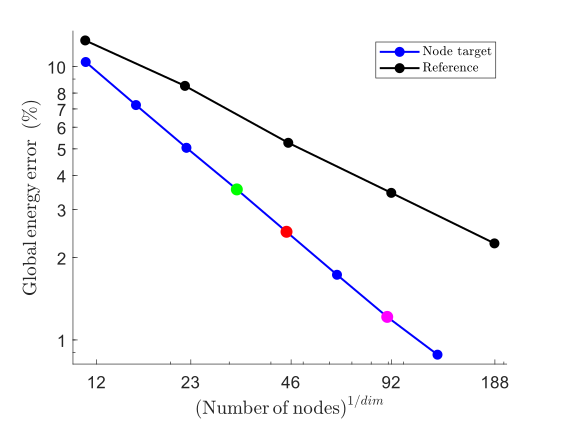}
			\caption{Various node targets}
		\end{subfigure}
		\caption{Energy error vs number of nodes for the L-shaped domain problem on Voronoi meshes of various initial densities with (a) ${n_{\text{v}}^{\text{targ}} = 1000}$, (b) ${n_{\text{v}}^{\text{targ}} = 2000}$, (c) ${n_{\text{v}}^{\text{targ}} = 8000}$, and (d) average final/adapted mesh error for various node targets.  
			\label{fig:LDomain_TargetNodes_ErrorConvergence_Voronoi}}
	\end{figure} 
	\FloatBarrier
	
	The error distribution during the mesh adaptation process for the L-shaped domain problem is depicted in Figure~\ref{fig:LDomain_TargetNodes_ErrorDistribution_Voronoi} for various node targets on Voronoi meshes.
	
	The nature of the convergence exhibited by the left column of figures is similar to that exhibited in Figure~\ref{fig:LDomain_TargetElements_ErrorDistribution_Structured} for element targets, albeit slightly more erratic due to the nature of Voronoi meshes.
	Additionally, the error distributions of the final adapted meshes, as indicated by the red markers and the box and whisker plots, are qualitatively very similar to those of Figure~\ref{fig:LDomain_TargetError_ErrorConvergence_Voronoi} for error targets and Figure~\ref{fig:LDomain_TargetElements_ErrorDistribution_Structured} for element targets.
	Specifically, the maximum errors fall within the upper bounds and the minimum errors are satisfactorily close to the lower bounds. Furthermore, the upper and lower quartiles indicate a narrow distribution of error around the average, and the average error is almost identical to the target error. Additionally, the average and median errors are almost identical which, again, emphasises the narrow distribution of the element-level errors.
	
	The results presented in Figure~\ref{fig:LDomain_TargetNodes_ErrorConvergence_Voronoi} demonstrated that the fully adaptive procedure was able to meet all specified node targets on Voronoi meshes. This is indicated by the red markers denoting the final adapted mesh lying exactly on the target node line. 
	The results presented in Figure~\ref{fig:LDomain_TargetNodes_ErrorDistribution_Voronoi} for Voronoi meshes demonstrated that the average element-level error almost exactly met the element-level target as the red markers denoting the average error strongly overlap the solid maroon target line. Furthermore, the element-level errors were satisfactorily equal as they (almost) all fell within the specified target error range.
	Thus, the fully adaptive procedure successfully generated quasi-optimal meshes for the specified node targets on Voronoi meshes.
	
	\FloatBarrier
	\FloatBarrier
	
	\subsubsection{Comparison of targets}
	\label{subsubsec:L_Domain_TargetComparison}
	
\begin{figure}[ht!]
		\centering
		\begin{subfigure}[t]{0.495\textwidth}
			\centering
			\includegraphics[width=0.95\textwidth]{{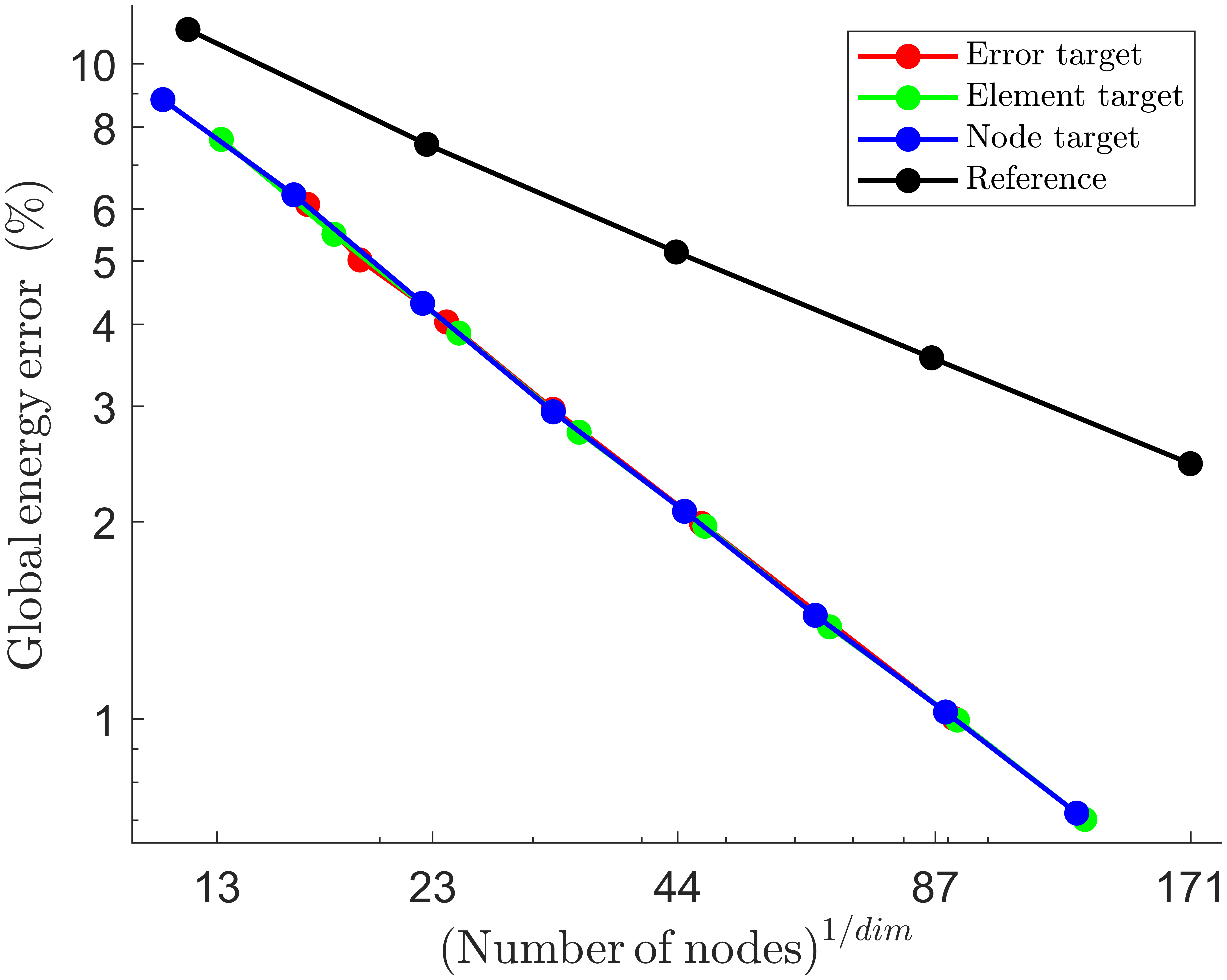}}
			\caption{Structured meshes}
		\end{subfigure}%
		\begin{subfigure}[t]{0.495\textwidth}
			\centering
			\includegraphics[width=0.95\textwidth]{{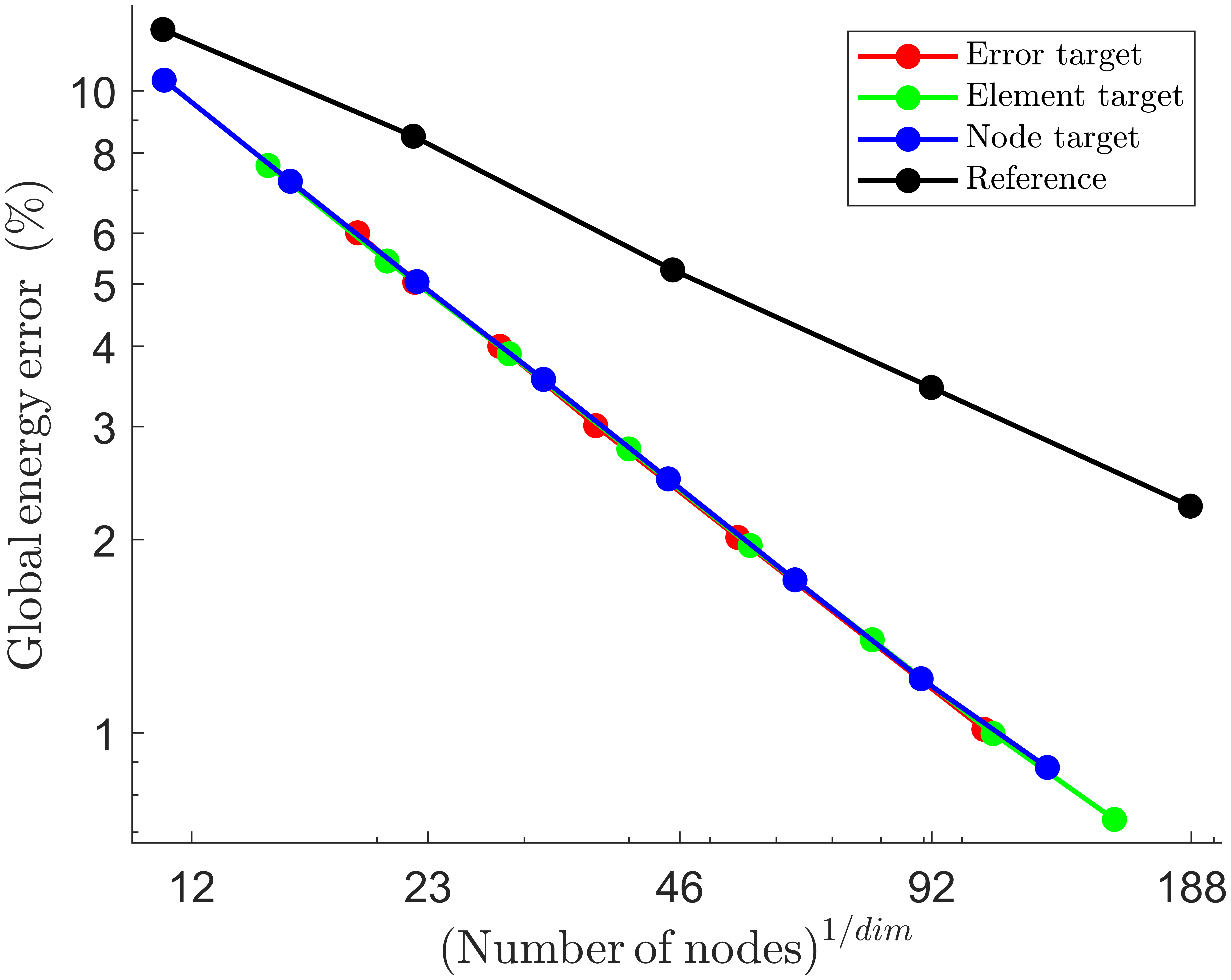}}
			\caption{Voronoi meshes}
		\end{subfigure}
		\caption{Energy error vs number of nodes for the L-shaped domain problem for various target types on (a) structured and (b) Voronoi meshes. 
			\label{fig:LDomain_TargetComparison}}
	\end{figure}

	The energy error approximation of the averaged final adapted mesh result for various remeshing target types is plotted against the number of nodes in Figure~\ref{fig:LDomain_TargetComparison} on a logarithmic scale for the L-shaped domain problem for (a) structured and (b) Voronoi meshes.
	From these results it is clear that the fully adaptive procedure is equally effective for error-, node-, or element-based targets on structured or Voronoi meshes.
	That is, the fully adaptive procedure can meet any of the prescribed targets while generating a quasi-optimal mesh.
	
	\FloatBarrier
	\FloatBarrier

	\subsection{Pseudo-dynamic punch}
	\label{subsec:DynamicPunch}
	In this section a pseudo-dynamic problem is presented to demonstrate the suitability of the proposed adaptive remeshing procedure for dynamic problems. 
	Adaptive remeshing procedures are particularly beneficial for dynamic problems because the zones of high stresses are continuously changing. An adaptive procedure allows for a high degree of refinement in a zone of high stress that can then be coarsened once the high stress has passed. This allows for a high degree of accuracy at a low computational cost.
	In this work, the remeshing requirements of a dynamic problem are mimicked by changing the boundary conditions of a static problem. 
	
\begin{figure}[ht!]
		\centering
		\begin{subfigure}[t]{0.5\textwidth}
			\centering
			\includegraphics[width=0.9\textwidth]{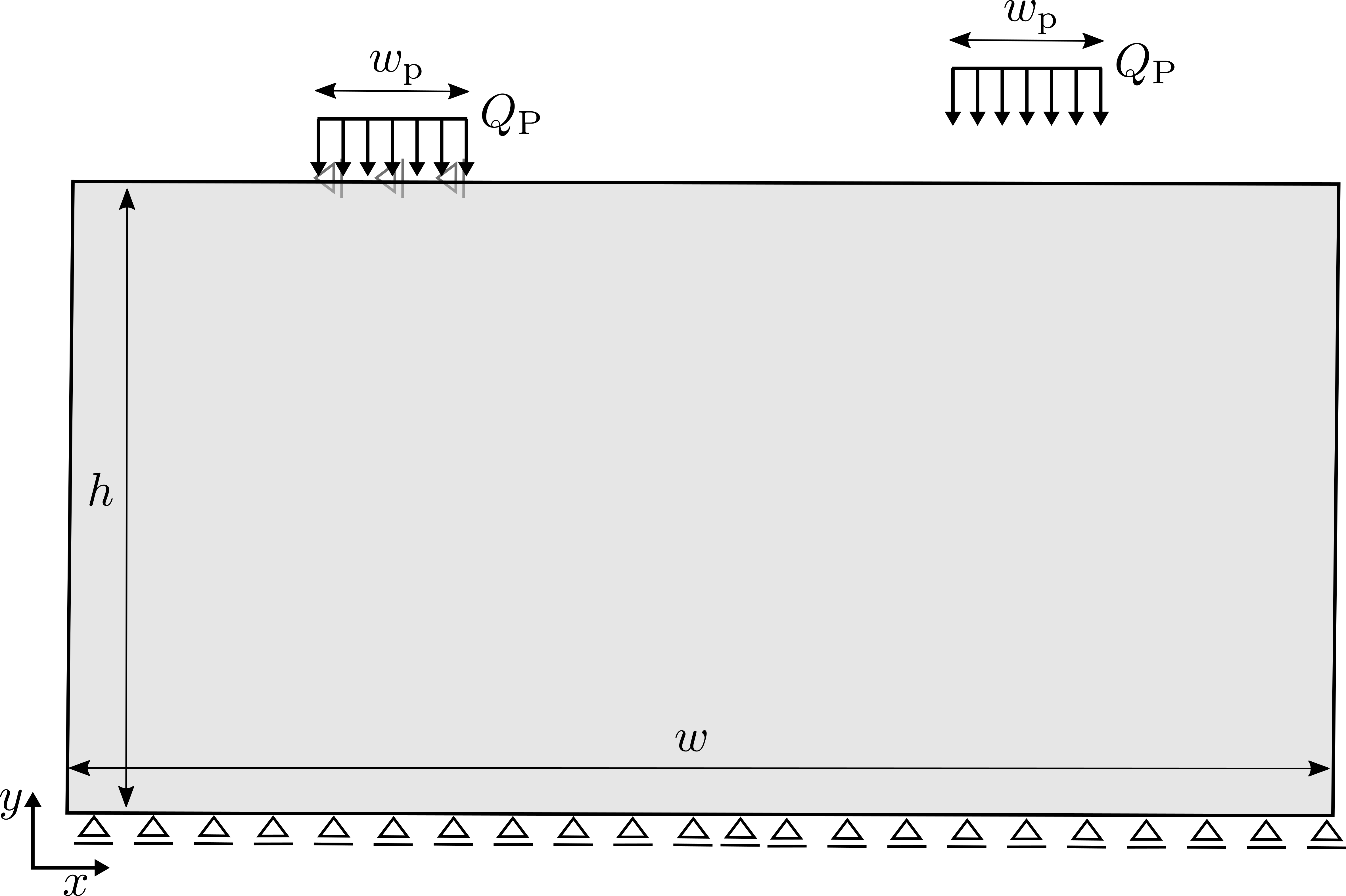}
			\caption{Problem geometry - Odd cycles}
		\end{subfigure}%
		\begin{subfigure}[t]{0.5\textwidth}
			\centering
			\includegraphics[width=0.9\textwidth]{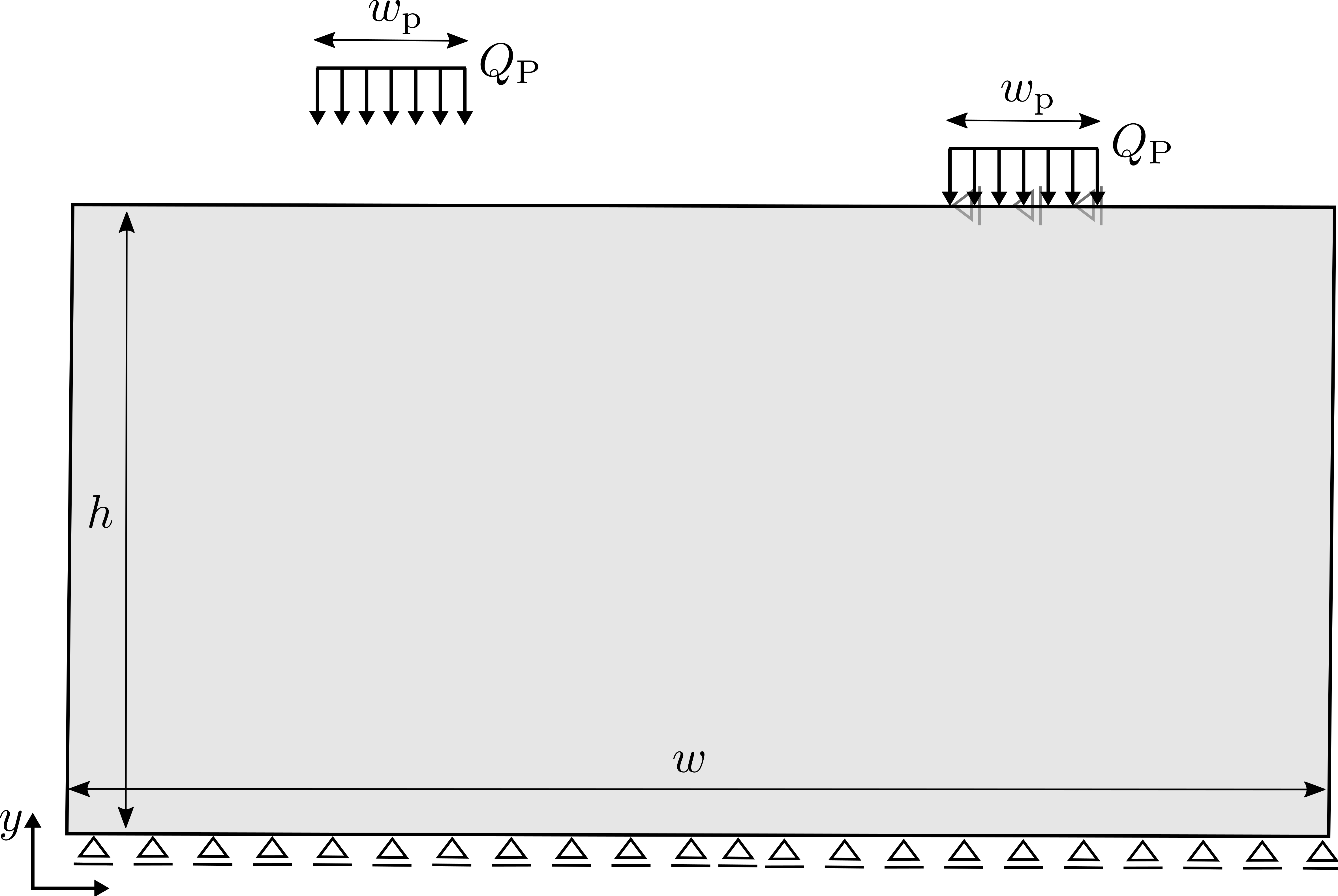}
			\caption{Problem geometry - Even cycles}
		\end{subfigure}
		\vskip \baselineskip 
		\vspace*{-5mm}
		\begin{subfigure}[t]{0.5\textwidth}
			\centering
			\includegraphics[width=0.85\textwidth]{{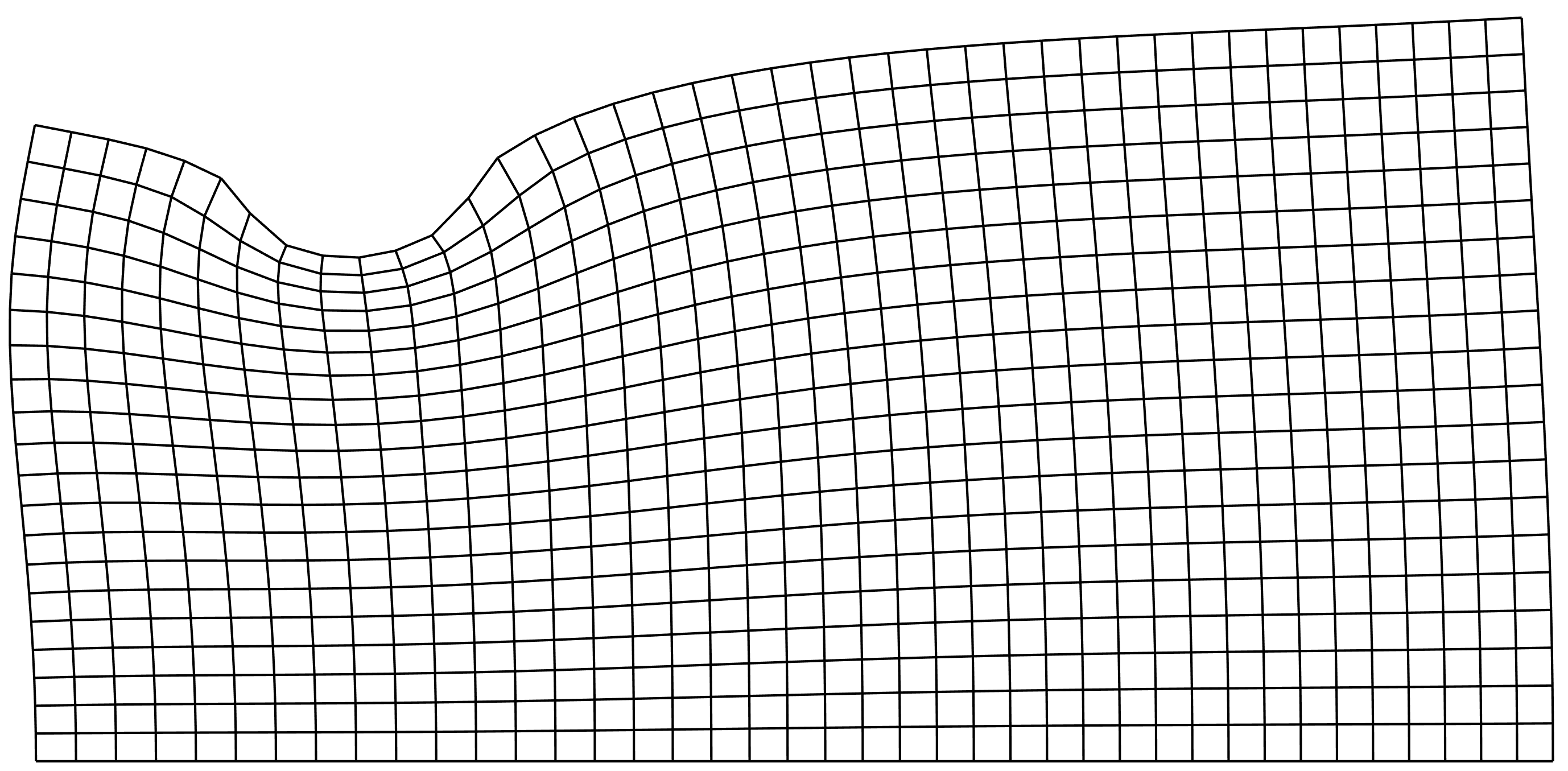}}
			\caption{Deformed configuration - Odd cycles}
		\end{subfigure}%
		\begin{subfigure}[t]{0.5\textwidth}
			\centering
			\includegraphics[width=0.85\textwidth]{{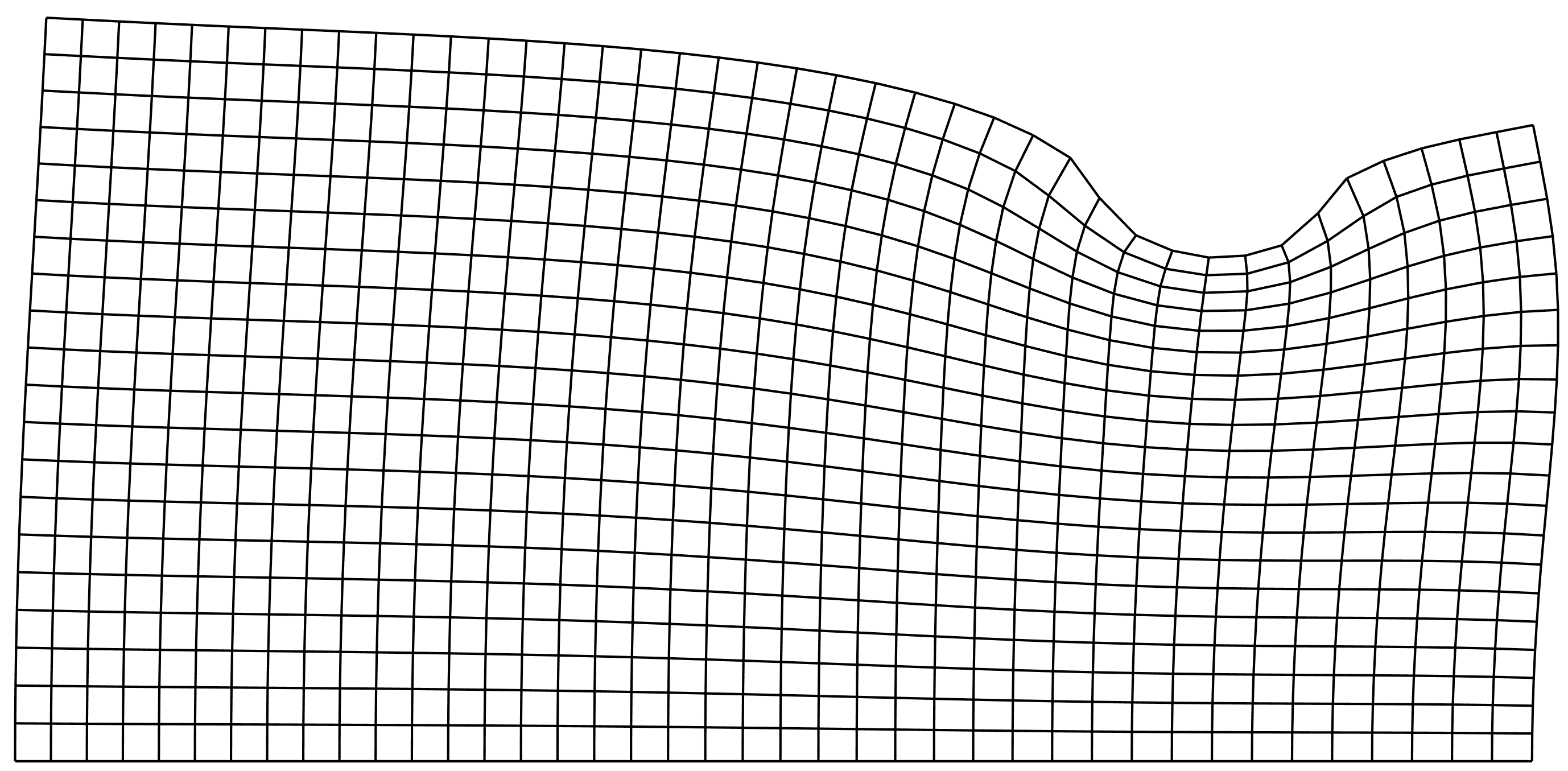}}
			\caption{Deformed configuration - Even cycles}
		\end{subfigure}
		\caption{Pseudo-dynamic punch problem geometry for (a) odd and (b) even numbered cycles and deformed configuration for (c) odd and (d) even numbered cycles.
			\label{fig:DynamicPunchGeometry}}
	\end{figure}

	The pseudo-dynamic punch problem comprises a rectangular body of width ${w=2~\text{m}}$ and height ${h=2~\text{m}}$ vertically constrained along its bottom edge (see Figure~\ref{fig:DynamicPunchGeometry}). The body is subjected to alternating punches of width ${w_{\text{p}}=0.3~\text{m}}$ modelled as uniformly distributed loads with a magnitude of ${Q_{\text{p}}}=0.675~\frac{\text{N}}{\text{m}}$. The centres of the punches are ${1.1~\text{m}}$ apart and ${0.45~\text{m}}$ from the left- and right-hand edges of the body. During application of the punch the region of the body experiencing the distributed load is horizontally constrained. 
	During odd numbered load cycles only the left-hand punch is active and during even numbered load cycles only the right-hand punch is active as illustrated in Figures~\ref{fig:DynamicPunchGeometry}(a) and (b) respectively. Figures~\ref{fig:DynamicPunchGeometry}(c) and (d) respectively depict the deformed configuration of the body for odd and even load cycles.
	
	\FloatBarrier
	\FloatBarrier
	
\begin{figure}[ht!]
		\centering
		\begin{subfigure}[t]{0.33\textwidth}
			\centering
			\includegraphics[width=0.95\textwidth]{{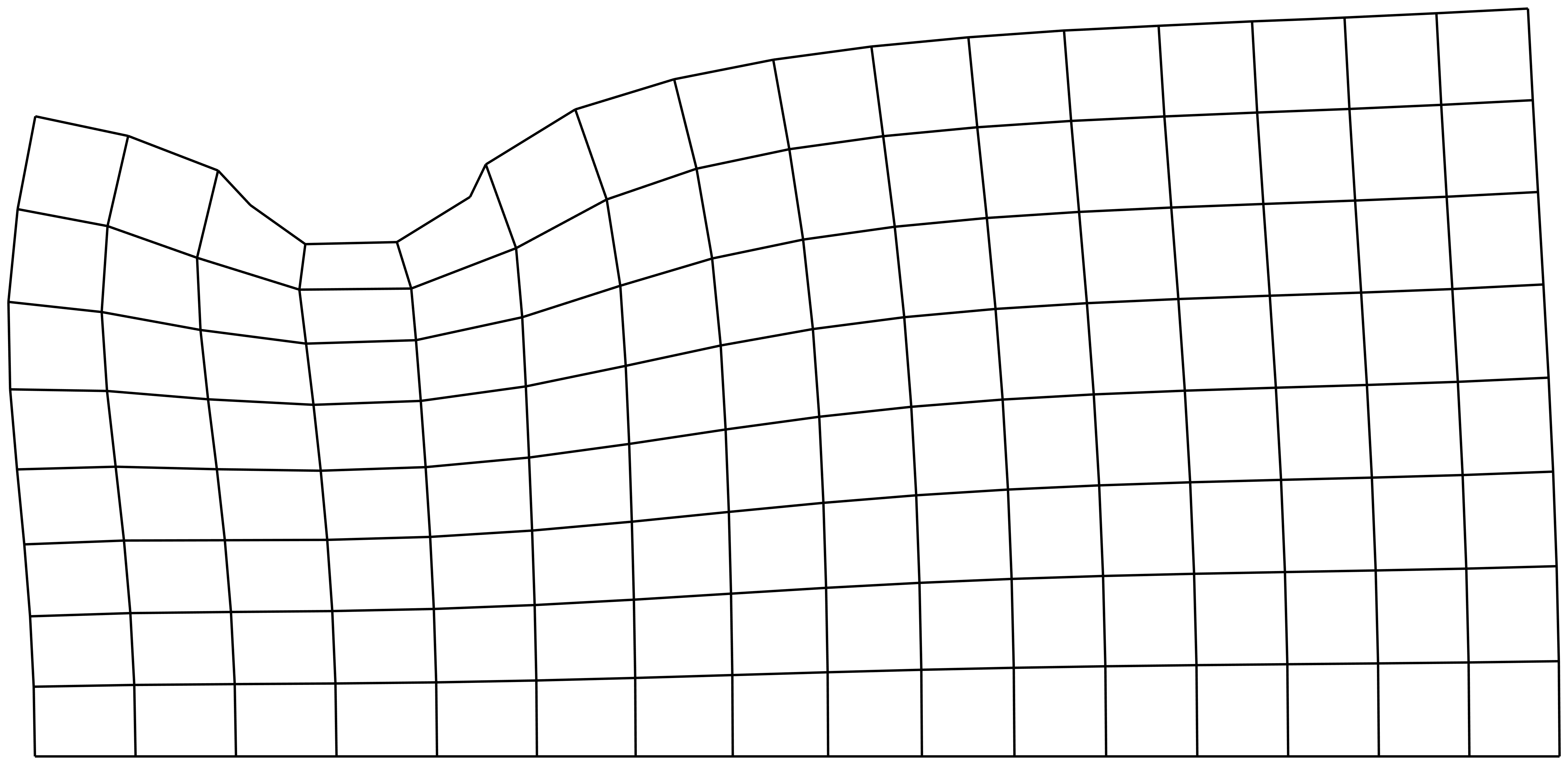}}
			\caption{Cycle 1 - Step 1 (Initial)}
		\end{subfigure}%
		\begin{subfigure}[t]{0.33\textwidth}
			\centering
			\includegraphics[width=0.95\textwidth]{{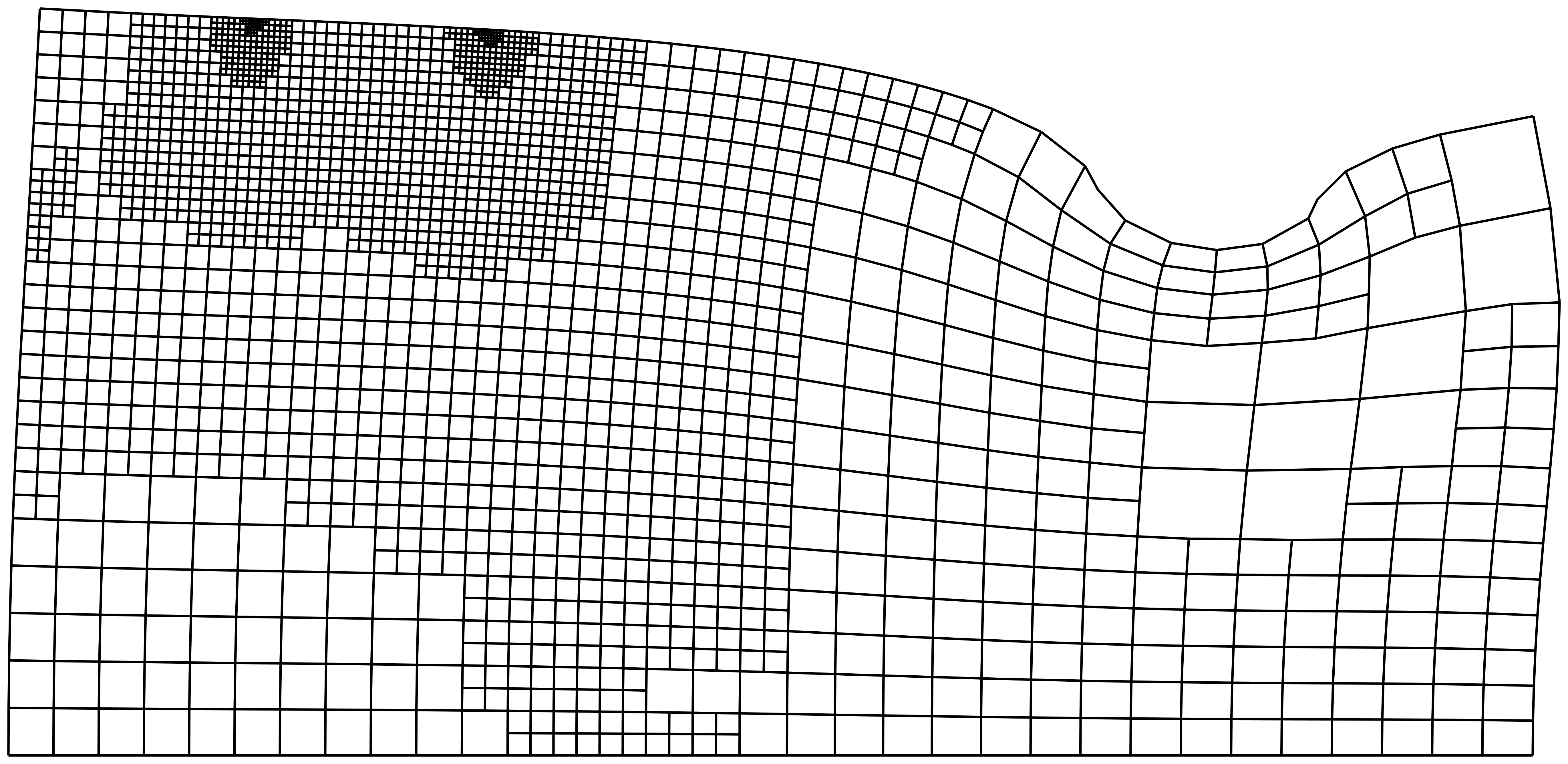}}
			\caption{Cycle 2 - Step 1 (Initial)}
		\end{subfigure}%
		\begin{subfigure}[t]{0.33\textwidth}
			\centering
			\includegraphics[width=0.95\textwidth]{{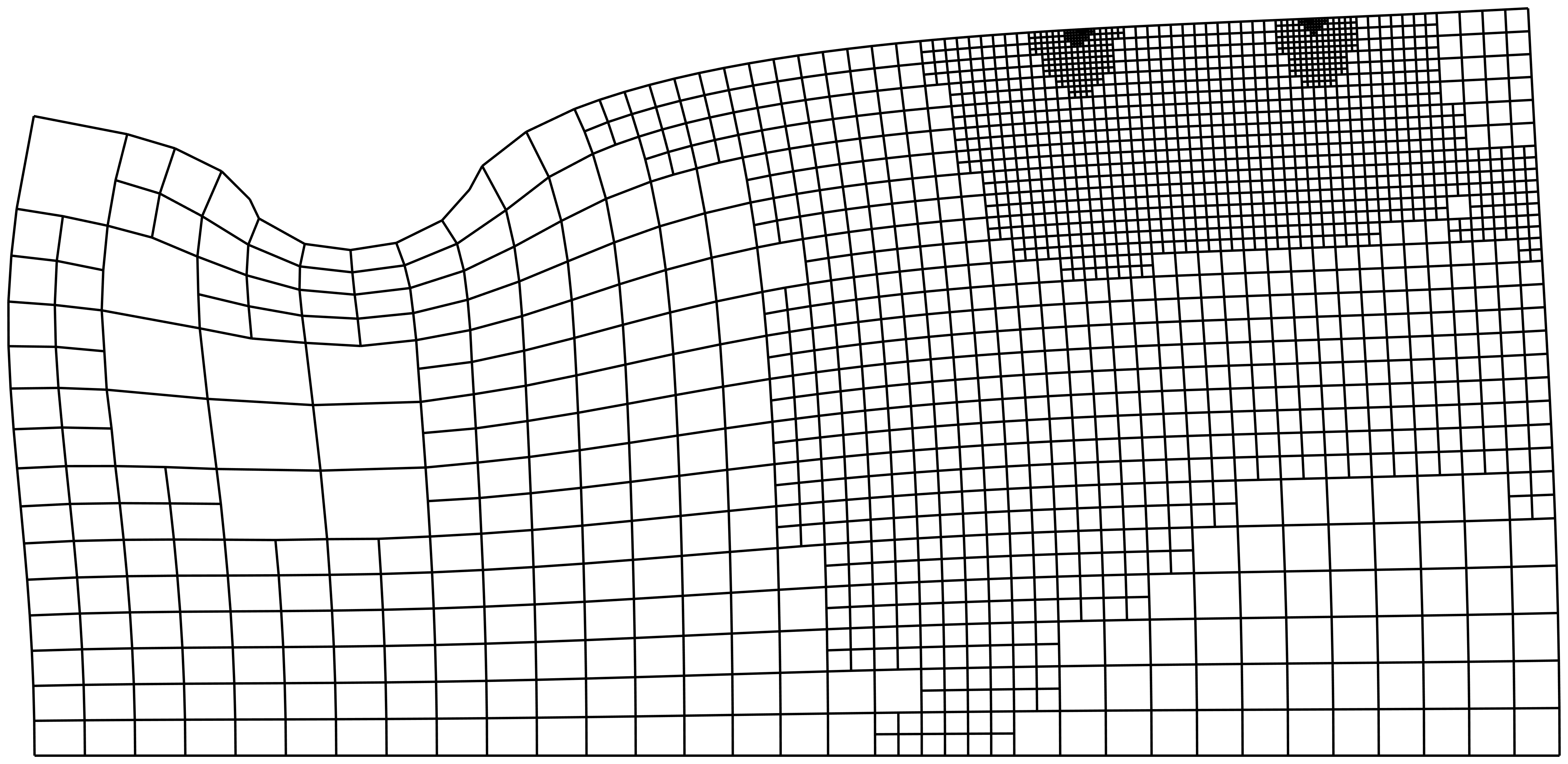}}
			\caption{Cycle 3 - Step 1 (Initial)}
		\end{subfigure}
		\vskip \baselineskip 
		\vspace*{-4mm}
		\begin{subfigure}[t]{0.33\textwidth}
			\centering
			\includegraphics[width=0.95\textwidth]{{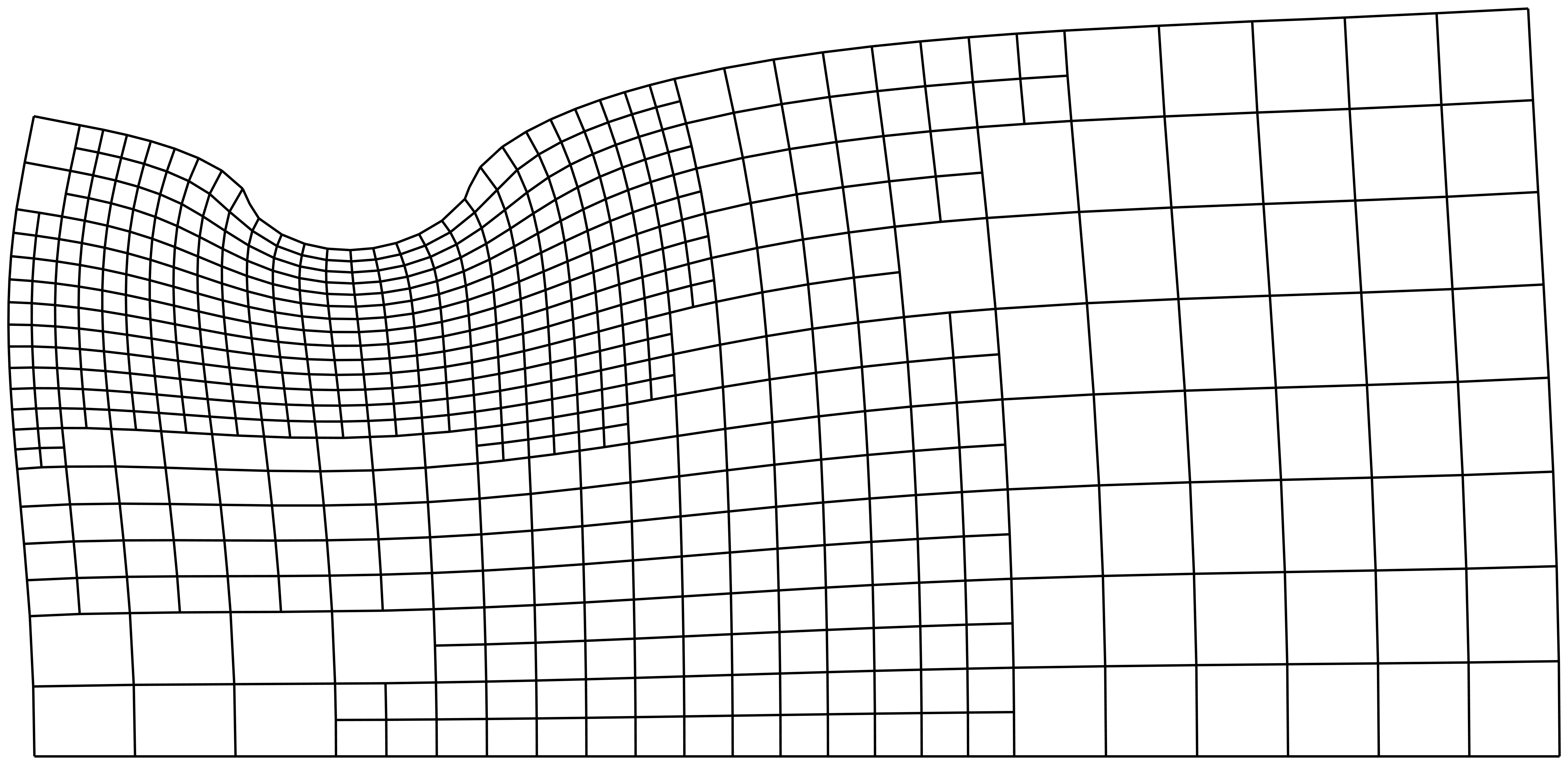}}
			\caption{Cycle 1 - Step 3}
		\end{subfigure}%
		\begin{subfigure}[t]{0.33\textwidth}
			\centering
			\includegraphics[width=0.95\textwidth]{{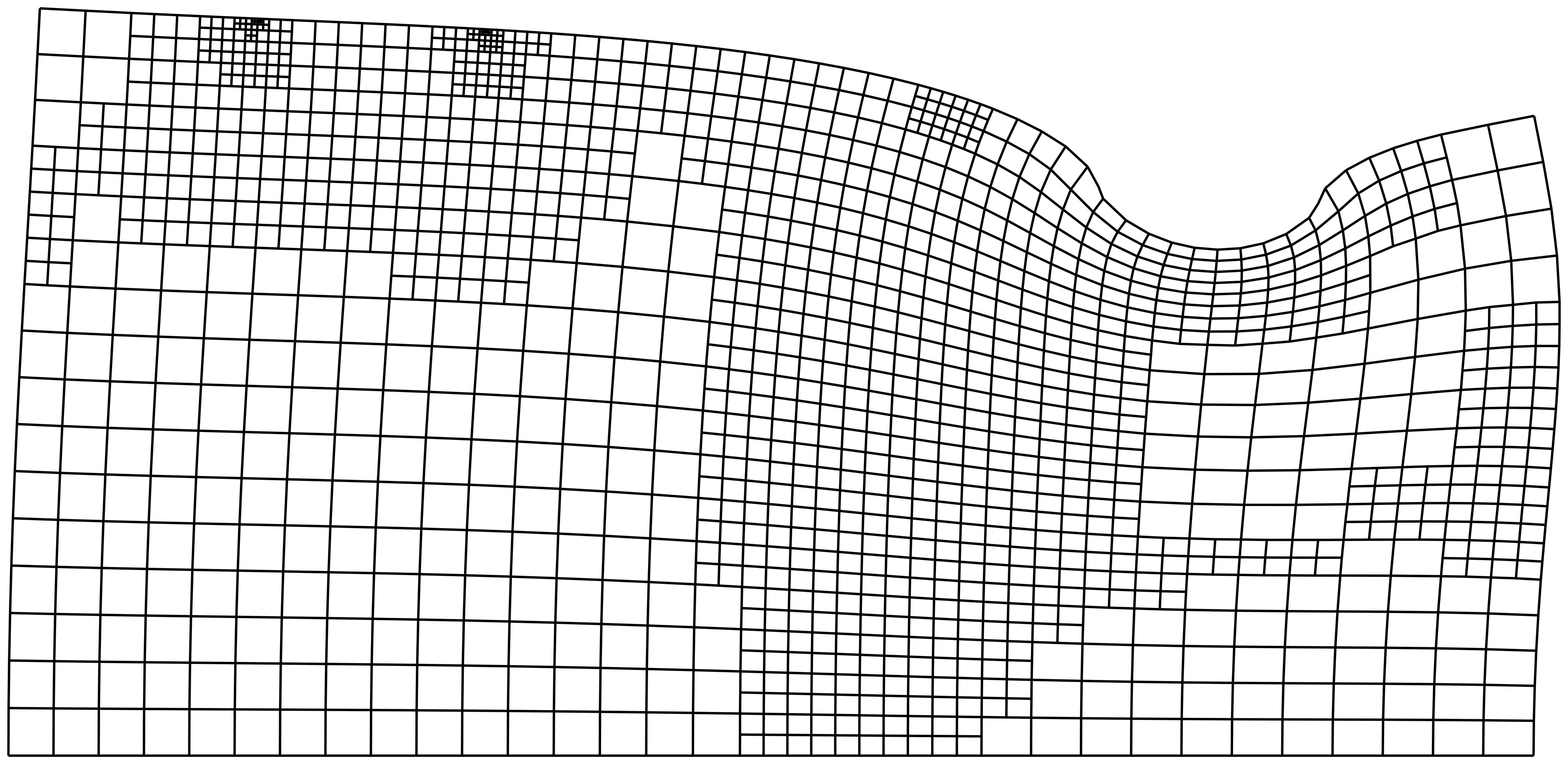}}
			\caption{Cycle 2 - Step 2}
		\end{subfigure}%
		\begin{subfigure}[t]{0.33\textwidth}
			\centering
			\includegraphics[width=0.95\textwidth]{{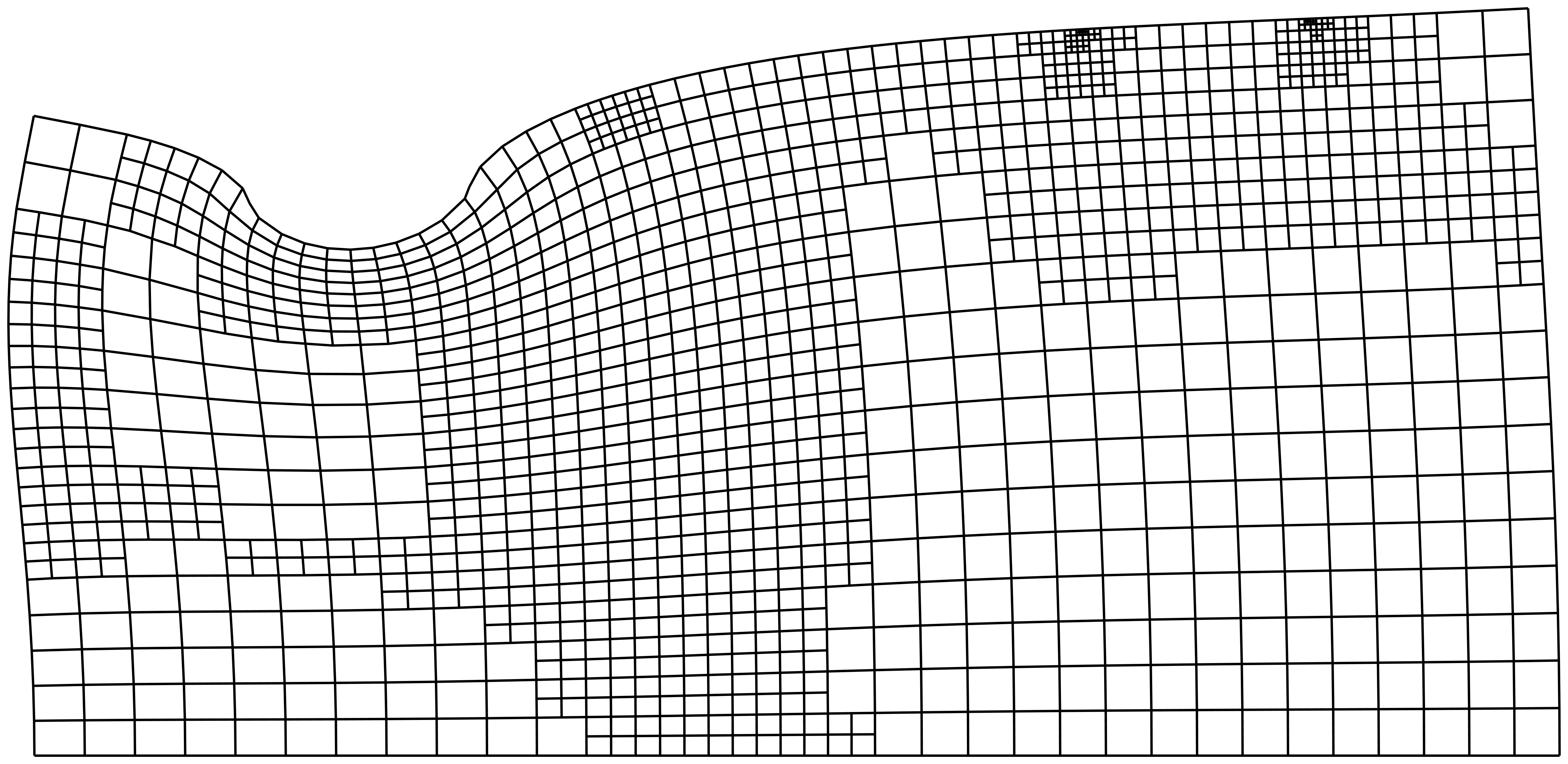}}
			\caption{Cycle 3 - Step 2}
		\end{subfigure}
		\vskip \baselineskip 
		\vspace*{-4mm}
		\begin{subfigure}[t]{0.33\textwidth}
			\centering
			\includegraphics[width=0.95\textwidth]{{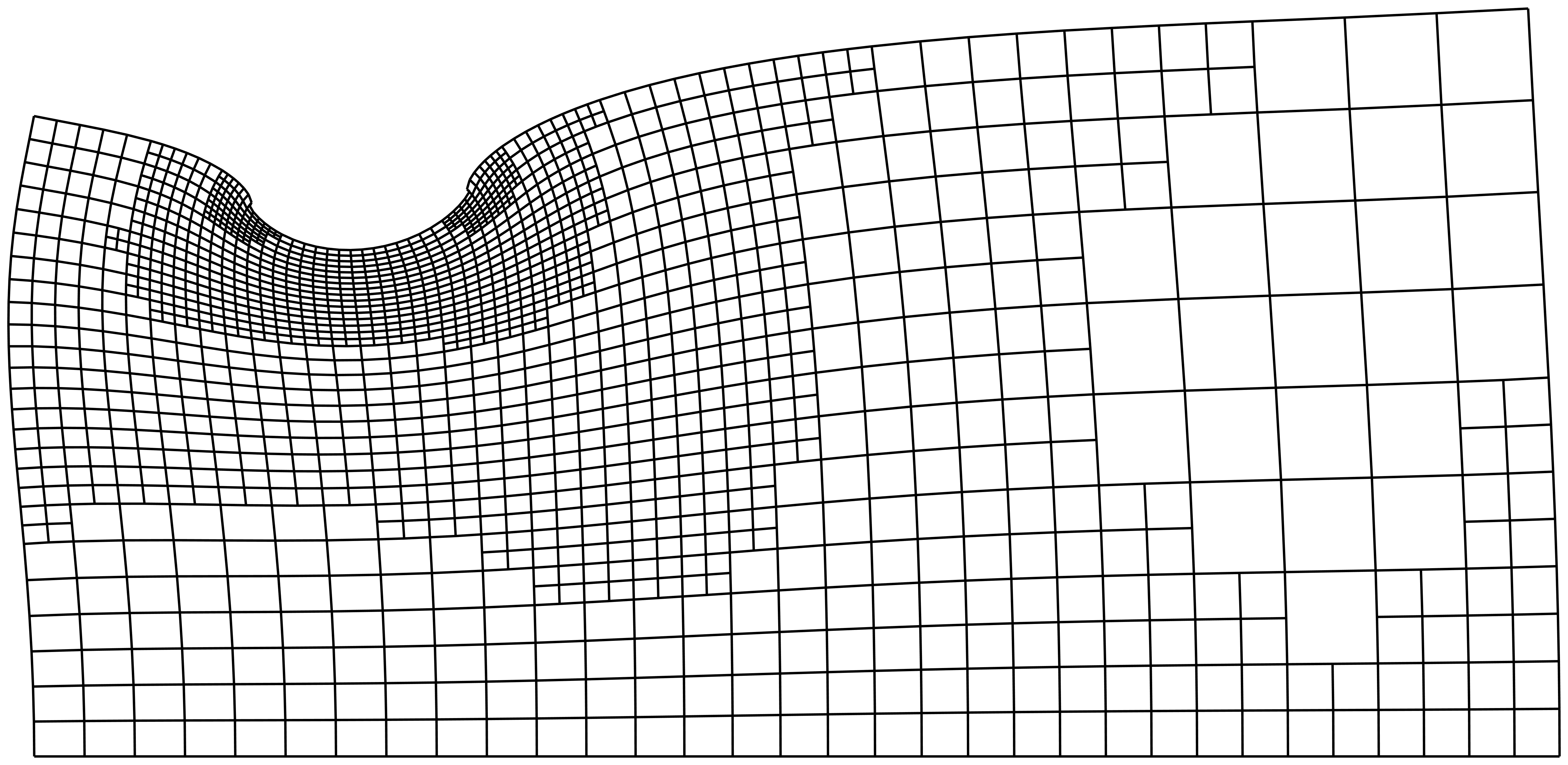}}
			\caption{Cycle 1 - Step 5}
		\end{subfigure}%
		\begin{subfigure}[t]{0.33\textwidth}
			\centering
			\includegraphics[width=0.95\textwidth]{{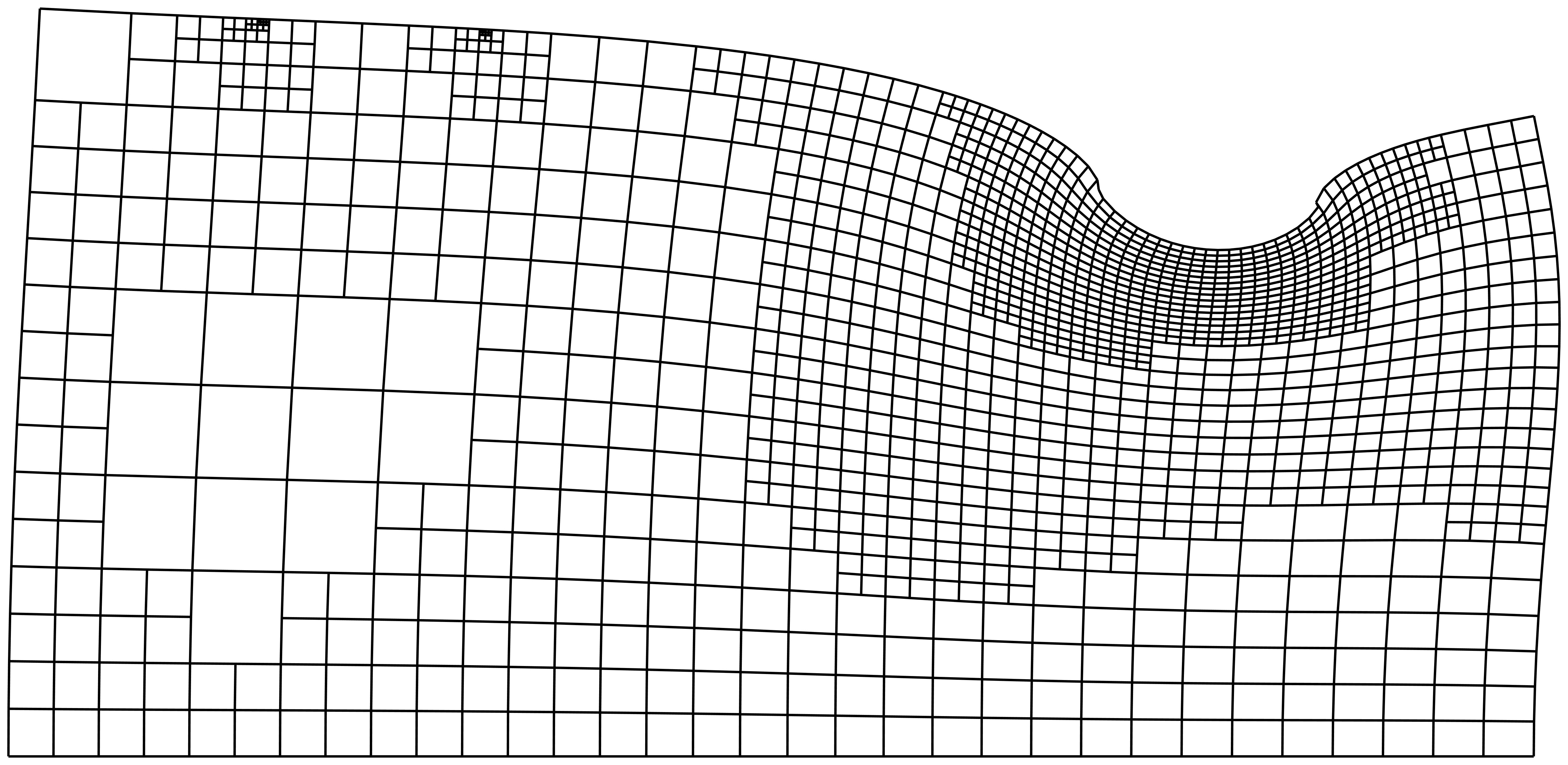}}
			\caption{Cycle 2 - Step 3}
		\end{subfigure}%
		\begin{subfigure}[t]{0.33\textwidth}
			\centering
			\includegraphics[width=0.95\textwidth]{{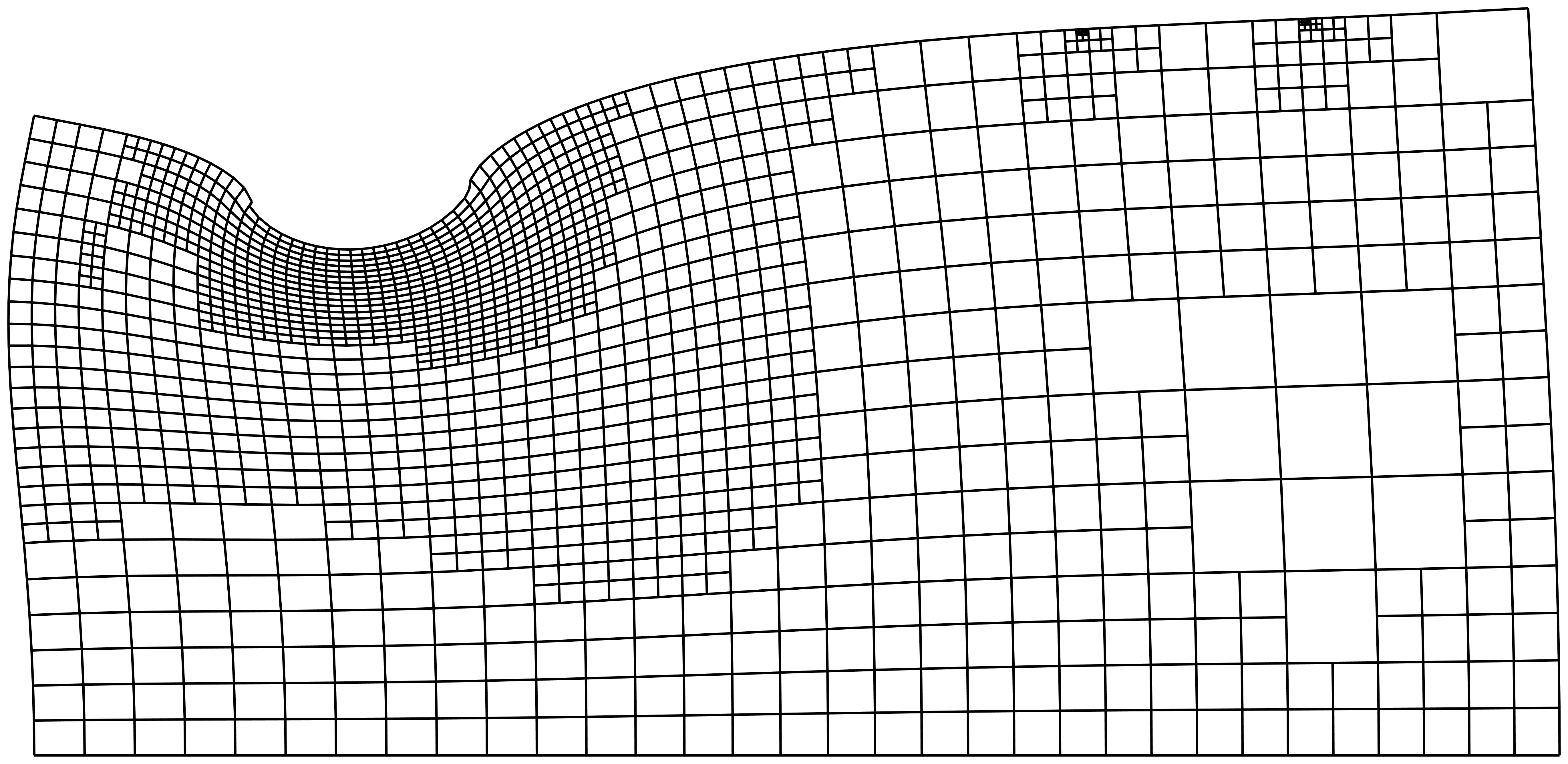}}
			\caption{Cycle 3 - Step 3}
		\end{subfigure}
		\vskip \baselineskip 
		\vspace*{-4mm}
		\begin{subfigure}[t]{0.33\textwidth}
			\centering
			\includegraphics[width=0.95\textwidth]{{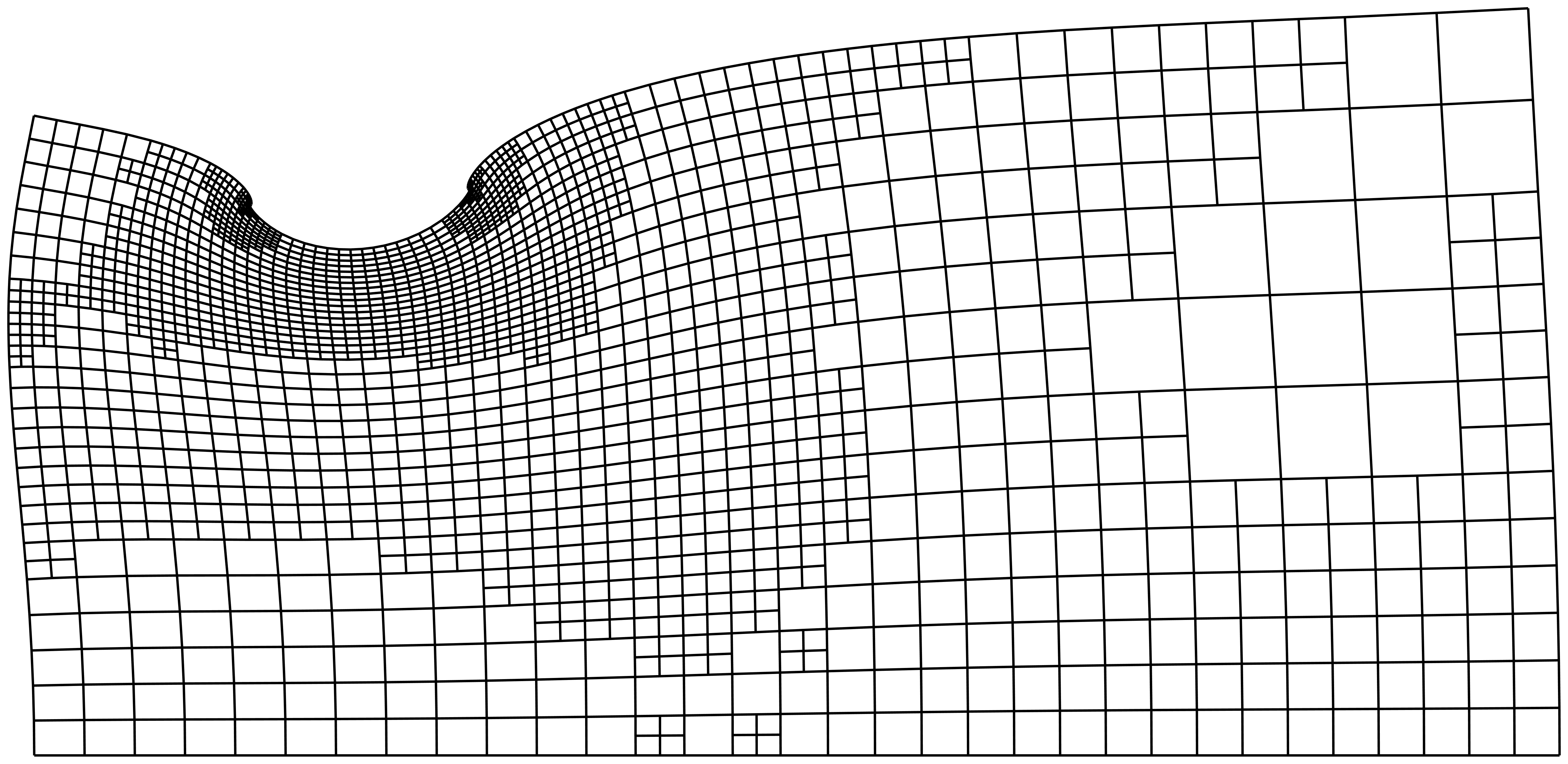}}
			\caption{Cycle 1 - Step 7}
		\end{subfigure}%
		\begin{subfigure}[t]{0.33\textwidth}
			\centering
			\includegraphics[width=0.95\textwidth]{{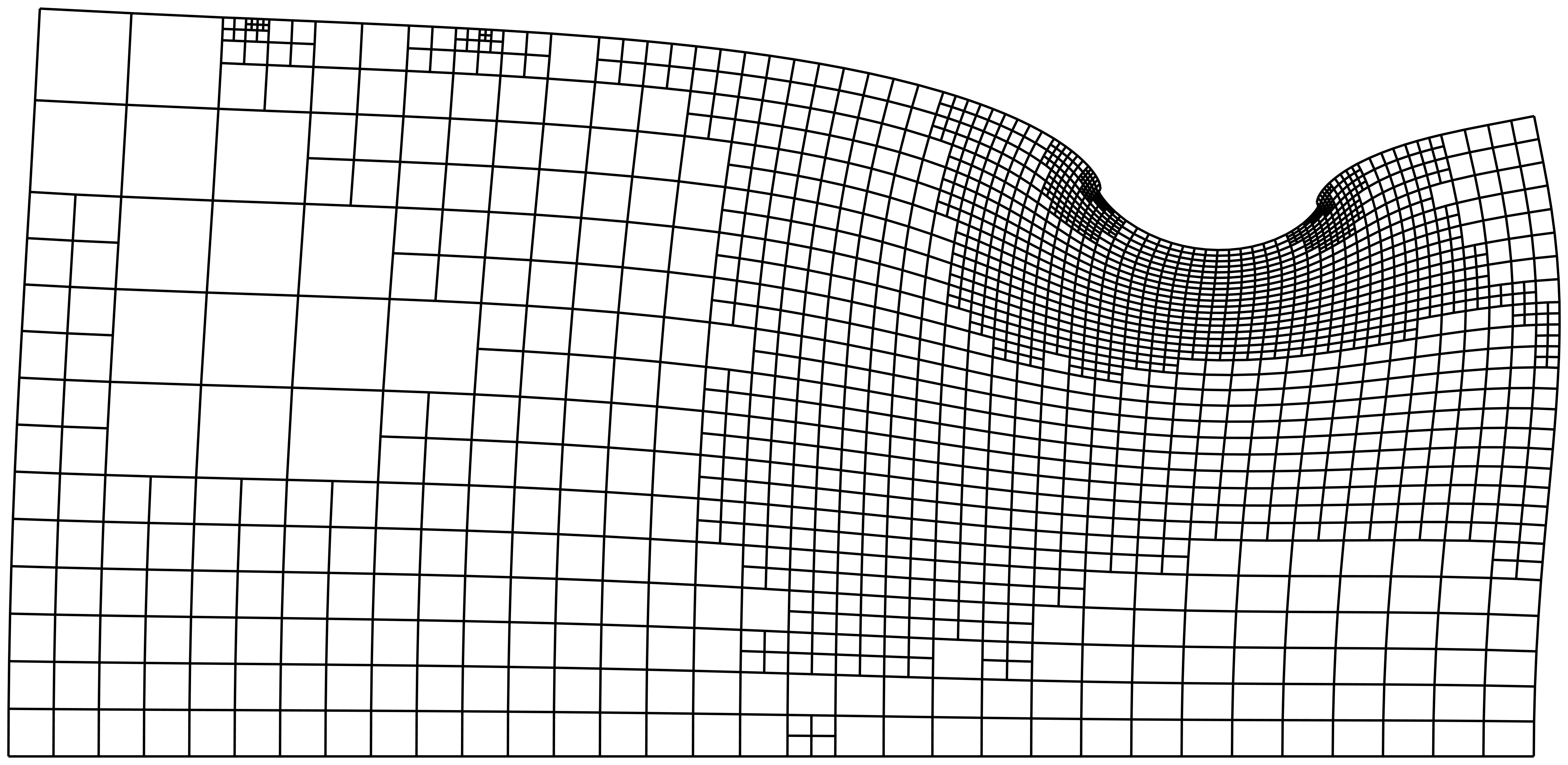}}
			\caption{Cycle 2 - Step 5}
		\end{subfigure}%
		\begin{subfigure}[t]{0.33\textwidth}
			\centering
			\includegraphics[width=0.95\textwidth]{{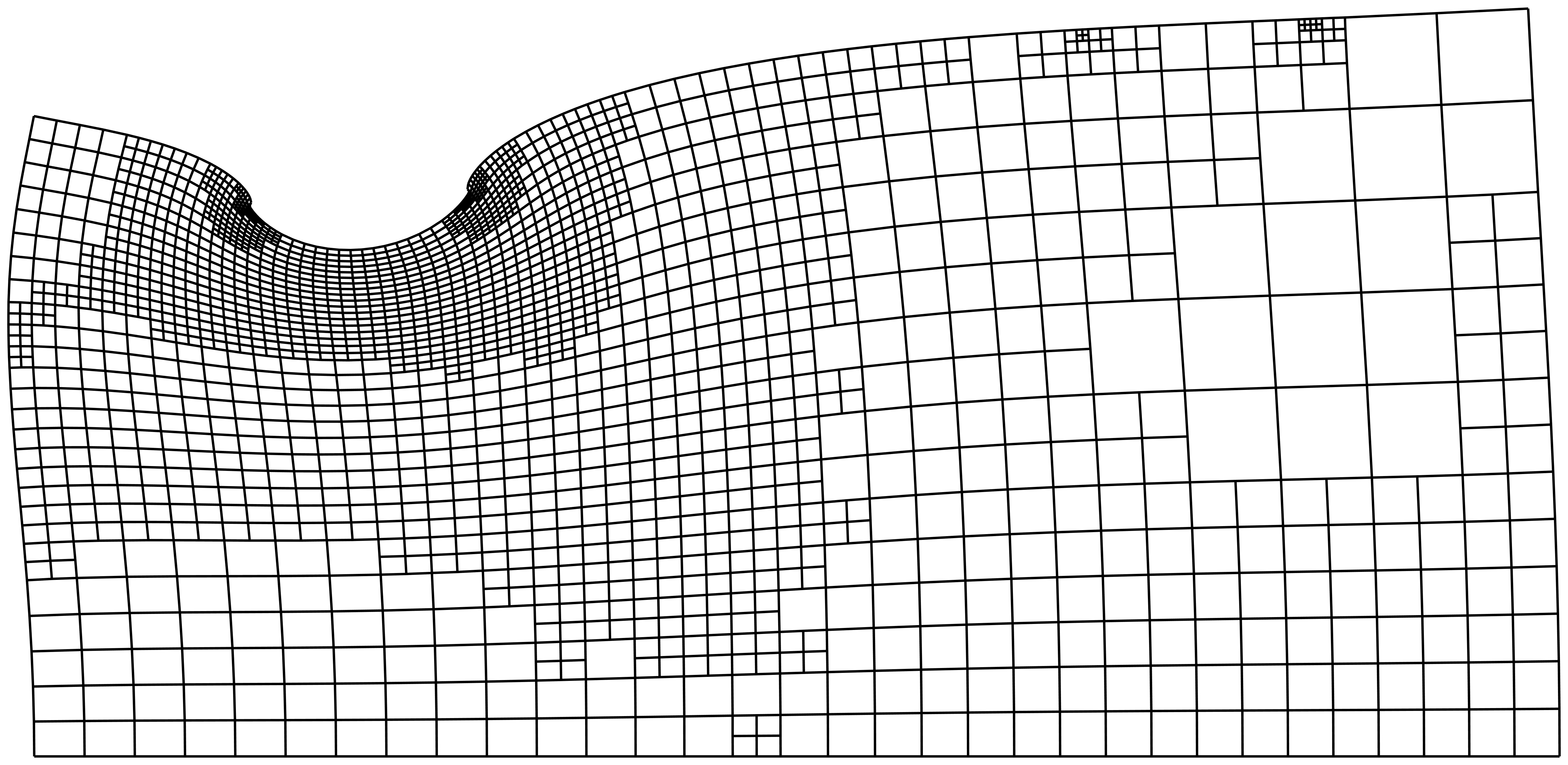}}
			\caption{Cycle 3 - Step 5}
		\end{subfigure}
		\vskip \baselineskip 
		\vspace*{-4mm}
		\begin{subfigure}[t]{0.33\textwidth}
			\centering
			\includegraphics[width=0.95\textwidth]{{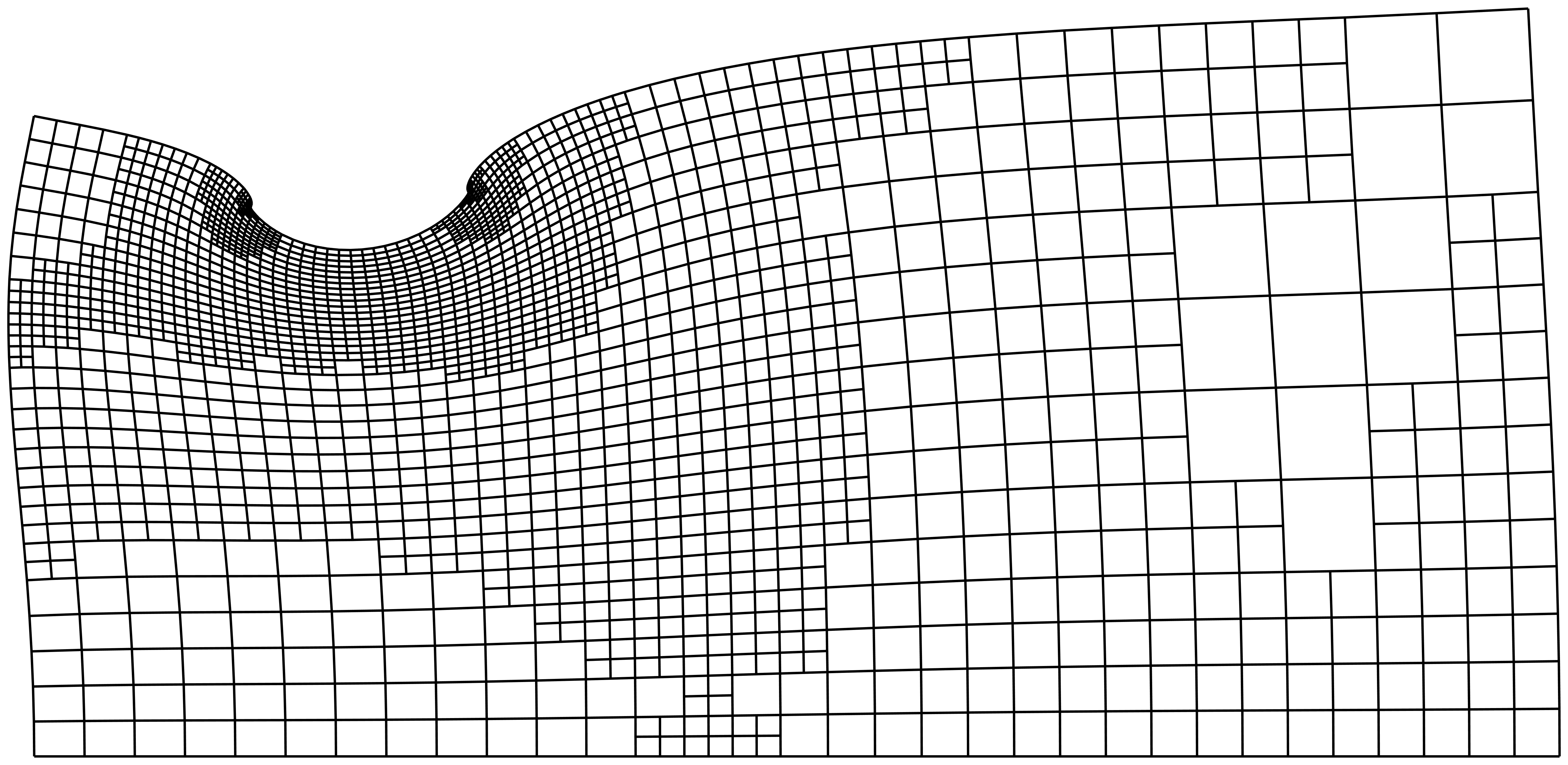}}
			\caption{Cycle 1 - Step 9}
		\end{subfigure}%
		\begin{subfigure}[t]{0.33\textwidth}
			\centering
			\includegraphics[width=0.95\textwidth]{{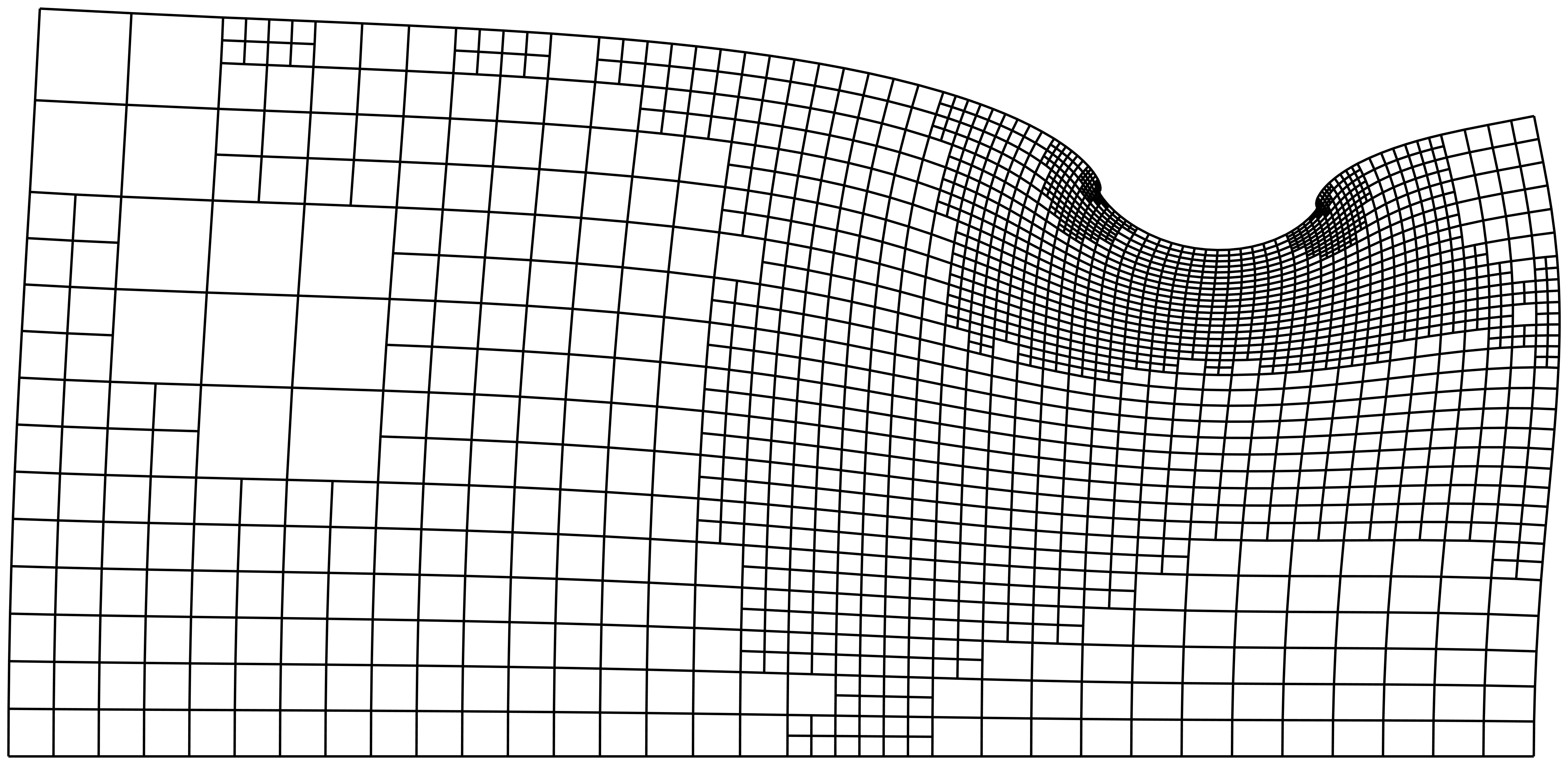}}
			\caption{Cycle 2 - Step 7}
		\end{subfigure}%
		\begin{subfigure}[t]{0.33\textwidth}
			\centering
			\includegraphics[width=0.95\textwidth]{{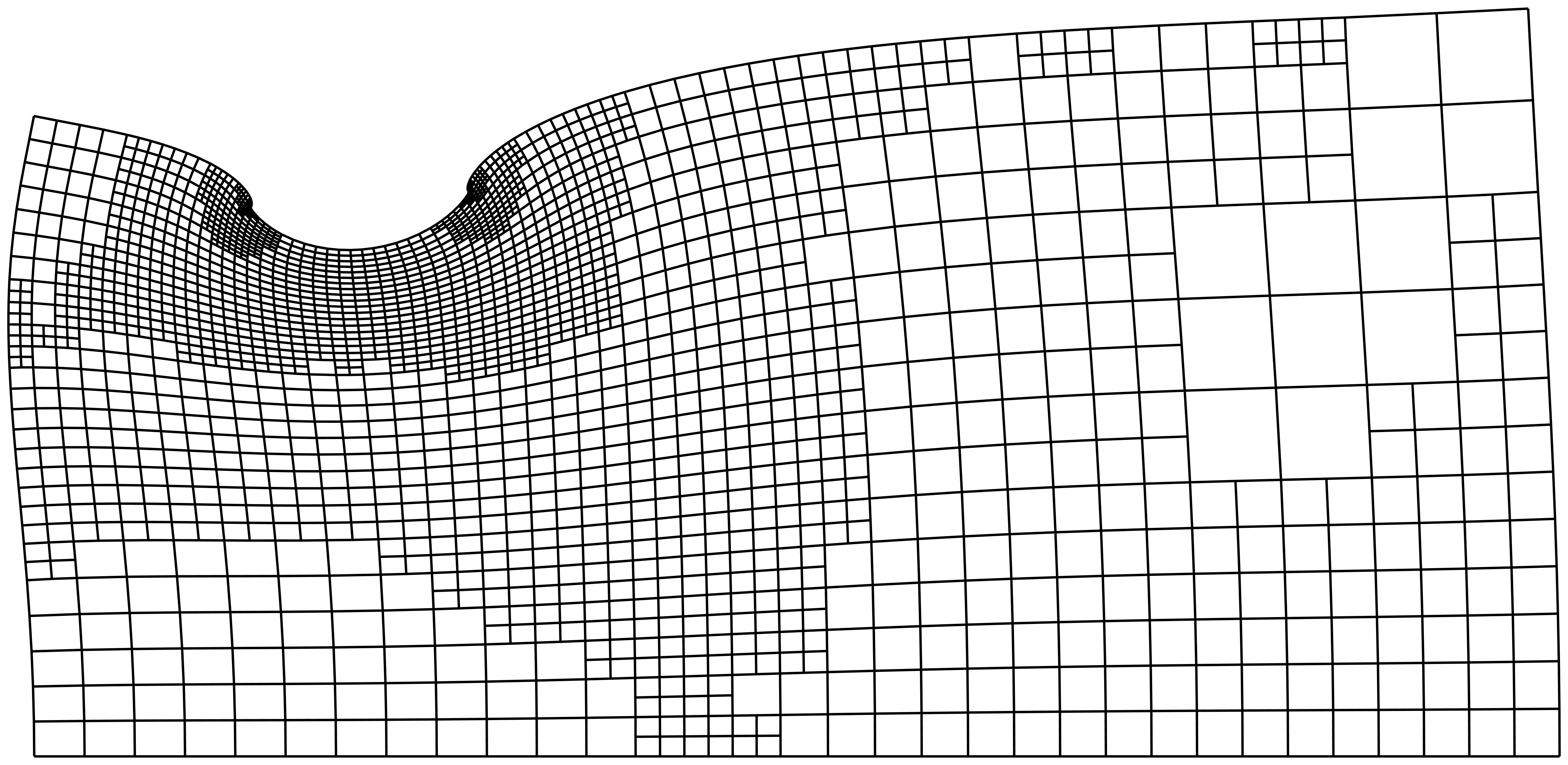}}
			\caption{Cycle 3 - Step 7}
		\end{subfigure}
		\vskip \baselineskip 
		\vspace*{-4mm}
		\begin{subfigure}[t]{0.33\textwidth}
			\centering
			\includegraphics[width=0.95\textwidth]{{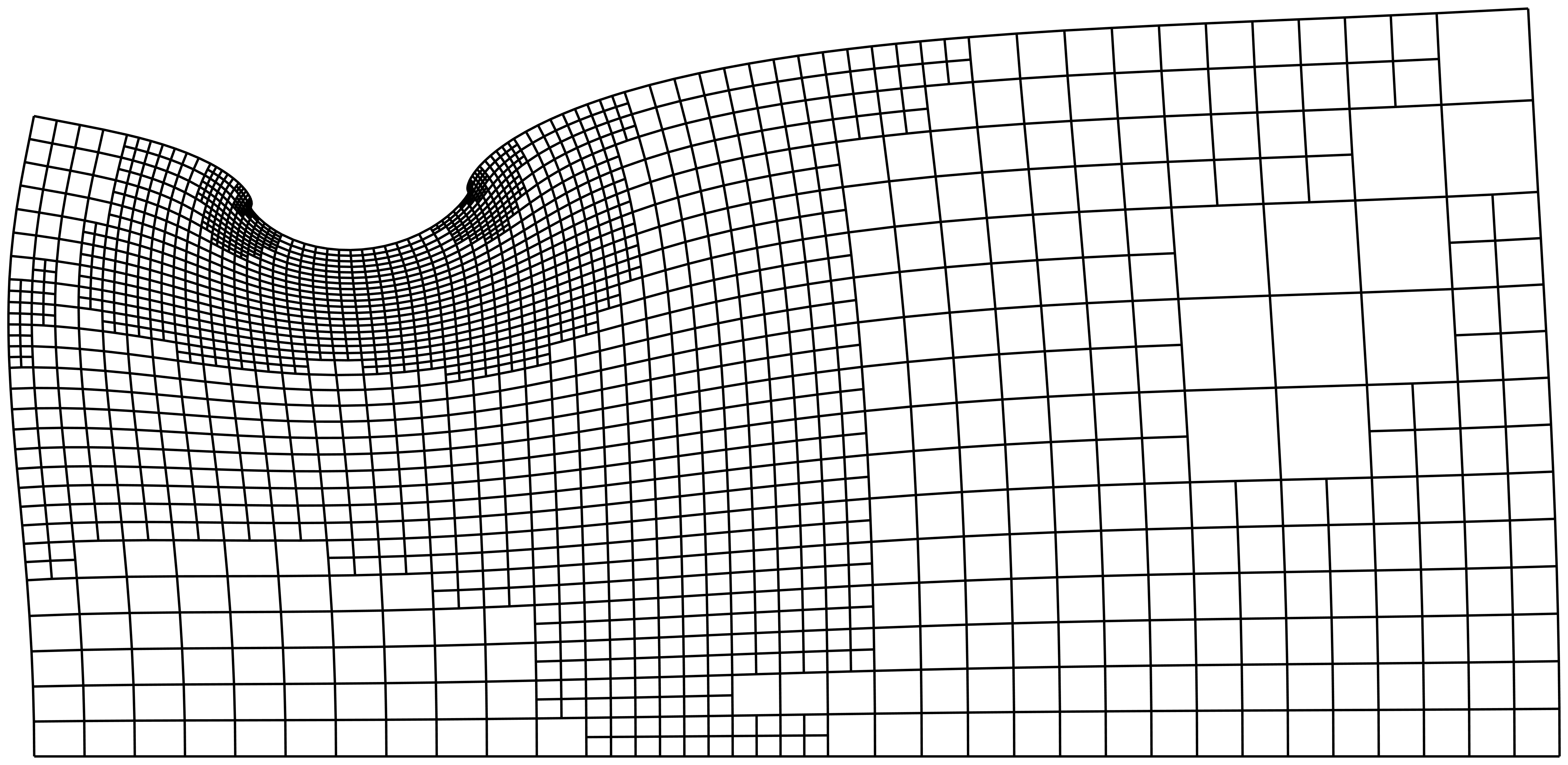}}
			\caption{Cycle 1 - Step 12 (Final)}
		\end{subfigure}%
		\begin{subfigure}[t]{0.33\textwidth}
			\centering
			\includegraphics[width=0.95\textwidth]{{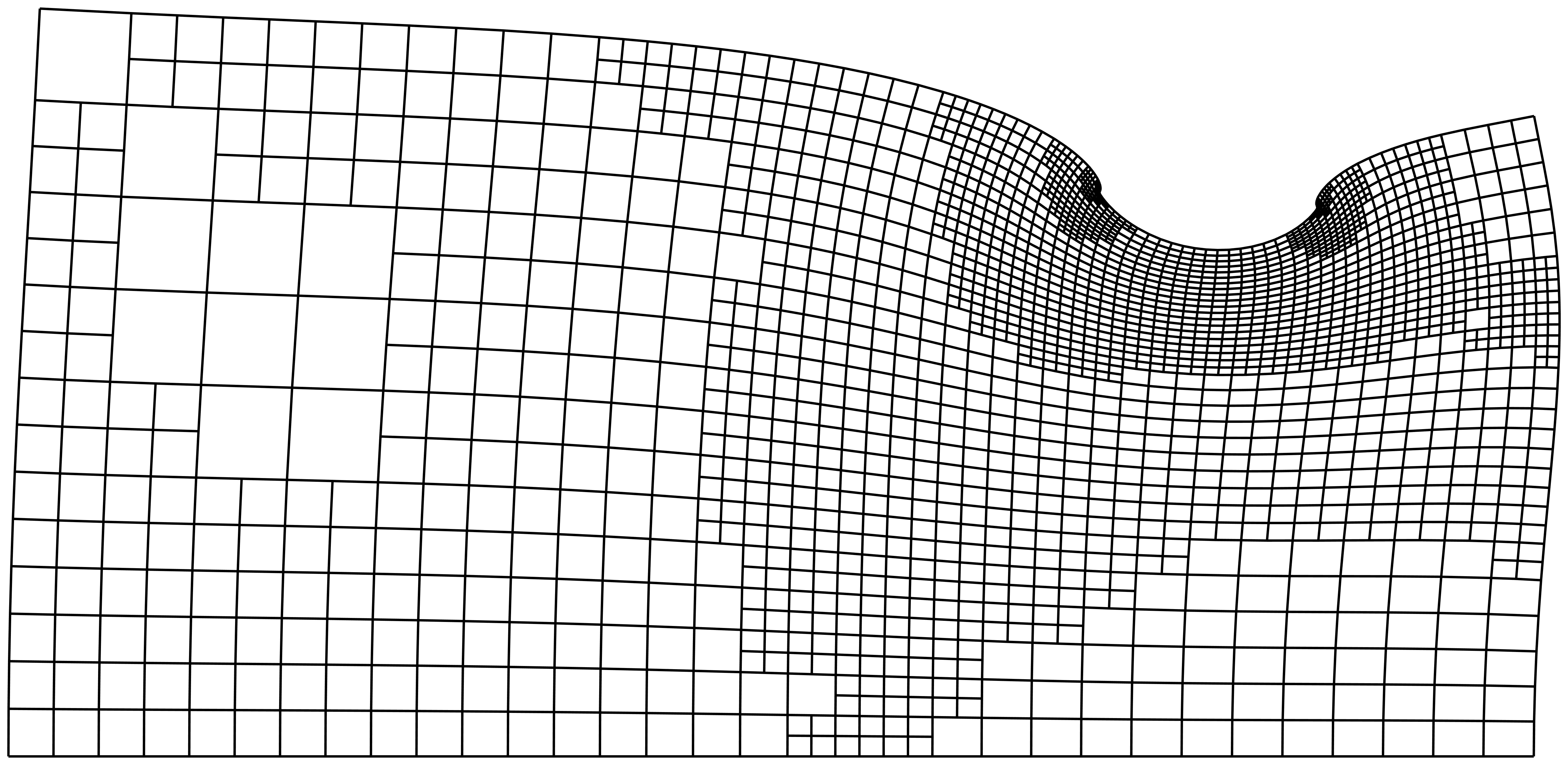}}
			\caption{Cycle 2 - Step 9 (Final)}
		\end{subfigure}%
		\begin{subfigure}[t]{0.33\textwidth}
			\centering
			\includegraphics[width=0.95\textwidth]{{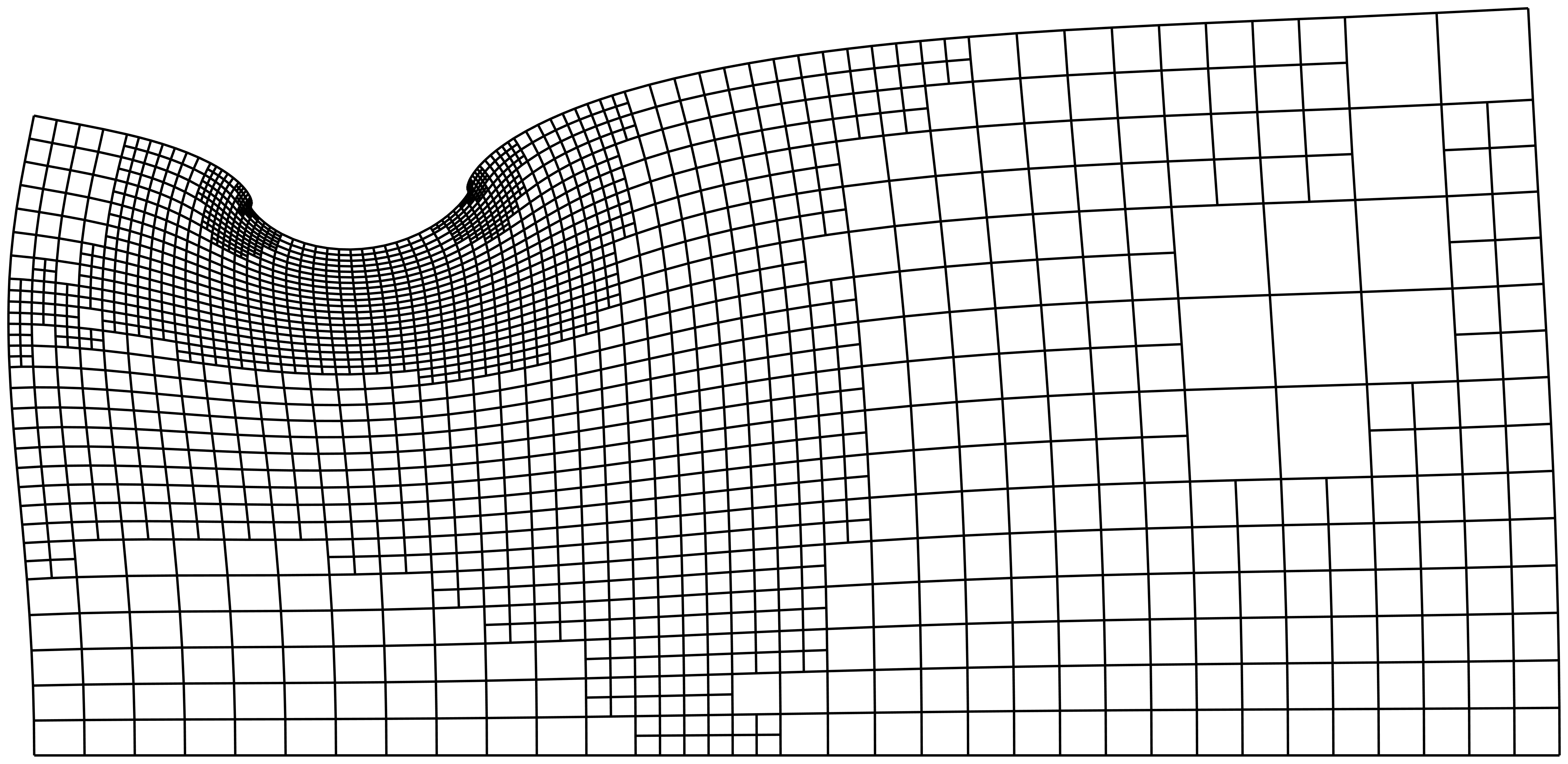}}
			\caption{Cycle 3 - Step 9 (Final)}
		\end{subfigure}
		\caption{Mesh evolution for the pseudo-dynamic punch problem from an initial uniform structured mesh with ${\|e\|_{\text{rel}}^{\text{targ}} = 5\%}$ for three load cycles.
			\label{fig:DynamicPunch_TargetError_MeshEvolution}}
	\end{figure}

	The mesh evolution during the fully adaptive remeshing process for the pseudo-dynamic punch problem is depicted in Figure~\ref{fig:DynamicPunch_TargetError_MeshEvolution} for an initially uniform structured mesh with an error target of ${\|e\|_{\text{rel}}^{\text{targ}} = 5\%}$. Here, three load cycles are considered with each column of figures corresponding to a particular load cycle. Each load cycle begins with the application of the corresponding boundary conditions and ends once the mesh adaptation is complete, i.e. the global error target and the termination criteria described in Section~\ref{subsec:ElementSelectionTargetError} have been met. 
	The mesh evolution is sensible and intuitive for this problem with the region around the active punch being most highly refined and the rest of the domain remaining relatively coarse. 
	Most notably, through this problem the reversibility of the proposed fully adaptive remeshing procedure, and thus its suitability for dynamic problems, is demonstrated. 
	
	\FloatBarrier
	\FloatBarrier
	
\begin{figure}[ht!]
		\centering
		\begin{subfigure}[t]{0.495\textwidth}
			\centering
			\includegraphics[width=0.95\textwidth]{{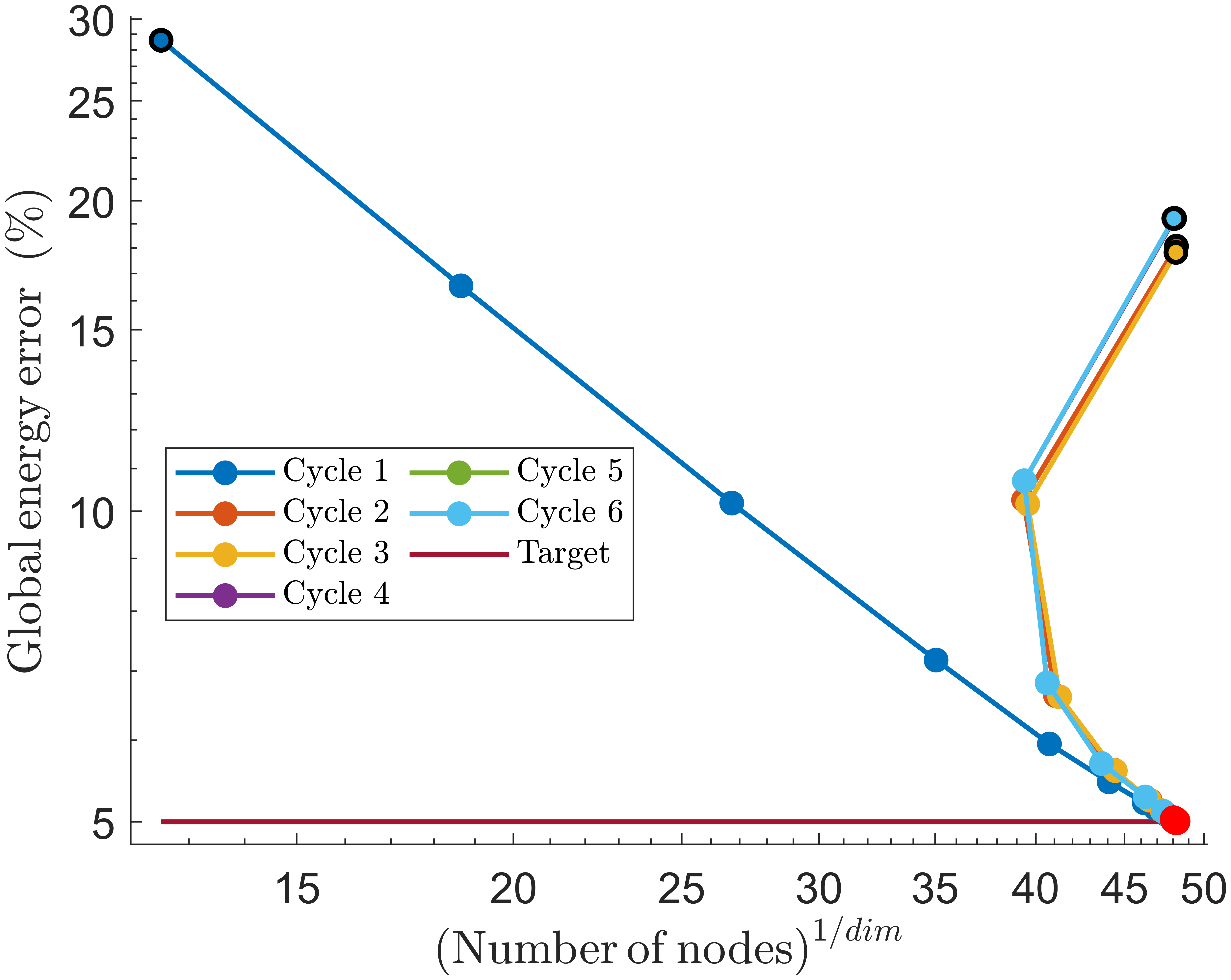}}
			\caption{Error convergence}
		\end{subfigure}%
		\begin{subfigure}[t]{0.495\textwidth}
			\centering
			\includegraphics[width=0.95\textwidth]{{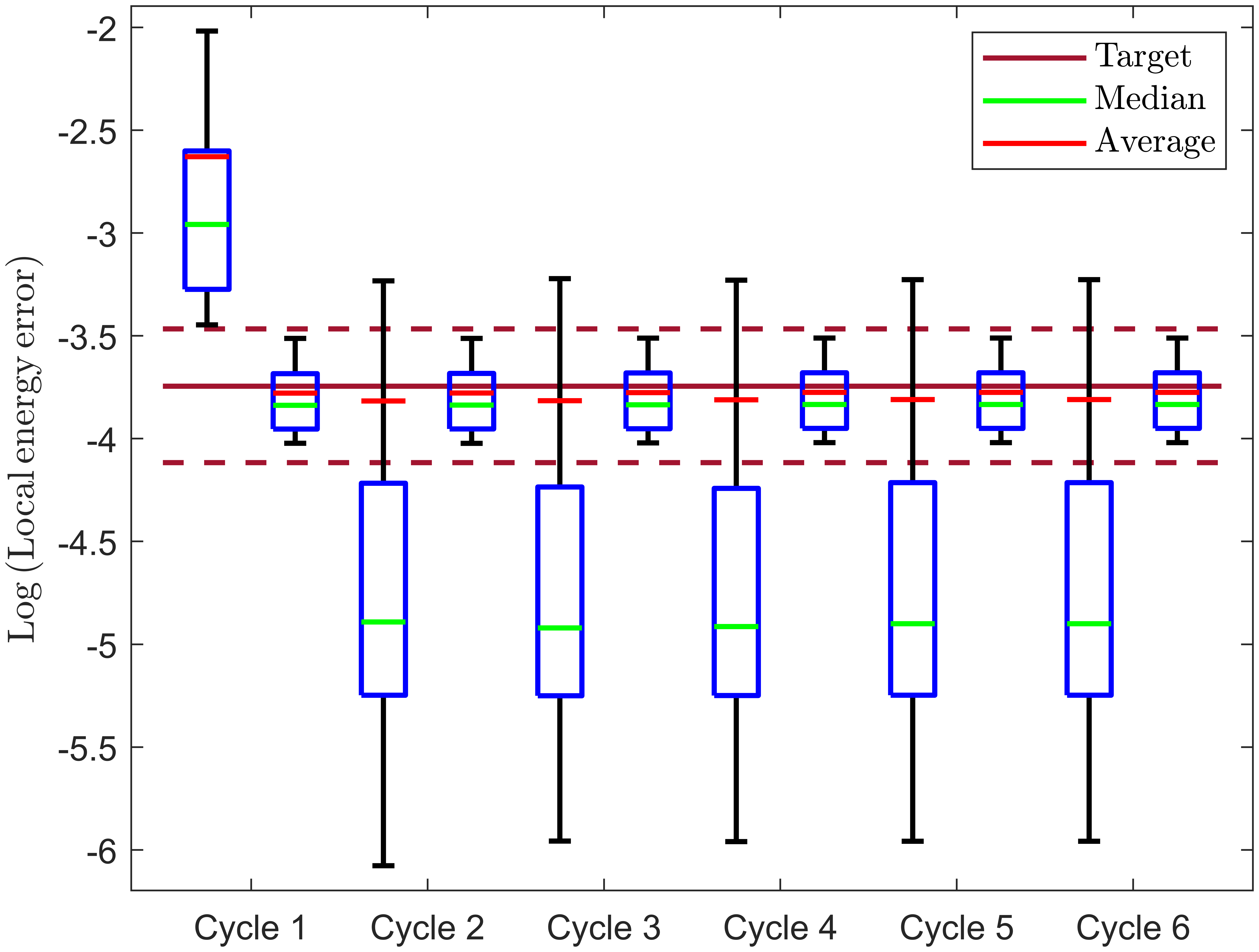}}
			\caption{Error distribution}
		\end{subfigure}
		\caption{Energy error (a) convergence and (b) box and whisker distribution for the pseudo-dynamic punch problem on structured meshes with ${\|e\|_{\text{rel}}^{\text{targ}} = 5\%}$ for six load cycles.
			\label{fig:DualPunch_ErrorAnalysis}}
	\end{figure}

	The energy error convergence and distribution for the pseudo-dynamic punch problem with an initially uniform structured mesh with ${\|e\|_{\text{rel}}^{\text{targ}} = 5\%}$ are plotted in Figures~\ref{fig:DualPunch_ErrorAnalysis}(a) and (b) respectively for six load cycles.
	
	In Figure~\ref{fig:DualPunch_ErrorAnalysis}(a) the energy error approximation is plotted against the number of nodes on a logarithmic scale. The convergence curve for each load cycle is plotted in a different colour and the outline of the marker denoting the first step in each cycle is indicated in black. Additionally, the final adapted mesh result is indicated by a red marker. 
	During the first cycle typical convergence behaviour is exhibited as the initially uniform coarse mesh is increasingly refined in the region around the punch until the target error and termination criteria are met.
	Thereafter, the boundary conditions are changed, while the mesh is not, thus the mesh is completely unsuitable for the new load conditions and a very high error is exhibited for the first step of the second cycle. The mesh adaptation process then iteratively improves the mesh until the error target and termination criteria are met. This process is repeated four more times with almost identical accuracy and efficiency exhibited by the final fully adapted meshes, as indicated by the almost identical locations of the red markers. 
	
	In Figure~\ref{fig:DualPunch_ErrorAnalysis}(b) the distribution of the element-level energy error is illustrated through a box and whisker plot where pairs of results correspond to the error evolution for a particular load cycle. 
	Here, the maximum and minimum element-level errors of the final adapted meshes fall within the prescribed error bounds for each load cycle. Furthermore, the upper and lower quartiles indicate a narrow distribution of error around the average, and the average error is almost identical to the target error. 
	
	The results presented in Figure~\ref{fig:DualPunch_ErrorAnalysis}(a) demonstrated that the fully adaptive procedure was able to meet all specified error target for every load cycle.
	The results presented in Figure~\ref{fig:DualPunch_ErrorAnalysis}(b) demonstrated that the element-level errors, on average, met the element-level target and were approximately equal as they fell within the specified target error range.
	Therefore, the fully adaptive procedure successfully generated quasi-optimal meshes for the specified error target for every load cycle, thus, demonstrating its suitability for dynamic problems. 
	
	\FloatBarrier
	\FloatBarrier

	\section{Discussion and conclusion} 
	\label{sec:Conclusion}
	In this work a novel fully adaptive remeshing procedure has been proposed for the virtual element method.
	The remeshing procedure comprises the novel combination of refinement and coarsening procedures for structured and unstructured/Voronoi meshes as well as novel procedures for the selection of elements to refine and element patches to coarsen.
	
	Three procedures for the selection of elements to refine and element patches to coarsen have been proposed for various adaptivity targets. Specifically, procedures were proposed for a target global error, a target number of elements, and a target number of nodes. In this work error was measured through an approximation of the well-known energy error. Additionally, all of the proposed element selection procedures were constructed to meet their respective targets while creating a mesh in which all elements have an approximately equal error. Thus, creating a quasi-even error distribution over the elements corresponding to a quasi-optimal mesh for a specific target.  
	
	The proposed fully adaptive procedures were studied numerically on a well-known benchmark problem.
	For each of the target types the mesh evolution during remeshing was analysed along with analysis of the error convergence and the distribution of error over the problem domain.
	In terms of the mesh evolution, the efficacy of the proposed procedures was evident as the adapted meshes had increased refinement in the most critical regions of the domain and remained relatively coarse, or were coarsened further, elsewhere. This efficacy was demonstrated on both structured and unstructured/Voronoi meshes and was demonstrated to be independent of the initial mesh.
	Furthermore, the efficacy of the proposed fully adaptive remeshing procedures was studied in the energy error norm. Here, the performance of the procedures was investigated on structured and unstructured/Voronoi meshes and was compared to a reference approach comprising meshes of uniform discretization density. Additionally, the influence of the initial mesh density on the performance of the adaptive procedures was investigated.
	The numerical results demonstrated the high degree of efficacy of all proposed adaptive procedures. The procedures were able to meet their specified targets from all initial meshes and on both mesh types.
	Furthermore, it was demonstrated that the meshes generated by the adaptive procedures had a quasi-even error distribution and thus represented quasi-optimal meshes for their respective targets.
	Additionally, the suitability of the proposed fully adaptive procedures for dynamic problems was demonstrated through the presentation of a novel pseudo-dynamic problem in which the challenges faced during a dynamic problem were mimicked.
	
	The good performance exhibited by the proposed fully adaptive procedures over a range of target types on both structured and unstructured/Voronoi meshes demonstrates its versatility, efficacy and suitability for application to the analysis of elastic problems using the virtual element method. 
	
	Future work of interest would be the extension to non-linear problems, higher-order formulations, and problems in three-dimensions would be of great interest.
	
	\section*{Conflict of interest}
	
	The authors declare that they have no known competing financial interests or personal relationships that could have appeared to influence the work reported in this paper.
	
	\section*{Acknowledgements}
	
	This work was carried out with support from the German Science Foundation (DFG) and the National Scientific and Technical Research Council of Argentina (CONICET) through project number DFG 544/68-1 (431843479).
	The authors acknowledge with thanks this support.
	
	\bibliographystyle{elsarticle-num}
	\bibliography{VEM_References}

\begin{thebibliography}{10}
\expandafter\ifx\csname url\endcsname\relax
  \def\url#1{\texttt{#1}}\fi
\expandafter\ifx\csname urlprefix\endcsname\relax\def\urlprefix{URL }\fi
\expandafter\ifx\csname href\endcsname\relax
  \def\href#1#2{#2} \def\path#1{#1}\fi

\bibitem{wriggers1970comparison}
P.~Wriggers, A.~Reiger, Comparison of different error measures for adaptive
  finite element techniques applied to contact problems involving large elastic
  strains, WIT Transactions on Engineering Sciences 24 (1970).

\bibitem{oliveira1971optimization}
E.~R.~A. Oliveira, Optimization of finite element solutions, in: Proceedings,
  Wright-Patterson Air Force Base Ohio, 1971, p. 423.

\bibitem{Zhu1988}
J.~Z. Zhu, O.~C. Zienkiewicz, Adaptive techniques in the finite element method,
  Communications in Applied Numerical Methods 4~(2) (1988) 197--204.
\newblock \href {https://doi.org/10.1002/cnm.1630040210}
  {\path{doi:10.1002/cnm.1630040210}}.

\bibitem{ONATE1993}
E.~O\~{n}ate, G.~Bugeda, A study of mesh optimality criteria in adaptive finite
  element analysis, Engineering Computations 10~(4) (1993) 307--321.
\newblock \href {https://doi.org/10.1108/eb023910}
  {\path{doi:10.1108/eb023910}}.

\bibitem{Li1995}
L.~Li, P.~Bettess, Notes on mesh optimal criteria in adaptive finite element
  computations, Communications in Numerical Methods in Engineering 11~(11)
  (1995) 911--915.
\newblock \href {https://doi.org/10.1002/cnm.1640111105}
  {\path{doi:10.1002/cnm.1640111105}}.

\bibitem{babuvvska1978error}
I.~Babu{\v{s}}ka, W.~C. Rheinboldt, Error estimates for adaptive finite element
  computations, SIAM Journal on Numerical Analysis 15~(4) (1978) 736--754.
\newblock \href {https://doi.org/10.1137/0715049} {\path{doi:10.1137/0715049}}.

\bibitem{Babuska1979}
I.~Babu{\v{s}}ka, W.~C. Rheinboldt, Adaptive approaches and reliability
  estimations in finite element analysis, Computer Methods in Applied Mechanics
  and Engineering 17-18 (1979) 519--540.
\newblock \href {https://doi.org/10.1016/0045-7825(79)90042-2}
  {\path{doi:10.1016/0045-7825(79)90042-2}}.

\bibitem{Zienkiewicz1987}
O.~C. Zienkiewicz, J.~Z. Zhu, A simple error estimator and adaptive procedure
  for practical engineering analysis, Int J Numer Methods Eng 24~(2) (1987)
  337--357.
\newblock \href {https://doi.org/10.1002/nme.1620240206}
  {\path{doi:10.1002/nme.1620240206}}.

\bibitem{Zienkiewicz1992}
O.~C. Zienkiewicz, J.~Z. Zhu, The superconvergent patch recovery and a
  posteriori error estimates. part 1: The recovery technique, International
  Journal for Numerical Methods in Engineering 33~(7) (1992) 1331--1364.
\newblock \href {https://doi.org/10.1002/nme.1620330702}
  {\path{doi:10.1002/nme.1620330702}}.

\bibitem{Zienkiewicz1992a}
O.~C. Zienkiewicz, J.~Z. Zhu, The superconvergent patch recovery and a
  posteriori error estimates. part 2: Error estimates and adaptivity,
  International Journal for Numerical Methods in Engineering 33~(7) (1992)
  1365--1382.
\newblock \href {https://doi.org/10.1002/nme.1620330703}
  {\path{doi:10.1002/nme.1620330703}}.

\bibitem{Zienkiewicz1995}
O.~C. Zienkiewicz, J.~Z. Zhu, Superconvergence and the superconvergent patch
  recovery, Finite Elements in Analysis and Design 19~(1-2) (1995) 11--23.
\newblock \href {https://doi.org/10.1016/0168-874x(94)00054-j}
  {\path{doi:10.1016/0168-874x(94)00054-j}}.

\bibitem{Arndt2022}
D.~Arndt, W.~Bangerth, M.~Feder, M.~Fehling, R.~Gassm{\"{o}}ller, T.~Heister,
  L.~Heltai, M.~Kronbichler, M.~Maier, P.~Munch, J.~P. Pelteret, S.~Sticko,
  B.~Turcksin, D.~Wells, The deal.{II} library, version 9.4, Journal of
  Numerical Mathematics 30~(3) (2022) 231--246.
\newblock \href {https://doi.org/10.1515/jnma-2022-0054}
  {\path{doi:10.1515/jnma-2022-0054}}.

\bibitem{Kirk2006}
B.~S. Kirk, J.~W. Peterson, R.~H. Stogner, G.~F. Carey, {libMesh} : a c++
  library for parallel adaptive mesh refinement/coarsening simulations,
  Engineering with Computers 22~(3-4) (2006) 237--254.
\newblock \href {https://doi.org/10.1007/s00366-006-0049-3}
  {\path{doi:10.1007/s00366-006-0049-3}}.

\bibitem{Park2012}
K.~Park, G.~H. Paulino, W.~Celes, R.~Espinha, Adaptive mesh refinement and
  coarsening for cohesive zone modeling of dynamic fracture, International
  Journal for Numerical Methods in Engineering 92~(1) (2012) 1--35.
\newblock \href {https://doi.org/10.1002/nme.3163}
  {\path{doi:10.1002/nme.3163}}.

\bibitem{VEIGA2012}
L.~{Beir\~{a}o da Veiga}, F.~Brezzi, A.~Cangiani, G.~Manzini, L.~D. Marini,
  A.~Russo, Basic principles of virtual element methods, Mathematical Models
  and Methods in Applied Sciences 23~(01) (2012) 199--214.
\newblock \href {https://doi.org/10.1142/s0218202512500492}
  {\path{doi:10.1142/s0218202512500492}}.

\bibitem{Veiga2014}
L.~{Beir\~{a}o da Veiga}, F.~Brezzi, L.~D. Marini, A.~Russo, The hitchhikers
  guide to the virtual element method, Mathematical Models and Methods in
  Applied Sciences 24~(08) (2014) 1541--1573.
\newblock \href {https://doi.org/10.1142/s021820251440003x}
  {\path{doi:10.1142/s021820251440003x}}.

\bibitem{VEMContactWriggers2016}
P.~Wriggers, W.~T. Rust, B.~D. Reddy, A virtual element method for contact,
  Computational Mechanics 58~(6) (2016) 1039--1050.
\newblock \href {https://doi.org/10.1007/s00466-016-1331-x}
  {\path{doi:10.1007/s00466-016-1331-x}}.

\bibitem{Wriggers2019}
P.~Wriggers, W.~T. Rust, A virtual element method for frictional contact
  including large deformations, Engineering Computations 36~(7) (2019)
  2133--2161.
\newblock \href {https://doi.org/10.1108/ec-02-2019-0043}
  {\path{doi:10.1108/ec-02-2019-0043}}.

\bibitem{Sorgente2021}
T.~Sorgente, S.~Biasotti, G.~Manzini, M.~Spagnuolo, The role of mesh quality
  and mesh quality indicators in the virtual element method, Advances in
  Computational Mathematics 48~(1) (2021).
\newblock \href {https://doi.org/10.1007/s10444-021-09913-3}
  {\path{doi:10.1007/s10444-021-09913-3}}.

\bibitem{Sorgente2022}
T.~Sorgente, S.~Biasotti, G.~Manzini, M.~Spagnuolo, Polyhedral mesh quality
  indicator for the {V}irtual {E}lement {M}ethod, Computers {\&} Mathematics
  with Applications 114 (2022) 151--160.
\newblock \href {https://doi.org/10.1016/j.camwa.2022.03.042}
  {\path{doi:10.1016/j.camwa.2022.03.042}}.

\bibitem{Huyssteen2020}
D.~{van Huyssteen}, B.~D. Reddy, A virtual element method for isotropic
  hyperelasticity, Computer Methods in Applied Mechanics and Engineering 367
  (2020) 113134.
\newblock \href {https://doi.org/10.1016/j.cma.2020.113134}
  {\path{doi:10.1016/j.cma.2020.113134}}.

\bibitem{Huyssteen2021}
D.~{van Huyssteen}, B.~D. Reddy, A virtual element method for transversely
  isotropic hyperelasticity, Computer Methods in Applied Mechanics and
  Engineering 386 (2021) 114108.
\newblock \href {https://doi.org/10.1016/j.cma.2021.114108}
  {\path{doi:10.1016/j.cma.2021.114108}}.

\bibitem{ReddyvanHuyssteen2021}
B.~D. Reddy, D.~{van {H}uyssteen}, Alternative approaches to the stabilization
  of virtual element formulations for hyperelasticity, in: F.~Aldakheel,
  B.~Hudobivnik, M.~Soleimani, H.~Wessels, C.~Weissenfels, M.~Marino (Eds.),
  Current Trends and Open Problems in Computational Mechanics, Springer, 2021.

\bibitem{DvH_BDR_MeshQuality}
D.~{van Huyssteen}, B.~D. Reddy, The incorporation of mesh quality in the
  stabilization of virtual element methods for nonlinear elasticity, Computer
  Methods in Applied Mechanics and Engineering 392 (2022) 114720.
\newblock \href {https://doi.org/10.1016/j.cma.2022.114720}
  {\path{doi:10.1016/j.cma.2022.114720}}.

\bibitem{WriggersIsotropic2017}
P.~Wriggers, B.~D. Reddy, W.~Rust, B.~Hudobivnik, Efficient virtual element
  formulations for compressible and incompressible finite deformations,
  Computational Mechanics 60~(2) (2017) 253--268.
\newblock \href {https://doi.org/10.1007/s00466-017-1405-4}
  {\path{doi:10.1007/s00466-017-1405-4}}.

\bibitem{Wriggers2023}
P.~Wriggers, A locking free virtual element formulation for {T}imoshenko beams,
  Computer Methods in Applied Mechanics and Engineering (2023) 116234\href
  {https://doi.org/10.1016/j.cma.2023.116234}
  {\path{doi:10.1016/j.cma.2023.116234}}.

\bibitem{Tang2020}
X.~Tang, Z.~Liu, B.~Zhang, M.~Feng, A low-order locking-free virtual element
  for linear elasticity problems, Computers {\&} Mathematics with Applications
  80~(5) (2020) 1260--1274.
\newblock \href {https://doi.org/10.1016/j.camwa.2020.04.032}
  {\path{doi:10.1016/j.camwa.2020.04.032}}.

\bibitem{Reddy2019}
B.~D. Reddy, D.~{van Huyssteen}, A virtual element method for transversely
  isotropic elasticity, Computational Mechanics 64~(4) (2019) 971--988.
\newblock \href {https://doi.org/10.1007/s00466-019-01690-7}
  {\path{doi:10.1007/s00466-019-01690-7}}.

\bibitem{Cangiani2017}
A.~Cangiani, E.~H. Georgoulis, T.~Pryer, O.~J. Sutton, A posteriori error
  estimates for the virtual element method, Numer. Math. 137~(4) (2017)
  857--893.
\newblock \href {https://doi.org/10.1007/s00211-017-0891-9}
  {\path{doi:10.1007/s00211-017-0891-9}}.

\bibitem{Veiga2015a}
L.~{Beir\~{a}o da Veiga}, G.~Manzini, Residual a posteriori error estimation
  for the virtual element method for elliptic problems, Mathematical Modelling
  and Numerical Analysis 49~(2) (2015) 577--599.
\newblock \href {https://doi.org/10.1051/m2an/2014047}
  {\path{doi:10.1051/m2an/2014047}}.

\bibitem{Berrone2017}
S.~Berrone, A.~Borio, A residual a posteriori error estimate for the virtual
  element method, Mathematical Models and Methods in Applied Sciences 27~(08)
  (2017) 1423--1458.
\newblock \href {https://doi.org/10.1142/s0218202517500233}
  {\path{doi:10.1142/s0218202517500233}}.

\bibitem{Mora2017}
D.~Mora, G.~Rivera, R.~Rodr{\'{i}}guez, A posteriori error estimates for a
  {V}irtual {E}lement {M}ethod for the {S}teklov eigenvalue problem, Computers
  and Mathematics with Applications 74~(9) (2017) 2172--2190.
\newblock \href {https://doi.org/10.1016/j.camwa.2017.05.016}
  {\path{doi:10.1016/j.camwa.2017.05.016}}.

\bibitem{Chi2019}
H.~Chi, L.~{Beir\~{a}o da Veiga}, G.~H. Paulino, A simple and effective
  gradient recovery scheme and a posteriori error estimator for the virtual
  element method ({VEM}), Computer Methods in Applied Mechanics and Engineering
  347 (2019) 21--58.
\newblock \href {https://doi.org/10.1016/j.cma.2018.08.014}
  {\path{doi:10.1016/j.cma.2018.08.014}}.

\bibitem{Guo2019}
H.~Guo, C.~Xie, R.~Zhao, Superconvergent gradient recovery for virtual element
  methods, Mathematical Models and Methods in Applied Sciences 29~(11) (2019)
  2007--2031.
\newblock \href {https://doi.org/10.1142/s0218202519500386}
  {\path{doi:10.1142/s0218202519500386}}.

\bibitem{Wei2023}
H.~Wei, Y.~Deng, F.~Wang, Gradient recovery type a posteriori error estimates
  of virtual element method for an elliptic variational inequality of the
  second kind, Nonlinear Analysis: Real World Applications 73 (2023) 103903.
\newblock \href {https://doi.org/10.1016/j.nonrwa.2023.103903}
  {\path{doi:10.1016/j.nonrwa.2023.103903}}.

\bibitem{NguyenThanh2018}
V.~M. Nguyen-Thanh, X.~Zhuang, H.~Nguyen-Xuan, T.~Rabczuk, P.~Wriggers, A
  virtual element method for {2D} linear elastic fracture analysis, Computer
  Methods in Applied Mechanics and Engineering 340 (2018) 366--395.
\newblock \href {https://doi.org/10.1016/j.cma.2018.05.021}
  {\path{doi:10.1016/j.cma.2018.05.021}}.

\bibitem{vanHuyssteen2022}
D.~{van Huyssteen}, F.~L. Rivarola, G.~Etse, P.~Steinmann, On mesh refinement
  procedures for the virtual element method for two-dimensional elastic
  problems, Computer Methods in Applied Mechanics and Engineering 393 (2022)
  114849.
\newblock \href {https://doi.org/10.1016/j.cma.2022.114849}
  {\path{doi:10.1016/j.cma.2022.114849}}.

\bibitem{Huyssteen2022}
D.~{van Huyssteen}, F.~L. Rivarola, G.~Etse, P.~Steinmann, On mesh refinement
  procedures for polygonal virtual elements, arXiv.org (Jul. 2022).
\newblock \href {http://arxiv.org/abs/2207.03792} {\path{arXiv:2207.03792}},
  \href {https://doi.org/10.48550/ARXIV.2207.03792}
  {\path{doi:10.48550/ARXIV.2207.03792}}.

\bibitem{Berrone2021}
S.~Berrone, A.~Borio, A.~D{\textquotesingle}Auria, Refinement strategies for
  polygonal meshes applied to adaptive {VEM} discretization, Finite Elements in
  Analysis and Design 186 (2021) 103502.
\newblock \href {https://doi.org/10.1016/j.finel.2020.103502}
  {\path{doi:10.1016/j.finel.2020.103502}}.

\bibitem{Choi2020}
H.~Choi, H.~Chi, K.~Park, G.~H. Paulino, Computational morphogenesis:
  Morphologic constructions using polygonal discretizations, International
  Journal for Numerical Methods in Engineering 122~(1) (2020) 25--52.
\newblock \href {https://doi.org/10.1002/nme.6519}
  {\path{doi:10.1002/nme.6519}}.

\bibitem{Huyssteen2024}
D.~van Huyssteen, F.~L. Rivarola, G.~Etse, P.~Steinmann, On adaptive mesh
  coarsening procedures for the virtual element method for two-dimensional
  elastic problems, Computer Methods in Applied Mechanics and Engineering 418
  (2024) 116507.
\newblock \href {https://doi.org/10.1016/j.cma.2023.116507}
  {\path{doi:10.1016/j.cma.2023.116507}}.

\bibitem{Artioli2017}
E.~Artioli, L.~{Beir\~{a}o da Veiga}, C.~Lovadina, E.~Sacco, Arbitrary order
  {2D} virtual elements for polygonal meshes: part {I}, elastic problem,
  Computational Mechanics 60~(3) (2017) 355--377.
\newblock \href {https://doi.org/10.1007/s00466-017-1404-5}
  {\path{doi:10.1007/s00466-017-1404-5}}.

\bibitem{Gain2014}
A.~L. Gain, C.~Talischi, G.~H. Paulino, On the virtual element method for
  three-dimensional linear elasticity problems on arbitrary polyhedral meshes,
  Computer Methods in Applied Mechanics and Engineering 282 (2014) 132--160.
\newblock \href {https://doi.org/10.1016/j.cma.2014.05.005}
  {\path{doi:10.1016/j.cma.2014.05.005}}.

\bibitem{Veiga2015}
L.~{Beir{\~{a}}o da Veiga}, C.~Lovadina, D.~Mora, A virtual element method for
  elastic and inelastic problems on polytope meshes, Computer Methods in
  Applied Mechanics and Engineering 295 (2015) 327--346.
\newblock \href {https://doi.org/10.1016/j.cma.2015.07.013}
  {\path{doi:10.1016/j.cma.2015.07.013}}.

\bibitem{PolyMesher}
C.~Talischi, G.~H. Paulino, A.~Pereira, I.~F.~M. Menezes, {PolyMesher}: a
  general-purpose mesh generator for polygonal elements written in {M}atlab,
  Structural and Multidisciplinary Optimization 45~(3) (2012) 309--328.
\newblock \href {https://doi.org/10.1007/s00158-011-0706-z}
  {\path{doi:10.1007/s00158-011-0706-z}}.

\end{thebibliography}
	
\end{document}